\documentclass[12pt]{amsart}

\usepackage[T1]{fontenc}
\usepackage{lmodern}
\usepackage{textcomp}

\usepackage{amsmath,amsthm,amssymb}
\usepackage{mathtools}
\usepackage{bm}
\usepackage{mathrsfs}
\usepackage{dsfont}
\usepackage{bbm}
\usepackage{leftidx}

\usepackage{graphicx}
\usepackage{tikz}
\usepackage{pdfpages}
\usepackage{xcolor}

\usepackage{fancyhdr}
\usepackage{calc}

\usepackage{enumitem}

\usepackage{comment}
\usepackage{upref}

\usepackage[centering,width=6in,truedimen]{geometry}
\usepackage[colorlinks,pdfstartview=FitH]{hyperref}

\hypersetup{
    colorlinks,
    citecolor=blue,
    filecolor=black,
    linkcolor=red,
    urlcolor=black
}

\theoremstyle{plain}
\newtheorem{theorem}{Theorem}[section]
\newtheorem{proposition}[theorem]{Proposition}
\newtheorem{lem}[theorem]{Lemma}
\newtheorem{Mlem}{Main Lemma}
\newtheorem{cor}[theorem]{Corollary}

\newtheorem*{ques*}{Question}

\theoremstyle{definition}
\newtheorem{definition}[theorem]{Definition}

\theoremstyle{remark}
\newtheorem{rem}[theorem]{Remark}
\DeclareMathOperator{\dist}{dist}
\newcommand{\D}{\mathop{}\!\mathrm{d}}
\newcommand{\N}{\mathbb N}
\newcommand{\Q}{\mathbb Q}
\newcommand{\R}{\mathbb R}
\newcommand{\Z}{\mathbb Z}
\newcommand{\C}{\mathbb C}

\newcommand{\E}{\mathbb E}
\newcommand{\F}{\mathbb F}
\newcommand{\G}{\mathbb G}
\newcommand{\K}{\mathbb K}
\newcommand{\T}{\mathbb T}
\newcommand{\V}{\mathbb V}
\newcommand{\W}{\mathbb W}
\newcommand{\CC}{\mathcal C}

\newcommand{\EE}{\mathcal E}
\newcommand{\HH}{\mathcal H}
\newcommand{\LL}{\mathbb L}

\newcommand{\VV}{\mathcal V}
\newcommand{\WW}{\mathcal W}
\newcommand{\XX}{\mathcal X}
\newcommand{\YY}{\mathcal Y}

\newcommand{\Sb}{\mathbb{S}}

\newcommand{\eps}{\varepsilon}
\newcommand{\1}{\vec{\bm 1}}
\newcommand{\vt}{{\vec{t}\,}}
\newcommand{\one}{{\vec{1}}}
\newcommand{\Gr}{\mathrm{Gr}}

\newcommand{\dir}{\mathrm{dir}}

\newcommand{\Span}{\mathrm{span}}

\newcommand{\norm}[1]{\left\|#1\right\|}

\numberwithin{equation}{section}

\newcommand{\on}{\operatorname}

\newcommand{\smblkcircle}{\bullet}

\begin{document}

\title[The Kakeya problem for codimension-two planes in $\R^n$]{The Kakeya problem for \\ codimension-two planes in $\R^n$}

\author{Chun-Kit Lai}
\address{Department of Mathematics, San Francisco State University, 1600 Holloway Ave, San Francisco, CA 94132
}
\email{cklai@sfsu.edu}

\author{Yixuan Pang}
\address{Department of Mathematics, University of Pennsylvania, Philadelphia, PA 19104, United States}
\email{pyixuan@sas.upenn.edu}

\author{Rahul Sethi}
\address{School of Mathematics, Georgia Institute of Technology, Atlanta, GA 30332, United States}
\email{rahul.sethi@math.gatech.edu}

\author{Chenjian Wang}
\address{Department of Mathematics, University of British Columbia, Vancouver, BC V6T 1Z2, Canada}
\email{chjwang@math.ubc.ca}


\begin{abstract}
We study the $(n-2)$-plane Kakeya problem in $\R^n$ ($n\geq4$) subject to convex Wolff axioms. In other words, we study finite families of $\delta$-neighborhoods of $(n-2)$-plane segments ($\delta$-flakes) in $\R^n$, with the property that not too many of them can be contained inside a common convex set. Our main result is that the union of $\delta$-flakes from such a set must have
almost maximal volume. This can be viewed as a higher-dimensional generalization of Wang and Zahl's work for $\delta$-tubes in $\R^3$. The latter implies the Kakeya set
conjecture in $\R^3$. As a consequence, we resolve Rogers' $2$-plane Kakeya problem in $\R^4$ subject to certain curvature conditions on two-dimensional families of direction subspaces in the Grassmannian $\Gr(4,2)$.
\end{abstract}

\maketitle

\setcounter{tocdepth}{2}
\tableofcontents

\section{Introduction}\label{Sec:introduction}

Let $\Gr(n,k)$ be the Grassmannian consisting of all $k$-subspaces (\textit{$k$-planes}) in $\R^n$, and let $\{ t+V: t\in \R^n, V \in \Gr(n,k) \}$ be the collection of all \textit{affine $k$-planes} in $\R^n$. To avoid triviality, we will  assume $n \geq 3$ throughout the paper.  

A \textit{$(\delta,k)$-flake} is the $\delta$-neighborhood of a unit cube in some affine $k$-plane. When $k=2$, we may simply call it a \textit{$\delta$-plate}. Suppose $\F$ is a set of $(\delta,k)$-flakes in $\R^n$. We write
\begin{align*}
    U(\F) \coloneqq \bigcup_{F\in\F} F.
\end{align*}
The recent solution of the Kakeya conjecture in $\R^3$ by the groundbreaking work of Wang and Zahl concerns $|U(\F)|$, the volume of $U(\F)$, when $n = 3$ and $k=1$. Roughly speaking,  if the distribution of $\F$ satisfies suitable nonconcentration conditions, then $|U(\F)|$ is almost as large as $|\F||F|$.  In this paper, our goal is to show that the same type of assumption also implies a lower bound for the volume of the union of flakes of codimension $2$ (i.e., $k = n-2$). 

To state our main results, we need to first fix our terminology.  We will  write $|\F|$ for the cardinality of $\F$, and $|F|$ for the volume of a $(\delta,k)$-flake $F \in \F$.\footnote{To make a distinction, we use $\norm{\cdot}$ to denote the operator norm or the Euclidean norm. For simplicity, we will often pretend that $|F| = \delta^{n-k}$ for a $(\delta,k)$-flake $F$, which is entirely harmless.} { For clarity, we sometimes write $F_\delta$ instead of $F$ to emphasize the scale.} In this case, we may simply say ``\textit{$\delta$-flake}'' instead of ``$(\delta, n-2)$-flake''. 

As in \cite{GWZ26streamlined}, we will describe the nonconcentration conditions for convex sets as follows. Suppose that $\W$ is a finite set of convex sets in $\R^n$ and suppose that $K$ is another convex set with $|K|>0$. We define
\begin{align*}
    \W[K] \coloneqq \{W\in\W: W \subset K\},
\end{align*}
and
\begin{align*}
    \Delta(\W, K) \coloneqq \frac{\sum_{W\in \W[K]}|W|}{|K|}.
\end{align*}
The maximal density measuring how the convex sets in $\W$ cluster into coarser convex sets  is defined  to be 
\begin{align*}
    \Delta_{max}(\W) \coloneqq \max_{K\text{ convex}} \Delta(\W,K).
\end{align*}

For a set of convex sets $\W$, we call a collection of Lebesgue measurable subsets $\{Y(W): W\in\W\}$ a \textit{shading} of $\W$ if $Y(W) \subset W$ for each $W\in\W$. We define 
\begin{align*}
    U(\W,Y) \coloneqq \bigcup_{W\in\W} Y(W).
\end{align*}
If the shading $Y$ is clear from the context, then we may also write $U(\W)$ for $U(\W,Y)$. We define the \textit{averaged density} of the shading $Y$ as 
\begin{align*}
    \lambda(\W,Y) \coloneqq \frac{\sum_{W \in \W} |Y(W)|}{\sum_{W \in \W} |W|}.
\end{align*}
For clarity, we also sometimes write $Y_\W$ and $\lambda(Y_\W)$ in place of $Y$ and $\lambda(\W,Y)$.

With these definitions, we can now state our main result:
\begin{theorem}[$(n-2)$-plane Kakeya in $\R^n$]\label{Thm:Main}
    For every $\eps>0$, there is some $\eta>0$ such that the following holds for all sufficiently small $\delta>0$.
    
    If $\F$ is a set of $\delta$-flakes in $B_1 \subset \R^n$ with $\Delta_{max}(\F) \leq \delta^{-\eta}$, and $Y$ is a shading on $\F$ with $\lambda(\F,Y) \geq \delta^\eta$, then 
    \begin{align*}
        |U(\F,Y)|
        \geq \delta^\eps |\F||F|.
    \end{align*}
\end{theorem}

  Wang and Zahl \cite{WZ26sticky,WZ25assouad,WZ25volume} settled the case of Theorem \ref{Thm:Main} when $n = 3$. We record it below for a comparison.
\begin{theorem}[Kakeya set conjecture in $\R^3$]\label{Thm:R3_Kakeya}
    For every $\eps>0$, there is some $\eta>0$ such that the following holds for all sufficiently small $\delta>0$.
    
    If $\T$ is a set of $\delta$-tubes in $B_1 \subset \R^3$ with $\Delta_{max}(\T) \leq \delta^{-\eta}$, and $Y$ is a shading of $\T$ with $\lambda(\T,Y) \geq \delta^\eta$, then 
    \begin{align*}
        |U(\T,Y)|
        \geq \delta^\eps |\T||T|.
    \end{align*}
    In particular, every Besicovitch set in $\R^3$ has Minkowski and Hausdorff dimension $3$. 
\end{theorem}

It is not hard to see that, by a straightforward tensor-product argument, Theorem \ref{Thm:Main} for $n>3$ implies Theorem \ref{Thm:R3_Kakeya}. Theorem~\ref{Thm:Main} tells us that as long as the ``codimension'' is $2$, nonconcentration with respect to convex sets continues to be an effective notion for the Kakeya problem. Furthermore, we notice that the analogous problem of Theorem \ref{Thm:Main} for $\Gr(n,n-1)$ has been resolved fairly easily in \cite[Lemma 5.7]{WZ25volume}\footnote{\color{black} See Lemma~\ref{Lem:slab_0} for the exact statement in our notation.} via a C\'ordoba-type $L^2$ argument, which is standard in resolving the Kakeya maximal conjecture in $\R^2$ (see \cite{Cordoba1977Kakeya}). Moreover, the special case for $\Gr(3,2)$ was used by Wang and Zahl in the same paper where they resolve the Kakeya conjecture in $\R^3$. This suggests that the resolution of the Kakeya conjecture in higher dimensions may also benefit from studying analogous statements in $\Gr(n,k)$ for all $2\le k\le n-1$. Our Theorem \ref{Thm:Main} resolves the case when $k = n-2$ without introducing any new nonconcentration conditions. 

Polynomial Wolff axioms first introduced by Guth and Zahl \cite{GuthZahl2018PolynomialWolff} and later established by Katz and Rogers \cite{KatzRogers2018PolynomialWolff} imply that direction-separated tubes cannot be concentrated nearby semialgebraic varieties. This nonconcentration condition is possibly needed in higher-codimensional cases (i.e., $1\le k<n-2$) (see \cite[Section~9]{guth2025introduction} for some related comments). {\color{black} It might be surprising that the Kakeya set conjecture in $\R^4$ for sticky sets of tubes was recently proved by Wang and Zakharov \cite{wang_zak2026_R4_kakeya} by only imposing the convex Wolff axioms. So it would be interesting to see whether new nonconcentration conditions are needed in non-sticky cases.}

Our result is undoubtedly related to different generalizations of higher-dimensional Kakeya-Besicovitch sets. A natural one to recall here is the \textit {$(n,k)$-Besicovitch set}, which is a measure zero set in $\R^n$ such that it contains a translate of any subspaces in $\Gr(n,k)$. Currently, it is still an open problem to decide if an $(n,k)$-Besicovitch set exists for $k>1$. Falconer \cite{Falconer1980Besicovitch} proved that no such set exists if $k>n/2$. Bourgain \cite{Bourgain1991Besicovitch} proved that no such set exists if 
$$
2^{k-1}+k\ge n.
$$
Oberlin \cite{Oberlin2010Xray} later extended the range to $(1+\sqrt{2})^{k-1} + k > n$. In particular, combining the results of Falconer and Bourgain, $(n,n-2)$-Besicovitch sets have been known to have positive Lebesgue measure. However, this type of Besicovitch sets fails our nonconcentration assumption in Theorem \ref{Thm:Main} whenever $n \geq 4$. Indeed, a direct counting shows that $\Delta_{max} (\F) \gtrsim \delta^{6-2n}$ for any discretization of $(n,n-2)$-Besicovitch sets, where $\F$ is a (direction-separated) $\delta$-discretized  version of the $(n-2)$-planes.

\begin{rem}\label{Rem:k-plane_history}
For more results on $(n,k)$-Besicovitch sets, we refer the reader to \cite{oberlin2007} for their Hausdorff dimension, \cite{FraserHarrisKroon2022Fourier} for their Fourier dimension, and Mattila's book \cite[Chapter~24]{mattila2015} for a general exposition. It is worth mentioning that the $k$-plane problem can also be formulated in terms of the dimension of the family in the affine Grassmannian, see \cite{gan2026}.
\end{rem} 

Because of the non-existence conjecture, another  higher-dimensional variant was introduced by Rogers \cite{rogers2006}, which turns out to be more relevant to our study. More precisely, he introduced a \textit{restricted} variant of $(n,2)$-Besicovitch sets, which only involves translates of every member of some $(n-2)$-dimensional ``curved'' submanifold in $\Gr(n,2)$. When $n=4$, it was proved that such Kakeya-type sets have Hausdorff dimension at least $7/2$. In Section \ref{Sec:applications}, we will show that the ``curved'' condition naturally extends to $\Gr(n,n-2)$, and indeed guarantees our nonconcentration condition on the $\delta$-flakes. Theorem \ref{Thm:Main} therefore resolves Rogers' Kakeya problem for $\Gr(4,2)$, and actually $\Gr(n,n-2)$ in general, i.e., such Kakeya-type sets have full Hausdorff dimension (see Section \ref{Sec:applications} for the precise definitions and details).

In fact, Rogers' formulation also works for general $k$-planes, i.e., the direction subspaces of them come from some $(n-k)$-dimensional submanifold $\mathcal{M} \subset \Gr(n,k)$. And a natural problem is to impose suitable conditions on $\mathcal{M}$ so that the $k$-planes are essentially disjoint, i.e., we can expect full dimension for the Kakeya-type sets.  So if we are interested in the Hausdorff dimension of the Kakeya-type sets that can possibly be Lebesgue measure zero, then the dimension complementarity $k + \dim \mathcal{M} = n$ is in some sense a more appropriate condition than considering the whole Grassmannian. For $(n,k)$-Besicovitch sets, we are in the ``supercritical'' regime when $k>1$ as $k + \dim \Gr(n,k) > n$. Recall from the last paragraph that our Theorem~\ref{Thm:Main} provides a satisfactory answer when $k=n-2$. 

For general $k$, although we have seen the necessity of the dimension condition $\dim \mathcal{M} = n-k$, it is still unclear what kind of curvature condition on $\mathcal{M}\subset \Gr(n,k)$ is suitable for Kakeya-type problems, and also for the stronger Kakeya maximal estimate associated to $\mathcal{M}$. A popular way to approach this is to make the additional assumption that $\mathcal{M}$ is induced by the normal bundle of a surface $\Sigma$ of codimension $k$ in $\R^n$, which allows us to encode the curvature property of $\mathcal{M}$ by that of $\Sigma$. See \cite{DMV} for partial results on boundedness of the restricted $2$-plane transform related to ``well-curved'' quadratic surfaces of codimension $2$. On the contrary, our main results only provide information about dimension bounds.

Finally, in view of the classical $k=1$ case, it is reasonable to believe that the study of restricted $k$-plane transforms, or even the restricted $k$-plane Kakeya sets, may facilitate that of higher-codimensional Fourier restriction. We will not discuss this in detail, but we hope to see progress in this direction in future work. Interested readers  may consult \cite[Section~1]{CMP2026} for references on Fourier restriction of higher-codimensional quadratic surfaces, which is the most extensively investigated in the literature.

\subsection{Overview of the Proof} 
Motivated by the work of Rogers \cite{rogers2006}, in Section~\ref{Sec:geometry}, we will introduce the concept of \textit{principal angles} to study the intersection behavior of the $\delta$-flakes. This serves as the basic language to describe the overlapping patterns. Our argument for $\Gr(n,n-2)$ will be based on induction on $n$, whose base case ($n = 3$) has been established by Theorem \ref{Thm:R3_Kakeya}. However, as the geometry of $\delta$-flakes in $\R^n$ ($n\geq4$) is more complicated than that of $\delta$-tubes in $\R^3$, there are many additional difficulties in executing the induction. In fact, we will need both the proof and the statement of Theorem~\ref{Thm:R3_Kakeya}, combined with some new geometric ideas. In this sense, our Theorem~\ref{Thm:Main} should be regarded as a nontrivial extension of Theorem~\ref{Thm:R3_Kakeya} instead of a brand new way to prove it.

From now on, we may assume $n \geq 4$. Our proof strategy of Theorem \ref{Thm:Main} is to divide the proof into sticky and non-sticky cases. 
 
We will study the sticky case first in this paper, as it contains most of the novelties. One obvious technical obstacle is that Wang-Zahl's proof of the sticky Kakeya conjecture for tubes in $\R^3$ \cite{WZ26sticky} relies heavily on existing deep results from projection theory, whose higher-dimensional versions seem to be needed here but unfortunately not yet systematically established in the literature. It is possible that one can still pursue along the projection-theoretic approach with additional work, but we decided to adopt a slicing approach to get around the technical projection theory parts. 
 
Our main idea and contribution here is to slice the $\delta$-flakes by parallel affine $3$-planes in $\R^n$, such that we have a set of $\delta$-tubes satisfying assumptions in Theorem~\ref{Thm:R3_Kakeya} inside each of the slices. If this step  is achievable, then we can simply apply Theorem~\ref{Thm:R3_Kakeya} in each slice and then use a Fubini-type argument to finish the proof. The difficulty is that we need to somehow find a way to relate the nonconcentration condition for $\delta$-flakes in $\R^n$ to that for $\delta$-tubes in each slice, which seems to be unclear in general due to the complicated geometry. To overcome this difficulty, we will build a novel induction scheme by applying a factoring lemma due to Wang and Zahl \cite{WZ25volume} inside each slice. This would induce a \textit{medium grains} structure in each slice, which, combined with the local grains structure provides favourable rigidity. We hope that the scheme initiated here will also be instructive to many other Kakeya-type problems where the slicing philosophy potentially applies.

Once the sticky case is settled, we will then reduce any general configurations of the $\delta$-flakes to the {sticky} case. To this end, we will mainly follow the streamlined presentation in \cite{GWZ26streamlined}. It turns out that there are not many additional difficulties in adapting their argument to $\delta$-flakes. However, due to the change of the dimensions, the parameters in the argument need to be modified accordingly. Moreover, the multiscale decomposition and the two-ends reduction need to be clarified more carefully. We decide to include the whole proof of the non-sticky part in this paper, so that readers can identify the difference between our proof and the original proof in \cite{GWZ26streamlined}.

As a historical side remark, the idea of general-to-sticky reduction was first systematically proposed by Katz and Tao \cite{TaoBlogKakeya}, and it has been developed carefully by many authors in other problems such as the Furstenberg sets conjecture \cite{OS2026, ren2023furstenberg}. One may consult \cite[Section~1.4]{WZ25volume} and \cite[Section~4]{guth2025introduction} for a more detailed history. It was not until recently that Wang and Zahl \cite{WZ25volume} successfully fleshed out the original Katz--Tao program \cite{TaoBlogKakeya} and proved Theorem~\ref{Thm:R3_Kakeya}. It is unclear how much this methodology can be adapted to other Kakeya-type problems, and our Theorem~\ref{Thm:Main} is one step in this direction.

One may wonder whether the slicing argument can be applied directly to a general set of $\delta$-flakes (not necessarily sticky), so that there is no need to first reduce to the sticky case and re-derive the proof in \cite{GWZ26streamlined}. We do not know the answer. The current proof relies heavily on the rigid graininess structure enjoyed by the sticky case, similarly to what happens in \cite{WZ26sticky}. Therefore, searching for a direct slicing proof without the stickiness assumption would be very interesting. One additional advantage of our approach is that Kakeya estimates for certain planks in $\R^n$ (satisfying convex Wolff axioms) are also obtained as by-products of the proof of the general-to-sticky reduction (Section~\ref{subsec:apply_plank_Rn}). These estimates generalize Theorem~\ref{Thm:Main} in a natural way.

\subsection{Organization of the paper}
Sections~\ref{Sec:setup}--\ref{Sec:refinements} provide the preliminaries. Section~~\ref{Sec:setup} introduces notation and definitions, including uniformity (Section~\ref{subsec:uniform}). Section~\ref{Sec:geometry} presents the principal-vector decomposition (Theorem~\ref{Thm:Principal_angle}) and develops geometric lemmas for $k$-planes and $(\delta,k)$-flakes. Section~\ref{Sec:organize} introduces the Katz--Tao and Frostman nonconcentration conditions (Definitions~\ref{Def:Katz_Tao}~and~\ref{Def:Frostman}), records the maximal density factoring lemma (Lemma~\ref{Lem:factoring}), and gives a codimension-one analogue of Theorem~\ref{Thm:Main} (Lemma~\ref{Lem:slab_0}). Section~\ref{Sec:refinements} develops refinement lemmas compatible with factoring and multiplicity estimates (Proposition~\ref{Prop:fact_mult} and Lemma~\ref{Lem:fact_mult}). Readers familiar with these preliminaries may skim these sections and return to them as needed.

Sections~\ref{Sec:stronger_sticky}--\ref{Sec:finish_sticky} establish the sticky Kakeya theorem for $\delta$-flakes (Theorem~\ref{Thm:2-plane_sticky}), with the main new ingredients appearing in Sections~\ref{Sec:plaininess}--\ref{Sec:medium_grains_end}. Section~\ref{Sec:stronger_sticky} introduces the relevant notions of stickiness (Definitions~\ref{Def:sticky}~and~\ref{Def:strong_sticky}) and formulates Theorem~\ref{Thm:2-plane_sticky} and its special case, Theorem~\ref{Thm:Main_Sticky}. Section~\ref{Sec:self_similar} sets up the extremal-family framework (Definitions~\ref{Def:sigma}~and~\ref{Def:eps_extremize}) and establishes self-similarity properties (Proposition~\ref{Prop:self-similar}), which will be used throughout the sticky argument. Section~\ref{Sec:plaininess} proves the planiness and the local grains structure (Proposition~\ref{Prop:local_grains}). The proof of Theorem~\ref{Thm:Main_Sticky} is then organized into ten steps, outlined at the beginning of Section~\ref{Sec:slicing_set_up}. Section~\ref{Sec:slicing_set_up} develops the slicing argument that produces medium grains (Steps~(i)--(vi)), while Section~\ref{Sec:medium_grains_end} combines medium grains with local grains to complete the proof of Theorem~\ref{Thm:Main_Sticky} (Steps~(vii)--(x)). Section~\ref{Sec:finish_sticky} then deduces Theorem~\ref{Thm:2-plane_sticky} from Theorem~\ref{Thm:Main_Sticky}.

Sections~\ref{Sec:partial_estimates}--\ref{Sec:Main_Lemma_2} establish Theorem~\ref{Thm:Main} through a reduction to the sticky case. Section~\ref{Sec:partial_estimates} introduces the assertions $K_{KT}(\beta)$ and $K_F(\beta)$ (Definitions \ref{Def:Katz_Tao_beta}~and~\ref{Def:Frostman_beta}) and reduces Theorem~\ref{Thm:Main} to Main Lemmas~\ref{Main_lem_1}~and~\ref{Main_lem_2}, whose iteration yields the desired conclusion. A roadmap for their proofs appears immediately before Section~\ref{subsec:large_CWA}. Sections~\ref{Sec:slab_plank}~and~\ref{Sec:reduce_to_sticky} are preparatory. Section~\ref{Sec:slab_plank} extends the multiplicity estimates to planks (Lemmas~\ref{Lem:plank_Katz_Tao}~and~\ref{Lem:plank_Frostman}), and to families of flakes factored by planks (Proposition~\ref{Prop:plank_factor_improve}). Section~\ref{Sec:reduce_to_sticky} establishes the multiscale decomposition (Lemma~\ref{Lem:multi-scale_decomposition}), which provides the sticky and structured non-sticky alternatives used in the reduction. Sections~\ref{Sec:Main_Lemma_1}~and~\ref{Sec:Main_Lemma_2} prove Main Lemmas~\ref{Main_lem_1}~and~\ref{Main_lem_2}, respectively, completing the proof of Theorem~\ref{Thm:Main}.

Section~\ref{Sec:applications} presents applications of the main results. Section~\ref{subsec:apply_Rogers_R4} applies Theorem~\ref{Thm:Main} to Rogers-type $(n-2)$-plane Kakeya sets in $\R^n$ under the stated curvature assumptions (Theorem~\ref{thm:rogers-application}), including Rogers' $2$-plane Kakeya problem in $\R^4$. Section~\ref{subsec:apply_plank_Rn} records the resulting Kakeya estimates for planks (Theorems~\ref{Thm:Rn_plank_Katz_Tao}~and~\ref{Thm:Rn_plank_Frostman}). Appendix~\ref{Appendix:Probability} collects the probabilistic tools and proves Lemmas~\ref{Lem:random_translate},~\ref{Lem:random_sub_(F)}~and~\ref{Lem:random_rigid_F}.

\section{Notation and Definitions}\label{Sec:setup}

Our convention basically follows that of \cite{GWZ26streamlined}.

\subsection{Linear algebra} 

For a finite set of vectors $\{v_j\}_{j=1}^k$ in $\R^n$, we use $\Span(v_1, \ldots, v_k)$ to denote the linear space they span. For simplicity, if there is a single vector $v \in \R^n$, then we define $\langle v \rangle \coloneqq \Span(v)$. This should not be confused with inner products $\langle \cdot, \cdot \rangle$.

For $1 \leq k_1 \leq k_2 \leq n$ and $V \in \Gr(n, k_2)$, we define $$\Gr(V, k_1) \coloneqq \{W \in \Gr(n,k_1): W \subset V\}.$$

Let $\pi_V$ be the orthogonal projection from $\R^n$ onto its subspace $V$. Let $\mathcal{I}$ be the identity operator on $\R^n$. Let $|v_1 \wedge v_2 \wedge \ldots \wedge v_k|$ be the $k$-dimensional volume of the parallelepiped generated by $\{v_j\}_{j=1}^k \subset \R^n$ ($k \leq n$).

\subsection{Convex sets and shadings}

Given a family of convex sets $\V$, we will often abbreviate expressions like $\sum_{V \in \V}f(V)$ as $\sum_{\V}f(V)$ for simplicity. Also, when all convex sets in $\V$ have comparable volume, we will write $|V|$ to denote the volume of one of them.

Let $n \geq 3$. The \textit{dimensions} of a convex set $V \subset \R^n$ are the lengths of the axes of its inner John
ellipsoid, listed in increasing order. We will denote this by $v_1 \times v_2 \times \cdots \times v_n$. {\color{black} To ease notation, we do not distinguish comparable scales when writing $v_1 \times v_2 \times \cdots \times v_n$, which is harmless. For simplicity, we use the boldface vector ``$\1$'' to denote the product of $(n-2)$ many ``$1$'', and ``$v\1$'' to denote the product of $(n-2)$ many ``$v$''.} A \textit{$\delta$-tube} is a convex set of dimension $\delta\1 \times \delta \times 1$. A \textit{$\delta$-plate} is a convex set of dimension $\delta\1 \times 1 \times 1$. A \textit{$\delta$-flake} is a convex set of dimension $\delta\times\delta\times\1$. A \textit{slab} is a convex set of dimension $\theta\times1\times\1$ for some $0\leq \theta\leq 1$. A \textit{plank} is a convex set of dimension $a\times b\times\1$ for some $0 < a \leq b \leq 1$. {\color{black} In general, we use $N_r(E)$ to denote the $r$-neighborhood of a set $E \subset \R^n$.}

Let $B^n(x,r)$ denote a ball of radius $r$ centered at $x$ in $\R^n$. For simplicity, we may also use the notation $B_r^n$ to denote a general ball of radius $r$ in $\R^n$, whose center is either $0$ or unimportant\footnote{It can be read from the context which of the two options holds.}. We omit the superscript $n$ when it is clear from the context. An additional convention is that $B_1 \coloneqq B^n(0,1)$. 

Let \(C>0\) and let \(K\) be a symmetric convex body in \(\mathbb{R}^n\) centered at \(x\) (i.e., if \(x+y\in K\) for some \(y\in\mathbb{R}^n\), then \(x-y\in K\)). For example, \(K\) may be a tube, a flake, a plank, a cube, or a ball. Then we use \(CK\) to denote the dilation of \(K\) with respect to \(x\) by a factor \(C\), i.e.,
\[
CK \coloneqq \{y : (y-x)/C \in K\}.
\]
Here \(C\) is allowed to depend on \(\delta\) (e.g., \(\log(1/\delta)\)), and $x$ will often be implicit when we write $CK$.

{\color{black} In Section~\ref{Sec:slicing_set_up}~and~\ref{Sec:medium_grains_end}, due to more delicate geometry, we use the vector ``$\one$'' to denote the product of $(n-3)$ many ``$1$'', and ``$v\one$'' to denote the product of $(n-3)$ many ``$v$''.}

A \textit{shading} of a set $V$ is a measurable set $Y(V) \subset V$. We will use $(\V, Y)$ to denote a family
of convex sets and their associated shadings $\{Y(V): V\in\V\}$. For brevity, we will sometimes say
``let $(\V, Y)$ be a family of convex subsets of $\R^n$''. Define the \textit{multiplicity} of a family of convex sets $(\W, Y)$ as 
\begin{align*}
    \mu(\W,Y) \coloneqq \frac{\sum_{W\in\W}|Y(W)|}{|U(\W,Y)|}.
\end{align*}
{\color{black} When $Y$ is clear, we may write $\mu(\W)$ instead.} Note that the conclusion of Theorem~\ref{Thm:Main} can be equivalently replaced by $\mu(\F) \leq \delta^{-\eps}$.

For $x\in\R^n$, we define $\V_Y(x) \coloneqq \{V\in\V: x\in Y(V)\}$ and $\mu(\V,Y)(x) \coloneqq |\V_Y(x)|$. We say that $(\V,Y)$ has \textit{constant multiplicity} if $\mu(\V,Y)(x)$ is comparable for all $x\in U(\V,Y)$. We say that $(\V',Y')$ is a \textit{refinement} of $(\V,Y)$ if $\V'\subset \V$ and $Y'(V) \subset Y(V)$ for each $V \in \V'$. For $c>0$, we say that $(\V',Y')$ is a \textit{$c$-refinement} of $(\V,Y)$ if $\sum_{\V'}|Y'(V)| \geq c \sum_{\V}|Y(V)|$ additionally holds. Finally, if $Y$ is a shading on $V$ and $\V'$ is a subset of $\V$, then we will often abuse notation and continue to use $Y$ to refer to the restriction of the shading to $\V'$.

\subsection{Pigeonholing and \texorpdfstring{$\lessapprox$}{⪅} notation}
There are many different types of dyadic pigeonholing steps in the proof, and one may consult \cite[Section~2.2]{GWZ26streamlined} for a list of them. These steps often involve harmless quantities of the form $\log(1/\delta)$ or $\log|\V|$. Also, in practice, any shading $Y(V)$ with small measure (i.e., $|Y(V)| \leq \delta^n$) can be replaced by $Y(V) = \varnothing$ without affecting our
arguments. To streamline our arguments, if $X$ and $Y$ are functions of a scale parameter $\rho$, then we use the notation ``$X\lessapprox Y$'' to mean that for any $\eps>0$, there is a constant $C = C(\eps)$ such that $X \leq C \rho^{-\eps} Y$. So $\rho$ will often be implicit. One can easily check that such notation reasonably applies to all pigeonholing steps, see \cite[Section~2.2]{GWZ26streamlined} for details. {\color{black} Throughout the paper, unless otherwise specified, we will use $C$ to denote a constant $\lesssim 1$, whose exact meaning may vary from line to line.}

\subsection{Uniform sets of flakes}\label{subsec:uniform} {\color{black} We say two flakes $F,F'$ are \textit{essentially distinct} if $|F \cap F'| < c_0|F|$ for some fixed small constant $0 < c_0 < 1/2$}, and a family of flakes is \textit{essentially distinct} if each pair of its elements is.  It will be useful to consider $\delta$-flakes $\F$ that are uniform in
the following sense.
\begin{definition}\label{Def:uniform_flakes}
    Let $\delta>0$, and $\F$ be a set of $\delta$-flakes in $B_1 \subset \R^n$. For $\rho \in [\delta,1]$, we say that $\F$ is \textit{uniform at scale $\rho$} if there is a set of $\rho$-flakes $\F_\rho$ so that     \begin{enumerate}[label={\rm (\roman*)}]
        \item ${\F} = \bigcup_{F_\rho \in \F_{\rho}} {\F}[F_{\rho}]$.
        \item The flakes ${\F}_{\rho}$ are essentially distinct.
        \item $|{\F}[F_{\rho}]|$ is constant up to a factor $\sim 1$ as $F_{\rho}$ varies in ${\F}_{\rho}$.
    \end{enumerate}
    
    Moreover, let $M = \lceil \log\log(1/\delta) \rceil$. We say that $\F$ 
    is \textit{uniform} if $\F$ is uniform at scale $\rho_k = \delta^{k/M}$ for each $k = 0,\ldots, M$. In this case, we call the numbers $N_k = |\F[F_{\rho_k}]|$ the \textit{branching numbers} of $\F$. 
\end{definition}
Note that by definition, if $\F$ is uniform, then it must be essentially distinct. Item (ii) above implies that each $F \in \F$ can only lie in $F_\rho$ for $\sim1$ choices of $F_\rho \in \F_\rho$. {\color{black} For $\F$ not necessarily uniform, we may still use $\F_\rho$ to denote a set of $\rho$-flakes with the smallest cardinality that satisfies Item~(i) above. We may call $\F_\rho$ a set of \textit{thickened flakes} at scale $\rho$.}

\begin{definition}\label{Def:uniform_shading}
    Let $\delta>0$, $M = \lceil \log\log(1/\delta) \rceil$, and $(\F,Y)$ be a set of $\delta$-flakes in $B_1 \subset \R^n$ and their associated shading. We say that $(\F,Y)$ 
    is \textit{uniform} if the following holds:
    \begin{itemize}
        \item $\F$ is uniform.
        \item For each point $x \in U(\F,Y)$, the sets $\F_Y(x)$ are uniform, and have roughly the same branching numbers.
    \end{itemize}
\end{definition}

The reason for setting $M \sim \log \log(1/\delta)$ in our above definitions is as follows. First, since $(\log|\F|)^{\log\log(1/\delta)} \lessapprox 1$, by pigeonholing, if $\F$ is an essentially distinct family of $\delta$-flakes in $B_1 \subset \R^n$, then it has a subset $\F'$ of cardinality $|\F'| \gtrapprox |\F|$ that is uniform, and every pair $(\F,Y)$ has a $\approx 1$-refinement that is uniform. 

Second, if $\F$ is uniform, then for \textit{every} scale $\rho \in [\delta,1]$, there is a set of $\rho$-flakes $\F_\rho$ that satisfies Item (i) and (ii) from Definition~\ref{Def:uniform_flakes}, and also the following slightly
weaker analogue of Item (iii): there is a number $N_\rho$ so that $|\F[F_\rho]| \approx N_\rho$ for all $F_\rho \in \F_\rho$. Similarly for
Definition~\ref{Def:uniform_shading}; the second bullet point is weakened to: the branching numbers $N_k$ are roughly
the same up to a factor $\approx 1$.

In Section~\ref{Sec:finish_sticky}, we will extend the meaning of ``uniform'' so that certain convex Wolff axiom constants are also uniformized (Proposition~\ref{Prop:CWA_uniform}).

\section{Geometry of \texorpdfstring{$k$-planes}{k-planes}}\label{Sec:geometry}
In this section, we give some basic geometric lemmas for $k$-planes and $(\delta,k)$-flakes in $\R^n$, {\color{black} which will be used throughout the paper}. In our proof of Theorem~\ref{Thm:Main} (especially the sticky case) we need to describe and analyze the overlapping pattern between two $k$-planes. The main tool is called the \textit{principal-vector decomposition}, which has been systematically studied by many people, see \cite{BG73, Halmos1969, MBI1992, wedin2006, MCspark25}.

\subsection{Principal-vector decomposition and the major angle}

Here we provide a version of such a decomposition that will be useful for our purposes. It is not presented in exactly the same form as any of the results in the references as mentioned, but is indeed an organic combination of them. 
\begin{theorem}[principal-vector decomposition]\label{Thm:Principal_angle}
    For any $V, W \in \Gr(n,k)$, there exist orthonormal bases $$\mathcal{V} \coloneqq \{v_1, \dots, v_k\} \subset V, \quad \mathcal{W} \coloneqq \{w_1, \dots, w_k\} \subset W,$$ and a set of angles $0 \leq \theta_1 \leq \ldots \leq \theta_k \leq \pi/2$, such that $\langle v_i, w_j \rangle = \delta_{ij}\cos\theta_i$, where $\delta_{ij}$ is the Kronecker delta. In particular, $\VV$ and $\WW$ correspond to the SVD (singular value decomposition) of $\pi_W|_V: V \rightarrow W$ and $\pi_V|_W: W \rightarrow V$, in the sense that $\pi_W(v_i) = (\cos\theta_i) w_i$ and $\pi_V(w_i) = (\cos\theta_i) v_i$ for each $i$. 
    
    We will call $\theta_1, \dots, \theta_k$ the \textit{principal angles} between $V$ and $W$, and call $\VV$ and $\WW$ \textit{principal vectors} of $V$ and $W$. The former is uniquely determined by $V$ and $W$, while the latter may not be. The whole structure is called a \textit{principal-vector decomposition}, which also enjoys the following properties:
    \begin{enumerate}[label=\rm{(\roman*)}]
        \item 
        $\norm{\pi_V - \pi_W} = \norm{(\mathcal{I}-\pi_W)|_V} = \norm{(\mathcal{I}-\pi_V)|_W} = \sin\theta_k$.

        \item 
        If $\dim (V \cap W) = l$, then $\theta_i = 0$ for $i \leq l$ and $\theta_i > 0$ for $i > l$. Moreover, $|(v_{l+1}\wedge \ldots \wedge v_k) \wedge (w_{l+1}\wedge \ldots \wedge w_k)| = \sin\theta_{l+1} \cdots \sin\theta_k$.
        
        \item The nonzero principal angles between $V^{\perp}$ and $W^{\perp}$ are exactly the same as the nonzero principal angles between $V$ and $W$.
    \end{enumerate}
    We may write $\theta_i$ as $\theta_i(V,W)$ to emphasize the subspaces. 
\end{theorem}

\begin{rem}\label{Rem:major_equiv}
Geometrically, the term $\norm{(\mathcal{I}-\pi_W)|_V}$ in Item~(i) is the largest distance from a point in $V \cap B_1$ to $W$, which is equal to the smallest $\rho$ such that $V\cap B_1 \subset N_\rho(W) \cap B_1$, or $V\cap B_1 \subset N_\rho(W\cap B_1)$. And similarly for $\norm{(\mathcal{I}-\pi_V)|_W}$. Therefore, Item~(i) means that if $V\cap B_1 \subset N_\rho(W) \cap B_1$, then $W\cap B_1 \subset N_\rho(V) \cap B_1$, and vice versa. This observation allows us to describe the inclusion between $(\delta,k)$-flakes at different scales.
\end{rem}

\begin{definition}\label{Def:flake_direction}
    For a $(\delta,k)$-flake $F$ in $\R^n$, we call the unique $V \in \Gr(n,k)$ parallel to its defining affine $k$-plane the \textit{direction subspace} of $F$, and we may write it as $\dir(F)$.
\end{definition}

\begin{definition}\label{Def:major_angle}
    In the setting of Theorem~\ref{Thm:Principal_angle}, we call $\angle(V,W) \coloneqq \sin(\theta_k(V,W))$ the \textit{major angle} between $V$ and $W$. Given two $(\delta,k)$-flakes $F_1, F_2$ in $\R^n$, we define their \textit{major angle} $\angle(F_1,F_2)$ by that between $\dir(F_1)$ and $\dir(F_2)$.
\end{definition} 

\begin{rem}\label{Rem:major_metric}
    Since the operator norm $\norm{\cdot}$ has the triangle inequality, we know from Item~(i) in Theorem~\ref{Thm:Principal_angle} that $\angle(V,W) = \norm{\pi_V - \pi_W}$ is a metric on $\Gr(n,k)$, which is a $k(n-k)$-dimensional smooth manifold. The normalized Hausdorff measure induced by this metric is comparable to the normalized Haar measure/Riemannian measure on $\Gr(n,k)$.
\end{rem}

\begin{lem}\label{Lem:flake_contain}
    If a $(\delta,k)$-flake $F_\delta$ is contained in a $(\rho,k)$-flake $F_\rho$ for some $\rho \in [\delta,1]$, then $\angle(F_\delta,F_\rho) \lesssim \rho$. Also, if two $(\delta,k)$-flakes $F_1$ and $F_2$ are contained in the same $(\rho,k)$-flake, then $\angle(F_1,F_2) \lesssim \rho$. Conversely, if two $(\delta,k)$-flakes $F_1$ and $F_2$ have nontrivial intersection, then $F_1$ is contained in the $(\rho,k)$-flake $N_\rho(F_2)$ with $\rho \lesssim \max\{\delta, \angle(F_1,F_2)\}$, and vice versa. All constants in these bounds depend only on $n$.
\end{lem}

Lemma~\ref{Lem:flake_contain} is self-evident from Remark~\ref{Rem:major_equiv}, and so we omit the proof. It is important because it justifies many intuitions about the geometry of $(\delta,k)$-flakes. In this sense, the major angle can be used to quantify how ``close'' the direction subspaces of flakes are and thus be a natural analogue of the classical angle between two lines/$\delta$-tubes in $\R^n$. Also, Lemma~\ref{Lem:flake_contain} means that the major angle between two $(\delta,k)$-flakes is well defined up to error $\delta$: In Definition~\ref{Def:major_angle}, we can use any affine $k$-plane $t+V$ with $t + V\cap B_1 \subset F$ to define the direction subspace of $F$. Roughly speaking, two $(\delta,k)$-flakes $F_1$ and $F_2$ in $B_1 \subset \R^n$ are essentially distinct if either $\angle(F_1,F_2) \gtrsim \delta$, or $\angle(F_1,F_2) \lesssim \delta$ but the translation vector connecting their centers has component $\gtrsim \delta$ in the common normal direction space.

\subsection{Bilinear Kakeya-type estimate and the minor angle}

The major angle is not enough to accurately reflect all geometric information, such as the volume of the intersection between two $(\delta,k)$-flakes when $k \geq 2$. This is a vital difference from the case for $\delta$-tubes. We need all principal angles.
\begin{lem}\label{Lem:flake_volume}
    Suppose $V, W \in \Gr(n,k)$ have principal angles $0 \leq \theta_1 \leq \ldots \leq \theta_k \leq \pi/2$. Then we have        \begin{align}\label{Eq:volume_local}
        |N_\delta(V) \cap N_\delta(W) \cap B_1| \sim \delta^{n}/\prod_{i=1}^k \max\{\delta, \sin\theta_i\}.
    \end{align}
    In particular, for any two $(\delta,k)$-flakes $F_1$ and $F_2$ in $\R^n$, we have    \begin{align}\label{Eq:volume_local'}
        |F_1 \cap F_2| \lesssim \delta^{n}/\prod_{i=1}^k \max\{\delta, \sin(\theta_i(F_1,F_2))\}.
    \end{align}
\end{lem}
\begin{proof}
    Note that \eqref{Eq:volume_local'} is an immediate corollary of \eqref{Eq:volume_local} by translation. 
    
    To prove \eqref{Eq:volume_local}, let us first do some reduction by using Fubini-type arguments. Let $L:=V\cap W$, $j:=\dim L$, and $S:=(V+W)\cap L^\perp$. Then $\theta_1, \ldots, \theta_j =0$. By Fubini, up to constants depending only on $n$, integration in the $L$-directions contributes a factor $\sim 1$, while integration in the $(V+W)^\perp$-directions contributes a factor $\delta^{n-2k+j}$.    
    It therefore suffices to prove \eqref{Eq:volume_local} for the complementary subspaces $V\cap L^\perp$ and $W\cap L^\perp$ in $S$ (in place of $V$ and $W$ in $\R^n$), whose principal angles are $\theta_{j+1},\ldots,\theta_k$, because 
    \[
    \delta^{n-2k+j}
    \frac{\delta^{2(k-j)}}
    {\prod_{i=j+1}^{k}\max\{\delta,\sin\theta_i\}}
    =
    \frac{\delta^n}
    {\prod_{i=1}^{k}\max\{\delta,\sin\theta_i\}}.
    \] 
    In other words, it suffices to only prove \eqref{Eq:volume_local} when $\R^n = V \oplus W$, in which case we have $n=2k$ (thus $k=n-k$).
    
    Take a principal-vector decomposition of $V$ and $W$ as in Theorem~\ref{Thm:Principal_angle}.  For each $1 \leq i \leq k$, let $v_i^*$ and $w_i^*$ be a basis of $H_i \coloneqq \Span(v_i,w_i)$ such that $\norm{v_i^*} = \norm{w_i^*} = 1$, $v_i^* \perp v_i$, $w_i^* \perp w_i$, and $v_i^* \cdot w_i^* = \cos\theta_i$. It is easy to verify that $\VV^\perp \coloneqq \{v_{1}^*, \dots, v_k^*\}$ and $\WW^\perp \coloneqq \{w_{1}^*, \dots, w_k^*\}$ are principal vectors of $V^\perp$ and $W^\perp$. This is in fact how one proves Item~(iii) of Theorem~\ref{Thm:Principal_angle}.

    Define the map
    \begin{align*}
        \Phi: \R^n &\longrightarrow \R^n\\
        x &\longmapsto (x \cdot v_1^*, \dots, x \cdot v_k^*, x \cdot w_1^*, \dots, x \cdot w_k^*).
    \end{align*}
    Then by Item~(ii) of Theorem~\ref{Thm:Principal_angle} (applied to the pair $V^\perp$ and $W^\perp$), we have
    \begin{align*}
        |\det\Phi| \equiv 
        |(v_{1}^*\wedge \ldots \wedge v_k^*) \wedge (w_{1}^*\wedge \ldots \wedge w_k^*)| 
        =        \prod_{i=1}^k \sin\theta_i \neq 0.
    \end{align*}
    Note that $N_\delta(V) \cap N_\delta(W) = \Phi^{-1}(B_\delta^k \times B_\delta^k)$
    , and so     \begin{align}\label{Eq:volume_global}
        |N_\delta(V) \cap N_\delta(W)| = |B_\delta^k \times B_\delta^k|/|\det\Phi| \sim \delta^{n}/\prod_{i=1}^k \sin\theta_i,
    \end{align}    
    which can be regarded as a simple global version of \eqref{Eq:volume_local}. 

    Although $\theta_i$ are all nonzero, they could still be fairly close to $0$, which makes the right-hand side of \eqref{Eq:volume_global} blow up. In \eqref{Eq:volume_local}, since we have localized the intersection to $B_1$, one should not expect the volume to blow up, and so we introduce a ``$\max$'' relative to $\delta$. To estimate 
    $$|N_\delta(V) \cap N_\delta(W) \cap B_1^n| = |\Phi^{-1}(B_\delta^k \times B_\delta^k)\cap B_1^n|,$$ 
    note that the $2$-planes $H_i$ are orthogonal to each other, and so $\Phi$ is decomposed into the product of $2$-dimensional transformations:
    \begin{align*}
        \Phi_i: H_i &\longrightarrow \R^2\\
        x &\longmapsto (x\cdot v_i^*, x\cdot w_i^*).
    \end{align*}

    Since $B_\delta^k \times B_\delta^k \sim B_\delta^n \sim B_\delta^2 \times \cdots \times B_\delta^2$ ($k$-th product), we have
    \begin{align*}
        \prod_{i=1}^k 
        |\Phi_i^{-1}(B_{c\delta}^2)\cap B_c^2|
        \leq 
        |\Phi^{-1}(B_\delta^k \times B_\delta^k)\cap B_1^n|
        \leq
        \prod_{i=1}^k 
        |\Phi_i^{-1}(B_{C\delta}^2)\cap B_1^2|
    \end{align*}
    for some constants $1 \lesssim c < C \lesssim 1$. Therefore, as long as we could show
    \begin{align*}
        |\Phi_i^{-1}(B_{c\delta}^2)\cap B_c^2| \gtrsim \delta^2/\max\{\delta, \sin\theta_i\}\gtrsim
        |\Phi_i^{-1}(B_{C\delta}^2)\cap B_1^2|
    \end{align*}
    for each $i$, then we are done. However, these $2$-dimensional estimates follow easily from planar geometry.
\end{proof}

Lemma~\ref{Lem:flake_volume} has the following immediate corollary, which will be useful in Section~\ref{Sec:plaininess}:
\begin{cor}\label{Cor:bilinear_flakes}
    Let $\{F_{i,j_i}\}_{j_i=1}^{J_i}$ \rm{(}$i=1,2$\rm{)} be two families of $(\delta, k)$-flakes in $\R^n$. Then
    \begin{align*}
        \int_{\R^n} \sum_{j_1=1}^{J_1} \sum_{j_2=1}^{J_2} \chi_{F_{1,j_1}} \chi_{F_{2,j_2}} \prod_{i=1}^k \max\{\delta, \sin(\theta_i(\dir(F_{1,j_1}), \dir(F_{2,j_2})))\} 
        \lesssim
        \delta^n J_1J_2.
    \end{align*}
\end{cor}
\begin{proof}
    Rewrite \eqref{Eq:volume_local'} in Lemma~\ref{Lem:flake_volume} as
    \begin{align*}
        \int_{\R^n} \chi_{F_1} \chi_{F_2} \prod_{i=1}^k \max\{\delta, \sin(\theta_i(F_1,F_2))\} \lesssim \delta^{n}.
    \end{align*}
    Apply this to each pair of $F_{1,j_1}$ and $F_{2,j_2}$, and then sum the resulting estimates up.
\end{proof}

Intuitively, Lemma~\ref{Lem:flake_volume} means that for the volume of intersection to be small, we need all principal angles to be large. Roughly speaking, if $\theta_j \leq \delta$ for some $j$, then the intersection is large as it at least contains a $(\delta,j)$-flake, which means the two $(\delta,k)$-flakes should be compressed into a common $(\delta,2k-j)$-flake. This intuition is made rigorous by the following lemma:
\begin{lem}\label{Lem:form_planes}
    Suppose $V, W \in \Gr(n,k)$ have principal angles $0 \leq \theta_1 \leq \ldots \leq \theta_k \leq \pi/2$ with $\theta_j \leq \rho$ for some $2k-n \leq j \leq k$ and $\rho \in [\delta,1]$.\footnote{Note that $\dim (V \cap W) \geq 2k-n$, so $\theta_j = 0$ for all $j \leq 2k-n$, and ``$j \geq 2k-n$'' is a natural requirement.} Then there exists some constant $C \lesssim 1$ and $H \in \Gr(n,2k-j)$ such that $N_{C\rho}(H)$ contains $(N_\delta(V) \cup N_\delta(W)) \cap B_1$. This $H$ contains the internal bisector $($line$)$ of a pair of principal vectors defining $\theta_i$ for each $1 \leq i \leq j$. In general, if two $(\delta,k)$-flakes intersect, then there exists some constant $C \lesssim 1$ and a $(C\rho, 2k-j)$-flake $F_\rho$ containing both of them.
\end{lem}
\begin{proof}
    By translation invariance and Lemma~\ref{Lem:flake_contain}, the second claim about flakes is an immediate corollary of the first claim about subspaces. And so it suffices to prove the first claim. 

    Take a principal-vector decomposition of $V$ and $W$ as in Theorem~\ref{Thm:Principal_angle}. For each $1\leq i \leq j$, let $h_i \coloneqq (v_i + w_i)/\norm{v_i + w_i}$ (the unit vector along the internal bisector of $v_i$ and $w_i$). Define 
    $$H' \coloneqq \Span(\{h_i\}_{i=1}^j, \{v_i\}_{i=j+1}^k, \{w_i\}_{i=j+1}^k).$$ 
    Then clearly $\dim H' \leq 2k-j$. Since $2k-j \leq n$, we can find some $H \in \Gr(n,2k-j)$ with $H' \subset H$. We need to verify that $(N_\delta(V) \cup N_\delta(W)) \cap B_1 \subset N_{C\rho}(H)$. For this purpose, note that for each $1 \leq i \leq j$, since $v_i$ and $w_i$ make an angle $\theta_i \leq \rho$, we have $\norm{h_i - v_i} = \norm{h_i - w_i} \leq \rho$, which implies that $v_i$ and $w_i$ both lie in $N_\rho(H)$. On the other hand, it is clear that $\{v_i\}_{i=j+1}^k, \{w_i\}_{i=j+1}^k \subset H \subset N_\rho(H)$, and so $\VV, \WW \subset N_\rho(H)$, which means that $N_\delta(V) \cap B_1 \subset N_{C\rho}(H)$ and $N_\delta(W) \cap B_1 \subset N_{C\rho}(H)$ for some $C \lesssim 1$. And so the proof is completed.
\end{proof}

The final geometric notion we need is as follows:
\begin{definition}\label{Def:minor_angle}
    In the setting of Theorem~\ref{Thm:Principal_angle}, if $k \geq 2$, then we call $\angle'(V,W) \coloneqq \sin(\theta_{k-1}(V,W))$ the \textit{minor angle} between $V$ and $W$. Given two $(\delta,k)$-flakes $F_1, F_2$ in $\R^n$, we define their \textit{minor angle} $\angle'(F_1,F_2)$ by that between $\dir(F_1)$ and $\dir(F_2)$.
\end{definition} 
In this paper, we care mostly about $(n-2)$-planes ($k=n-2$, $n \geq 4$). Note that for a pair of $\delta$-flakes, $\theta_{n-3}$ and $\theta_{n-2}$ are the only possibly nonzero principal angles, and their relative rotational position is completely determined by their major angle and minor angle.

\begin{rem}\label{Rem:Rogers_angle}
    The major angle forms a metric by Remark~\ref{Rem:major_metric}, but the minor angle does not. In \cite{rogers2006}, Rogers also considered major and minor angles in the study of $2$-plane Kakeya (especially in $\R^4$), although his definition is superficially different from Definition~\ref{Def:major_angle}~and~\ref{Def:minor_angle}. More precisely, for a pair of $\delta$-plates $P_1,P_2$ in $\R^n$ with $\dir(P_1) = V$ and $\dir(P_2) = W$, he defines $\angle(P_1,P_2)$ as $\norm{\pi_V - \pi_W}$, and $\angle'(P_1,P_2)$ as $$\min_{\substack{\ell,\ell' \in \Gr(n,1)\\ \ell\subset V, \ell' \subset W}} \norm{\pi_\ell - \pi_{\ell'}}.$$
    However, by the structural properties of the principal-vector decomposition in Theorem~\ref{Thm:Principal_angle}, the two definitions are actually exactly the same. Equivalence of $\angle(P_1,P_2)$ comes directly from Item~(i) there. For $\angle'(P_1,P_2)$, simply note that
\begin{align*}
    \min_{\substack{\ell,\ell' \in \Gr(n,1)\\ \ell\subset V, \ell' \subset W}} \norm{\pi_\ell - \pi_{\ell'}} &= \min_{\substack{\ell,\ell' \in \Gr(n,1)\\ \ell\subset V, \ell' \subset W}} \sin(\theta_1(\ell,\ell')) = \sqrt{1- \max_{\substack{\ell,\ell' \in \Gr(n,1)\\ \ell\subset V, \ell' \subset W}} \cos(\theta_1(\ell,\ell'))^2}\\
    &= \sqrt{1- \cos(\theta_1(V,W))^2} = \sin(\theta_1(V,W)).
\end{align*}
\end{rem}

Remark~\ref{Rem:Rogers_angle} is the historical reason why we have insisted on putting a ``$\sin$'' in the definition of major and minor angles, although $\sin\theta \sim \theta$ whenever $\theta \in [0, \pi/2]$. In the $(n-2)$-plane case ($n\geq4$), we can summarize what we get so far as follows:
\begin{itemize}
    \item 
    By Item~(i) in Theorem~\ref{Thm:Principal_angle}, $\norm{\pi_V - \pi_W} = \angle(V,W)$.
    \item 
    By Item~(ii) in Theorem~\ref{Thm:Principal_angle}, $|(v_{n-3}\wedge v_{n-2}) \wedge (w_{n-3}\wedge w_{n-2})| = \angle'(V,W)\angle(V,W)$. Note that this is true even if $\dim V \cap W = n-3$, i.e., $\angle'(V,W) = 0$. 
    \item 
    In Corollary~\ref{Cor:bilinear_flakes}, the Jacobian factor becomes $$\delta^{n-4} \max\{\delta, \angle'(F_{1,j_1}, F_{2,j_2})\} \max\{\delta, \angle(F_{1,j_1},F_{2,j_2})\}.$$ This together with the prior bullet point recovers the special case of the multilinear $k_j$-plane estimate in \cite[Theorem~1.5]{zhang2018} when $d=4$ and $k_1=k_2=2$. However, our Corollary~\ref{Cor:bilinear_flakes} does not necessarily require the two families to be ``totally transverse'', which will be more flexible for our purposes.

    \item In Lemma~\ref{Lem:form_planes}, we can roughly say that if $\angle'(F_1,F_2)$ is small, then $F_1$ and $F_2$ are compressed into a suitable neighborhood of an affine $(n-1)$-plane.
\end{itemize}

\subsection{Additional geometric lemmas} In Section~\ref{Sec:plaininess}--\ref{Sec:medium_grains_end}, to analyze the sticky case, we need more results on $k$-plane geometry. Still, only some special cases will be used, but we keep the statements as general as possible for potential future applications.

\begin{lem}\label{Lem:total_transversality}
Let $V, W \in \Gr(n,k)$ with $\angle(V,W) \leq 1/10$, then $N_{\delta}(V) \cap N_{\delta}(W^\perp) \subset B(0, 2\delta)$ for all $\delta > 0$. In particular, $V \cap W^\perp = \{0\}$.
\end{lem}
\begin{proof}
Take $x \in N_{\delta}(V) \cap N_{\delta}(W^\perp)$, then $\norm{\pi_{V^\perp}(x)}, \norm{\pi_{W}(x)} \leq \delta$. By Remark~\ref{Rem:major_metric}, we have $\norm{\pi_V - \pi_W} = \angle(V,W) \leq 1/10$, and so $\norm{\pi_V(x) - \pi_W(x)} \leq (1/10)\norm{x}$. Thus by the Pythagorean theorem and triangle inequality, we get 
\begin{align*}
\norm{x}^2 &= \norm{\pi_{V^\perp}(x)}^2 + \norm{\pi_V(x)}^2\\
&\leq 
\norm{\pi_{V^\perp}(x)}^2 + (\norm{\pi_W(x)} + \norm{\pi_V(x) - \pi_W(x)})^2\\
&\leq \delta^2 + (\delta + (1/10)\norm{x})^2,
\end{align*}
rearranging which yields $\norm{x} \leq \frac{10(1+\sqrt{199})}{99}\delta < 2\delta$. So $N_{\delta}(V) \cap N_{\delta}(W^\perp) \subset B(0, 2\delta)$. 

Finally, if $V \cap W^\perp \neq \{0\}$, then $\dim V \cap W^\perp \geq 1$, and so $N_\delta(V \cap W^\perp) \subset N_\delta(V) \cap N_\delta(W^\perp)$ is unbounded, contradicting $N_{\delta}(V) \cap N_{\delta}(W^\perp) \subset B(0, 2\delta)$. Therefore, we must have $V \cap W^\perp = \{0\}$.
\end{proof}

\begin{lem}\label{Lem:faithful_proj}
Let $1 \leq k_1 \leq k_2 \leq n$, $W \in \Gr(n, k_1)$, and $V \in \Gr(n, k_2)$. Suppose that $W \cap B_1 \subset N_{\rho}(V)$ for some $\rho \leq 1/10$, then $\dim \pi_V(W) = k_1$, $\angle(W, \pi_V(W)) \leq \rho$, and $W \cap V^\perp = \{0\}$. Moreover, $V^\perp \cap B_1 \subset N_\rho(W^\perp)$.
\end{lem}
\begin{proof}
Clearly $\dim \pi_V(W) \leq \dim W = k_1$. Suppose to the contrary that $\dim \pi_V(W) \leq k_1-1$, then there exists a unit vector $w_0 \in W \cap \ker \pi_V$. Note that $w_0 \in W \cap B_1$, but $w_0 \not\in N_{1/10}(V)$, contradicting the assumption that $W \cap B_1 \subset N_\rho(V)$ ($\rho \leq 1/10$). Thus we must have $\dim \pi_V(W) = k_1$.

Now for each $w \in W \cap B_1$, we have $\pi_V(w) \in \pi_V(W)$ and $\norm{w - \pi_V(w)} \leq \rho$,  so $w \in N_{\rho}(\pi_V(W))$. Thus $W\cap B_1 \subset N_{\rho}(\pi_V(W))$, which by Remark~\ref{Rem:major_equiv} implies that $\angle(W, \pi_V(W)) \leq \rho$. 

Since $\angle(W, \pi_V(W)) \leq \rho$, by Lemma~\ref{Lem:total_transversality}, we get $W \cap \pi_V(W)^\perp = \{0\}$. As $\pi_V(W) \subset V$, we must have $W \cap V^\perp = \{0\}$.

Finally, by Item~(iii) of Theorem~\ref{Thm:Principal_angle}, we have $\angle(W^\perp, \pi_V(W)^\perp) = \angle(W, \pi_V(W)) \leq \rho$, which by Remark~\ref{Rem:major_equiv} implies $\pi_V(W)^\perp \cap B_1 \subset N_\rho(W^\perp)$. Note that $\pi_V(W) \subset V$, and so $V^\perp \cap B_1 \subset \pi_V(W)^\perp \cap B_1 \subset N_\rho(W^\perp)$, as desired.
\end{proof}

\begin{lem}\label{Lem:slice_transversality}
Let $1\leq k_1 \leq k_2 \leq n$, $W \in \Gr(n,k_1)$, and $V \in \Gr(n,k_2)$. Suppose that $W \cap B_1 \subset N_{1/10}(V)$, then $\dim (V \cap W^\perp) = k_2-k_1$, $V = \pi_{V}(W) \oplus (V \cap W^\perp)$, and $N_{\delta}(V) \cap N_{\delta}(W^\perp) \subset N_{2\delta}(V \cap W^\perp)$ for all $\delta>0$.
\end{lem}
\begin{proof}
By Lemma~\ref{Lem:faithful_proj}, we have $\dim \pi_{V}(W) = k_1$ and $\angle(W, \pi_{V}(W)) \leq 1/10$. Suppose $\pi_{V}(W) \cap W^\perp \neq \{0\}
$, then there exists a unit vector $v_0 \in \pi_{V}(W) \cap W^\perp$. Note that $v_0 \in \pi_{V}(W) \cap B_1 \subset N_{1/10}(W)$ by Remark~\ref{Rem:major_equiv}, contradicting the fact that $v_0 \in W^\perp \cap \mathbb{S}^{n-1}$. Therefore, we must have $\pi_{V}(W) \cap (V \cap W^\perp) = \{0\}$, which implies 
$$\dim(V \cap W^\perp) \leq \dim V - \dim \pi_V(W) = k_2 - k_1.$$ 
On the other hand, 
$$\dim(V \cap W^\perp) \geq \dim V + \dim W^\perp - \dim \R^n = k_2 + (n-k_1) - n = k_2 - k_1,$$
so we must have $\dim(V \cap W^\perp) = k_2-k_1$ and $V = \pi_{V}(W) \oplus (V \cap W^\perp)$.

Now we can project everything along $V \cap W^\perp$ to the $(n-k_2+k_1)$-plane $(V \cap W^\perp)^\perp$, and it remains to show the following claim: Let $V', W' \in \Gr(n-k_2+k_1, k_1)$ and suppose that $\angle(V', W') \leq 1/10$, then $N_{\delta}(V') \cap N_{\delta}(W'^\perp) \subset B(0, 2\delta)$. By Lemma~\ref{Lem:total_transversality}, the claim immediately holds true.
\end{proof}

\begin{lem}\label{Lem:geo_fact_1}
    For $i=1,2$, let $V_i, V_i' \in \Gr(n,k_i)$ with $\angle(V_i,V_i') \leq \rho$, where $\rho \leq 1/100$. Suppose that $V_1 \cap B_1 \subset N_{1/20}(V_2^\perp)$. Then $V_1+V_2$ and $V_1'+V_2'$ are both direct sums, and $(V_1' \oplus V_2')\cap B_1 \subset N_{3\rho}(V_1 \oplus V_2)$. In particular, $\angle(V_1 \oplus V_2, V_1' \oplus V_2') \leq 3\rho$.
\end{lem}
\begin{proof}
By Item~(iii) of Theorem~\ref{Thm:Principal_angle}, we know that $\angle(V_2^\perp, V_2'^\perp) = \angle(V_2, V_2') \leq \rho$. So by Remark~\ref{Rem:major_equiv}, we have 
$$V_1' \cap B_1 \subset N_{\rho}(V_1 \cap B_1) \subset N_{\rho + 1/20}(V_2^\perp \cap B_1) \subset N_{2\rho + 1/20}(V_2'^\perp \cap B_1) \subset N_{1/10}(V_2'^\perp).$$ 
Thus Lemma~\ref{Lem:faithful_proj} implies $V_1 \cap V_2 = \{0\}$ and $V_1' \cap V_2' = \{0\}$. Therefore, we can define $V = V_1 \oplus V_2$ and $V' = V_1' \oplus V_2'$, which satisfy $\dim V = \dim V' = k_1 + k_2$. By Remark~\ref{Rem:major_equiv}, it remains to show that $V' \cap B_1 \subset N_{3\rho}(V)$.

For each unit vector $v_0 \in V_2'^\perp \cap V'$, since $V' = V_1' \oplus V_2'$, there exist unique $v_1 \in V_1'$ and $v_2 \in V_2'$ such that $v_0 = v_1 + v_2$. We claim that $\norm{v_1} < 1.2$ and $\norm{v_2} < 0.12$. Note that $v_1 = v_0 + (-v_2)$ is subject to the orthogonal direct sum $V' = (V_2'^\perp \cap V') \oplus V_2'$, we must have $v_1 \neq 0$. Also, by taking inner products on both sides with $v_2$, we get $|\langle v_1, v_2 \rangle| = \norm{v_2}^2$. The key observation is that $V_1' \cap B_1 \subset N_{1/10}(V_2'^\perp)$ implies $\norm{\pi_{V_2'}(v_1/\norm{v_1})} \leq 1/10$, which in turn implies $|\langle v_1, v_2 \rangle| \leq \norm{v_2} \norm{\pi_{V_2'}(v_1)} \leq (1/10) \norm{v_1}\norm{v_2}$. If $v_2 \neq 0$, then by $|\langle v_1, v_2 \rangle| = \norm{v_2}^2$, we get $\norm{v_2} \leq (1/10)\norm{v_1}$. Thus $\norm{v_1} \leq \norm{v_0} + \norm{v_2} \leq 1 + (1/10)\norm{v_1}$, which implies $\norm{v_1} \leq 10/9 < 1.2$ and $\norm{v_2} \leq (1/10)\norm{v_1} < 0.12$. On the other hand if $v_2 = 0$, then $\norm{v_1} = \norm{v_0} = 1$. So the claim always holds true.

We know that $v_0 = \norm{v_1}(v_1/\norm{v_1}) + \norm{v_2}(v_2/\norm{v_2})$ with $v_1/\norm{v_1} \in V_1' \cap B_1 \subset N_\rho(V_1) \subset N_{\rho}(V)$ and $v_2/\norm{v_2} \in V_2' \cap B_1 \subset N_\rho(V_2) \subset N_{\rho}(V)$. So by linearity, we have 
\begin{align*}
\norm{\pi_{V^\perp}(v_0)} \leq \norm{v_1}\norm{\pi_{V^\perp}(v_1/\norm{v_1})} + \norm{v_2}\norm{\pi_{V^\perp}(v_2/\norm{v_2})}
< 1.2\rho + 0.12\rho < 1.4\rho.
\end{align*}

Now for each unit vector $v \in V'$, since $V' = (V_2'^\perp \cap V') \oplus V_2'$, there exist unique $v_1' \in V_2'^\perp \cap V'$ and $v_2' \in V_2'$ such that $v = v_1' + v_2'$. By orthogonality, we must have $\norm{v_1'}, \norm{v_2'} \leq \norm{v} = 1$. Writing $v = \norm{v_1'}(v_1'/\norm{v_1'}) + \norm{v_2'}(v_2'/\norm{v_2'})$ and taking the $v_0$ above as $v_1'/\norm{v_1'}$, we obtain 
\begin{align*}
\norm{\pi_{V^\perp}(v)} \leq \norm{v_1'}\norm{\pi_{V^\perp}(v_0)} + \norm{v_2'}\norm{\pi_{V^\perp}(v_2'/\norm{v_2'})}
< 1.4\rho + \rho < 3\rho.
\end{align*}
Therefore, we can conclude that $V' \cap B_1 \subset N_{3\rho}(V)$, as desired.
\end{proof}

\begin{lem}\label{Lem:geo_fact_2}
    For $i=1,2$, let $V_i\in \Gr(n,k_i)$ and $V_i' \in \Gr(n, k_i')$ with $k_i' \leq k_i$ and $V_i' \cap B_1 \subset N_\rho(V_i)$, where $\rho \leq 1/100$. Suppose that $V_1 \cap B_1 \subset N_{1/20}(V_2^\perp)$. Then $V_1+V_2$ and $V_1'+V_2'$ are direct sums, and $(V_1' \oplus V_2')\cap B_1 \subset N_{3\rho}(V_1 \oplus V_2)$. 
\end{lem}
\begin{proof}
    By applying Lemma~\ref{Lem:faithful_proj} to $V_1$ and $V_2^\perp$, we know that $V_1 \cap V_2 = \{0\}$, which means $V_1 + V_2 = V_1 \oplus V_2$. For $i = 1,2$, let $\bar{V}_i' \coloneqq \pi_{V_i}(V_i')$, then by applying Lemma~\ref{Lem:faithful_proj} again to $V_i'$ and $V_i$, we must have $\dim \bar{V}_i' = k_i'$ and $\angle(\bar{V}_i', V_i') \leq \rho$. Note that $\bar{V}_1' \cap B_1 \subset V_1 \cap B_1 \subset N_{1/20}(V_2^\perp) \subset N_{1/20}(\bar{V}_2'^\perp)$. So we can apply Lemma~\ref{Lem:geo_fact_1} to conclude that $\bar{V}_1' + \bar{V}_2'$ and $V_1'+ V_2'$ are direct sums, and $(V_1' \oplus V_2') \cap B_1 \subset N_{3\rho}(\bar{V}_1' \oplus \bar{V}_2') \subset N_{3\rho}(V_1 \oplus V_2)$. 
\end{proof}

\begin{lem}\label{Lem:geo_fact_3}
Let $V \in \Gr(n,k)$ and $W \in \Gr(n,n-k)$ satisfy $\angle(V,W^\perp) \leq \rho$ for some $\rho \leq 1/100$. Then for any $\ell \in \Gr(n,1)$ with $\ell \cap B_1 \subset N_\rho(V)$, we have $W + \ell$ intersects $V$ at a unique line $\ell'$. Moreover, $\angle(\ell, \ell') \leq 2\rho$. 
\end{lem}
\begin{proof}
By Lemma~\ref{Lem:total_transversality}, $\angle(V,W^\perp) \leq \rho \leq 1/10$ implies $V \cap W = \{0\}$, and so $\R^n = V \oplus W$. By Remark~\ref{Rem:major_equiv}, we know that $\ell \cap B_1 \subset N_\rho(V \cap B_1) \subset N_{2\rho}(W^\perp) \subset N_{1/10}(W^\perp)$, which by Lemma~\ref{Lem:faithful_proj} implies $\ell \cap W= \{0\}$. Thus $\dim (W+\ell) = n-k+1$, and so 
$$\dim (W+ \ell)\cap V = \dim (W+ \ell) + \dim V - \dim (W+\ell+V) = (n-k+1)+k-n = 1.$$ 
This justifies the first claim.

It remains to prove $\angle(\ell, \ell') \leq 2\rho$. Fix a unit vector $v_0 \in \ell \cap \mathbb{S}^{n-1}$. Since $\R^n = V \oplus W$, we have a unique decomposition $v_0 = v_1 + w_1$ with $v_1 \in V$ and $w_1 \in W$. If $w_1 = 0$, then $v_0 = v_1 \in V$, which means $\ell = \ell' \subset V$, and so $\angle(\ell,\ell')=0< \rho$. Thus we may assume $w_1 \neq 0$. Since $\R^n = V \oplus V^\perp$, we can further uniquely decompose $w_1$ as $w_1 = v_2 + v_3$ with $v_2 \in V$ and $v_3 \in V^\perp$. Note that $v_0 = (v_1 + v_2) + v_3$ with $v_1 + v_2 \in V$ and $v_3 \in V^\perp$, and so it must be the unique decomposition of $v_0$ under $\R^n = V \oplus V^\perp$. In particular, as $v_0 \in \ell \cap B_1 \subset N_\rho(V)$, we get $\norm{v_3} = \norm{\pi_{V^\perp}(v_0)} \leq \rho$ by definition of neighborhoods. 

By Item~(iii) of Theorem~\ref{Thm:Principal_angle}, we have $\angle(V^\perp, W) = \angle(V, W^\perp) \leq \rho$, which by Remark~\ref{Rem:major_equiv} implies $W \cap B_1 \subset N_\rho(V^\perp)$. In particular, as $w_1/\norm{w_1} \in W \cap B_1 \subset N_\rho(V^\perp)$, we get $\norm{v_2} = \norm{\pi_{V}(w_1)} = \norm{w_1} \cdot \norm{\pi_{V}(w_1/\norm{w_1})} \leq \rho\norm{w_1}$ by definition of
neighborhoods. Now by the triangle inequality, we get $\norm{w_1} \leq \norm{v_2} + \norm{v_3} \leq \rho\norm{w_1} + \rho$, and so $\norm{w_1} \leq \rho/(1-\rho) < 2\rho$. Also, $\norm{v_1} > \norm{v_0} - \norm{w_1} > 1- 2\rho > 0$, which means $v_1 \neq 0$. Note that $\ell = \langle v_0 \rangle$ and $\ell' = \langle v_1 \rangle$. Therefore, we obtain $\angle(\ell, \ell') = \dist(v_0, \ell') = \min_{v \in \ell'} \norm{v_0 - v} \leq \norm{v_0 - v_1} = \norm{w_1} < 2\rho$, as desired. 
\end{proof}

\section{Organizing convex sets}\label{Sec:organize}

In this section, we introduce and analyze two types of nonconcentration conditions for convex sets, called Katz--Tao conditions and
Frostman conditions. We record them from \cite[Section~3]{GWZ26streamlined} for the reader's convenience.

\begin{definition}\label{Def:Katz_Tao}
    A family of convex sets $\V$ is called \textit{$C$-convex-Katz--Tao} (or \textit{$C$-Katz--Tao} for short) if $\Delta_{max}(\V) \leq C$. Informally, we say $\V$ is \textit{Katz--Tao} if $\Delta_{max}(\V) \lessapprox 1$.
\end{definition}

\begin{definition}\label{Def:Frostman}
    If $\V$ is a set of convex sets in another convex set $K$, then we define
    \begin{align*}
        C_F(\V,K) \coloneqq \frac{\sup_{K'\subset K}\Delta(\V,K')}{\Delta(\V,K)},
    \end{align*}
where the supremum is taken over all convex subsets of $K$. {\color{black} Note that we may equivalently replace $\sup_{K'\subset K}\Delta(\V,K')$ with $\Delta_{max}(\V)$ in the numerator.}

We say that $\V$ is \textit{$C$-convex-Frostman}
(or \textit{$C$-Frostman} for short) \textit{in $K$} if $C_F(\V, K) \leq C$. Informally, we say that $\V$ is \textit{Frostman in $K$}
if $C_F(\V, K) \lessapprox 1$. {\color{black} For simplicity, if $\V$ has the form $\W[K]$ for some $\W$, then we abbreviate $C_F(\W[K],K)$ as  $C_F(\W[K])$.} In particular, if $K = B_1$ (thus $\V = \V[K]$ trivially holds), then we abbreviate $C_F(\V,K)$ as $C_F(\V)$, and we say that $\V$ is
\textit{$C$-Frostman}. 
\end{definition}

The two conditions enjoy nice heredity properties amenable to induction on scales. This is basically \cite[Remark~3.3]{GWZ26streamlined}, and we record them here for the reader’s convenience.

\begin{proposition}[heredity properties]\label{Prop:inherit}
    \phantom{x}    \begin{enumerate}[label={\rm (\Alph*)}]
    \item The property of being Frostman is \textit{inherited upwards}: Suppose $\V = \bigsqcup_W \V_W$ is $C$-Frostman inside $K$, where each $W\in\W$ is contained in $K$ and $V \subset W$ for each $V \in \V_W$. Then there is a subset $\W' \subset \W$ that is $\lessapprox C$-Frostman in $K$. If in addition the quantity $\sum_{\V_W}|V|/|W|$ is approximately the same for each $W \in \W$, then we may take $\W' = \W$, and in fact $\W$ is $O(C)$-Frostman in $K$.

    \item The property of being Katz--Tao is \textit{inherited downwards}: If $\V$ is $C$-Katz--Tao, then every subset $\V' \subset \V$ is also $C$-Katz--Tao.
    \end{enumerate}
\end{proposition}
One may consult \cite[Section~3]{GWZ26streamlined} for a detailed proof of Proposition~\ref{Prop:inherit}.

\begin{rem}\label{Rem:inherit}
    In Proposition~\ref{Prop:inherit}(A), if we assume that $\sum_{\V_W}|V|/|W|$ is constant up to a factor $\leq \Bar{C}$ instead, then $\W$ is $\Bar{C}C$-Frostman in $K$.
\end{rem}

If $\V$ is a family of convex sets with $\Delta_{max}(\V) \gg 1$, it turns out to be useful to consider the convex sets $W$ so that $\Delta(\V, W)$ is near maximal. These convex sets W can be used to organize
$\V$. A key tool in this philosophy is a \textit{factoring lemma} first proved by Wang and Zahl \cite[Proposition~4.6]{WZ25volume}. Here we record the version in \cite{GWZ26streamlined}, and one can find a detailed proof therein.

\begin{lem}[Maximal density factoring lemma, {\cite[Lemma~4.1]{GWZ26streamlined}}]\label{Lem:factoring}
    Let $\V$ be a finite family of convex subsets of $\R^n$. Then there is a set $\V' \subset \V$ with $\sum_{\V'} |V| \gtrapprox \sum_{\V} |V|$; a family $\W$ of convex sets in $\R^n$, each of approximately the same dimensions; and a partition
    \begin{align*}
        \V' = \bigsqcup_{W\in\W} \V_W,
    \end{align*}
    where each set $\V_W$ is contained in $W$. This partition has the following properties.
    \begin{enumerate}[label={\rm (\roman*)}]
        \item $\W$ is Katz--Tao, i.e. $\Delta_{max}(\W) \sim 1$.

        \item $\V_W$ is Frostman in $W$, i.e. $\Delta(\V_W, W) \sim \Delta_{max}(\V_W)$ for each $W \in \W$. Furthermore, $\Delta(\V_W, W) \sim \Delta_{max}(\V')$.
    \end{enumerate}
\end{lem}

\begin{definition}\label{Def:factor}
    Let $\V$ and $\W$ be families of convex subsets of $\R^n$, and suppose there is a partition $\V = \bigsqcup_\W \V_W$. We say that \textit{$\W$ factors $\V$} if the following is true:
    \begin{itemize}
        \item $\W$ is Katz--Tao, i.e. $\Delta_{max}(\W) \sim 1$.

        \item $\V_W$ is Frostman in $W$, i.e. $\Delta(\V_W, W) \sim \Delta_{max}(\V_W)$ for each $W \in \W$. Furthermore, $\Delta(\V_W, W) \sim \Delta_{max}(\V)$.

        \item Each $W \in \W$ has approximately the same dimensions, and each set in $\V_W$ is contained in $W$.
    \end{itemize}
\end{definition}

Thus Lemma~\ref{Lem:factoring} says that $\W$ factors $\V'$. We call the set $\W$ coming from Lemma~\ref{Lem:factoring} the
\textit{maximal density factoring} of $\V$. Thanks to the factoring lemma, we can often reduce the problem of understanding
$\V$ to the problem of understanding the Frostman sets $\V[W]$ for $W \in \W$ and the problem of
understanding the Katz--Tao set $\W$. This will facilitate certain unusual induction scheme and allows reduction to the sticky case. Wang and Zahl \cite{WZ25volume} successfully carried it out for $\delta$-tubes in $\R^3$, while we care about general $\delta$-flakes in $\R^n$. 

We end this section with the $\Gr(n,n-1)$ version of Theorem~\ref{Thm:Main}, which is essentially \cite[Lemma~5.7]{WZ25volume}, but is recorded here for the reader's convenience:
\begin{lem}\label{Lem:slab_0}
    Suppose $\mathbb{S}$ is a set of $\delta \times 1 \times \1$ slabs in $B_1 \subset \R^n$ (not necessarily essentially distinct) with a shading $Y$. Then \begin{align}\label{Eq:slab_0}
        |U(\mathbb{S},Y)| \gtrsim 
        \log(2/\delta)^{-1} \Delta_{max}(\mathbb{S})^{-1} \lambda(\mathbb{S},Y)^2 |\mathbb{S}||S|.
    \end{align}
    Equivalently, we have     \begin{align}\label{Eq:slab_0'}
        |U(\mathbb{S},Y)| \gtrsim 
        \log(2/\delta)^{-1} C_F(\mathbb{S})^{-1} \lambda(\mathbb{S},Y)^2.
    \end{align}
\end{lem}
\begin{proof}
    The first inequality \eqref{Eq:slab_0} follows immediately from \cite[Lemma~5.7]{WZ25volume} (which relies on a simple C\'ordoba-type $L^2$ argument) by noting that ``$\lambda$'' therein corresponds to $\lambda(\mathbb{S},Y)$ here, while ``$m=C_{KT-CW}(\mathcal{S})$'' therein corresponds to $\Delta_{max}(\mathbb{S})$ here. The second inequality \eqref{Eq:slab_0'} is clearly equivalent to \eqref{Eq:slab_0}, as $C_F(\mathbb{S}) \sim \Delta_{max}(\mathbb{S})/(|\mathbb{S}||S|)$ by Definition~\ref{Def:Frostman}.
\end{proof}

Lemma~\ref{Lem:slab_0} will be used in many places throughout the paper.

\section{Multiplicities, refinements, and factoring}\label{Sec:refinements}

In this section, we deal with technical issues related to shadings. One may consult \cite[Section~5]{GWZ26streamlined} for motivations, intuitions, and remarks. We will generalize things to $\delta$-flakes in $\R^n$.

\begin{proposition}[Factoring and multiplicity]\label{Prop:fact_mult}
    Let \( \V = \bigsqcup_{\W} \V_W \), where \( \V \) and \( \W \) are families of convex sets in \( \mathbb{R}^n \), and each family \( \V_W \) is \( C \)-Frostman in \( W \). Suppose that the sets in \( \V \) have dimensions {\color{black} $v_1 \times v_2 \times v_3\1$, and the sets in \( \W \) have dimensions $w_1 \times w_2 \times \cdots \times w_n$}. Let \( Y \) be a shading on \( \V \).

    Then there is a refinement \( (\V', Y') \) and a subset \( \W' \subset \W \) with a shading \( Y_{\W'} \) so that the following holds. Define \( \V'_W = \V_W \cap \V' \) so that \( \V' = \bigsqcup_{W \in \W'} \V'_W \).

\begin{enumerate}[label={\rm (\arabic*)}]
    \item \( (\V', Y') \) is a \( \gtrapprox 1 \) refinement of \( (\V, Y) \).
    \item \( \lambda(\W', Y_{\W'}) \gtrapprox C^{-1} \lambda(\V, Y)^2 \).
    \item \( (\W', Y_{\W'}) \) has constant multiplicity.
    \item For each \( W \in \W' \), \( (\V'_W, Y') \) has constant multiplicity. This multiplicity is also the same for all \( W \in \W' \).
    \item For each \( W \in \W' \), we have
    \[
    \mu(\V, Y) \lessapprox \mu(\W', Y_{\W'}) \, \mu(\V_W', Y').
    \]
    \item For every \( x \in U(\V', Y') \), we have 
    \[
    \V'_{Y'}(x) \subset \bigcup_{W \in \W'} (\V'_W)_{Y'}(x).
    \]
    \item \( |U(\V', Y') \cap B(x, {\color{black} w_1})| \) is roughly the same for every \( x \in U(\W', Y_{\W'}) \), up to a factor \( \lesssim 1 \).
\end{enumerate}
\end{proposition}

\begin{rem}\label{Rem:weaker_fact_mult}
    The hypothesis that $\V_W$ is $C$-Frostman in $W$ can be replaced by the following weaker requirement. Let $W \in \W$ have dimensions $w_1 \times w_2 \times \cdots \times w_n$. Then instead of requiring that $\V_W$ is $C$-Frostman in $W$, we require that the $w_1$-thickening of each $V \in \V_W$ is $C$-Frostman in W.
\end{rem}

The shading $Y_{\W'}$ on $\W'$ is compatible with the shading $Y'$ on $\V$ according to \cite[Definition~5.7]{GWZ26streamlined}, and is called the \textit{induced shading}. The proof of Proposition~\ref{Prop:fact_mult} is mostly the same as that of \cite[Proposition~5.1]{GWZ26streamlined}, and we only need to establish the following analogue of \cite[Lemma~5.9]{GWZ26streamlined}:
\begin{lem}\label{Lem:Frostman_density}
    Let $\V = \bigsqcup_\W \V_W$ be a collection of convex sets of dimensions $v_1 \times v_2 \times v_3\1$ in $\R^n$, where the sets of $\V_W$ are contained in $W$. Let $Y$ be a shading on $\V$, and let $Y_\W$ be the induced shading on $\W$. Suppose that for each $W \in \W$, the set $\V_W$ is $C$-Frostman in $W$. Then for each $W \in \W$, we have 
    \begin{align*}
        |Y_\W(W)| \gtrapprox 
        C^{-1} \lambda(\V_W, Y)^2 |W|.
    \end{align*}
\end{lem}
\begin{proof}
    The argument is similar to Córdoba’s $L^2$ proof for slabs (Lemma~\ref{Lem:slab_0}). After dyadic pigeonholing, we may select a set $\V_W' \subset \V_W$ and a number $\lambda \gtrapprox \lambda(\V_W,Y)$ so that $|Y(V)| \sim \lambda|V|$ for each $V \in \V_W'$, and $\sum_{\V_W'}|Y(V)| \gtrapprox \sum_{\V_W}|Y(V)|$. We thus have
\[
\lambda |\V_W'| \gtrsim |V|^{-1} \sum_{\V_W'} |Y(V)| \gtrapprox |V|^{-1} \sum_{\V_W} |Y(V)| \gtrsim \lambda (\V_W, Y) |\V_W|.
\]

In particular, since \( \V_W \) is \( C \)-Frostman relative to \( W \), we have that \( \V_W' \) is \( \lessapprox \frac{\lambda}{\lambda (\V_W, Y)} C \)-Frostman relative to \( W \).

Let \( w_1 \times w_2 \times \cdots \times w_n \) be the dimensions of the set \( W \) that we fixed above. We have \( v_i \leq w_i \) for \( i = 1, 2, \dots, n \). Define \( a_i = \max\{v_i, w_1\} \), and define 
\[
\K = \{ N_{w_1}(V) \cap W : V \in \V_W' \};
\]
each \( K \in \K \) has dimensions comparable to \( a_1 \times a_2 \times a_3\1 \). (We note that the sets \( K \in \K \) need not be essentially distinct. We have \( |\K| = |\V_W'|. \)) Define 
\[
Y_{\K}(K) = K \cap \bigcup_{\V_W'} N_{w_1}(Y(V)).
\]

Each \( K \in \K \) contains a set \( V \in \V_W'[K] \) with \( |Y(V)| \gtrsim \lambda |V| \), and \( v_1 \leq w_1 \), so we have that the \( w_1 \)-covering number of \( Y(V) \) satisfies

\[
\mathcal{E}_{w_1}(Y(V)) \gtrsim 
\begin{cases} 
\lambda w_1^{-(n-1)} v_2 v_3^{n-2}, & \text{if } v_2 \geq w_1, \\[4pt]
\lambda w_1^{-(n-2)} v_3^{n-2}, & \text{if } v_3 \geq w_1 > v_2, \\[4pt]
1, & \text{if } w_1 > v_3.
\end{cases}
\]
In each of the above scenarios, we have 
\[
|Y_{\K}(K)| \sim w_1^n \mathcal{E}_{w_1}(Y(V)) \gtrsim \lambda |K|
\]
for each \( K \in \K \).

By Proposition~\ref{Prop:inherit}(A), there is a subset \( \K' \subset \K \) that is \( \lessapprox \frac{\lambda}{\lambda(\V_{W},Y)} C \)-Frostman relative to \( W \). Define \( \theta_0 = a_2/a_3 \). Observe that if \( K, K' \in \K' \) with \( |K \cap K'| \geq \theta |K| \) for some \( \theta \in [\theta_0, 1] \) ($\theta_0$ is the smallest possible nontrivial value), then \( K' \) must be contained in the rectangular prism of dimensions comparable to \( a_1 \times \theta^{-1} a_2 \times a_3\1 \) given by \( N_{\theta^{-1} a_2}(K) \cap W \). This prism has volume roughly \( \theta^{-1} |K| \), and thus contains \( \lessapprox \frac{\lambda}{\lambda(\V_{W},Y)} C \frac{\theta^{-1} |K|}{|W|} |\K'| \) sets from \( \K' \). We now compute

\[
\begin{aligned}
\Big\| \sum_{K \in \K'} \chi_K \Big\|_2^2 
&= \sum_{K \in \K'} 
   \sum_{\substack{\theta \in [\theta_0, 1]\\
   \theta \text{ dyadic}}} 
   \sum_{\substack{K' \in \K' \\ 
    |K \cap K'| \sim \theta |K|}} 
   |K \cap K'| \\
&\lessapprox \sum_{K \in \K'}    \sum_{\substack{\theta \in [\theta_0, 1]\\ \theta \text{ dyadic}}} 
   \left( \frac{\lambda}{\lambda(\V_{W},Y)} C \frac{\theta^{-1} |K|}{|W|} |\K'| \right) (\theta |K|) \\
&\lessapprox \frac{\lambda}{\lambda(\V_{W},Y)} C \frac{|K|^2|\K'|^2}{|W|}.
\end{aligned}
\]

Finally, by Cauchy–Schwarz and the fact that \( \lambda \gtrapprox \lambda(\V_{W},Y) \), we have

\[
\begin{aligned}
|Y_\W(W)| &\gtrsim |U(\K', Y_\K)| \\
&\gtrsim \frac{ \bigl( \lambda |K| |\K'| \bigr)^2 }{ \bigl\| \sum_{K \in \K'} \chi_{Y_\K(K)} \bigr\|_2^2 } \\
&\gtrapprox C^{-1} \lambda \, \lambda(\V_{W},Y) |W| \\
&\gtrapprox C^{-1} \lambda(\V_{W},Y)^2 |W|,
\end{aligned}
\]
as desired.
\end{proof}

We also state an easier variant of Proposition~\ref{Prop:fact_mult} when $\V$ and $\W$ are both sets of flakes. We will only use the case when $k=n-2$ in this paper.

\begin{lem}\label{Lem:fact_mult}
Suppose \( \F \) is a set of \( (\delta,k) \)-flakes in $B_1 \subset \R^n$ and \( Y \) is a shading on \( \F \).

Let \( \rho \in [\delta, 1] \) and let \( \F = \bigsqcup_{F_\rho \in \F_\rho} \F[F_\rho] \). Let \( Y \) be a shading on \( \F \) with \( |Y(F)|/|F| \gtrapprox \lambda \) for each \( F \in \F \).

Then there is a subset \( \F'_\rho \subset \F_\rho \), a shading \( Y_{\F_\rho'}\) on \( \F_\rho' \), and a shading \( Y' \) on \( \F' = \bigsqcup_{F_\rho\in\F_\rho'} \F[F_\rho] \) so that $|\F[F_\rho]|$ are comparable for all $F_\rho \in \F_\rho'$, and 
\begin{enumerate}[label={\rm (\arabic*)}]
    \item \( (\F', Y') \) is a \( \gtrapprox 1 \) refinement of \( (\F, Y) \).
    \item \( \lambda(\F_\rho', Y_{\F_\rho'}) \gtrapprox \lambda \).
    \item \( (\F_\rho', Y_{\F_\rho'}) \) has constant multiplicity.
    \item For each \( F_\rho \in \F_\rho' \), \( (\F'[F_\rho], Y') \) has constant multiplicity. This multiplicity is also the same for all \( F_\rho \in \F'_\rho \).
    \item For each \( F_\rho \in \F'_\rho \) we have
    \[
    \mu(\F, Y) \lessapprox \mu(\F_\rho', Y_{\F_\rho'}) \, \mu(\F'[F_\rho], Y').
    \]
    \item For every \( x \in U(\F', Y') \), we have 
    \[
    \F_{Y'}(x) \subset \bigcup_{\substack{F_\rho \in \F_\rho' \\ x \in Y_{\F_\rho'}(F_\rho)}} (\F'[F_\rho])_{Y'}(x).
    \]
    \item For every \( x \in U(\F_\rho', Y_{\F_\rho'}) \) we have 
    \[
    |U(\F', Y')| {\color{black} \sim} |U(\F_\rho', Y'_{\F_\rho})| \,
    \frac{|U(\F', Y') \cap B(x, \rho)|}{|B(x, \rho)|}.
    \]
\end{enumerate}
\end{lem}

The proof of Lemma~\ref{Lem:fact_mult} is the same as the proof of Proposition~\ref{Prop:fact_mult}, except that in place of Lemma~\ref{Lem:Frostman_density} we use the following fact: If \( F_\rho \supset F \), then \( |F_\rho \cap N_\rho(Y(F))| \gtrsim \frac{|Y(F)|}{|F|} |F_\rho| \). Thus we do not need to assume that $k=n-2$ and \( \F_{F_\rho} \) is $C$-Frostman in \( F_\rho \), and we obtain the estimate \( \lambda(\F_\rho', Y_{\F_\rho'}) \gtrapprox \lambda \), which is stronger than its analogue (Item~(2)) in Proposition~\ref{Prop:fact_mult}.

\begin{rem}\label{Rem:factor_only_shading}
In Lemma~\ref{Lem:fact_mult}, if $|\F[F_\rho]|$ is constant up to a factor $\approx 1$ for all $F_\rho \in \F_\rho$, then we can take $\F' = \F$ and $\F_\rho' = \F_\rho$. In other words, we can just refine the shadings, while keeping the original set of flakes unchanged. This can be easily seen from the proof of \cite[Proposition~5.1]{GWZ26streamlined}.
\end{rem}

\section{Sticky Kakeya for flakes}\label{Sec:stronger_sticky}

{\color{black} Now we are in the right place to formally start the proof of Theorem~\ref{Thm:Main}, which is set up by induction on the dimension $n$. By Theorem~\ref{Thm:R3_Kakeya}, the $n=3$ case is already known, which serves as the base case. From now on, we always assume $n\geq 4$, and it suffices to prove Theorem~\ref{Thm:Main} under the induction hypothesis that it holds when $n$ is replaced with $n-1$. Throughout the proof, this hypothesis is only used at one place, i.e., Step (ii) in Section~\ref{Sec:slicing_set_up}.}

Let us begin with the sticky case. In this section, we will introduce versions of the $(n-2)$-plane sticky Kakeya theorem, and indicate the relationship between them.

Let \( \F \) be a uniform set of \( \delta \)-flakes. Let \( \rho \in [\delta, 1] \) and let \( \F_\rho \) be the thickening of \( \F \) at scale \( \rho \). We begin with the following two observations:
\begin{enumerate}[label={\rm (\Alph*)}]
\item If \( \F \) is Frostman, then by Proposition~\ref{Prop:inherit}(A), \( \F_\rho \) is Frostman. In general, however, the flakes in \( \F[F_\rho] \) will not be Frostman inside \( F_\rho \).

\item If \( \F \) is Katz--Tao, then by Proposition~\ref{Prop:inherit}(B), the flakes in \( \F[F_\rho] \) will be Katz--Tao. In general, however, the flakes in \( \F_\rho \) will not be Katz--Tao.
\end{enumerate}

It will be helpful to consider the special cases where the sets \( \F[F_\rho] \) are Frostman in \( F_\rho \), or the set \( \F_\rho \) is Katz--Tao.

\begin{definition}\label{Def:sticky}
Let \( \F \) be a set of \( \delta \)-flakes in \( \mathbb{R}^n \).
\begin{enumerate}[label={\rm (\Alph*)}]
\item
We say that \( \F \) is \emph{$C$ Frostman at every scale} if for every scale \( \delta \leq \rho \leq 1 \), we have
\[
C_F(\F[F_\rho]) \leq C \quad \text{for every } F_\rho \in \F_\rho.
\]

\item 
We say that \( \F \) is \emph{$C$ Katz--Tao at every scale} if for every scale \( \delta \leq \rho \leq 1 \), we have
\[
\Delta_{max}(\F_\rho) \leq C.
\]
\end{enumerate}
\end{definition}

We will usually consider sets of flakes that are uniform, in the sense of Definition~\ref{Def:uniform_flakes}, but this is not required in the above definition.

We say that \( \F \) is \emph{sticky} if it is both Frostman at every scale and Katz--Tao at every scale. In that case, we expect an $(n-2)$-plane sticky Kakeya theorem, which says that \( |U(\F)| \gtrapprox 1 \) and \( \mu(\F) \lessapprox 1 \). A mild generalization says that if \( \F \) is Frostman at every scale, then \( |U(\F)| \gtrapprox 1 \); while if \( \F \) is Katz--Tao at every scale, then \( \mu(\F) \lessapprox 1 \). The precise statement is as follows:
\begin{theorem}\label{Thm:2-plane_sticky}
For all \( \eps > 0 \), there exists \( \eta, \delta_0 > 0 \) so that the following holds for all \( \delta \in (0, \delta_0] \). Let \( (\F, Y) \) be a set of {\color{black} essentially distinct} \( \delta \)-flakes in \( B_1 \subset \mathbb{R}^n \), with \( \lambda(\F, Y) \geq \delta^\eta \).
\begin{enumerate}[label={\rm (\Alph*)}]
\item
If \( \F \) is $\delta^{-\eta}$ Frostman at every scale, then  
\[
|U(\F, Y)| \geq \delta^\eps.
\]

\item 
If \( \F \) is $\delta^{-\eta}$ Katz--Tao at every scale, then  
\[
\mu(\F, Y) \leq \delta^{-\eps}.
\]
\end{enumerate}
\end{theorem}

{\color{black} Similarly to Theorem~\ref{Thm:Main}, by a straightforward tensor-product argument, the special case when $n=3$ solved by Wang and Zahl \cite{WZ26sticky, WZ25assouad, WZ25volume} can be deduced from Theorem \ref{Thm:2-plane_sticky}.} However, the geometry of $\delta$-flakes is more complicated than $\delta$-tubes due to
multiple principal angles (see Section~\ref{Sec:geometry}), and so Theorem~\ref{Thm:2-plane_sticky} does not come immediately from Wang and Zahl's work.

We will first prove Theorem~\ref{Thm:2-plane_sticky} in the special case when $|\F| \approx \delta^{-2}$, under which Parts~(A)~and~(B) are essentially equivalent. This motivates the following stronger version of stickiness:
\begin{definition}\label{Def:strong_sticky}
    We say that a set of $\delta$-flakes ${\mathbb F}$ in $B_1 \subset \R^n$ is \textit{$\delta^{-\varepsilon}$-sticky} if 
\begin{enumerate}[label={\rm(\alph*)}]
    \item $\F$ is essentially distinct.
    \item $\Delta_{max}(\F) \leq \delta^{-\eps}$.
    \item For each $\delta \leq \rho \leq 1$ and $F_\rho \in \F_\rho$, we have
    $$    \delta^{\varepsilon}\left(\frac{\rho}{\delta}\right)^2 \le |{\mathbb F}[F_{\rho}]| \le \delta^{-\varepsilon}\left(\frac{\rho}{\delta}\right)^2.
    $$ 
\end{enumerate}
\end{definition}

\begin{rem}\label{Rem:strong_sticky}
    Item~(c) implies $\delta^{-2+\eps} \leq |\F| \leq \delta^{-2-\eps}$, and so for each $\delta \leq \rho \leq 1$, we have $\delta^{2\eps} \rho^{-2} \leq |\F_\rho| \leq \delta^{-2\eps} \rho^{-2}$. Item~(b) implies $\Delta_{max}(\F[F_\rho]) \leq \delta^{-\eps}$ by Proposition~\ref{Prop:inherit}(B). However, Item (c) also implies $\delta^\eps \leq \Delta(\F[F_\rho], F_\rho) \leq \delta^{-\eps}$, and so $C_F(\F[F_\rho]) \leq \delta^{-2\eps}$, i.e., $\F$ is $\delta^{-2\eps}$ Frostman at every scale.
\end{rem}

We know from Remark~\ref{Rem:strong_sticky} that Theorem~\ref{Thm:2-plane_sticky}, once established, should immediately imply:
\begin{theorem}\label{Thm:Main_Sticky}
For all \( \eps > 0 \), there exists \( \eta, \delta_0 > 0 \) so that the following holds for all \( \delta \in (0, \delta_0] \). Let \( (\F, Y) \) be a set of \( \delta \)-flakes in \( B_1 \subset \mathbb{R}^n \) such that $\F$ is $\delta^{-\eta}$-sticky and \( \lambda(\F, Y) \geq \delta^\eta \). Then  
\[
|U(\F, Y)| \geq \delta^\eps.
\]
\end{theorem}
\begin{rem}
Standard pigeonholing arguments allows us to only consider Theorem~\ref{Thm:Main_Sticky} in the case when $(\F,Y)$ is uniform.
\end{rem}

In Section~\ref{Sec:finish_sticky}, we will show that the converse is also true, i.e., Theorem~\ref{Thm:Main_Sticky} actually implies Theorem~\ref{Thm:2-plane_sticky}. More precisely, we will first show that Part~(A) of Theorem~\ref{Thm:2-plane_sticky} implies Part~(B), and then show that Theorem~\ref{Thm:Main_Sticky} implies Part~(A). These arguments are kind of standard and only relies on certain probabilistic methods, so we choose to postpone them until later on. Therefore, for the following few sections (Section~\ref{Sec:self_similar}--\ref{Sec:medium_grains_end}), we will focus on how to prove Theorem~\ref{Thm:Main_Sticky}, presents the principal challenge and embodies our key technical novelty.

\section{Extremal families of flakes and multiscale structure}\label{Sec:self_similar}

In this section, we will set up the induction framework for Theorem~\ref{Thm:Main_Sticky} and provide some useful self-similarity properties. We will always work in $\R^n$.

\begin{definition}\label{Def:sigma}
    For $\eps,\delta \in (0,1)$, define 
    \begin{equation*}
        M(\eps,\delta) = \inf_{\F, Y} |U(\F,Y)|.
    \end{equation*}
    Here, the infimum is taken over all sets of $\delta$-flakes $(\F, Y)$ with the following properties:
    \begin{enumerate}[label=(\alph*)]
        \item\label{item_a} $\F$ is $\delta^{-\eps}$-sticky.
        \item\label{item_b} $\lambda(\F,Y) \geq \delta^\eps$.
    \end{enumerate}
    Next, let
    \begin{equation*}
        N(\eps) = \limsup_{\delta \rightarrow 0^+} \frac{\log M(\eps,\delta)}{\log \delta},
    \end{equation*}
    which is defined for all $\eps \in (0,1)$. Observe that for each $\delta\in (0,1)$, if $\eps \leq \eps'$, then $M(\eps,\delta) \geq M(\eps',\delta)$, and hence $N(\eps) \leq N(\eps')$. In particular, since $0 \leq N(\eps) \leq n$, we have that 
    \begin{equation*}
        \lim_{\eps \rightarrow 0^+} N(\eps) = \inf_{\eps\in (0,1)} N(\eps).
    \end{equation*}
    Denote this common value by $\sigma$.
\end{definition}

We may regard $\sigma \geq 0$ as fixed from now on. Observe that $\sigma=0$ immediately implies Theorem~\ref{Thm:Main_Sticky}. Also, since we only care about small $\eps$ in practice, we may always assume that $0 < \eps < 1/100$. To prove $\sigma=0$, we will {\color{black} start by assuming $\sigma>0$} and study families of flakes that nearly extremize the quantity
$M$ from Definition~\ref{Def:sigma} to get a contradiction.

\begin{definition}\label{Def:eps_extremize}
    Let $\eps, \delta > 0$. We say a set of $\delta$-flakes $(\F,Y)$ is \textit{$\eps$-extremal} if it satisfies Items~\ref{item_a} and \ref{item_b} from Definition~\ref{Def:sigma}, and furthermore $|U(\F,Y)| \leq \delta^{\sigma-\eps}$.
\end{definition}
\begin{rem}\label{Rem:extreme_mult}
    Note that any $\eps$-extremal set of $\delta$-flakes has $\mu(\F,Y) \geq \delta^{-\sigma+3\eps}$. On the other hand, if we only assume $(\F, Y)$ to satisfy Items~\ref{item_a}~and~\ref{item_b}, then for any given $\eps_0>0$, by choosing $\eps$ small enough depending on $\eps_0$, we have $|U(\F,Y)| \geq \delta^{\sigma + \eps_0}$ (thus $\mu(\F,Y) \leq \delta^{-\sigma-2\eps_0}$) for small $\delta$. {\color{black} Finally, if $(\F',Y')$ is a refinement of it with $\lambda(\F',Y') \geq \delta^{\eps'}$ for some $\eps' > \eps$, then $(\F,Y')$ is $\eps'$-extremal by letting $Y'(F) = \varnothing$ for $F \in \F \setminus \F'$. And abusing notation, we may sometimes say that $(\F',Y')$ is $\eps'$-extremal. This does not affect the arguments as $U(\F,Y') = U(\F',Y')$.}
\end{rem}
From the definition of $\sigma$, we immediately conclude that extremal collections of
flakes always exist for small $\delta$. More precisely, we have the following.

\begin{lem}\label{Lem:extremizer_exist}
    Let $\eps, \delta_0 > 0$. Then there exists an $\eps$-extremal set of $\delta$-flakes $(\F, Y)$ for some $\delta \in (0, \delta_0)$.
\end{lem}

We can now state the main result of this section, which can be regarded as an upgraded version of Lemma~\ref{Lem:fact_mult}. It says that extremal collections of $\delta$-flakes look self-similar at all scales.
\begin{proposition}\label{Prop:self-similar}
    For all $\eps>0$, there exist $\eta>0$ and $\delta_0>0$ so that the following holds for all $\delta \in (0,\delta_0]$. Let $(\F, Y)$ be an $\eta$-extremal set of $\delta$-flakes, and let $\rho \in [\delta^{1-\eps}, \delta^\eps]$. Then there is a $\approx 1$ refinement $(\F', Y')$ of $(\F, Y)$, a subset $\F_\rho' \subset \F_\rho$, and a shading $Y_{\F_\rho'}$ on $\F_\rho'$ with $\lambda(\F_\rho', Y_{\F_\rho'}) \gtrapprox \lambda(\F, Y)$, such that the following properties hold: 
    \begin{enumerate}[label={\rm (\roman*)}]
        \item $(\F_\rho', Y_{\F_\rho'})$ is $\eps$-extremal.

        \item 
        For each $x \in \R^n$, $\mu(\F_\rho', Y_{\F_\rho'})(x) \leq \rho^{-\sigma-\eps}$. In fact, on $U(\F_\rho', Y_{\F_\rho'})$, $\mu(\F_\rho', Y_{\F_\rho'})(x)$ is constant and satisfies $$\rho^{-\sigma + \eps} \leq \mu(\F_\rho',Y_{\F_\rho'})(x) \leq \rho^{-\sigma-\eps}.$$
        \item 
        For each $F_\rho \in \F_\rho'$ and $x \in \R^n$, $\mu(\F'[F_\rho], Y')(x) \leq (\delta/\rho)^{-\sigma-\eps}$. In fact, on $U(\F'[F_\rho], Y')$, $\mu(\F'[F_\rho], Y')(x)$ is constant and satisfies $$(\delta/\rho)^{-\sigma+\eps} \leq \mu(\F'[F_\rho], Y')(x) \leq (\delta/\rho)^{-\sigma-\eps}.$$ The constant is also the same for all $F_\rho \in \F_\rho'$.
        \item $\F' = \bigsqcup_{\F_\rho'} \F'[F_\rho]$, and for each $F_\rho \in \F_\rho'$, $U(\F'[F_\rho], Y') \subset Y_{\F_\rho'}(F_\rho)$.
        \item 
        $|U(\F',Y') \cap B(x, \rho)|$ is roughly the same for every $x \in U(\F_\rho', Y_{\F_\rho'})$. In other words,
        \[
        |U(\F', Y')| \sim |U(\F_\rho', Y'_{\F_\rho})| \,
        \frac{|U(\F', Y') \cap B(x, \rho)|}{|B(x, \rho)|}.
        \]
    \end{enumerate}
\end{proposition}

{\color{black} By the convention in \cite{WZ26sticky}, we say that \textit{$(\F_\rho', Y_{\F_\rho'})$ covers $(\F', Y_{\F'})$} if Item~(iv) above holds, and call such a cover \textit{balanced} if Item~(v) additionally holds.} 
The proof of Proposition~\ref{Prop:self-similar} is very similar to that for \cite[Proposition~3.2]{WZ26sticky} (see also \cite[Proposition~6.2]{WZ25assouad}). Here we provide an exposition using (the proof of) \cite[Proposition~5.1]{GWZ26streamlined}. Before that, we need some technical lemmas:
\begin{lem}\label{Lem:retain_sticky}
    For all $\eta_1 > 0$, there exist $\delta_0>0$ so that the following holds for $\delta \in (0,\delta_0]$. Let $\F$ be a $\delta^{-\eta_1}$-sticky set of $\delta$-flakes, and let $\F'$ be an arbitrary subset of $\F$ with $|\F'| \gtrapprox \delta^{\eta_2} |\F|$, then $\F'$ contains a uniform subset $\F''$ with $|\F''| \gtrapprox |\F'|$ that is $\delta^{-(4\eta_1+\eta_2)}$-sticky. 
\end{lem}
\begin{proof}
    Let $M$ and $\rho_k$ be as in Definition~\ref{Def:uniform_flakes}. By pigeonholing for $\F'$ at scale $\rho_1$, there is $\widetilde{\F}_{\rho_1} \subset \F_{\rho_1}$ and a constant $C$ such that 
    \begin{itemize}
        \item $\widetilde{\F}_{\rho_1}$ is \textit{spaced} in the sense that $\{(\log(1/\delta))F_{\rho_1}\}_{F_{\rho_1} \in \widetilde{\F}_{\rho_1}}$ is essentially distinct. In particular, for any $\tau \leq \rho_1$, the sets $\{F_\tau \in \F_{\tau}: \F[F_\tau] \cap \F[F_{\rho_1}] \neq \varnothing\}_{F_{\rho_1} \in \widetilde{\F}_{\rho_1}}$ are pairwise disjoint.

        \item Let $\F'_{(\rho_1)} \coloneqq \bigsqcup_{F_{\rho_1} \in \widetilde{\F}_{\rho_1}} \F'[F_{\rho_1}]$. Then $|\F'_{(\rho_1)}| \gtrsim (\log(1/\delta))^{-C} |\F'|$.
     \end{itemize} 
     Then by further pigeonholing for $\F'_{(\rho_1)} \coloneqq \bigsqcup_{F_{\rho_1} \in \widetilde{\F}_{\rho_1}} \F'[F_{\rho_1}]$ at scale $\rho_2$, there is $\widetilde{\F}_{\rho_2} \subset \F_{\rho_2}$ such that 
    \begin{itemize}
        \item $\widetilde{\F}_{\rho_2}$ is \textit{spaced} in the sense that $\{(\log(1/\delta))F_{\rho_2}\}_{F_{\rho_2} \in \widetilde{\F}_{\rho_2}}$ is essentially distinct. In particular, for any $\tau \leq \rho_2$, the sets $\{F_\tau \in \F_{\tau}: \F[F_\tau] \cap \F[F_{\rho_2}] \neq \varnothing\}_{F_{\rho_2} \in \widetilde{\F}_{\rho_2}}$ are pairwise disjoint.

        \item         Let $\F'_{(\rho_2)} \coloneqq \bigsqcup_{F_{\rho_2} \in \widetilde{\F}_{\rho_2}} \F_{(\rho_1)}'[F_{\rho_2}]$. Then $|\F'_{(\rho_2)}| \gtrsim (\log(1/\delta))^{-C} |\F'_{(\rho_1)}|$.
     \end{itemize} 
     Iterating such pruning procedure down to $\rho_M = \delta$. Since $(\log(1/\delta))^{CM} \lessapprox 1$, we have $|\F_{(\rho_M)}'|\gtrapprox |\F'|$. The good thing is that $\F_{(\rho_M)}'$ is spaced at all scales $\rho_k$, and so is uniquely encoded by a tree with $N$ levels: For any $k$, there is no ambiguity in talking about its thickened flakes $\F'_{(\rho_M), \rho_k}$ at each scale $\rho_k$ (exactly the elements in $\F_{\rho_k}''$ that intersect $\F'_{(\rho_M)}$). 
     
     Now we do another round of pigeonholing starting from the leaves of the tree. To begin with, there is $\F''_{(\rho_M), \rho_{M-1}} \subset \F'_{(\rho_M), \rho_{M-1}}$ such that $|\F'_{(\rho_M)}[F_{\rho_{M-1}}]|$ is roughly the same for all $F_{\rho_{M-1}} \in \F''_{(\rho_M), \rho_{M-1}}$, and 
     $$\F_{(\rho_{M-1})}'' \coloneqq \bigsqcup_{F_{\rho_{M-1}} \in \F''_{(\rho_M), \rho_{M-1}}} \F'_{(\rho_M)}[F_{\rho_{M-1}}]$$ 
     satisfies $|\F_{(\rho_{M-1})}''| \gtrsim (\log(1/\delta))|\F'_{(\rho_M)}|$. Next, there is $\F''_{(\rho_M), \rho_{M-2}} \subset \F'_{(\rho_M), \rho_{M-2}}$ such that $|\F''_{(\rho_{M-1})}[F_{\rho_{M-2}}]|$ is roughly the same for all $F_{\rho_{M-2}} \in \F''_{(\rho_M), \rho_{M-2}}$, and 
     $$\F_{(\rho_{M-2})}'' \coloneqq \bigsqcup_{F_{\rho_{M-2}} \in \F''_{(\rho_M), \rho_{M-2}}} \F'_{(\rho_{M-1})}[F_{\rho_{M-2}}]$$ 
     satisfies $|\F_{(\rho_{M-2})}''| \gtrsim (\log(1/\delta))|\F''_{(\rho_{M-1})}|$. Iterating such pruning procedure up to $\rho_1$. Define $\F'' \coloneqq \F''_{(\rho_1)}$. Since $(\log(1/\delta))^{CM} \lessapprox 1$, we have $|\F''| \gtrapprox |\F'_{(\rho_M)}| \gtrapprox |\F'|$. Note that $\F''$ is still spaced at all scales $\rho_k$, and it is also uniform in the sense of Definition~\ref{Def:uniform_flakes}. Let $N_k$ be its branching numbers. Then
     \begin{align*}
         N_k \cdot |\F''_{\rho_k}| \gtrsim |\F''| \gtrapprox |\F'|.
     \end{align*}
     
     By tracing the entire pruning process, we see that $\F''_{\rho_k}$ is a subset of $\F_{\rho_k}$. By Remark~\ref{Rem:strong_sticky}, since $\F$ is $\delta^{-\eta_1}$-sticky, we know that $\delta^{-2+\eta_1} \leq |\F| \leq \delta^{-2-\eta_1}$, and for each $\delta \leq \rho \leq 1$, we have $\delta^{2\eta_1} \rho^{-2} \leq |\F_\rho| \leq \delta^{-2\eta_1} \rho^{-2}$. In particular, $|\F''_{\rho_k}| \leq \delta^{-2\eta_1} \rho_k^{-2}$, and $|\F'| \gtrapprox \delta^{\eta_2}|\F| \geq \delta^{-2+\eta_1+\eta_2}$. Combining these two inequalities with the above displayed relation, we get $N_k \gtrapprox \delta^{3\eta_1+\eta_2}(\rho_k/\delta)^{2}$. On the other hand, for each $F_{\rho_k} \in \F''_{\rho_k}$,  $\F''[F_{\rho_k}] \subset \F[F_{\rho_k}]$, so we trivially have $N_k \sim |\F''[F_{\rho_k}]| \leq |\F[F_{\rho_k}]| \leq \delta^{-\eta_1}(\rho_k/\delta)^{2}$.

     We also need to check scales other than $\rho_k$. This is straightforward due to our spaced condition. For each $\delta \leq \rho \leq 1$, there exists some $1\leq k \leq M$ such that $\rho_k \leq \rho \leq \rho_{k-1}$. Note that $\F''_\rho = \bigsqcup_{F_{\rho_{k-1}}\in \F_{\rho_{k-1}}} (\F''_{\rho_k}[F_{\rho_{k-1}}])_\rho$. So we immediately have
     \begin{align*}
         \delta^{3\eta_1+\eta_2}(\rho_k/\delta)^{2} \lessapprox N_k \leq |\F''[F_\rho]| \leq N_{k-1} \lesssim \delta^{-\eta_1}(\rho_{k-1}/\delta)^{2}.
     \end{align*}
     Since $\rho_k/\rho_{k-1} = \delta^{1/M} \gtrapprox 1$, this immediately implies
     \begin{align*}     \delta^{4\eta_1+\eta_2}(\rho/\delta)^{2} \leq |\F''[F_\rho]| \leq \delta^{-2\eta_1}(\rho/\delta)^{2},
     \end{align*}
     as long as $\delta$ is sufficiently small compared with $\eta_1$. 

     Finally, by Proposition~\ref{Prop:inherit}(B), we have $\Delta_{max}(\F'') \leq \Delta_{max}(\F) \leq \delta^{-\eta_1}$. Therefore, we can conclude that $\F''$ is $\delta^{-(4\eta_1+\eta_2)}$-sticky. 
\end{proof}

By slightly modifying the proof of Lemma~\ref{Lem:retain_sticky}, we actually have:
\begin{lem}\label{Lem:split_sticky}
    For all $\eta_1 > 0$, there exist $\delta_0>0$ so that the following holds for $\delta \in (0,\delta_0]$. Let $\eps > 0$ and $\rho \in [\delta^{1-\eps}, \delta^\eps]$. Let $\F$ be a $\delta^{-\eta_1}$-sticky set of $\delta$-flakes, and let $\F'$ be an arbitrary subset of $\F$ with $|\F'| \gtrapprox \delta^{\eta_2} |\F|$, then $\F'$ contains a uniform subset $\F''$ with $|\F''| \gtrapprox |\F'|$ that is $\delta^{-(4\eta_1+\eta_2)}$-sticky. Moreover, $\F''$ is uniform at scale $\rho$, and satisfies
    \begin{itemize}
        \item $\F_\rho''$ is a $\rho^{-(5\eta_1+\eta_2)/\eps}$-sticky set of $\rho$-flakes.
        \item 
        For each $F_\rho \in \F_\rho''$, let $\widehat{\F}$ be the image of the flakes in $\F''[F_\rho]$ under the rescaling taking $F_\rho$ to $B_1$, then $\widehat{\F}$ is a $(\delta/\rho)^{-(4\eta_1+\eta_2)/\eps}$-sticky set of $\delta/\rho$-flakes.
    \end{itemize}
\end{lem}
\begin{proof}
    Recall that in the proof of Lemma~\ref{Lem:retain_sticky}, we did two rounds of multiscale pigeonholing with respect to the scales $\rho_k$ ($k=1, \dots, M$). Now we do everything with respect to the scales $\rho_k$ as well as $\rho$. Let $\F''$ be the resulting set of $\delta$-flakes. Note that $\F''$ enjoys all the properties as described in Lemma~\ref{Lem:retain_sticky}, and it is also uniform at scale $\rho$. We need to check $\F''$ satisfies the two bullet points.

    Let $F_\rho \in \F_\rho''$. For any $\tau \in [\delta,\rho]$, and $F_\tau \in (\F''[F_\rho])_\tau$, we have $(\F''[F_\rho])[F_\tau] = \F''[F_\tau]$ by the spaced condition. We know from the proof of Lemma~\ref{Lem:retain_sticky} that $\delta^{4\eta_1+\eta_2}(\tau/\delta)^{2} \leq |\F''[F_\tau]| \leq \delta^{-2\eta_1}(\tau/\delta)^{2}$, and so
    \begin{align*}    \delta^{4\eta_1+\eta_2}(\tau/\delta)^{2} \leq |(\F''[F_\rho])[F_\tau]| 
    \leq \delta^{-2\eta_1}(\tau/\delta)^{2}.
    \end{align*}
    After rescaling, this implies that for any $\theta \in [\delta/\rho ,1]$, $\tau \coloneqq \rho\,\theta \in [\delta,\rho]$ and $F_\theta \in \widehat{\F}_\theta$,
    \begin{align*}    (\delta/\rho)^{(4\eta_1+\eta_2)/\eps}(\theta/(\delta/\rho))^{2}
    \leq \delta^{4\eta_1+\eta_2}(\tau/\delta)^{2} \leq |\widehat{\F}[F_\theta]|,\\
    |\widehat{\F}[F_\theta]|
    \leq \delta^{-2\eta_1}(\tau/\delta)^{2}
    \leq
    (\delta/\rho)^{-2\eta_1/\eps}(\theta/(\delta/\rho))^{2}.
    \end{align*}
    Here we used $\delta \geq (\delta/\rho)^{1/\eps}$ due to $\rho \geq \delta^{1-\eps}$. Also, 
    $$\Delta_{max}(\widehat{\F}) = \Delta_{max}(\F''[F_\rho]) \leq \Delta_{max}(\F'') \leq \delta^{-\eta_1} \leq (\delta/\rho)^{-\eta_1/\eps}.$$ Thus, the second bullet point follows.

    Now we turn to the first bullet point. For any $\rho_k \in [\rho, 1]$, and $F_{\rho_k} \in (\F_\rho'')_{\rho_k}$, we have $ |\F_\rho''[F_{\rho_k}]| \sim |\F''[F_{\rho_k}]|/|\F''[F_\rho]|$ by the uniform and spaced condition. However, recall that
     \begin{align*}
         \delta^{3\eta_1+\eta_2}(\rho_k/\delta)^{2} \lessapprox |\F''[F_{\rho_k}]| \lesssim \delta^{-\eta_1}(\rho_{k}/\delta)^{2}
     \end{align*}
     and 
     \begin{align*}
         \delta^{3\eta_1+\eta_2}(\rho/\delta)^{2} \lessapprox |\F''[F_{\rho}]| \lesssim \delta^{-\eta_1}(\rho/\delta)^{2},
     \end{align*}
     so we have
     \begin{align*}
         \delta^{4\eta_1 + \eta_2}(\rho_k/\rho)^2 \lessapprox |\F_\rho''[F_{\rho_k}]| \lessapprox
         \delta^{-(4\eta_1 + \eta_2)}(\rho_k/\rho)^2.
     \end{align*}
     We also need to check scales other than $\rho_k$. This is straightforward due to our spaced condition. For each $\rho \leq \theta \leq 1$, there exists some $1\leq k \leq M$ such that $\rho_k \leq \theta \leq \rho_{k-1}$.\footnote{If $\rho_{k_0}$ is the smallest $\rho_k > \rho$, then we replace $\rho_k \leq \theta \leq \rho_{k-1}$ with $\rho \leq \theta \leq \rho_{k_0}$, and the arguments are exactly the same.} Note that $(\F''_\rho)_\theta = \bigsqcup_{F_{\rho_{k-1}}\in \F_{\rho_{k-1}}} (\F''_{\rho_k}[F_{\rho_{k-1}}])_\theta$. So we immediately have
     \begin{align*}
         \delta^{4\eta_1 + \eta_2}(\rho_k/\rho)^2 \lessapprox |\F_\rho''[F_{\rho_k}]| \leq |\F_\rho''[F_{\theta}]| \leq |\F_\rho''[F_{\rho_{k-1}}]| \lessapprox
         \delta^{-(4\eta_1 + \eta_2)}(\rho_{k-1}/\rho)^2.
     \end{align*}
     Since $\rho_k/\rho_{k-1} = \delta^{1/M} \gtrapprox 1$, this immediately implies
     \begin{align*}     \delta^{5\eta_1+\eta_2}(\theta/\rho)^{2} \leq |\F_\rho''[F_{\theta}]|
     \leq \delta^{-(5\eta_1+\eta_2)}(\theta/\rho)^{2},
     \end{align*}
     as long as $\delta$ is sufficiently small compared with $\eta_1$. By using $\delta \geq \rho^{1/\eps}$ due to $\rho \leq \delta^{\eps}$, we get
     \begin{align*}     \rho^{(5\eta_1+\eta_2)/\eps}(\theta/\rho)^{2} \leq |\F_\rho''[F_{\theta}]|
     \leq \rho^{-(5\eta_1+\eta_2)/\eps}(\theta/\rho)^{2},
     \end{align*}

     Finally, by Remark~\ref{Rem:strong_sticky}, we know that $|\F_\rho''| \leq |\F_\rho| \leq \delta^{-2\eta_1}\rho^{-2}$, which means $\Delta(\F_\rho'', B_1) \leq \delta^{-2\eta_1}$. Also, $C_F(\F) \leq \delta^{-2\eta_1}$ implies $C_F(\F'') \lessapprox \delta^{-2\eta_1-\eta_2}$ since $|\F''| \gtrapprox |\F'| \gtrapprox \delta^{\eta_2}|\F|$, which in turn implies $C_F(\F_\rho'') \lessapprox \delta^{-2\eta_1-\eta_2}$ by uniformity of $\F''$ at scale $\rho$ and Proposition~\ref{Prop:inherit}(A). Therefore, $\Delta_{max}(\F_\rho'') = C_F(\F_\rho'') \Delta(\F_\rho'', B_1) \lessapprox \delta^{-4\eta_1-\eta_2}$ by Definition~\ref{Def:Frostman}. And so $\Delta_{max}(\F_\rho'') \leq \delta^{-5\eta_1-\eta_2}$ as long as $\delta$ is sufficiently small compared with $\eta_1$. By using $\delta \geq \rho^{1/\eps}$, this becomes $\Delta_{max}(\F_\rho'') \leq \rho^{-(5\eta_1+\eta_2)/\eps}$. This completes the proof of the first bullet point.
\end{proof}

\begin{proof}[Proof of Proposition~\ref{Prop:self-similar}] 
We start by largely repeating the proof of \cite[Proposition~5.1]{GWZ26streamlined} with $(\F,Y)$ in place of $(\V,Y)$ and $\F_\rho$ in place of $\W$, but with some modifications.

By pigeonholing, we can select a subset $\F'' \subset \F$ with $|\F''| \gtrapprox \delta^{\eta}|\F|$ such that $(\F'',Y)$ is a $\gtrapprox 1$ refinement of $(\F,Y)$, and $|Y(F)|/|F| \sim \lambda(\F'',Y) \gtrapprox \lambda(\F,Y) \geq \delta^\eta$ for each $F \in \F''$. Since $\F$ is $\delta^{-\eta}$-sticky, we can apply Lemma~\ref{Lem:split_sticky} to $\F$ with $\eta$ in place of $\eta_1, \eta_2$, and $\F''$ in place of $\F'$ there. Let $\F'$ be the output of Lemma~\ref{Lem:split_sticky} ($\F''$ therein). Then $\F'$ is uniform with $|\F'| \gtrapprox |\F''|$ and is $\delta^{-5\eta}$-sticky. Moreover, $\F'$ is uniform at scale $\rho$, and satisfies
\begin{itemize}
    \item $\F_\rho'$ is a $\rho^{-6\eta/\eps}$-sticky set of $\rho$-flakes.
    \item 
    For each $F_\rho \in \F_\rho'$, let $\widehat{\F}$ be the image of the flakes in $\F'[F_\rho]$ under the rescaling taking $F_\rho$ to a flake of thickness $1$, then $\widehat{\F}$ is a $(\delta/\rho)^{-5\eta/\eps}$-sticky set of $\delta/\rho$-flakes.
\end{itemize}
Note that $\F' = \bigsqcup_{\F_\rho'} \F'[F_\rho]$, and $(\F',Y)$ is still a $\gtrapprox 1$ refinement of $(\F,Y)$.

Since $\F'$ is uniform at scale $\rho$ and $|Y(F)|/|F| \gtrapprox \lambda(\F, Y)$ for each $F \in \F'$, we can apply Lemma~\ref{Lem:fact_mult} (Remark~\ref{Rem:factor_only_shading}) to $(\F', Y)$. The output is a refinement $(\F',Y')$, whose properties can be summarized as follows:
\begin{enumerate}[label={\rm (\arabic*)}]
    \item \( (\F', Y') \) is a \( \gtrapprox 1 \) refinement of \( (\F, Y) \) with $\F' = \bigsqcup_{\F_\rho'} \F'[F_\rho]$.
    \item \( \lambda(\F_\rho', Y_{\F_\rho'}) \gtrapprox \lambda(\F, Y) \).
    \item \( (\F_\rho', Y_{\F_\rho'}) \) has constant multiplicity $\mu_{\rm coarse}$.
    \item For each \( F_\rho \in \F_\rho' \), \( (\F'[F_\rho], Y') \) has constant multiplicity. This multiplicity is also the same for all \( F_\rho \in \F_\rho' \), and we denote it by $\mu_{\rm fine}$. 
    \item $\mu(\F, Y) \lessapprox \mu_{\rm coarse} \mu_{\rm fine}$.
    \item For every \( x \in U(\F', Y') \),
    \[
    \F'_{Y'}(x) \subset \bigcup_{F_\rho \in \F_\rho'} (\F'[F_\rho])_{Y'}(x).
    \]
    \item \( |U(\F', Y') \cap B(x, \rho)| \) is roughly the same for every \( x \in U(\F_\rho', Y_{\F_\rho'}) \). 
    \item 
    For each $F_\rho \in \F_\rho'$, $U(\F'[F_\rho], Y') \subset Y_{\F_\rho'}(F_\rho)$. Thus $U(\F',Y') \subset U(\F_\rho', Y_{\F_\rho'})$. 
\end{enumerate}
The last item above is not explicitly listed in \cite[Proposition~5.1]{GWZ26streamlined}, but is clear from the proof. So far we have finished the construction of $(\F',Y')$ and $(\F_\rho', Y_{\F_\rho'})$, as well as the proof of Item~(iv) and (v) in Proposition~\ref{Prop:self-similar}. It suffices to show that Items~(i)(ii)(iii) hold as long as we choose $\eta$ small depending on $\eps$. 

By Definition~\ref{Def:eps_extremize}, if $(\F,Y)$ is $\eta$-extremal, then $|U(\F,Y)| \leq \delta^{\sigma-\eta}$. Thus
\begin{align}\label{Eq:total_mult_extreme}
    \mu(\F,Y) = \frac{\sum_{\F} |Y(F)|}{|U(\F,Y)|} = \lambda(\F,Y) \frac{\sum_{\F} |F|}{|U(\F,Y)|} \geq \delta^{\eta} \frac{\delta^{-2+\eta}\delta^2}{\delta^{\sigma-\eta}} = \delta^{-\sigma+3\eta}.
\end{align}

Since Item~(2) above says \( \lambda(\F_\rho', Y_{\F_\rho'}) \gtrapprox \lambda(\F, Y) \geq \delta^\eta\), for $\delta$ small depending on $\eta$ we have $\lambda(\F_\rho', Y_{\F_\rho'}) \geq \delta^{2\eta}$. This implies $\lambda(\F_\rho', Y_{\F_\rho'}) \geq \rho^{2\eta/\eps}$ by our assumption that $\rho \leq \delta^\eps$. Recall that $\F_\rho'$ is $\rho^{-6\eta/\eps}$-sticky. So by choosing $\eta < \eps^3$ small enough depending on $\eps$ and using Definition~\ref{Def:sigma} for $\sigma$ at scale $\rho$, we have 
\begin{enumerate}[label=(\alph*)]
    \item $\F_\rho'$ is $\rho^{-\eps^3}$-sticky.
    \item $\lambda(\F_\rho', Y_{\F_\rho'}) \geq \rho^{\eps^3}$.
    \item $|U(\F_\rho', Y_{\F_\rho'})| \geq \rho^{\sigma + \eps^3}$ for $\delta$ small depending on $\eps$. Thus 
    \begin{align*}
    \mu(\F_\rho', Y_{\F_\rho'}) = \frac{\sum_{\F_\rho'} |Y(F_\rho)|}{|U(\F_\rho', Y_{\F_\rho'})|} \leq \frac{\sum_{\F_\rho'} |F_\rho|}{|U(\F_\rho', Y_{\F_\rho'})|} \leq \frac{\rho^{-\eps^3} \rho^{-2} \rho^2}{\rho^{\sigma + \eps^3}} = \rho^{-\sigma - 2\eps^3}.
    \end{align*}
    The second-to-last step uses the fact that $\F_\rho'$ is $\rho^{-\eps^3}$-sticky. In view of Item~(3), we conclude that    \begin{align}\label{Eq:rho_coarse}
        \mu_{\rm coarse} = \mu(\F_\rho', Y_{\F_\rho'}) \leq \rho^{-\sigma - 2\eps^3}.
    \end{align}
    This verifies the first half of Item~(ii) in Proposition~\ref{Prop:self-similar}. 
\end{enumerate}

On the other hand, since $\lambda(\F', Y') \gtrapprox \lambda(\F',Y) \sim \lambda(\F'',Y) \gtrapprox \delta^\eta$ and $\F' = \bigsqcup_{\F_\rho'} \F'[F_\rho]$, we know the existence of one $F_\rho \in \F_\rho'$ with $\lambda(\F'[F_\rho],Y') \gtrapprox \delta^\eta$. For $\delta$ small depending on $\eta$, we have $\lambda(\F'[F_\rho],Y') \geq \delta^{2\eta}$. This implies $\lambda(\F'[F_\rho],Y') \geq (\delta/\rho)^{2\eta/\eps}$ by our assumption that $\rho \geq \delta^{1-\eps}$. Let $(\widehat{\F},\widehat{Y})$ be the image of $(\F'[F_\rho], Y')$ under the rescaling taking $F_\rho$ to $B_1$, then $(\widehat{\F},\widehat{Y})$ is a set of $\delta/\rho$-flakes with $\lambda(\widehat{\F},\widehat{Y}) = \lambda(\F'[F_\rho],Y') \geq (\delta/\rho)^{2\eta/\eps}$. Recall also that $\widehat{\F}$ is $(\delta/\rho)^{-5\eta/\eps}$-sticky. So by choosing $\eta$ small depending on $\eps$ and using Definition~\ref{Def:sigma} for $\sigma$ at scale $\delta/\rho$, we have 
\begin{enumerate}[label=(\alph*)]
    \item $\widehat{\F}$ is $(\delta/\rho)^{-\eps^3}$-sticky.
    \item $\lambda(\widehat{\F}, \widehat{Y}) \geq (\delta/\rho)^{\eps^3}$.
    \item $|U(\widehat{\F}, \widehat{Y})| \geq (\delta/\rho)^{\sigma + \eps^3}$ for $\delta$ small depending on $\eps$. Thus
    \begin{align*}
    \mu(\widehat{\F}, \widehat{Y}) = \frac{\sum_{\widehat{\F}} |\widehat{Y}(\widehat{F})|}{|U(\widehat{\F}, \widehat{Y})|} \leq \frac{\sum_{\widehat{\F}} |\widehat{F}|}{|U(\widehat{\F}, \widehat{Y})|} \leq \frac{(\delta/\rho)^{-\eps^3} (\delta/\rho)^{-2} (\delta/\rho)^2}{(\delta/\rho)^{\sigma + \eps^3}} = (\delta/\rho)^{-\sigma - 2\eps^3}.
    \end{align*}
    The second-to-last step uses the fact that $\widehat{\F}$ is $(\delta/\rho)^{-\eps^3}$-sticky. In view of Item~(4) and the rescaling invariance of multiplicities, we conclude that \begin{align}\label{Eq:rho_fine}
        \mu_{\rm fine} = \mu(\F'[F_\rho], Y') = \mu(\widehat{\F}, \widehat{Y}) \leq (\delta/\rho)^{-\sigma - 2\eps^3}.
    \end{align}
    This verifies the first half of Item~(iii) in Proposition~\ref{Prop:self-similar}.
\end{enumerate}

Now combining \eqref{Eq:total_mult_extreme}, \eqref{Eq:rho_coarse}, \eqref{Eq:rho_fine} with Item~(5), we can further conclude that
\begin{align*}
    \mu_{\rm coarse} \gtrapprox \mu(\F,Y)/\mu_{\rm fine} &\geq \delta^{3\eta + 2\eps^3} \rho^{-\sigma - 2\eps^3}\\
    \Longrightarrow
    \mu(\F_\rho', Y_{\F_\rho'}) &\geq \delta^{4\eta + 2\eps^3} \rho^{-\sigma - 2\eps^3} \geq \rho^{-\sigma + 4\eps^2}
\end{align*}
by using $\eta < \eps^3$ and $\delta \geq \rho^{1/\eps}$, 
and that
\begin{align*}
    \mu_{\rm fine} \gtrapprox \mu(\F,Y)/\mu_{\rm coarse} 
    &\geq \delta^{3\eta+2\eps^3} (\delta/\rho)^{-\sigma - 2\eps^3}\\
    \Longrightarrow
    \mu(\F'[F_\rho], Y') &\geq \delta^{4\eta+2\eps^3} (\delta/\rho)^{-\sigma - 2\eps^3} \geq 
    (\delta/\rho)^{-\sigma + 4\eps^2}
\end{align*}
by using $\eta < \eps^3$ and $\delta \geq (\delta/\rho)^{1/\eps}$. These lower bounds confirm the second halves of Items~(ii)~and~(iii). Moreover, we have
\begin{align*}
|U(\F_\rho', Y_{\F_\rho'})| = \frac{\sum_{\F_\rho'} |Y(F_\rho)|}{\mu(\F_\rho', Y_{\F_\rho'})} \leq \frac{\sum_{\F_\rho'} |F_\rho|}{\mu(\F_\rho', Y_{\F_\rho'})} \leq \frac{\rho^{-\eps^3} \rho^{-2} \rho^2}{\rho^{-\sigma + 4\eps^2}} \leq \rho^{\sigma - 5\eps^2},
\end{align*}
which together with our previous Items~(a)~and~(b) for $(\F_\rho', Y_{\F_\rho'})$ gives Item~(i) in Proposition~\ref{Prop:self-similar}. This completes the proof of the whole Proposition~\ref{Prop:self-similar}. 

Finally, although we will not need it, we remark that 
\begin{align*}
|U(\widehat{\F}, \widehat{Y})| = \frac{\sum_{\widehat{\F}} |\widehat{Y}(\widehat{F})|}{\mu(\widehat{\F}, \widehat{Y})} \leq \frac{\sum_{\widehat{\F}} |\widehat{F}|}{\mu(\widehat{\F}, \widehat{Y})} \leq \frac{(\delta/\rho)^{-\eps^3} (\delta/\rho)^{-2} (\delta/\rho)^2}{(\delta/\rho)^{-\sigma + 4\eps^2}} 
\leq
(\delta/\rho)^{\sigma - 5\eps^2},
\end{align*}
which together with our previous Items~(a)~and~(b) for $(\widehat{\F}, \widehat{Y})$ shows that $(\widehat{\F}, \widehat{Y})$ is also $\eps$-extremal.
\end{proof}
In the proof of Proposition~\ref{Prop:self-similar}, $Y_{\F_\rho'}$ is essentially induced by $Y'$, so we can further assume that for any $F \in \F'$ and $F_\rho \in \F_\rho'$ containing $F$, we have $N_{\rho/10n}(Y'(F)) \subset Y_{\F_\rho'}(F_\rho)$. This regularity of $Y_{\F_\rho'}$ is technically useful in the next section, but does not affect the main ideas.

\section{Planiness and graininess}\label{Sec:plaininess}

In this section, we will use the geometric lemmas established in Section~\ref{Sec:geometry} to find a plane map and then construct the local grains as in \cite[Section~4]{WZ26sticky}. The main difficulty is that flakes have both major and minor angles, which makes the geometry much more complicated than tubes, as we will see in the proof. 

Before that, let us first show a reduction about robust transversality: If $(\F,Y)$ is an extremal collection of flakes, then the direction subspaces of the flakes passing through a typical point cannot focus too tightly. This can be viewed as a natural analogue of \cite[Lemma~3.4]{WZ26sticky}.

\begin{lem}\label{Lem:robust_transverse}
For all $\eps>0$, there exist $\eta>0$ and $\delta_0>0$ so that the following
holds for all $\delta \in (0,\delta_0]$. Let $(\F,Y)$ be an $\eta$-extremal set of $\delta$-flakes. Then there is a $\approx 1$ refinement $(\F',Y')$ of $(\F,Y)$ so that for each $x \in \R^n$ and $V \in \Gr(n,n-2)$, we have
\begin{align}\label{Eq:robust_transverse}
    |\{F\in {\mathbb F}'_{Y'}(x): \angle(V,\dir(F)) <\delta^{\varepsilon}\}| \le \delta^{-\sigma+\varepsilon\sigma/2 }.
\end{align}
\end{lem}
\begin{proof}
    Apply Proposition~\ref{Prop:self-similar} with $\eps\sigma/2$ in place of $\eps$, and $\rho = \delta^\eps$. By Lemma~\ref{Lem:flake_contain}, the $\delta$-flakes contributing to \eqref{Eq:robust_transverse} are all contained in the $\rho$-flake $N_\rho(V) + x$. Since the resulting set $(\F_\rho', Y_{\F_\rho'})$ of $\rho$-flakes is essentially distinct, these $\delta$-flakes can be further covered by $O(1)$ $\rho$-flakes from $\F_\rho'$. Now by Item~(iii) in Proposition~\ref{Prop:self-similar}, we can bound the left hand side of \eqref{Eq:robust_transverse} by 
    \begin{align*}
        \mu(\F'[F_\rho], Y')(x)
        \leq
        (\delta/\rho)^{-\sigma-\eps\sigma/2}
        =
        \delta^{-\sigma + \eps\sigma/2 + \eps^2\sigma/2},
    \end{align*}
    which yields the desired bound. 
\end{proof}
Roughly speaking, Lemma~\ref{Lem:robust_transverse} allows us to always assume the major angles between flakes to be $\approx 1$.  We will see the existence of extremal sets of flakes $(\F,Y)$ with additional structural
properties, the first of which is the existence of a Lipschitz plane map. 

\begin{definition}\label{Def:plane_map}
    Let $(\F,Y)$ be a set of $\delta$-flakes. A \textit{plane map} for $(\F,Y)$ is a function $\HH: U(\F,Y) \rightarrow \Gr(n,1)$ so that for each $x\in U(\F,Y)$, we have 
    \begin{align*}
        \dir(F) \cap B_1 \subset N_{2\delta}(\HH(x)^\perp) \quad \text{for all } F \in \F_Y(x).
    \end{align*}
\end{definition}
We expect $U(\F,Y)$ to arrange itself into a union of thin slab segments, which we call \textit{local grains}. They look like an Ahlfors-David regular set (\cite[Definition~4.2]{WZ26sticky}) in the normal direction. We record the specific definition for the reader's convenience:
\begin{definition}\label{Def:AD_set}
    For $n \geq 1$, $\alpha\in [0,n]$, $C \geq 1$ and $\delta>0$, we say a set $E \subset \R^n$ is a $(\delta, \alpha, C)_n$-ADset if for all $x \in \R^n$, all $\rho \geq \delta$, and all $r \geq \rho$, we have
    \begin{align*}
        \EE_\rho(E \cap B(x,r)) \leq C(r/\rho)^\alpha.
    \end{align*}
    {\color{black} Here $\EE_\rho(\cdot)$ denotes the smallest number of $\rho$-balls needed to cover a set, which is comparable to the largest number of $\rho$-separated points in a set.}
\end{definition}
One may consult \cite[Remark~4.1]{WZ26sticky} for some intuitive explanation of this definition. We can now state the main result of this section.  
\begin{proposition}\label{Prop:local_grains}
     For all $\eps, \delta_0>0$, there exist $\delta \in (0,\delta_0]$ and an $\eps$-extremal set of $\delta$-flakes $(\F,Y)$ with a $1$-Lipschitz\footnote{By Remark~\ref{Rem:major_metric}, the major angle yields a natural metric on $\Gr(n,1)$. And \textit{$C$-Lipschitz} means $\angle(\HH(x), \HH(y)) \leq C \norm{x-y}$ for all admissible pairs of $x,y$.} plane map $\HH$ for it, such that for all $\rho \in [\delta,1]$ and $x\in U(\F,Y)$, $\pi_{\HH(x)}(U(\F,Y) \cap B(x,\rho^{1/2}))$ is a $(\rho, 1-\sigma, \delta^{-\eps})_1$-ADset\footnote{Here we are fixing an arbitrary isometry $\HH(x) \cong \R$ to make sense of the terminology. Equivalently, since $\HH(x) \subset \R^n$, one can also just write ``$(\rho, 1-\sigma, \delta^{-\eps})_n$-ADset'' instead.}. 
\end{proposition}
\begin{rem}
    We do not know a priori whether $1 - \sigma \geq 0$ or not, but one can see from the proof that this must be true. Intuitively, this is because we will use local hairbrush-type arguments to grow grains somewhere in the proof, which roughly tells that the ``dimension'' of the union is at least $n-1$, and so $\sigma \leq 1$.
\end{rem}
Roughly speaking, Proposition~\ref{Prop:local_grains} implies that each $N_\rho(U(\F,Y)) \cap B(x,\rho^{1/2})$ is contained in a union of at most $\rho^{\sigma/2 - 1/2}$ parallel rectangular slabs of dimensions roughly $\rho \times \rho^{1/2} \times \rho^{1/2}\1$. Since $U(\F,Y)$ can be covered by $\delta^{-O(\eps)}\rho^{\sigma/2-n/2}$ balls of radius $\rho^{1/2}$, this implies that $|N_\rho(U(\F,Y))| \lesssim \delta^{-O(\eps)}\rho^\sigma$. On the other hand, this containment is nearly sharp, since $|N_\rho(U(\F,Y))| \gtrsim \delta^{O(\eps)}\rho^\sigma$. 
\begin{rem}
    Our Proposition~\ref{Prop:local_grains} provides the so-called \textit{local grains} structure for $(n-2)$-plane Kakeya in $\R^n$, whose statement is similar to that in \cite[Proposition~4.1(2)]{WZ26sticky}, but the proof here is more complicated, as one will soon see. It is conceivable that a \textit{global grains} structure similar to that in \cite[Proposition~4.1(1)]{WZ26sticky} could also be established, but we will not need it for our purposes. In fact, we only need to apply Proposition~\ref{Prop:local_grains} at some special scales later on, but we think that it is beneficial to present the full version of it, which does not incur additional difficulties.
\end{rem}

{\color{black} The remaining parts of this section will be devoted to the proof of Proposition~\ref{Prop:local_grains}, which will be divided into three subsections. In Section~\ref{subsec:find_plane_map}, we will show that every extremal set
of flakes has a refinement with a plane map. In Section~\ref{subsec:Lips_reg_plane_map}, we will show
that if an extremal set of flakes has a plane map, then after restricting to a refinement, the plane map must be Lipschitz with controlled Lipschitz norm. In Section~\ref{subsec:local_grains}, we will show that if an extremal set of flakes has a Lipschitz
plane map, then after restricting to a refinement, the flakes arrange themselves into local grains. }

\subsection{Finding a plane map}\label{subsec:find_plane_map}

In this section, we will show that every extremal set
of flakes has a refinement with a plane map. Our main tools will be the bilinear Kakeya-type estimate Corollary~\ref{Cor:bilinear_flakes} proved in Section~\ref{Sec:geometry}, as well as the multilinear Kakeya estimate due to Carbery and Valdimarsson (\cite[Theorem~1]{CV2013}, see also \cite[Theorem~4.2]{WZ26sticky}). We recall the former in the special case when $k=n-2$, and the two families coincide:
\begin{cor}\label{Cor:bilinear_flakes'}
    Let $\F$ be a set of $\delta$-flakes in $\R^n$, then 
    \begin{align*}
        \int_{\R^n} \sum_{F_1,F_2 \in \F} \chi_{F_1} \chi_{F_2} \angle'(F_1, F_2) \angle(F_1, F_2) 
        \lesssim
        (\delta^2 |\F|)^2.
    \end{align*}
\end{cor}

First, we prove that every extremal set of flakes has a refinement that “weakly” agrees with a plane map. Here comes the difference from the tube case in \cite{WZ26sticky}: If we use Corollary~\ref{Cor:bilinear_flakes'} as in \cite[Section~4.1]{WZ26sticky}, then we are no longer immediately reduced to the plany case. Instead, there is an additional case when a generic point $x$ in $U(\F,Y)$ has all $\F_Y(x)$ sharing a common {\color{black} $(n-3)$-plane} through $x$. However, by using certain slicing arguments and applying trilinear Kakeya for tubes in each slice, this additional case can actually be further reduced to the plany case. Once the weak planiness has been established, it is not hard to upgrade it to a strong one.

\begin{lem}\label{Lem:plane_or_pivot}
    There exist $\eta,\delta_0>0$ so that the following holds for all $\delta \in (0,\delta_0]$. Let $(\F,Y)$ be an $\eta$-extremal set of $\delta$-flakes. Then at least one of the following cases holds true:
    \begin{enumerate}[label=\rm{(\roman*)}]
        \item There is a $\gtrapprox 1$ refinement $(\F',Y')$ and a function $\HH: U(\F',Y') \rightarrow \Gr(n,1)$ such that for each $x \in U(\F',Y')$, we have
        \begin{align*}
            \dir(F) \cap B_1 \subset N_{\delta^{\sigma/5}}(\HH(x)^\perp) \quad \text{for all } F \in \F'_{Y'}(x).
        \end{align*}
        
        \item There is a $\gtrapprox 1$ refinement $(\F',Y')$ and a function $\XX: U(\F',Y') \rightarrow \Gr(n,n-3)$ such that for each $x \in U(\F',Y')$, we have
        \begin{align*}
            \XX(x) \cap B_1 \subset N_{\delta^{\sigma/2}}(\dir(F)) \quad \text{for all } F \in \F'_{Y'}(x). 
        \end{align*}
    \end{enumerate}
\end{lem}
The proof of this lemma relies on two additional ingredients, which we now introduce. The first is about the maximum number of triangles sharing a common edge, which is usually called the \textit{booksize} of a graph. The version we present here is by no means optimal, but suffices for our purposes.
\begin{lem}\label{Lem:booksize_graph}
    Let $\mathbf{G}$ be a simple graph with $n$ vertices and $m$ edges such that $m \geq \lambda n^2$ for some $\lambda > 1/4$. Then there exists $uv \in \mathbf{E}$ {\rm(}$\mathbf{E}$ is the set of edges in $\mathbf{G}${\rm)} such that the number of triangles containing $uv$ is at least $(4\lambda-1)n$. 
\end{lem}
\begin{proof}
    Let $\mathcal{N}(x)$ be the set of neighbors of a vertex $x$ in $\mathbf{G}$. Let $\mathcal{T}$ be the number of triangles in $\mathbf{G}$, and for each $uv \in \mathbf{E}$, let $\mathcal{Q}(uv) \coloneqq |\mathcal{N}(u) \cap \mathcal{N}(v)|$ be the number of triangles containing $uv \in \mathbf{E}$. Since each triangle is contributed by three sides, we have 
    \begin{align*}
        \sum_{uv\in E} \mathcal{Q}(uv) = 3\mathcal{T}, 
    \end{align*}
    which guarantees the existence of an edge $uv \in \mathbf{E}$ with $\mathcal{Q}(uv) \geq 3\mathcal{T}/m$. 

    By a standard Goodman-type supersaturation inequality (see \cite[Exercise~1.1.4]{zhao2023}), we have 
    \begin{align*}
        \mathcal{T} \geq \frac{4m}{3n}\Big( m - \frac{n^2}{4} \Big).
    \end{align*}
    So we can conclude that
    \begin{align*}
        \max_{uv\in \mathbf{E}}|\mathcal{N}(u) \cap \mathcal{N}(v)| \geq \frac{3\mathcal{T}}{m} \geq \frac{4}{n}\Big( m - \frac{n^2}{4} \Big) \geq (4\lambda - 1)n,
    \end{align*}
    which is exactly what we want.
\end{proof}

The second ingredient is geometric and basically says the following: If each pair of three flakes intersects transversely nearby an affine $(n-3)$-plane and the three resulting planes differ significantly, then they must be contained in a small neighborhood of an affine $(n-1)$-plane.
\begin{lem}\label{Lem:flake_triangle}
    Let $\sigma > 0$, and $V_i \in \Gr(n,n-2)$ {\rm(}$i=1,2,3${\rm)} satisfy $\angle(V_i,V_j) \geq \delta^{\sigma/10}$ and $\angle'(V_i,V_j) \leq \delta^{9\sigma/10}$ {\rm(}$\forall\,\{i,j\} \subset \{1,2,3\}${\rm)}. Let $H$ be the output $(n-1)$-plane of Lemma~\ref{Lem:form_planes} for the pair $V_1$ and $V_2$. Recall that $H$ contains the internal bisector $\ell_{12}$ of a pair of principal vectors defining $\angle'(V_1,V_2)$. Let $X_{12} \coloneqq (V_1\cap V_2) + \ell_{12}$. For $\delta>0$ small depending on $\sigma$, if $X_{12} \cap B_1 \not\subset N_{\delta^{\sigma/2}}(V_3)$, then $V_3 \cap B_1 \subset N_{\delta^{\sigma/5}}(H)$.
\end{lem}
\begin{proof}
    Similarly to $X_{12}$, we define $X_{13} \coloneqq (V_1 \cap V_3) + \ell_{13}$, and $X_{23} \coloneqq (V_2 \cap V_3) + \ell_{23}$, where $\ell_{13}$ and $\ell_{23}$ bisect $\angle'(V_1,V_3)$ and $\angle'(V_2,V_3)$ respectively. They are all $(n-3)$-planes.
    
    By our assumption of minor angles, we have 
    $$X_{13} \cap B_1 \subset N_{\delta^{9\sigma/10}}(V_1) \cap N_{\delta^{9\sigma/10}}(V_3) \cap B_1 \subset N_{C\delta^{9\sigma/10}}(H),$$ 
    and 
    $$X_{23} \cap B_1 \subset N_{\delta^{9\sigma/10}}(V_2) \cap N_{\delta^{9\sigma/10}}(V_3) \cap B_1 \subset N_{C\delta^{9\sigma/10}}(H).$$ 
    We claim that $\angle(X_{13}, X_{23}) \geq \delta^{2\sigma/3}$. If this is not the case, then $\angle(X_{13}, X_{23}) \leq \delta^{2\sigma/3}$, and so $X_{13} \cap B_1 \subset N_{\delta^{2\sigma/3}}(X_{23})$. This together with $X_{23} \cap B_1 \subset N_{\delta^{9\sigma/10}}(V_2) \cap B_1$ implies that $X_{13} \cap B_1 \subset N_{2\delta^{2\sigma/3}}(V_2)  \cap B_1$. On the other hand, we know that $X_{13} \cap B_1 \subset N_{\delta^{9\sigma/10}}(V_1) \cap B_1$, and so $X_{13} \cap B_1 \subset N_{\delta^{9\sigma/10}}(V_1) \cap N_{2\delta^{2\sigma/3}}(V_2)  \cap B_1$. This in turn implies $X_{13} \cap B_1 \subset N_{\delta^{8\sigma/15}}(X_{12}) \cap B_1$. Here for ``in turn implies'', we are using the geometric picture that $\angle(V_1,V_2) \geq \delta^{\sigma/10}$ and $\angle'(V_1,V_2) \leq \delta^{9\sigma/10}$ mean that for any $\tau \geq \delta^{9\sigma/10}$ (and in particular $\tau = 2\delta^{2\sigma/3}$), we have $N_{\tau}(V_1) \cap N_{\tau}(V_2) \cap B_1$ is a prism of dimension $\tau \times \tau \times (\tau/\delta^{\sigma/10})  \times 1 \times \cdots \times 1$ ($(n-3)$ many $1$) whose central plane is $X_{12}$. Therefore, by Remark~\ref{Rem:major_equiv}, we can conclude that $\angle(X_{13}, X_{12}) \leq \delta^{8\sigma/15}$, and so $X_{12} \cap B_1 \subset N_{\delta^{8\sigma/15}}(X_{13}) \cap B_1 \subset N_{\delta^{\sigma/2}}(V_3)$, contradicting our assumption that $X_{12} \cap B_1 \not\subset N_{\delta^{\sigma/2}}(V_3)$. Thus the claim must hold true. 

    Let $X_{13}'\coloneqq \pi_{V_3}(X_{13})$, and $X_{23}'\coloneqq \pi_{V_3}(X_{23})$. By Lemma~\ref{Lem:faithful_proj}, $X_{13} \cap B_1 \subset N_{\delta^{9\sigma/10}}(V_3)$ implies $\dim X_{13}' = n-3$ and $\angle(X_{13}, X_{13}') \leq \delta^{9\sigma/10}$, while $X_{23} \cap B_1 \subset N_{\delta^{9\sigma/10}}(V_3)$ implies $\dim X_{23}' = n-3$ and $\angle(X_{23}, X_{23}') \leq \delta^{9\sigma/10}$. In particular, we also have $X_{13}'  \cap B_1, X_{23}' \cap B_1 \subset N_{C\delta^{9\sigma/10}}(H)$. 
    Now since $\angle(X_{13}, X_{23}) \geq \delta^{2\sigma/3}$ and the major angle is a metric by Remark~\ref{Rem:major_metric}, we have 
    $$\angle(X_{13}', X_{23}') \geq \angle(X_{13}, X_{23}) - \angle(X_{13}, X_{13}') - \angle(X_{23}, X_{23}') \geq \delta^{2\sigma/3}/2.$$ 
    This allows us to find some unit vector $w_0 \in X_{13}'$ with $\dist(w_0, X_{23}') \geq \delta^{2\sigma/3}/2$. In other words, let $w \coloneqq w_0 - \pi_{X_{23}'}(w_0)$, then it is normal to $X_{23}'$ and has $\norm{w} \geq \delta^{2\sigma/3}/2$. Let $\HH$ be the unit normal vector of $H$, then 
    $$|\HH \cdot (w/\norm{w})| \leq (|\HH \cdot w_0| + |\HH \cdot \pi_{X_{23}'}(w_0)|)/\norm{w} \lesssim \delta^{9\sigma/10}/\delta^{2\sigma/3} = \delta^{7\sigma/30}.$$ 
    Finally, any vector $v$ in $V_3 \cap B_1$ satisfies 
    $$|\HH \cdot v|\leq|\HH \cdot(w/\norm{w})|+|\HH \cdot \pi_{X_{23}'}(v)| \leq \delta^{\sigma/5},$$ 
    and so the conclusion holds true.
\end{proof}

\begin{proof}[Proof of Lemma~\ref{Lem:plane_or_pivot}]
    Define $\eps = \sigma/5$. We will select $\eta$ much smaller than some $\eta_1 \leq \min\{\eps/44,\eps\sigma/32\}$. Apply Lemma~\ref{Lem:robust_transverse} to $(\F,Y)$ with $\eps/4$ in place of $\eps$ to obtain a $\gtrapprox 1$ refinement $(\F_1,Y_1)$ so that for each $x\in \R^n$ and $V \in \Gr(n,n-2)$, we have
\begin{align}\label{Eq:robust_transverse'}
|\{F\in ({\mathbb F}_1)_{Y_1}(x): \angle(V,\dir(F)) \leq \delta^{\varepsilon/4}\}| \le \delta^{-\sigma+(\eps/4)\sigma/2 } \leq \delta^{-\sigma +4\eta_1}.
\end{align}

Let $(\F_2,Y_2)$ be a constant multiplicity refinement of $(\F_1,Y_1)$, with multiplicity $\mu \gtrapprox \delta^{-\sigma + 3\eta_1}$ (recall Remark~\ref{Rem:extreme_mult}). By Corollary~\ref{Cor:bilinear_flakes'}, we have
\begin{align}\label{Eq:bilinear_flake}
    \int_{\R^n} \sum_{F_1,F_2 \in \F_2} \chi_{F_1} \chi_{F_2} \angle'(F_1, F_2) \angle(F_1, F_2) 
    \lesssim
    (\delta^2 |\F_2|)^2 \leq \delta^{-2\eta_1}.
\end{align}
We say a point $x \in U(\F_2,Y_2)$ is \textit{broad} if 
\begin{align}\label{Eq:def_broad_pt}
    |\{(F_1,F_2) \in ((\F_2)_{Y_2}(x))^2: \angle'(F_1, F_2) \angle(F_1, F_2) \geq \delta^{\sigma-11\eta_1} \}| \geq \delta^{-2\sigma +7\eta_1},
\end{align}
otherwise we say it is \textit{narrow}. Both the set of broad and the set of narrow points can be viewed as unions
of $\delta$-cubes. Observe that each broad point contributes at least $\delta^{-\sigma-4\eta_1}$ 
to the integrand on the LHS of \eqref{Eq:bilinear_flake}, and thus by \eqref{Eq:bilinear_flake} the measure of the set of broad
points $\lesssim \delta^{\sigma + 2\eta_1}$. On the other hand, $U(\F_2, Y_2)$ has measure $\geq \delta^{\sigma + \eta_1}$ by Definition~\ref{Def:sigma} and Remark~\ref{Rem:extreme_mult}. 

Define $Y_3(F) \coloneqq \{x \in Y_2(F): x \text{ narrow}\}$, and define $\F_3 = \F_2$. Then it is still true that $|(\F_3)_{Y_3}(x)| =\mu$\footnote{Rigorously speaking, we can only locate $|(\F_3)_{Y_3}(x)|$ in some fixed range, say $[\mu, (101/100)\mu]$. However, this does not harm our computations later on, so for simplicity, we pretend that the multiplicity is exactly $\mu$.} for each $x \in U(\F_3,Y_3)$, and $|U(\F_3,Y_3)| \geq (1/2)|U(\F_2,Y_2)|$. 

Now for each $x \in U(\F_3,Y_3)$, we have
\begin{align}\label{Eq:narrow_dominate}
    |\{(F_1,F_2) \in ((\F_3)_{Y_3}(x))^2: \angle'(F_1, F_2) \angle(F_1, F_2) < \delta^{\sigma-11\eta_1} \}|  \geq \mu^2 - \delta^{-2\sigma +7\eta_1} \geq \frac{99}{100}\mu^2.
\end{align}
On the other hand, by \eqref{Eq:robust_transverse'}, we have
\begin{align}\label{Eq:non_trans_small}
    |\{(F_1,F_2) \in ((\F_3)_{Y_3}(x))^2: \angle(F_1,F_2) < \delta^{\eps/4} \}| \leq \mu \delta^{-\sigma + 4\eta_1} \leq \frac{1}{100}\mu^2. 
\end{align}
Comparing \eqref{Eq:narrow_dominate} and \eqref{Eq:non_trans_small}, we conclude that there are at least $(98/100)\mu^2$ pairs $(F_1,F_2)$ with $\angle'(F_1, F_2) \angle(F_1, F_2) < \delta^{\sigma-11\eta_1} \leq \delta^{\sigma - \eps/4}$ and $\angle(F_1,F_2) \geq \delta^{\eps/4}$, which implies $\angle'(F_1,F_2) \leq \delta^{\sigma - \eps/2}$. By applying Lemma~\ref{Lem:form_planes} with $k=n-2$, $j=n-3$, and $\rho = \delta^{\sigma - \eps/2}$, we know that $F_1$ and $F_2$ are contained in a $\lesssim \delta^{\sigma - \eps/2}$(thus $\leq \delta^{\sigma/2}$)-neighborhood of a hyperplane orthogonal to $\HH(x,\{F_1,F_2\}) \in \Gr(n,1)$. Recall that $\HH(x,\{F_1,F_2\})$ is constructed by taking the internal bisector of a pair of principal vectors defining $\angle'(F_1, F_2)$, which we denote by $\ell(x, \{F_1,F_2\})$. This seems almost what we want for Case~(i), but the tricky issue is to identify a common hyperplane for each pair of flakes, which is not immediately possible and accounts for the appearance of Case~(ii). Now we turn to the details.

Fix any $x \in U(\F_3,Y_3)$. Construct a finite graph with vertices $(\F_3)_{Y_3}(x)$ and edges 
$$\mathbf{E}\coloneqq \{ \{F_1,F_2\} \subset (\F_3)_{Y_3}(x): \angle(F_1,F_2) \geq \delta^{\eps/4}, \angle'(F_1,F_2) \leq \delta^{\sigma - \eps/2} \}.$$ 
The previous discussions tell us $|(\F_3)_{Y_3}(x)| = \mu$ and $|\mathbf{E}| \geq (49/100) \mu^2$. So by applying Lemma~\ref{Lem:booksize_graph} to this graph with $\lambda = 49/100$, we know the existence of a pair $\{F_1,F_2\} \in \mathbf{E}$ such that the set 
$$\mathbf{Q}(x, \{F_1,F_2\}) \coloneqq \{F_3 \in (\F_3)_{Y_3}(x): \{F_1,F_3\}\in \mathbf{E}, \{F_2,F_3\}\in \mathbf{E}\}$$ 
has cardinality $\geq (24/25)\mu$. Fix such a pair for each $x\in U(\F_3,Y_3)$ and write them as $\{F_1^x,F_2^x\}$ to emphasize the dependence on $x$. Define $Y_4(F) \coloneqq \{x\in Y_3(F): F \in \mathbf{Q}(x, \{F_1^x,F_2^x\}) \}$, and define $\F_4 = \F_3$. Then $(\F_4,Y_4)$ is still a $\gtrapprox 1$ refinement of $(\F,Y)$.

Now there are two possibilities. For each $x\in U(\F_4,Y_4)$, let $$X_x \coloneqq (\dir(F_1^x) \cap \dir(F_2^x)) + \ell(x, \{F_1^x, F_2^x\}).$$ Consider 
$$\mathbf{R}(x) \coloneqq \{F \in \F_4: x\in Y_4(F), X_x \cap B_1 \not\subset N_{\delta^{\sigma/2}}(\dir(F)) \}.$$ 
If 
\begin{align*}
    \int_{U(\F_4,Y_4)} |\mathbf{R}(x)| \D x \geq \frac{1}{2} \sum_{\F_4} |Y_4(F)|,
\end{align*}
then we define $Y'(F) \coloneqq \{ x\in Y_4(F): F \in \mathbf{R}(x) \}$ and $\F' = \F_4$. Note that $(\F',Y')$ is a $\gtrapprox 1$ refinement of $(\F,Y)$. Moreover, for each $x \in U(\F',Y')$, define $\HH(x) \coloneqq \HH(x,\{F_1^x,F_2^x\})$. Then by Lemma~\ref{Lem:flake_triangle} (recall that $\eps = \sigma/5$), we know that $(\F',Y')$ and $\HH$ satisfy the desired property in Case~(i). And so Case~(i) holds. Otherwise, we must have
\begin{align*}
    \int_{U(\F_4,Y_4)} |(\F_4)_{Y_4}(x) \setminus \mathbf{R}(x)| \D x \geq \frac{1}{2} \sum_{\F_4} |Y_4(F)|,
\end{align*}
and we define $Y'(F) \coloneqq \{ x\in Y_4(F): F \in (\F_4)_{Y_4}(x) \setminus \mathbf{R}(x) \}$ and $\F' = \F_4$. Note that $(\F',Y')$ is a $\gtrapprox 1$ refinement of $(\F,Y)$. Moreover, for each $x \in U(\F',Y')$, define $\XX(x) \coloneqq X_x \in \Gr(n, n-3)$. Then by our definition of $Y'$, we know that $(\F',Y')$ and $\XX$ satisfy the desired property in Case~(ii). And so Case~(ii) holds. This completes the proof of Lemma~\ref{Lem:plane_or_pivot}.
\end{proof}

As we have said before, Case~(ii) of Lemma~\ref{Lem:plane_or_pivot} is only an intermediate case and can be further reduced to Case~(i) by slicing arguments.  In fact, we have the following result:
\begin{proposition}\label{Prop:planiness}
There exist $\eta,\delta_0>0$ so that the following holds for all $\delta \in (0,\delta_0]$. Let $(\F,Y)$ be an $\eta$-extremal set of $\delta$-flakes. Then there is a $\gtrapprox 1$ refinement $(\F'',Y'')$ and a function $\HH: U(\F'',Y'') \rightarrow \Gr(n,1)$ such that for each $x \in U(\F'',Y'')$, we have
\begin{align*}
    \dir(F) \cap B_1 \subset N_{\delta^{\sigma/5}}(\HH(x)^\perp) \quad \text{for all } F \in \F''_{Y''}(x).
\end{align*}
\end{proposition}

\begin{proof}
    {\color{black} We will select $\eta$ much smaller than some $\eta_1$, which is in turn much smaller than that $\eta_1$ chosen in the proof of Lemma~\ref{Lem:plane_or_pivot}.} It suffices to argue that Case~(ii) of Lemma~\ref{Lem:plane_or_pivot} is still plany up to further suitable refinements. We will argue by using a slicing argument and exploit multilinear Kakeya for tubes inside each slice. To set thing up, suppose that we are in Case~(ii) in the proof of Lemma~\ref{Lem:plane_or_pivot}, i.e., there is a $\gtrapprox 1$ refinement $(\F',Y')$ and a function $\XX: U(\F',Y') \rightarrow \Gr(n, n-3)$ such that for each $x \in U(\F',Y')$, we have
    \begin{align*}
        \XX(x) \cap B_1 \subset N_{\delta^{\sigma/2}}(\dir(F)) \quad \text{for all } F \in \F'_{Y'}(x). 
    \end{align*}
    
    To find a good direction subspace for slicing, let us first sort the points $x$ into several groups according to the rough ranges of $\XX(x)$. By Remark~\ref{Rem:major_metric}, the major angle yields a metric on the finite-dimensional manifold $\Gr(n, n-3)$, so we can cover $\Gr(n, n-3)$ with a finitely overlapping union of $O(1)$ small caps of diameter $1/100$. By pigeonholing, there exists (at least) one cap, which we denote by $\CC$, such that 
    \begin{align*}
        \int_{\XX^{-1}(\CC)} \mu(\F',Y')(x) \D x
        \gtrsim 
        \int_{U(\F',Y')} \mu(\F',Y')(x) \D x.
    \end{align*}
    Define $Y_1(F) \coloneqq Y'(F) \cap \XX^{-1}(\CC)$ and $\F_1 \coloneqq \F'$, then $(\F_1, Y_1)$ is a $\gtrsim 1$ refinement of $(\F', Y')$. Also, for each $x \in U(\F_1,Y_1)$, we have $\XX(x) \in \CC$. 
    
    Let $(\F_2,Y_2)$ be a constant multiplicity refinement of $(\F_1,Y_1)$, with multiplicity \begin{align}\label{Eq:slice_tube_mult}
        \delta^{-\sigma + 3\eta_1} \leq \mu \leq \delta^{-\sigma - 2\eta_1}
    \end{align}(recall Remark~\ref{Rem:extreme_mult}). Let $X_0 \in \Gr(n, n-3)$ be the center of $\CC$, in the sense that $\angle(W, X_0) \leq 1/100$ for each $W \in \CC$. In particular, we have $\angle(\XX(x), X_0) \leq 1/100$ for each $x \in U(\F_2, Y_2)$. We will slice $(\F_2,Y_2)$ by affine $3$-planes $\{X_0^\perp + \vt\}_{\vt\in X_0 \cap B_1}$ perpendicular to $X_0$. Without loss of generality we may assume that $Y_2(F) \neq \varnothing$ for each $F \in \F_2$. This means $\XX(x) \cap B_1 \subset N_{\delta^{\sigma/2}}(\dir(F))$ for some $x \in \XX^{-1}(\CC)$, and so $X_0 \cap B_1 \subset N_{1/100}(\XX(x) \cap B_1) \subset N_{1/50}(\dir(F))$. By Lemma~\ref{Lem:slice_transversality} with $V = \dir(F)$ and $W = X_0$, we have $\dim(\dir(F) \cap X_0^\perp) = 1$, $\dir(F) = \pi_{\dir(F)}(X_0) \oplus (\dir(F) \cap X_0^\perp)$, and $N_\delta(\dir(F)) \cap N_\delta(X_0^\perp) \subset N_{2\delta}(\dir(F)\cap X_0^\perp)$. 
    
    The geometric importance of the above observations is that each $F \in \F_2$ intersects our slices transversely, in the sense that $F \cap(X_0^\perp + \vt)$ is (comparable to) a $3$-dimensional $\delta$-tube in $X_0^\perp + \vt$. It is for such nondegeneracy that we choose to concentrate on a dominating cap $\CC$ at the beginning. Define $$\T_\vt \coloneqq \{F \cap(X_0^\perp + \vt) : F \in \F_2\}$$ to be the resulting set of $\delta$-tubes in $X_0^\perp + \vt$ (we do not require $\T_\vt$ to be essentially distinct), and define the shading $Y_\vt$ on it by $$Y_\vt(F \cap(X_0^\perp + \vt)) \coloneqq Y_2(F) \cap (X_0^\perp + \vt).$$ Note that $U(\T_\vt, Y_\vt) = U(\F_2,Y_2) \cap (X_0^\perp + \vt)$, and that by Fubini's theorem, we have
    \begin{align}\label{Eq:Fubini_plany}
        \int_{U(\F_2,Y_2)} \mu(\F_2,Y_2)(x) \D x
        = \int_{X_0 \cap B_1} \int_{U(\T_\vt, Y_\vt)} \mu(\T_\vt, Y_\vt)(y) \D y \D \vt,
    \end{align}
    where the inner integral on the RHS is with respect to the standard $3$-dimensional Lebesgue measure on $X_0^\perp + \vt$. 
    
    Since the LHS of \eqref{Eq:Fubini_plany} is $\geq \delta^{50}$, we may safely discard those slices on the RHS with $\int_{U(\T_\vt, Y_\vt)} \mu(\T_\vt, Y_\vt)(y) \D y \leq \delta^{100}$, while still keeping the two sides of \eqref{Eq:Fubini_plany} comparable. On the other hand, we clearly have $\int_{U(\T_\vt, Y_\vt)} \mu(\T_\vt, Y_\vt)(y) \D y \leq \delta^{-100}$ for all $\vt \in X_0 \cap B_1$, and so the values of the inner integrals are restricted to $[\delta^{100}, \delta^{-100}]$. This allows us to apply dyadic pigeonholing with respect to this value to find a subset $A \subset X_0 \cap B_1$ such that $\int_{U(\T_\vt, Y_\vt)} \mu(\T_\vt, Y_\vt)(y) \D y$ is comparable to some constant $M$ for all $\vt \in A$, and
    \begin{align*}
        \int_{A} \int_{U(\T_\vt, Y_\vt)} \mu(\T_\vt, Y_\vt)(y) \D y \D \vt \gtrapprox \int_{U(\F_2,Y_2)} \mu(\F_2,Y_2)(x) \D x.
    \end{align*}
    For this inequality, observe that the LHS is $\lesssim |A|M \lesssim M$, while the RHS is $$\sum_{\F_2} |Y_2(F)| \gtrapprox \sum_{\F'} |Y'(F)| \gtrapprox \sum_{\F} |Y(F)| = \lambda(\F,Y) |F| |\F| \geq \delta^{\eta} \delta^{2} \delta^{-2+\eta} = \delta^{2\eta}.$$ 
    And so we can conclude that $\int_{U(\T_\vt, Y_\vt)} \mu(\T_\vt, Y_\vt)(y) \D y \sim M \gtrapprox \delta^{2\eta}$ for all $\vt \in A$. Another key observation is that since $(\F_2, Y_2)$ has constant multiplicity $\mu$, any restricted shading system $(\T_\vt,Y_\vt)$ also has constant multiplicity $\mu$. By \eqref{Eq:slice_tube_mult}, we have $|U(\T_\vt,Y_\vt)| \sim (1/\mu) \int_{U(\T_\vt, Y_\vt)} \mu(\T_\vt, Y_\vt)(y) \D y \gtrapprox \delta^{\sigma + 2\eta_1 + 2\eta}$, and so  \begin{align}\label{Eq:slice_tube_union}
        |U(\T_\vt,Y_\vt)| \geq \delta^{\sigma + 4\eta_1}
    \end{align}
    for all $\vt \in A$. 

    Now we are in the right place to apply the arguments in the proof of \cite[Lemma~4.3]{WZ26sticky} to each $(\T_\vt,Y_\vt)$ in the slice $X_0^\perp + \vt$ for $\vt \in A$. Although the statement of \cite[Lemma~4.3]{WZ26sticky} involves concepts like ``$\eta$-extremal'', one can see from the proof that we actually only need the set of tubes to satisfy several resulting properties, which are much more flexible and easy to check in our setting:
    \begin{itemize}
        \item $(\T_\vt,Y_\vt)$ has constant multiplicity \eqref{Eq:slice_tube_mult}, and has union \eqref{Eq:slice_tube_union}.
        
        \item The robust transversality assumption for $(\T_\vt,Y_\vt)$ is inherited directly from that for $(\F_2,Y_2)$ by recalling \eqref{Eq:robust_transverse'} at the beginning of the proof of Lemma~\ref{Lem:plane_or_pivot}. More precisely, let $\eps = \sigma/5$, and we want to show that for each $y \in U(\T_\vt,Y_\vt)$ and each line $\ell \in \Gr(X_0^\perp,1)$, we have   \begin{align}\label{Eq:rob_trans_tube}
            |\{T \in (\T_\vt)_{Y_\vt}(y): \angle( \ell, \dir(T)) \leq \delta^{5\eps/16} \}| \leq \delta^{-\sigma + 3\eta_1}.
        \end{align}
        
        Suppose this is not true, then there exists some $y_0 \in U(\T_\vt,Y_\vt)$ and $\ell_0 \in \Gr(X_0^\perp,1)$ such that 
        $$|\{T \in (\T_\vt)_{Y_\vt}(y_0): \angle(\ell_0, \dir(T)) \leq \delta^{5\eps/16} \}| > \delta^{-\sigma + 3\eta_1}.$$ Let 
        $$\F_0 \coloneqq \{F \in (\F_2)_{Y_2}(y_0): \angle(\ell_0, \dir(F) \cap X_0^\perp) \leq \delta^{5\eps/16} \},$$ then $|\F_0| > \delta^{-\sigma + 3\eta_1}$. Let $V_0\coloneqq \XX(y_0) + \ell_0$. Fix any $F \in \F_0$. We want to apply Lemma~\ref{Lem:geo_fact_1} with 
        \begin{align*}
            V_1 = \XX(y_0), & \quad  V_1' = \pi_{\dir(F)}(\XX(y_0)), \\
            V_2 = \ell_0, 
            & \quad V_2' = \dir(F) \cap X_0^\perp.
        \end{align*} Let us check the assumptions. Since $\XX(y_0) \cap B_1 \subset N_{\delta^{\sigma/2}}(\dir(F))$, by Lemma~\ref{Lem:faithful_proj} with $\rho = \delta^{\sigma/2}$, we get $\dim V_1' = n-3$ and $\angle(V_1, V_1') \leq \delta^{\sigma/2} \leq \delta^{5\eps/16}$. On the other hand, the definition of $\F_0$ immediately yields $\angle(V_2, V_2') \leq \delta^{5\eps/16}$. Finally, recall that $\angle(\XX(y_0), X_0) \leq 1/100$, and so $V_1 \cap B_1 \subset N_{1/100}(X_0) = N_{1/100}((X_0^\perp)^\perp) \subset N_{1/100}(V_2^\perp)$. Therefore, we can apply Lemma~\ref{Lem:geo_fact_1} with $\rho = \delta^{5\eps/16}$ to conclude that $V_1 + V_2 = V_0$ and $V_1' + V_2' \subset \dir(F)$ are both direct sums, which implies $V_0 = V_1 \oplus V_2 \in \Gr(n,n-2)$ and $\dir(F) = V_1' \oplus V_2'$ by counting dimensions. Moreover, $\angle(V_0, \dir(F)) = \angle(V_1 \oplus V_2, V_1' \oplus V_2') \leq 3\delta^{5\eps/16} \leq \delta^{\eps/4}$. This holds for all $F\in \F_0$, and so we have \begin{align*}
            |\{F\in ({\mathbb F}_2)_{Y_2}(y_0): \angle(V_0,\dir(F)) \leq \delta^{\varepsilon/4}\}| \geq |\F_0| > \delta^{-\sigma +3\eta_1},
        \end{align*}
        which contradicts \eqref{Eq:robust_transverse'} ($(\F_2, Y_2)$ is a refinement of the flakes and shadings therein by trace back the chain of refinements). And so \eqref{Eq:rob_trans_tube} must hold true.      

        \item The version of the multilinear Kakeya inequality due to Carbery and Valdimarsson (\cite[Theorem~1]{CV2013}, see also \cite[Theorem~4.2]{WZ26sticky}) does not require any additional structural properties of tubes, and so applies equally well to each $\T_t$. (We do not even need the tubes $\T_t$ to be essentially distinct!)
    \end{itemize}
    Therefore, we will omit the remaining details, as they are essentially the same as those for \cite[Lemma~4.3]{WZ26sticky}. As a result of all these processes, we get a $\gtrapprox 1$ refinement $(\T_\vt', Y_\vt')$ of $(\T_\vt, Y_\vt)$ and a function $\YY_\vt: U(\T_\vt', Y_\vt') \rightarrow \Gr(X_0^\perp,2)$ for each $\vt \in A$, such that for each $y \in U(\T_\vt', Y_\vt')$, we have    \begin{align}\label{Eq:slice_tube_plany}
        \dir(T) \cap B_1 \leq N_{\delta^{\sigma/2}}(\YY_\vt(y)) \quad \text{for all } T \in (\T_\vt')_{Y_\vt'}(y).
    \end{align}
    We point out that the numerology regarding small parameters like $\eps$ and $\eta_1$ only depends on the properties listed above, and are therefore uniform over all slices. This yields some admissible requirements on our final choice of $\eta_1$. 
    
    It is time to put the refinements in all slices together and see how planiness for $\delta$-tubes inside each slice tells us about the planiness for the original set of $\delta$-flakes in $\R^n$. Define 
    $$Y''(F) \coloneqq Y_2(F) \cap (\bigsqcup_{\vt \in A} Y_\vt'(F \cap (X_0^\perp + \vt))) \footnote{There could be measurability issues, see Remark~\ref{Rem:measurability} for some technical explanations.},$$ 
    where we adopt the convention that $Y_\vt'(T) = \varnothing$ if $T \in \T_\vt\setminus\T_\vt'$. Let $\F'' = \F_2$. Note that 
    \begin{align*}
        \sum_{\F''} |Y''(F)| = \int_A \sum_{\T_\vt'} |Y_\vt'(T)| \D \vt \gtrapprox \int_A \sum_{\T_\vt} |Y_\vt(T)| \D \vt \gtrapprox \sum_{\F_2} |Y_2(F)|,
    \end{align*}
    and so $(\F'', Y'')$ is a $\gtrapprox 1$ refinement of $(\F_2,Y_2)$, and thus of $(\F,Y)$. Besides, for each $x \in U(\F'', Y'')$, there is a unique $\vt$, which is $\pi_{X_0}(x)$, such that $x \in X_0^\perp + \vt$. So we can equivalently rewrite \eqref{Eq:slice_tube_plany} as follows: For each $x \in U(\F'',Y'')$, we have     \begin{align}\label{Eq:slice_tube_plany'}
        (\dir(F) \cap X_0^\perp) \cap B_1 \subset N_{\delta^{\sigma/2}}(\YY_{\pi_{X_0}(x)}(x)) \quad \text{for all } F \in (\F'')_{Y''}(x).
    \end{align}
    
    It remains to verify that $(\F'', Y'')$ satisfies the desired planiness property. Fix any $x \in U(\F'',Y'')$ and any $F \in (\F'')_{Y''}(x)$. Let $\HH^\perp(x) \coloneqq \XX(x)+\YY_{\pi_{X_0}(x)}(x)$. We want to apply Lemma~\ref{Lem:geo_fact_2} with 
    \begin{align*}
        V_1 = \XX(x), & \quad V_1' = \pi_{\dir(F)}(\XX(x)), \\
        V_2 = \YY_{\pi_{X_0}(x)}(x), & \quad V_2'  = \dir(F) \cap X_0^\perp.
    \end{align*} 
    Let us check the assumptions. Since $\XX(x) \cap B_1 \subset N_{\delta^{\sigma/2}}(\dir(F))$, by Lemma~\ref{Lem:faithful_proj} with $\rho = \delta^{\sigma/2}$, we get $\dim V_1' = n-3$ and $\angle(V_1, V_1') \leq \delta^{\sigma/2}$, i.e., $V_1' \cap B_1 \subset N_{\delta^{\sigma/2}}(V_1)$. On the other hand, \eqref{Eq:slice_tube_plany'} immediately yields $V_2' \cap B_1 \subset N_{\delta^{\sigma/2}}(V_2)$. Finally, recall that $\angle(\XX(x), X_0) \leq 1/100$, and so $V_1 \cap B_1 \subset N_{1/100}(X_0) = N_{1/100}((X_0^\perp)^\perp) \subset N_{1/100}(V_2^\perp)$. Therefore, we can apply Lemma~\ref{Lem:geo_fact_2} with $\rho = \delta^{\sigma/2}$ to conclude that $V_1 + V_2 = \HH^\perp(x)$ and $V_1' + V_2' \subset \dir(F)$ are both direct sums, which implies $\HH^\perp(x) = V_1 \oplus V_2 \in \Gr(n,n-1)$ and $\dir(F) = V_1' \oplus V_2'$ by counting dimensions. Moreover, we have 
    $$\dir(F) \cap B_1 = (V_1' \oplus V_2') \cap B_1 \subset N_{3\delta^{\sigma/2}}(V_1 \oplus V_2) \subset N_{\delta^{\sigma/5}}(\HH^\perp(x)).$$ Define $\HH(x) \coloneqq (\HH^\perp(x))^\perp \in \Gr(n,1)$. Then $\HH(x)^\perp = \HH^\perp(x)$, and we get $\dir(F) \cap B_1 \subset N_{\delta^{\sigma/5}}(\HH(x)^\perp)$. This finishes the proof of Proposition~\ref{Prop:planiness}. 
\end{proof} 

\begin{rem}\label{Rem:measurability}
In general we want the shadings $Y''(F)$ in the proof to be measurable sets, but an uncountable union of measurable sets may no longer be measurable. This issue is purely technical, as we can first apply Lemma~\ref{Lem:fact_mult} to the $\delta$-flakes with $\rho=\delta$ to replace the shadings with their $\delta$-neighborhoods, and then do all the slicing arguments that we have presented. In fact, we can do the replacement at the very beginning of Section~\ref{Sec:self_similar} when we define ``$\sigma$'' and ``$\eps$-extremal'', i.e., limit the class of shadings to those that are unions of $\delta$-balls. In this way the distribution of the shadings is constant at scale $\delta$, which allows us to choose the same $(\T_\vt', Y_\vt')$ and $\YY_\vt$ whenever $\vt$ is in a cube of side length $\delta$. This local constancy guarantees the measurability of $Y''(F)$. We omit such details that are tedious and not insightful. Another way of fixing the issue is to use slices with thickness $\delta$ and apply a ``thickened'' version of the proof of \cite[Lemma~4.3]{WZ26sticky} inside each slice. We will see arguments in this spirit in Section~\ref{Sec:slicing_set_up}.
\end{rem}

Proposition~\ref{Prop:planiness} says that extremal sets of $\delta$-flakes must satisfy a weak type of planiness, which involves a $\delta^{\sigma/5}$ compression that is $\ll 1$ but $\gg \delta$. However, by moving to a coarser scale thanks to the stickiness assumption, we can ensure that a suitable thickening of the flakes agrees with their plane map $\HH$ in a strong sense. This is the content of the next lemma. {\color{black} In what follows, for any $\rho \in [\delta,1]$ a \textit{$\rho$-cube} is a set of the form $[0,\rho)^n + x$, where $x \in (\rho \Z)^n$.} 
\begin{lem}\label{Lem:strong_planiness}
    For all $\eps>0$, there exist $\eta, \delta_0 >0$ so that the following holds for all $\delta \in (0, \delta_0]$. Let $(\F,Y)$ be an $\eta$-extremal set of $\delta$-flakes, and let $\rho \in [\delta^{1-\eps}, \delta^\eps]$. Suppose there exists a function $\HH: U(\F,Y) \rightarrow \Gr(n,1)$ so that for each $x \in U(\F,Y)$, we have 
    \begin{align*}
        \dir(F) \cap B_1 \subset N_\rho(\HH(x)^\perp) \quad \text{for all } F \in \F_Y(x).
    \end{align*}
    Then there exists a $\gtrapprox 1$ refinement $(\F', Y')$ of $(\F,Y)$ and an $\eps$-extremal set of $\rho$-flakes $(\F_\rho', Y_\rho')$ with plane map $\HH'$ so that $U(\F', Y') \subset U(\F_\rho', Y_\rho')$. Furthermore, on each $(\rho/10n)$-cube $Q$ with $Q \cap U(\F_\rho', Y_\rho') \neq \varnothing$, $\HH'$ is constant with value $\HH'(Q)$, and $\HH$ and $\HH'$ are consistent in the sense that there exists $x \in Q \cap U(\F',Y')$ with $\HH'(Q) = \HH(x)$. 
\end{lem}
\begin{proof}
    Let $\eps_1>0$ be a small constant to be chosen below. We will select $\eta$ very small compared to $\eps_1$, and $\eps_1$ small compared to $\eps$. Apply Proposition~\ref{Prop:self-similar} to $(\F,Y)$ with $\eps_1$ in place of $\eps$ and with $\rho$ as stated in the lemma. We obtain a $\gtrapprox 1$ refinement $(\F'',Y'')$ of $(\F,Y)$ and an $\eps_1$-extremal set of $\rho$-flakes $(\F_\rho'',Y_\rho'')$ that covers $(\F'',Y'')$, such that all items in Proposition~\ref{Prop:self-similar} are satisfied with $\eps_1$ in place of $\eps$. 
    
    Since $|U(\F'',Y'')| \leq |U(\F,Y)| \leq \delta^{\sigma-\eta}$ and $\sum_{\F''}|Y''(F)| \gtrapprox \sum_{\F}|Y(F)| \geq \delta^{2\eta}$, after discarding those $x \in U(\F'',Y'')$ with $\mu(\F'',Y'')(x) \leq \delta^{-\sigma + 4\eta}$, we get a $\gtrsim 1$ refinement $(\F',Y')$ of $(\F'',Y'')$ such that $\mu(\F',Y')(x) \geq \delta^{-\sigma + 4\eta}$ for all $x \in U(\F',Y')$. In particular, for each $(\rho/10n)$-cube $Q$ with $Q \cap U(\F',Y') \neq \varnothing$, there exists a point $x_Q \in Q$ with $\mu(\F',Y')(x_Q) \geq \delta^{-\sigma + 4\eta}$. We say the flake $F_\rho \in \F_\rho''$ is \textit{associated} to the $(\rho/10n)$-cube $Q$ if $x_Q \in U(\F'[F_\rho], Y')$; this necessarily implies $Q \subset Y_\rho''(F_\rho)$. By Item~(iii) in Proposition~\ref{Prop:self-similar}, each $(\rho/10n)$-cube $Q$ with $Q \cap U(\F',Y') \neq \varnothing$ satisfies 
    \begin{align*}
        |\{F_\rho \in \F_\rho'': F_\rho \text{ is associated to } Q\}| \geq \frac{\delta^{-\sigma + 4\eta}}{(\rho/\delta)^{\sigma + \eps_1}} = \delta^{4\eta + \eps_1} \rho^{-\sigma}.
    \end{align*}
    On the other hand, by choosing $\eta$ small enough depending on $\eps_1$, we have $|U(\F',Y')| \geq \delta^{\sigma + \eps_1}$ and $|U(\F_\rho'', Y''_{\F_\rho})| \geq \rho^{\sigma + \eps_1}$ for small $\delta$ (recall Remark~\ref{Rem:extreme_mult}). Thus by Item~(v) in Proposition~\ref{Prop:self-similar}, for each $(\rho/10n)$-cube $Q$ with $Q \cap U(\F',Y') \neq \varnothing$, we have 
    \begin{align*}
        \frac{|U(\F', Y') \cap Q|}{|Q|} \lesssim \frac{|U(\F'', Y'') \cap B(x_Q, \rho)|}{|B(x_Q, \rho)|} \sim \frac{|U(\F'', Y'')|}{|U(\F_\rho'', Y''_{\F_\rho})|} \leq \frac{\delta^{\sigma-\eta}}{\rho^{\sigma + \eps_1}}.
    \end{align*}    
    This further implies that
    \begin{align*}
        |\{Q: Q \cap U(\F',Y') \neq \varnothing\}| \gtrsim \frac{\delta^{\sigma + \eps_1}}{|Q| \delta^{\sigma-\eta} / \rho^{\sigma + \eps_1}} = \delta^{\eta+\eps_1} \rho^{\sigma+\eps_1} |Q|^{-1}.
    \end{align*}
    Combining this lower bound for the number of $Q$'s with the previous lower bound for the number of $\rho$-flakes associated to a single $Q$, we can conclude that if we define $Y_\rho'(F_\rho)$ to be the union of $(\rho/10n)$-cubes $Q$ associated to $F_\rho$ and define $\F_\rho' = \F_\rho''$, then 
    \begin{align*}
        \sum_{\F_\rho'} |Y_\rho'(F_\rho)| 
        \gtrsim 
        \delta^{4\eta + \eps_1} \rho^{-\sigma} \cdot |Q| \cdot \delta^{\eta+\eps_1} \rho^{\sigma+\eps_1} |Q|^{-1} \geq \delta^{5\eta + 3\eps_1}.
    \end{align*}
    In particular, we have $\lambda(\F_\rho',Y_\rho') \gtrsim \delta^{5\eta + 3\eps_1}\rho^{\eps_1} \geq \rho^{(5\eta+3\eps_1)/\eps + \eps_1}$. 
    Since $\eps_1$ is sufficiently small depending on $\eps$, we have $\lambda(\F_\rho',Y_\rho') \geq \rho^\eps$, which means $(\F_\rho',Y_\rho')$ is $\eps$-extremal. By our construction, we know that $U(\F', Y') \subset U(\F_\rho', Y_\rho')$, $U(\F_\rho', Y_\rho')$ is a disjoint union of $(\rho/10n)$-cubes, and each $(\rho/10n)$-cube $Q$ in $U(\F_\rho', Y_\rho')$ must intersect $U(\F', Y')$ at $x_Q$.

    Next, we shall define the function $\HH': U(\F_\rho', Y_\rho') \rightarrow \Gr(n,1)$ as follows. For each $(\rho/10n)$-cube $Q \subset U(\F_\rho', Y_\rho')$, we define $\HH'(x) \equiv \HH(x_Q)$ for all $x \in Q$. Clearly $\HH'$ is constant on $(\rho/10n)$-cubes, and satisfies the consistency condition with $\HH$ as stated in Lemma~\ref{Lem:strong_planiness}. It remains to show that $\HH'$ is a plane map for $(\F_\rho', Y_\rho')$. Let $x \in U(\F_\rho', Y_\rho')$ be contained in a $(\rho/10n)$-cube $Q$, and let $F_\rho \in (\F_\rho')_{Y_\rho'}(x)$. Then there exists a $\delta$-flake $F \in \F'[F_\rho] \cap \F_Y(x_Q)$. By our assumption, we know that $\dir(F) \cap B_1 \subset N_{\rho}(\HH(x_Q)^\perp)$, and so $\dir(F_\rho) \cap B_1 \subset N_{\rho}(\dir(F) \cap B_1) \subset N_{2\rho}(\HH(x_Q)^\perp)$, which is exactly what we want.
\end{proof}

\subsection{Lipschitz regularity of the plane map}\label{subsec:Lips_reg_plane_map}
So far we have seen the existence of a plane map, but we still need to know about some of its properties. We will show that if an extremal set of flakes has a plane map, then after restricting to a
refinement, this plane map must be Lipschitz with controlled Lipschitz norm.

\begin{lem}\label{Lem:plane_map_mild}
    For all $\eps>0$, there exist $\eta, \delta_0 >0$ so that the following holds for all $\delta \in (0, \delta_0]$. Let $(\F,Y)$ be an $\eta$-extremal set of $\delta$-flakes with plane map $\HH$, and let $\rho \in [\delta^{1-\eps}, \delta^\eps]$. Then there exists an $\eps$-extremal refinement $(\F',Y')$, so that if $x,y \in U(\F',Y')$ are contained in a common $\rho$-cube, then $\angle(\HH(x), \HH(y)) \leq \rho$. 
\end{lem}
\begin{proof}
    Let $\eps_1, \eps_2 > 0$ be small parameters to be chosen below. We will select $\eta$ very small compared to $\eps_1$, $\eps_1$ very small compared to $\eps_2$, and $\eps_2$ very small compared to $\eps$.

    Let $(\F_1,Y_1)$ be the $\gtrapprox 1$ refinement obtained by applying Lemma~\ref{Lem:robust_transverse} to $(\F,Y)$ with $\eps_1$ in place of $\eps$. Next, apply Proposition~\ref{Prop:self-similar} to $(\F_1,Y_1)$ with $\eps_2$ in place of $\eps$, and $\rho$ as specified in the lemma. We obtain a $\gtrapprox 1$ refinement $(\F_2,Y_2)$ of $(\F_1,Y_1)$ and an $\eps_2$-extremal set of $\rho$-flakes $(\F_\rho,Y_\rho)$. 
    
    Since $|U(\F_2,Y_2)| \leq |U(\F,Y)| \leq \delta^{\sigma-\eta}$ and $\sum_{\F_2}|Y_2(F)| \gtrapprox \sum_{\F}|Y(F)| \geq \delta^{2\eta}$, after discarding those $x \in U(\F_2,Y_2)$ with $\mu(\F_2,Y_2)(x) \leq \delta^{-\sigma + 4\eta}$, we get a $\gtrsim 1$ refinement $(\F_3,Y_3)$ of $(\F_2,Y_2)$ such that $\mu(\F_3,Y_3)(x) \geq \delta^{-\sigma + 4\eta}$ for all $x \in U(\F_3,Y_3)$. If $\eta>0$ is selected sufficiently small (depending on $\eps_1$ and $\sigma$), then for each $x \in U(\F_3,Y_3)$, there are at least 
    $$\mu(\F_3,Y_3)(x)^2 - \mu(\F_3,Y_3)(x)\delta^{-\sigma + \eps_1\sigma/2} \geq \mu(\F_3,Y_3)(x)^2/2$$ 
    pairs $F, F' \in \F_3(x)$ with $\angle(F,F') \geq \delta^{\eps_1}$. Thus by Item~(iii) in Proposition~\ref{Prop:self-similar}, there are at least 
    $$(\mu(\F_3,Y_3)(x)^2/2)/((\rho/\delta)^{\sigma+\eps_2})^2 \geq \delta^{8\eta + 2\eps_2} \rho^{-2\sigma}$$ 
    pairs $F_\rho, F_\rho' \in \F_\rho$ with $x \in U(\F_3[F_\rho], Y_3)$, $x \in U(\F_3[F_\rho'], Y_3)$, and $\angle(F_\rho, F_\rho') \geq \delta^{\eps_1} - 2\rho$ (recall from Remark~\ref{Rem:major_metric} that the major angle induces a metric); we will choose $\eps_1$ sufficiently small so that $\delta^{\eps_1} - 2\rho \geq \delta^{\eps_1}/2$. 
    
    For each $\rho$-cube $Q$, let $\F_\rho(Q) \coloneqq \{F_\rho \in \F_\rho: Y_\rho(F_\rho) \cap Q \neq \varnothing\}$. Note that for each $F_\rho \in \F_\rho(Q)$, $Y_\rho(F_\rho)$ must contain a $(\rho/10n)$-cube living inside $10Q$. Since there are at most $O(1)$ such choices of $(\rho/10n)$-cubes, each of which is contained in at most $\rho^{-\sigma-\eps_2}$ $Y_\rho(F_\rho)$ by Item~(ii) in Proposition~\ref{Prop:self-similar}, we can conclude that $|\F_\rho(Q)| \lesssim \rho^{-\sigma-\eps_2}$. By the information we get so far, for each $\rho$-cube $Q$ intersecting $U(\F_3, Y_3)$, we have
    \begin{align*}
        &|\{ (x, F_\rho, F_\rho') \in (Q \cap U(\F_3,Y_3)) \times \F_\rho \times \F_\rho:\\ 
        &\quad
        x \in U(\F_3[F_\rho], Y_3) \cap U(\F_3[F_\rho'], Y_3), 
        \angle(F_\rho,F_\rho')\geq \delta^{\eps_1}/2 \}|\\
        & \geq 
        |Q \cap U(\F_3,Y_3)| \delta^{8\eta+2\eps_2} \rho^{-2\sigma}\\
        & \gtrsim 
        \delta^{8\eta+4\eps_2} |Q \cap U(\F_3,Y_3)| |\F_\rho(Q)|^2,
    \end{align*}
    where we use $|\cdot|$ on the LHS to denote the product of the Lebesgue measure in $\R^n$ and the counting measure. Note that we may equivalently replace $(Q \cap U(\F_3,Y_3)) \times \F_\rho \times \F_\rho$ on the LHS with $(Q \cap U(\F_3,Y_3)) \times \F_\rho(Q) \times \F_\rho(Q)$. Thus by pigeonholing, we can select two $\rho$-flakes from $\F_\rho(Q)$, which we denote by $F_{\rho,Q}$ and $F_{\rho,Q}'$, with $\angle(F_{\rho,Q},F_{\rho,Q}') \geq \delta^{\eps_1}/2$, so that 
    \begin{align}\label{Eq:bilinear_perfect_overlap}
        | Q \cap U(\F_3[F_\rho], Y_3) \cap U(\F_3[F_\rho'], Y_3) | 
        \gtrsim   \delta^{8\eta+4\eps_2} |Q \cap U(\F_3,Y_3)|.
    \end{align}

    Define $F_Q$ to be the set on the LHS of \eqref{Eq:bilinear_perfect_overlap}. Let $\F_4 = \F_3$, and for each $F \in \F_4$, define 
    \[
    Y_4(F) \coloneqq Y_3(F) \cap \bigcup_{Q: Q \cap U(\F_3,Y_3) \neq \varnothing} F_Q.
    \]
    Then $$|U(\F_4, Y_4)| = \sum_{Q} |F_Q| \gtrsim \sum_{Q} (\delta^{8\eta+4\eps_2} |Q \cap U(\F_3,Y_3)|) = \delta^{8\eta+4\eps_2} |U(\F_3,Y_3)|.$$
    Also, recall that $\mu(\F_3,Y_3)(x) \geq \delta^{-\sigma + 4\eta}$ for each $x \in U(\F_3,Y_3)$, and observe that this property is inherited by $(\F_4,Y_4)$, i.e., $\mu(\F_4,Y_4)(x) \geq \delta^{-\sigma + 4\eta}$ for each $x \in U(\F_4,Y_4)$. Therefore, we have $$\sum_{\F_4}|Y_4(F)| \gtrsim \delta^{8\eta+4\eps_2} |U(\F_3,Y_3)| \cdot \delta^{-\sigma + 4\eta} \geq \delta^{12\eta + 5\eps_2},$$ where we used the fact that $|U(\F_3,Y_3)| \geq \delta^{\sigma + \eps_2}$ if $\eta$ is chosen sufficiently small (recall Remark~\ref{Rem:extreme_mult}). 

    Let $Q$ be a $\rho$-cube intersecting $U(\F_4,Y_4)$, and let $x \in U(\F_4,Y_4) \cap Q$. Then there exist $F_x \in \F_4[F_{\rho,Q}]$ and $F_x' \in \F_4[F_{\rho,Q}']$ so that $x \in Y_4(F_x) \cap Y_4(F_x')$, and hence 
    \[
    \dir(F_{\rho,Q}) \cap B_1 \subset N_\rho(\dir(F_x) \cap B_1) \subset N_{2\delta+\rho}(\HH(x)^\perp)
    \]
    and 
    \[
    \dir(F_{\rho,Q}') \cap B_1 \subset N_\rho(\dir(F_x') \cap B_1) \subset N_{2\delta+\rho}(\HH(x)^\perp).
    \]
    Since $\angle(F_{\rho,Q},F_{\rho,Q}') \geq \delta^{\eps_1}/2$, we can find $v_0 \in \dir(F_{\rho,Q})~\cap~\mathbb{S}^{n-1}$ such that $\dist(v_0, \dir(F_{\rho,Q}')) \geq \delta^{\eps_1}/2$. Let $\{v_j\}_{j=1}^{n-2}$ be an orthonormal basis of $\dir(F_{\rho,Q}')$, and let $H_{\rho,Q}$ be the hyperplane spanned by $v_0$ and $\dir(F_{\rho,Q}')$. Then for any $v \in H_{\rho,Q} \cap B_1$, there is a unique decomposition $v = \sum_{j=0}^{n-2} a_j v_j$. Note that the inverse of the $(n-1) \times (n-1)$ matrix formed by $v_0, v_1, \ldots, v_{n-2}$ (coordinate representation inside $H_{\rho,Q}$) is  its adjugate matrix multiplied by $|\wedge_{j=0}^{n-2} v_j|^{-1} \lesssim \delta^{-\eps_1}$, and so has all entries $\lesssim \delta^{-\eps_1}$. This means $|a_j| \lesssim \delta^{-\eps_1}$ for each $0 \leq j \leq n-2$, and so we can deduce from $\{v_j\}_{j=0}^{n-2} \subset N_{2\delta+\rho}(\HH(x)^\perp)$ that $v \in N_{C\delta^{-\eps_1}\rho}(\HH(x)^\perp)$ for some $C \lesssim 1$. In other words, we have $H_{\rho,Q} \cap B_1 \subset N_{C\delta^{-\eps_1}\rho}(\HH(x)^\perp)$, which by Item~(iii) of Theorem~\ref{Thm:Principal_angle} and Remark~\ref{Rem:major_equiv} means $\angle(\HH(x), H_{\rho, Q}^\perp) = \angle(\HH(x)^\perp, H_{\rho, Q}) \leq C\delta^{-\eps_1}\rho$. This happens for all $x \in U(\F_4,Y_4) \cap Q$. 
    
    By a further obvious pigeonholing step with respect to $\HH(x)$ inside each $Q$, we can find a refinement $(\F', Y')$ of $(\F_4,Y_4)$ such that (i) for each $\rho$-cube intersecting $(\F', Y')$, we have $\angle(\HH(x), \HH(y)) \leq \rho$ for all $x,y\in Q \cap U(\F', Y')$, and (ii) 
    $$\sum_{\F'} |Y'(F)| \gtrsim \delta^{(n-1)\eps_1} \sum_{\F_4} |Y_4(F)| \gtrsim \delta^{12 \eta + (n-1)\eps_1 + 5\eps_2}.$$ 
    If we choose $\eta, \eps_1, \eps_2$ sufficiently small, then $(\F', Y')$ will be $\eps$-extremal. 
\end{proof}

Applying Lemma~\ref{Lem:plane_map_mild} and then deleting a constant fraction of the $\rho$-cubes to
obtain a $\sqrt{n}\rho$-separated set of cubes, we obtain the following result:

\begin{cor}\label{Cor:plane_map_mild}
For all $\eps>0$, there exist $\eta, \delta_0 >0$ so that the following holds for all $\delta \in (0, \delta_0]$. Let $(\F,Y)$ be an $\eta$-extremal set of $\delta$-flakes with plane map $\HH$, and let $\rho \in [\delta^{1-\eps}, \delta^\eps]$. Then there exists an $\eps$-extremal refinement $(\F',Y')$, so that if $x,y \in U(\F',Y')$ satisfy $\norm{x-y} \leq \rho$, then $\angle(\HH(x), \HH(y)) \leq \rho$.
\end{cor}

Corollary~\ref{Cor:plane_map_mild} says that after restricting to a refinement, the plane map $\HH$ will
be Lipschitz at a specific scale $\rho$. Lemma~\ref{Lem:eps_Lipschitz} iterates this result at many scales.

\begin{lem}\label{Lem:eps_Lipschitz}
For all $\eps>0$, there exist $\eta, \delta_0 >0$ so that the following holds for all $\delta \in (0, \delta_0]$. Let $(\F,Y)$ be an $\eta$-extremal set of $\delta$-flakes with plane map $\HH$ that is constant on $(\delta/10n)$-cubes. Then there exists an $\eps$-extremal refinement $(\F',Y')$, so that $\HH$ is $\delta^{-\eps}$-Lipschitz on $U(\F',Y')$.
\end{lem}
\begin{proof}
After deleting a constant fraction of the $(\delta/10n)$-cubes to
obtain a $\sqrt{n}\delta$-separated set of cubes, we can suppose that if $x, y \in U(\F,Y)$ satisfy $\norm{x-y} \leq \delta$, then $x$ and $y$ are contained in a common $(\delta/10n)$-cube, and hence $\HH(x) = \HH(y)$. Let $N = \lceil 2/\eps \rceil$. If $\eta, \delta_0$ are selected sufficiently small, then we can apply Corollary~\ref{Cor:plane_map_mild} iteratively to $(\F,Y)$ for $\rho = \rho_j = \delta^{j/N}$, $j = 1, \dots, N-1$, and the resulting refinement $(\F', Y')$ will be $\eps$-extremal. We claim that $\HH$ is $\delta^{-1/N}$-Lipschitz on $(\F',Y')$. To see this, let $x,y \in U(\F',Y')$. If $\norm{x - y} \leq \delta$, then $\HH(x) = \HH(y)$ and we are done. If $\norm{x-y} > \delta^{1/N}$, then $\angle(\HH(x), \HH(y)) \leq 1$ and again we are done. Otherwise, there is an index $1\leq k \leq N-1$ so that $\delta^{(k+1)/N} < \norm{x-y} \leq \delta^{k/N}$, and 
\begin{align*}
\angle(\HH(x), \HH(y)) \leq \delta^{k/N} \leq \delta^{-1/N} \norm{x-y}.
\end{align*}
This verifies the claim, and thus the conclusion of the lemma.
\end{proof} 

Combining the results we obtained so far with Lemma~\ref{Lem:extremizer_exist}, we have the following conclusion:
\begin{lem}\label{Lem:exist_eps_Lips}
    For all $\eps, \delta_0 > 0$, there exist $\delta \in (0, \delta_0]$ and an $\eps$-extremal set of $\delta$-flakes that has a $\delta^{-\eps}$-Lipschitz plane map. 
\end{lem}
\begin{proof}
By Lemma~\ref{Lem:extremizer_exist}, we can take an $\eta$-extremal set of $\delta$-flakes for a very small $\eta$. Then Proposition~\ref{Prop:planiness} gives the weak planiness, which is further upgraded to a strong one by Lemma~\ref{Lem:strong_planiness} with $\rho = \delta^{\sigma/5}$. Finally, Lemma~\ref{Lem:eps_Lipschitz} provides the desired Lipschitz property of the plane map. 
\end{proof}

\subsection{Local grains}\label{subsec:local_grains}
With Lemma~\ref{Lem:exist_eps_Lips} in hand, we are finally in the position to study the graininess structure. We will show that if an extremal set of flakes has a Lipschitz plane map, then after restricting to a refinement, the flakes arrange
themselves into parallel slab segments, which are called \textit{local grains}. In what follows, $\HH$ will be a
plane map for a set of flakes. 

Our goal is to show that the set of local grains inside a cube arrange themselves into an AD-regular set, in the sense of Definition~\ref{Def:AD_set}. As a first step, Lemma~\ref{Lem:grains_number_given_scale} controls
how many local grains can occur inside a cube at a given scale. {\color{black} The idea of the proof is similar to that of \cite[Lemma~4.9]{WZ26sticky}, except that we use a projection method to handle local hairbrushes, instead of directly executing a C\'ordoba-type $L^2$ argument. The reason is that the geometry of flakes is much more complicated than tubes and calls for a more robust argument. Also, we write things in a more detailed and rigorous way.} 

\begin{lem}\label{Lem:}\label{Lem:grains_number_given_scale}
For all $\eps>0$, there exist $\eta, \delta_0 >0$ so that the following holds for all $\delta \in (0, \delta_0]$. Let $(\F,Y)$ be an $\eta$-extremal set of $\delta$-flakes with $\delta^{-\eta}$-Lipschitz plane map $\HH$. Let $\rho \in [\delta^{1-\eps}, \delta^\eps]$ and let $\tau \in [\rho^{1-\eps}, \rho^{1/2}]$. For each $\tau$-cube $\Tilde{Q}$ intersecting $U(\F,Y)$, define $\HH(\Tilde{Q}) \coloneqq \HH(y_{\Tilde{Q}})$, where $y_{\Tilde{Q}}$ is any fixed point in $\Tilde{Q} \cap U(\F,Y)$. 

Then there exists an $\gtrapprox 1$ refinement $(\F',Y')$ which is $2\eta$-extremal, so that for each $\tau$-cube $\Tilde{Q}$ intersecting $U(\F,Y)$, we have
\begin{align}\label{Eq:E_rho_local_grain}  \EE_\rho\big( \pi_{\HH(\Tilde{Q})}(U(\F',Y') \cap \Tilde{Q})\big) \leq \rho^{-\eps}(\tau/\rho)^{1-\sigma}.
\end{align}
\end{lem}
\begin{proof}
It suffices to only consider the case when $\tau/\rho$ is an integer. Let $\eps_1, \eps_2, \eps_3 > 0$ be small constants to be chosen below. We will select $\eta$ very small compared to $\eps_1$, $\eps_1$ very small compared to $\eps_2$, $\eps_2$ very small compared
to $\eps_3$, and $\eps_3$ very small compared to $\eps$.

Applying Proposition~\ref{Prop:self-similar} to $(\F, Y)$ with $\eps_1$ in place of $\eps$, we obtain a $\gtrapprox 1$ refinement $(\F_1,Y_1)$ which is covered by an $\eps_1$-extremal set of $\rho$-flakes $(\F_\rho, Y_\rho)$. For each $F_\rho \in \F_\rho$, let $Y_{\rho}'(F_\rho)$ be the union of all $(\rho/10n)$-cubes intersecting $U(\F_1[F_\rho], Y_1)$, which we denote by $\Q(F_\rho)$, then $Y_\rho'(F_\rho) \subset Y_\rho(F_\rho)$ and $|Y_\rho'(F_\rho)| \gtrsim |Y_\rho(F_\rho)|$ as $Y_\rho$ is the induced shading of $Y_1$. This means $(\F_\rho, Y_\rho')$ is a $\gtrsim 1$ refinement of $(\F_\rho, Y_\rho)$. Furthermore, by pigeonholing with respect to $\{Q \cap U(\F_1[F_\rho], Y_1): F_\rho \in \F_\rho, Q \in \Q(F_\rho) \}$, we get a subset $\Q'(F_\rho) \subset \Q(F_\rho)$ for each $F_\rho \subset \F_\rho$ so that $|Q \cap U(\F_1[F_\rho], Y_1)|$ are comparable for all $F_\rho \in \F_\rho$ and $Q \in \Q'(F_\rho)$, and $$\sum_{F_\rho \in \F_\rho}|U(\Q'(F_\rho)) \cap U(\F_1[F_\rho], Y_1)| \gtrapprox \sum_{F_\rho \in \F_\rho}|U(\F_1[F_\rho], Y_1)|.$$ In particular, by Item~(iii) of Proposition~\ref{Prop:self-similar}, this implies that $\sum_{F \in \F_1[F_\rho]} |Q \cap Y_1(F)|$ is comparable for all $F_\rho \in \F_\rho$ and $Q \in \Q'(F_\rho)$, and if we define $Y_2(F) \coloneqq Y_1(F) \cap U(\Q'(F_\rho))$ for $F \in \F_1[F_\rho]$ ($F_\rho \in \F_\rho$), then $(\F_1, Y_2)$ is a $\gtrapprox 1$ refinement of $(\F_1, Y_1)$ (thus of $(\F,Y)$). Also, if we define $Y_\rho''(F_\rho) \coloneqq U(\Q'(F_\rho))$ for $F_\rho \in \F_\rho$, then $(\F_\rho, Y_\rho'')$ covers $(\F_1, Y_2)$ and is a $2\eps_1$-extremal refinement of $(\F_\rho,Y_\rho)$ for the following reasons. Recall that $\F_\rho$ is $\rho^{-\eps_1}$-sticky and note that $|U(\F_\rho, Y_\rho'')| \leq |U(\F_\rho, Y_\rho)| \leq \rho^{\sigma - \eps_1}$. It remains to bound $\lambda(\F_\rho, Y_\rho'')$ from below. For each $F_\rho \in \F_\rho$, pick a $F(F_\rho) \in \F_1[F_\rho]$ that maximizes $|Y_2(F)|/|F|$ (among all $F \in \F_1[F_\rho]$). Then $$|Y_2(F(F_\rho))| \geq (1/|\F_1[F_\rho]|) \sum_{F \in \F_1[F_\rho]}|Y_2(F)|\geq \delta^{\eta}(\rho/\delta)^{-2} \sum_{F \in \F_1[F_\rho]}|Y_2(F)|,$$ which implies 
$$ \sum_{F_\rho \in \F_\rho} |Y_\rho''(F_\rho)| \gtrsim \sum_{F_\rho \in \F_\rho} |Y_2(F(F_\rho))|\frac{|F_\rho|}{|F(F_\rho)|} \geq \delta^{\eta} \sum_{F \in \F_1} |Y_2(F)| \gtrapprox \delta^{\eta}\sum_{F\in \F} |Y(F)| \geq \delta^{3\eta}.$$ Since $|\F_\rho| \leq \rho^{-2-\eps_1}$, we have $\lambda(\F_\rho, Y_\rho'') \gtrapprox \delta^{3\eta} \rho^{\eps_1} \geq \rho^{3\eta/\eps + \eps_1}$, and so $\lambda(\F_\rho, Y_\rho'') \geq \rho^{2\eps_1}$ if $\eta$ is much smaller than $\eps_1$. To summarize, we get a $\gtrapprox 1$ refinement $(\F_1, Y_2)$ of $(\F, Y)$, and a $2\eps_1$-extremal cover $(\F_\rho, Y_\rho'')$ of it. By our construction, $(\F_\rho, Y_\rho'')$ is formed by unions of $(\rho/10n)$-cubes, and any of its $\gtrapprox 1$ refinement (that is still formed by unions of $(\rho/10n)$-cubes) will also naturally induce a $\gtrapprox 1$ refinement of $(\F_1, Y_2)$ (thus of $(\F, Y)$). This allows us to only focus on suitably refining $(\F_\rho, Y_\rho'')$.

Applying Proposition~\ref{Prop:self-similar} again to $(\F_\rho, Y_\rho'')$ with $\tau$ in place of $\rho$ and $\eps_2$ in place of $\eps$, we obtain a $c$-refinement ($c \gtrapprox 1$) $(\F_{\rho,1}, Y_{\rho, 1})$ which is covered by an $\eps_2$-extremal set of $\tau$-flakes $(\F_\tau, Y_\tau)$. By item~(v) of Proposition~\ref{Prop:self-similar}, for any $(\tau/10n)$-cube $\Tilde{Q}$ intersecting $U(\F_{\rho,1}, Y_{\rho, 1})$, we have (by noting that $\Tilde{Q} \subset B(x, \tau)$ for any $x \in U(\F_{\rho,1}, Y_{\rho, 1}) \cap \Tilde{Q}$) 
\begin{align*}
    |U(\F_{\rho,1}, Y_{\rho, 1}) \cap \Tilde{Q}| &\leq |U(\F_{\rho,1}, Y_{\rho, 1}) \cap B(x, \tau)| \\
    &\sim |B(x, \tau)| \cdot \frac{|U(\F_{\rho,1}, Y_{\rho, 1})|}{|U(\F_\tau, Y_\tau)|} \\
    & \lesssim \tau^n \rho^{\sigma - \eps_1} \tau^{-\sigma - \eps_3}
\end{align*}
by choosing $\eps_2$ small depending on $\eps_3$. However, the tricky issue is that $(\F_{\rho,1}, Y_{\rho, 1})$ may no longer be formed by unions of $(\rho/10n)$-cubes. This can be remedied as follows. For each $F_\rho \in \F_{\rho,1}$, we define $Y_{\rho,2}(F_\rho)$ to be the union of all $(\rho/10n)$-cubes $Q \subset Y_{\rho}''(F_\rho)$ with $|Y_{\rho, 1}(F_\rho) \cap Q| \geq (c/2)|Q|$. Then $|Y_{\rho,1}(F_\rho) \setminus Y_{\rho,2}(F_\rho)| \leq (c/2)|Y_\rho''(F_\rho)|$ for each $F_\rho \in \F_{\rho,1}$, and so 
\begin{align*}
\sum_{\F_{\rho,1}} |Y_{\rho,2}(F_\rho)| &\geq \sum_{\F_{\rho,1}} |Y_{\rho,2}(F_\rho) \cap Y_{\rho,1}(F_\rho)| \\
&= \sum_{\F_{\rho,1}} |Y_{\rho,1}(F_\rho)| - \sum_{\F_{\rho,1}} |Y_{\rho,1}(F_\rho) \setminus Y_{\rho,2}(F_\rho)| \\
&\geq c \sum_{\F_\rho}|Y_\rho''(F_\rho)| - (c/2)\sum_{\F_{\rho,1}}|Y_\rho''(F_\rho)| \geq (c/2)\sum_{\F_{\rho}}|Y_{\rho}''(F_\rho)|,
\end{align*} 
which means $(\F_{\rho,1}, Y_{\rho,2})$ is a $\gtrapprox 1$ refinement of $(\F_\rho, Y_\rho'')$. Note that for each $(\rho/10n)$-cube $Q \subset U(\F_{\rho,1}, Y_{\rho,2})$, we must have $|U(\F_{\rho,1}, Y_{\rho,1}) \cap Q| \geq (c/2)|Q| \gtrapprox |Q|$. Since $\tau/\rho$ is an integer, for any $(\tau/10n)$-cube $\Tilde{Q}$ intersecting $U(\F_{\rho,1}, Y_{\rho,2})$ (and thus $U(\F_{\rho,1}, Y_{\rho,1})$), we have 
\begin{equation}
\begin{split}
    |U(\F_{\rho,1}, Y_{\rho, 2}) \cap \Tilde{Q}| &= \sum_{(\rho/10n)\text{-cube } Q \subset U(\F_{\rho,1}, Y_{\rho, 2}) \cap  \Tilde{Q}} |Q| \\
    & \lessapprox \sum_{(\rho/10n)\text{-cube } Q \subset U(\F_{\rho,1}, Y_{\rho, 2}) \cap  \Tilde{Q}} |U(\F_{\rho,1}, Y_{\rho, 1}) \cap Q| \\
    & \leq |U(\F_{\rho,1}, Y_{\rho, 1}) \cap \Tilde{Q}| \\
    & \lesssim \rho^{\sigma - \eps_1} \tau^{n-\sigma -\eps_3}.
\end{split}\label{Eq:tau_cube_density}
\end{equation}
In fact, the density bound \eqref{Eq:tau_cube_density} is the only property that we need later on regarding scale $\tau$, and refine $(\F_{\rho,1}, Y_{\rho, 2})$ will not affect this property. So we may forget about $(\F_\tau, Y_\tau)$ and only aim at finding a suitable refinement of $(\F_{\rho,1}, Y_{\rho, 2})$ from now on. We point out that $(\F_{\rho,1}, Y_{\rho, 2})$ is $3\eps_1$-extremal.

Applying Lemma~\ref{Lem:robust_transverse} to $(\F_{\rho,1}, Y_{\rho, 2})$ with $\eps_2$ in place of $\eps$, we obtain a $\gtrapprox 1$ refinement $(\F_{\rho,3}, Y_{\rho,3})$ so that for each $x \in \R^n$ and $V \in \Gr(n,n-2)$, we have
\begin{align}\label{Eq:loc_const_trans}
    |\{F\in ({\mathbb F}_{\rho,3})_{Y_{\rho,3}}(x): \angle(V,\dir(F)) <\rho^{\eps_2}\}| 
    \le \rho^{-\sigma+\eps_2\sigma/2}.
\end{align}
For each $(\rho/10n)$-cube $Q$ intersecting $U(\F_{\rho,3}, Y_{\rho,3})$, we can choose an $x_Q \in Q$ so that $\mu(\F_{\rho,3}, Y_{\rho,3})(x_Q) = \max_{x\in Q} \mu(\F_{\rho,3}, Y_{\rho,3})(x) > 0$. If we define $Y_{\rho,3}'(F_\rho)$ to be the union of all $(\rho/10n)$-cubes $Q$ with $x_Q \in U(\F_{\rho,3}[F_\rho], Y_{\rho,3})$ for each $F_\rho \in \F_{\rho,3}$, then it is easy to verify that $(\F_{\rho,3}, Y_{\rho,3}')$ is still a $\gtrapprox 1$ refinement of $(\F_{\rho,1}, Y_{\rho, 2})$ and satisfies \eqref{Eq:loc_const_trans} with $Y_{\rho,3}'$ in place of $Y_{\rho,3}$. However, the good thing is that now $(\F_{\rho,3}, Y_{\rho,3}')$ is formed by unions of $(\rho/10n)$-cubes. Abusing notation, we may assume that $(\F_{\rho,3}, Y_{\rho,3})$ itself is formed by unions of $(\rho/10n)$-cubes.

Let $(\F_{\rho,4}, Y_{\rho,4})$ be a $\gtrsim 1$ refinement of $(\F_{\rho,3}, Y_{\rho,3})$ that is formed by unions of $(\rho/10n)$-cubes and has the following property ($*$): For each $F_\rho \in \F_{\rho,4}$ and each $\tau$-cube $\Tilde{Q}$ that intersects $Y_{\rho,4}(F_\rho)$, we have $|\Tilde{Q} \cap Y_{\rho,4}(F_\rho)| \geq \rho^{2+3\eps_1} \tau^{n-2}$. This can be done by recalling that $\tau/\rho$ is an integer and defining 
\begin{align*}
    Y_{\rho,4}(F_\rho) \coloneqq \bigsqcup_{\Tilde{Q}: |\Tilde{Q} \cap Y_{\rho,3}(F_\rho)| \geq \rho^{2+3\eps_1} \tau^{n-2}} \Tilde{Q} \cap Y_{\rho,3}(F_\rho),
\end{align*}
which ensures that the total mass of shadings we discarded is 
\begin{align*}
&\lesssim \tau^{-(n-2)} \rho^{2+3\eps_1} \tau^{n-2} \rho^{-2-\eps_1}\\ 
&\ll \rho^{\eps_1} \leq \rho^{-\eps_1} \sum_{\F_\rho} |Y_\rho(F_\rho)| \lessapprox \rho^{-\eps_1} \delta^{-\eta} \sum_{\F_{\rho,1}} |Y_{\rho,2}(F_\rho)| \lessapprox \sum_{\F_{\rho,4}}|Y_{\rho,4}(F_\rho)|.
\end{align*} 
Let $(\F_{\rho,5}, Y_{\rho,5})$ be a constant multiplicity refinement of $(\F_{\rho,4}, Y_{\rho,4})$. Clearly, we may assume that $(\F_{\rho,5}, Y_{\rho,5})$ is still formed by unions of $(\rho/10n)$-cubes. Finally, by the same procedure that we did for $(\F_{\rho,3}, Y_{\rho,3})$, let $(\F_{\rho,6}, Y_{\rho,6})$ be a $\gtrsim 1$ refinement of $(\F_{\rho,5}, Y_{\rho,5})$ that is formed by unions of $(\rho/10n)$-cubes and again has the following property ($*$). After all these refinements, if $F_{\rho,0} \in \F_{\rho,6}$ and if $\Tilde{Q}$ is a $\tau$-cube that intersects $Y_{\rho,6}(F_{\rho,0})$, then $|\Tilde{Q} \cap Y_{\rho,6}(F_{\rho,0})| \geq \rho^{2+3\eps_1}\tau^{n-2}$, and for any $x \in \Tilde{Q} \cap Y_{\rho,6}(F_{\rho,0})$, we have 
\begin{align*}
\mu(\F_{\rho,5}, Y_{\rho,5})(x) &= (\sum_{\F_{\rho, 5}}|Y_{\rho, 5}(F_\rho)|) / |U(\F_{\rho,5}, Y_{\rho, 5})| \\
&\gtrapprox (\sum_{\F_\rho}|Y_\rho''(F_\rho)|)/|U(\F_\rho, Y_\rho)| \geq \delta^{3\eta}/\rho^{\sigma - \eps_1} \geq \rho^{-\sigma + 2\eps_1}.
\end{align*} Because of the refinement we did in \eqref{Eq:loc_const_trans}, if $\eps_1$ is sufficiently small compared to $\eps_2$, then there is a $\rho$-flake $F_{\rho, x} \in (\F_{\rho,5})_{Y_{\rho,5}}(x) \subset (\F_{\rho,4})_{Y_{\rho,4}}(x)$ with $\angle(F_{\rho,0}, F_{\rho, x}) \geq \rho^{\eps_2}$.  

We claim that the union of the sets $\{\Tilde{Q} \cap Y_{\rho,4}(F_{\rho,x})\}_{x \in \Tilde{Q} \cap Y_{\rho,6}(F_{\rho,0})}$ is contained inside a slab segment of dimensions $\delta^{-\eta} \rho \times \tau \times \tau\1$. Indeed, since the plane map $\HH$ is $\delta^{-\eta}$-Lipschitz and $\Tilde{Q}$ has diameter $\leq \sqrt{n} \tau \leq \sqrt{n} \rho^{1/2}$, if $y \in \Tilde{Q} \cap Y_{\rho,6}(F_{\rho,0})$, then we have 
\begin{align}\label{Eq:brush_upper_vol}
    \Tilde{Q} \cap \bigcup_{x \in \Tilde{Q} \cap Y_{\rho,6}(F_{\rho,0})} F_{\rho,x} 
    \subset
    \Tilde{Q} \cap N_{C_n\delta^{-\eta}\rho}(y + \HH(\Tilde{Q})^\perp)
\end{align}
for some constant $C_n$ depending only on $n$. 
Here $\HH(\Tilde{Q})$ is well defined as $\Tilde{Q}$ must intersect $U(\F,Y)$. This verifies the claim. Note that $\Tilde{Q} \cap F_{\rho,0}$ is also contained in this same slab segment.

Next, we want to show that the union of the sets $\{\Tilde{Q} \cap Y_{\rho,4}(F_{\rho,x})\}_{x \in \Tilde{Q} \cap Y_{\rho,6}(F_{\rho,0})}$ has large volume inside the rectangular prism. For this purpose, we project these sets onto the central hyperplane $y + \HH(\Tilde{Q})^\perp$, which essentially turns things into slabs in $\R^{n-1}$, and then apply Lemma~\ref{Lem:slab_0}. Let $\pi_{\Tilde{Q}, y}$ be the orthogonal projection from $\R^n$ to $y + \HH(\Tilde{Q})^\perp$. For each $(\rho/10n)$-cube $Q \subset \Tilde{Q}$, let $Q_0 \subset Q$ be a $(\rho/10n^2)$-cube with the same center as $Q$ and one face parallel to $\HH(\Tilde{Q})^\perp$. For each $x \in \Tilde{Q} \cap Y_{\rho,6}(F_{\rho,0})$, define $\Tilde{F}_{\rho, x} \coloneqq \Tilde{N}_\rho(\pi_{\Tilde{Q},y}(\Tilde{Q} \cap F_{\rho,x}))$ and $\dir(\Tilde{F}_{\rho, x}) \coloneqq \pi_{\Tilde{Q}, y}(\dir({F}_{\rho, x}))$, where $\Tilde{N}$ denotes the neighborhood in $y + \HH(\Tilde{Q})^\perp \cong \R^{n-1}$. Then define $$\Tilde{Y}_{\rho,4}(\Tilde{F}_{\rho,x}) \coloneqq \cup_{Q \subset \Tilde{Q} \cap Y_{\rho,4}(F_{\rho,x})} \pi_{\Tilde{Q}, y}(Q_0).$$ Note that $\Tilde{F}_{\rho, 0}$, $\dir(\Tilde{F}_{\rho, 0})$, and $\Tilde{Y}_{\rho, 6}(\Tilde{F}_{\rho, 0})$ can also be defined in the same obvious way. Note that Lemma~\ref{Lem:faithful_proj} implies $\dim \dir(\Tilde{F}_{\rho,x}) = \dim \dir(F_{\rho,x}) = n-2$ and $\dim \dir(\Tilde{F}_{\rho,0}) = \dim \dir(F_{\rho,0}) = n-2$. Now inside $y + \HH(\Tilde{Q})^\perp$, $\Tilde{Y}_{\rho,4}$ defines a shading on each $\rho \times \tau\1$ slab $\Tilde{F}_{\rho,x}$ with 
$$|\Tilde{Y}_{\rho,4}(\Tilde{F}_{\rho,x})| \gtrsim \delta^{\eta} \rho^{-1} |\Tilde{Q} \cap Y_{\rho,4}(F_{\rho,x})| \geq \delta^{\eta} \rho^{1 + 3\eps_1}\tau^{n-2}.$$ 
Similarly, we have 
$$|\Tilde{Y}_{\rho,6}(\Tilde{F}_{\rho,0})| \gtrsim \delta^{\eta} \rho^{-1} |\Tilde{Q} \cap Y_{\rho,6}(F_{\rho,0})| \geq \delta^{\eta} \rho^{1 + 3\eps_1}\tau^{n-2},$$ 
which allows us to choose a set of points $z_1, \dots, z_k \in \Tilde{Y}_{\rho,6}(\Tilde{F}_{\rho,0})$ that are $\geq \rho$ separated, with $k \gtrsim \delta^\eta \rho^{3\eps_1}(\tau/\rho)^{n-2}$. By our construction of the shadings, each $z_i$ ($1 \leq i \leq k$) is contained in $\Tilde{Y}_{\rho,4}(\Tilde{F}_{\rho,x_i})$ for some $x_i \in \Tilde{Q} \cap Y_{\rho,6}(F_{\rho,0})$. Let $\Tilde{\F}_{\rho} \coloneqq \{\Tilde{F}_{\rho,x_i}\}_{i=1}^k$, then $\lambda(\Tilde{\F}_\rho, \Tilde{Y}_{\rho,4}) \gtrsim (\delta^{\eta} \rho^{1 + 3\eps_1}\tau^{n-2})/(\rho\tau^{n-2}) = \delta^\eta \rho^{3\eps_1}$. 

For each $i$, recall that $\angle(F_{\rho,0}, F_{\rho, x_i}) \geq \rho^{\eps_2}$, and note that $$\angle(\Tilde{F}_{\rho,x_i}, F_{\rho, x_i}) \lesssim (\delta^{-\eta} \rho)/(\rho^{3\eps_1}\tau) \leq \rho^{-\eta/\eps - 3\eps_1 + \eps}$$ by \eqref{Eq:brush_upper_vol} and (a rescaled version of) Lemma~\ref{Lem:flake_contain}, which means $$\angle(\Tilde{F}_{\rho,x_i}, F_{\rho, x_i}) \leq \rho^{\eps/2}$$ if $\eta$ and $\eps_1$ are much smaller than $\eps$. Similarly, $\angle(\Tilde{F}_{\rho,0}, F_{\rho, 0}) \leq \rho^{\eps/2}$. Therefore, by Remark~\ref{Rem:major_metric}, we have $$\angle(\Tilde{F}_{\rho,x_i}, \Tilde{F}_{\rho, 0}) \geq \rho^{\eps_2} - 2\rho^{\eps/2} \geq \rho^{2\eps_2}$$ for each $i$, if $\eps_2$ is much smaller than $\eps$. It is worth pointing out that since we are now in $\R^{n-1}$, the angle ``$\angle$'' between two $(n-2)$-planes can be interpreted in the classical sense. For any convex set $W \subset \pi_{\Tilde{Q},y}(\Tilde{Q})$, if it contains some $\Tilde{F}_{\rho, x_j} \in \Tilde{\F}_\rho$, then it has dimension $a \times \tau\1$ for some $a \in [2\rho, \tau]$. In particular, we have 
$$|W \cap \Tilde{Y}_{\rho,6}(\Tilde{F}_{\rho,0})| \leq |\pi_{\Tilde{Q},y}(\Tilde{Q}) \cap \Tilde{N}_a(\Tilde{F}_{\rho,x_j}) \cap \Tilde{F}_{\rho,0}| \lesssim \rho^{1-2\eps_2} a \tau^{n-3},$$
which by the separation property of $z_i$'s implies that
\begin{align*}
    |W \cap \{z_i\}_{i=1}^k| \lesssim 
    \rho^{-(n-2)-2\eps_2} a \tau^{n-3}.
\end{align*}
Since each $z_i \in \Tilde{F}_{\rho,x_i}$, we have $|\Tilde{\F}_\rho[W]|\leq |W \cap \{z_i\}_{i=1}^k|$, and so
\begin{align*}
    |\Tilde{\F}_\rho[W]| 
    \lesssim 
    \rho^{-(n-2)-2\eps_2} a \tau^{n-3} 
    \lesssim 
    \delta^{-\eta} \rho^{-3\eps_1-2\eps_2} \tau^{-(n-1)} |W| |\Tilde{\F}_\rho| \lesssim \delta^{-\eta} \rho^{-3\eps_1-2\eps_2}  \frac{|W|}{|\pi_{\Tilde{Q},y}(\Tilde{Q})|} |\Tilde{\F}_\rho|
\end{align*}
by recalling $|\Tilde{\F}_\rho| = k \gtrsim \delta^\eta \rho^{3\eps_1}(\tau/\rho)^{n-2}$. This means $\Tilde{\F}_\rho$ is $\lesssim \delta^{-\eta} \rho^{-3\eps_1-2\eps_2}$-Frostman in $\pi_{\Tilde{Q},y}(\Tilde{Q})$, and so by (a rescaled version of) Lemma~\ref{Lem:slab_0}\footnote{As pointed out therein, this lemma does not require $\Tilde{\F}_\rho$ to be essentially distinct.}, we have
\begin{align*}
    |U(\Tilde{\F}_\rho, \Tilde{Y}_{\rho,4})| 
    &\gtrapprox C_F(\Tilde{\F}_\rho, \pi_{\Tilde{Q},y}(\Tilde{Q}))^{-1}  \lambda(\Tilde{\F}_\rho, \Tilde{Y}_{\rho,4})^2 |\pi_{\Tilde{Q},y}(\Tilde{Q})|\\
    &\gtrsim
    \delta^{\eta} \rho^{3\eps_1+2\eps_2} (\delta^\eta \rho^{3\eps_1})^2 \tau^{n-1} 
    = \delta^{3\eta}\rho^{9\eps_1 + 2\eps_2} \tau^{n-1}.
\end{align*}
Finally, by our construction of the shadings, for each $z \in U(\Tilde{\F}_\rho, \Tilde{Y}_{\rho,4})$, we know that $\pi_{\Tilde{Q},y}^{-1}(z) \cap \Tilde{Q} \cap \bigcup_{x \in \Tilde{Q} \cap Y_{\rho,6}(F_{\rho,0})} F_{\rho,x}$ at least contains a line segment of length $\rho/10n^2$, and so we obtain
\begin{align}\label{Eq:brush_lower_vol}
    \Big|\Tilde{Q} \cap \bigcup_{x \in \Tilde{Q} \cap Y_{\rho,6}(F_{\rho,0})} F_{\rho,x}\Big| \gtrsim \rho |U(\Tilde{\F}_\rho, \Tilde{Y}_{\rho,4})| 
    \gtrapprox     \delta^{3\eta}\rho^{1 + 9\eps_1 + 2\eps_2} \tau^{n-1}.
\end{align}

Combining \eqref{Eq:brush_upper_vol} and \eqref{Eq:brush_lower_vol}, we can conclude that if $t \in \pi_{\HH(\Tilde{Q})}(\Tilde{Q} \cap U(\F_{\rho, 6}, Y_{\rho, 6}))$, then 
\begin{align}\label{Eq:grain_lower_vol}
|\Tilde{Q} \cap U(\F_{\rho, 4}, Y_{\rho, 4}) \cap N_{C_n\delta^{-\eta}\rho}(t+\HH(\Tilde{Q})^\perp)| \gtrapprox \delta^{3\eta}\rho^{1 + 9\eps_1 + 2\eps_2} \tau^{n-1}.
\end{align}
Informally, \eqref{Eq:grain_lower_vol} says that each $\rho$-separated point $t \in \pi_{\HH(\Tilde{Q})}(\Tilde{Q} \cap U(\F_{\rho, 6}, Y_{\rho, 6}))$ corresponds to a $\rho \times \tau \times \tau\1$ grain (i.e. rectangular prism whose short axis is parallel to $\HH(\Tilde{Q})$).

If $\norm{t-t'} \geq 2C_n \delta^{-\eta}\rho$, then the sets $N_{C_n\delta^{-\eta}\rho}(t+\HH(\Tilde{Q})^\perp)$ and $N_{C_n\delta^{-\eta}\rho}(t'+\HH(\Tilde{Q})^\perp)$ are disjoint up to a null set. Hence \eqref{Eq:grain_lower_vol} together with \eqref{Eq:tau_cube_density} (cover any $\tau$-cube by $O(1)$ $(\tau/10n)$-cubes) implies that for any $\tau$-cube $\Tilde{Q}$ intersecting $U(\F,Y)$, we have 
\begin{align*}
\EE_{2C_n\delta^{-\eta} \rho}\big( \pi_{\HH(\Tilde{Q})}(\Tilde{Q} \cap U(\F_{\rho, 6}, Y_{\rho, 6})) \big) &\lessapprox 
\frac{|\Tilde{Q} \cap U(\F_{\rho, 4}, Y_{\rho, 4})|}{\delta^{3\eta}\rho^{1 + 9\eps_1 + 2\eps_2} \tau^{n-1}} \\
&\lessapprox
\frac{\rho^{\sigma-\eps_1} \tau^{n-\sigma-\eps_3}}{\delta^{3\eta}\rho^{1 + 9\eps_1 + 2\eps_2} \tau^{n-1}}\\
&= \delta^{-3\eta} \rho^{-10\eps_1 - 2\eps_2 - \eps_3} (\tau/\rho)^{1-\sigma}.
\end{align*}
This immediately gives 
\begin{equation}
\begin{split}
    \EE_{\rho}\big( \pi_{\HH(\Tilde{Q})}(\Tilde{Q} \cap U(\F_{\rho, 6}, Y_{\rho, 6})) \big) 
    &\leq 
    2C_n\delta^{-\eta} \EE_{2C_n\delta^{-\eta} \rho}\big( \pi_{\HH(\Tilde{Q})}(\Tilde{Q} \cap U(\F_{\rho, 6}, Y_{\rho, 6})) \big)\\
    &\lessapprox 
    \delta^{-4\eta} \rho^{-10\eps_1 - 2\eps_2 - \eps_3} (\tau/\rho)^{1-\sigma}.
\end{split}\label{Eq:E_rho_local_grain'}
\end{equation}
To conclude the proof, select $\eps_1, \eps_2, \eps_3, \eta$ sufficiently small (depending on $\eps$) so that the LHS of \eqref{Eq:E_rho_local_grain'} is at most $\rho^{-\eps}(\tau/\rho)^{1-\sigma}$ (this is possible since $\rho \leq \delta^\eps$). Let $(\F', Y')$ be the $\gtrapprox 1$ refinement of $(\F_1, Y_2)$ (thus of $(\F, Y)$) induced by the $\gtrapprox 1$ refinement $(\F_{\rho, 6}, Y_{\rho, 6})$ of $(\F_\rho, Y_\rho'')$. Since $U(\F', Y') \subset U(\F_{\rho, 6}, Y_{\rho, 6})$, \eqref{Eq:E_rho_local_grain'} implies that $(\F', Y')$ satisfies \eqref{Eq:E_rho_local_grain}. Since $(\F, Y)$ is $\eta$-extremal, we know that $(\F', Y')$ is $2\eta$-extremal by Remark~\ref{Rem:extreme_mult}.
\end{proof}

Lemma~\ref{Lem:grains_number_given_scale} says that when restricted to a $\tau$-cube and examined at scale $\rho$, the set $U(\F',Y')$ looks like a union of about $(\tau/\rho)^{1-\sigma}$ parallel grains, whose normal directions all coincide and are given by the plane map $\HH$. Lemma~\ref{Lem:grains_single_rho_tau} asserts that a similar
statement continues to hold if we replace a $\tau$-cube with a (potentially much larger) $\tau \times \rho^{1/2} \times \rho^{1/2}\1$ 
thick slab segment, whose central hyperplane is orthogonal to the plane map $\HH$. To prove this lemma, we need to show that at scale $\rho$, the grains inside a $\tau$-cube are pieces of the grains inside a $\rho^{1/2}$-cube. 
{\color{black} The proof is similar to that in \cite[Lemma~4.10]{WZ26sticky}, but written in a fully rigorous way.}

\begin{lem}\label{Lem:grains_single_rho_tau}
For all $\eps>0$, there exist $\eta, \delta_0 >0$ so that the following holds for all $\delta \in (0, \delta_0]$. Let $(\F,Y)$ be an $\eta$-extremal set of $\delta$-flakes with $\delta^{-\eta}$-Lipschitz plane map $\HH$. Let $\rho \in [\delta^{1-\eps}, \delta^\eps]$ and let $\tau \in [\rho^{1-\eps}, \rho^{1/2}]$. For each $\rho^{1/2}$-cube $\bar{Q}$ intersecting $U(\F,Y)$, define $\HH(\bar{Q}) \coloneqq \HH(y_{\bar{Q}})$, where $y_{\bar{Q}}$ is any fixed point in $\bar{Q} \cap U(\F,Y)$. 

Then there exists an $\eps$-extremal refinement $(\F',Y')$, so that for each $\rho^{1/2}$-cube $\bar{Q}$ intersecting $U(\F,Y)$ and each line segment $J \subset \HH(\bar{Q})$ of length $\tau$, we have
\begin{align}\label{Eq:grains_single_rho_tau}
\EE_\rho\left( J \cap \pi_{\HH(\bar{Q})}(U(\F',Y') \cap \bar{Q}) \right) \leq \rho^{-\eps}(\tau/\rho)^{1-\sigma}.
\end{align}
\end{lem}
\begin{proof}
It suffices to only consider the case when both $\rho^{1/2}/\tau$ and $\tau/\rho$ are integers. Let $\eps_1, \eps_2 > 0$ be small constants to be chosen below. We will select $\eta$ very small compared to $\eps_1$, $\eps_1$ very small
compared to $\eps_{2}$, and $\eps_2$ very small compared to $\eps$. 

By the proof of Lemma~\ref{Lem:grains_number_given_scale}, we can first get a $\gtrapprox 1$ refinement $(\F_1, Y_1)$ of $(\F, Y)$, and a $2\eps_1$-extremal cover $(\F_{\rho,1}, Y_{\rho,1})$ of it, which is formed by unions of $(\rho/10n)$-cubes. Moreover, any $c$-refinement of $(\F_{\rho,1}, Y_{\rho,1})$ that is still formed by unions of $(\rho/10n)$-cubes will also naturally induce a $c$-refinement of $(\F_1, Y_1)$ (thus of $(\F, Y)$). This allows us to only focus on suitably refining $(\F_{\rho,1}, Y_{\rho,1})$. Then we apply the remaining parts of Lemma~\ref{Lem:grains_number_given_scale} twice to $(\F_{\rho,1}, Y_{\rho,1})$, first at the scale $\tau$ and next at the scale $\rho^{1/2}$. As a result, we obtain a $\gtrapprox 1$ refinement $(\F_{\rho,2}, Y_{\rho,2})$ (formed by unions of $(\rho/10n)$-cubes) of $(\F_{\rho,1}, Y_{\rho,1})$ with the following two properties: 
\begin{itemize}
    \item For each $\rho^{1/2}$ cube $\bar{Q}$ intersecting $(\F,Y)$, we have
    \begin{align}\label{Eq:rho_1/2_counting}
    \EE_{\rho}\big( \pi_{\HH(\bar{Q})}(U(\F_{\rho, 2}, Y_{\rho, 2})  \cap \bar{Q}) \big) 
    \leq 
    \rho^{-\eps_2} \rho^{(\sigma-1)/2}.
    \end{align}

    \item For each $\rho^{1/2}$ cube $\bar{Q}$, and each $\tau$-cube $\Tilde{Q} \subset \bar{Q}$ intersecting $(\F,Y)$, we have     \begin{align}\label{Eq:tau_counting}
    \EE_{\rho}\big( \pi_{\HH(\bar{Q})}(U(\F_{\rho, 2}, Y_{\rho, 2})  \cap \Tilde{Q}) \big) 
    \lesssim 
    \delta^{-\eta} 
    \rho^{-\eps_2} (\tau/\rho)^{1-\sigma}.
    \end{align}
\end{itemize}
Note that in \eqref{Eq:tau_counting} we have the unit vector $\HH(\bar{Q})$ rather than $\HH(\Tilde{Q})$. This requires justification. By the proof of Lemma~\ref{Lem:grains_number_given_scale}, we initially get the estimate $\EE_{\rho}\big( \pi_{\HH(\Tilde{Q})}(U(\F_{\rho, 2}, Y_{\rho, 2})  \cap \Tilde{Q}) \big) \leq \rho^{-\eps_2} (\tau/\rho)^{1-\sigma}$, where $\HH(\Tilde{Q}) \coloneqq \HH(y_{\Tilde{Q}})$ for a fixed point $y_{\Tilde{Q}} \in \Tilde{Q} \cap U(\F,Y)$. To obtain \eqref{Eq:tau_counting}, we bound 
\begin{equation}
\begin{split}
  \EE_{\rho}\big( \pi_{\HH(\bar{Q})}(U(\F_{\rho, 2}, Y_{\rho, 2})  \cap \Tilde{Q}) \big) & \leq 10n\delta^{-\eta} \EE_{10n\delta^{-\eta} \rho}\big( \pi_{\HH(\bar{Q})}(U(\F_{\rho, 2}, Y_{\rho, 2})  \cap \Tilde{Q}) \big) \\
  &\lesssim 
  \delta^{-\eta} \EE_{n \delta^{-\eta} \rho}\big( \pi_{\HH(\Tilde{Q})}(U(\F_{\rho, 2}, Y_{\rho, 2})  \cap \Tilde{Q}) \big) \\
  &\leq 
  \delta^{-\eta} \EE_{\rho}\big( \pi_{\HH(\Tilde{Q})}(U(\F_{\rho, 2}, Y_{\rho, 2})  \cap \Tilde{Q}) \big) \\
  &\leq 
  \delta^{-\eta} \rho^{-\eps_2} (\tau/\rho)^{1-\sigma}.
\end{split}
\end{equation}
The key step is the second inequality, which used the fact that $\on{diam}(\Tilde{Q}) \leq \sqrt{n}\rho^{1/2}$ and $\angle(\HH(\bar{Q}), \HH(\Tilde{Q})) \leq \sqrt{n} \delta^{-\eta}\rho^{1/2}$; the latter statement follows from the fact that $\Tilde{Q} \subset \bar{Q}$, $\on{diam}(\bar{Q}) \leq \sqrt{n} \rho^{1/2}$, and $\HH$ is $\delta^{-\eta}$-Lipschitz. 

Let $(\F_{\rho,2}', Y_{\rho,2}')$ be the constant multiplicity refinement of $(\F_{\rho,2}, Y_{\rho,2})$ formed by unions of $(\rho/10n)$-cubes. Then $(\F_{\rho,2}', Y_{\rho,2}')$ is $3\eps_1$-extremal, and so $\mu(\F_{\rho,2}', Y_{\rho,2}') = (\sum_{\F_{\rho,2}'}|Y_{\rho,2}'(F)|)/|U(\F_{\rho,2}', Y_{\rho,2}')| \leq \rho^{-2-3\eps_1}\rho^2/\rho^{\sigma + \eps_2} \leq \rho^{-\sigma - 2\eps_2}$ if $\eps_1$ is very small compared to $\eps_2$ (recall Remark~\ref{Rem:extreme_mult}). Applying Proposition~\ref{Prop:self-similar} to $(\F_{\rho,2}', Y_{\rho,2}')$ with $\eps_2$ in place of $\eps$ and $\rho^{1/2}$ in place of $\rho$, we obtain a $\gtrapprox 1$ refinement $(\F_{\rho,3}, Y_{\rho,3})$ of $(\F_{\rho,2}', Y_{\rho,2}')$, which is covered by an $\eps_2$-extremal set of $\rho^{1/2}$-flakes $(\F_{\rho^{1/2}}, Y_{\rho^{1/2}})$. 

Let us first show that $U(\F_{\rho^{1/2}}, Y_{\rho^{1/2}})$ is contained in a union of $\rho^{1/2}$-cubes of cardinality $\lesssim \rho^{(\sigma-n)/2 - \eps_2/2}$. Let $\{x_i\}_{i=1}^M$ be a maximal set of $\rho^{1/2}$-separated points in $U(\F_{\rho,3}, Y_{\rho,3})$, then $\{B(x_i, 2\rho^{1/2})\}_{i=1}^M$ is a finitely overlapping covering of $U(\F_{\rho^{1/2}}, Y_{\rho^{1/2}})$, which is induced by $(\F_{\rho,3}, Y_{\rho,3})$. Divide this covering into $j$ groups $\{B(x_{i,j}, 2\rho^{1/2})\}_{i=1}^{M_j}$ with $j \lesssim 1$ such that the balls inside each group are $10n\rho^{1/2}$ separated. Again, since $(\F_{\rho^{1/2}}, Y_{\rho^{1/2}})$ is induced by $(\F_{\rho,3}, Y_{\rho,3})$, for each $j$, the (unique) $(\rho^{1/2}/10n)$-cube $\bar{Q}_{i,j}'$ containing $x_{i,j}$ is wholly contained in both $B(x_{i,j}, 2\rho^{1/2})$ and $U(\F_{\rho^{1/2}}, Y_{\rho^{1/2}})$. On the other hand, each $B(x_{i,j},2\rho^{1/2})$ can be covered by $\rho^{1/2}$-cubes $\{\bar{Q}_{i,j,k}\}_{k=1}^{M_{i,j}}$ with $M_{i,j} \lesssim 1$ that are disjoint for fixed $j$. Note that $U(\F_{\rho,3}, Y_{\rho,3}) \subset \cup_{i,j,k} \bar{Q}_{i,j,k}$, and for each $j$, we have 
\begin{align*}
|\sqcup_{i,k} \bar{Q}_{i,j,k}| 
= \sum_{i} |\sqcup_k \bar{Q}_{i,j,k}| 
\lesssim
\sum_{i} (\rho^{1/2})^n \lesssim \sum_{i} |\bar{Q}_{i,j}'|
= |\sqcup_i \bar{Q}_{i,j}'| 
\leq 
|U(\F_{\rho^{1/2}}, Y_{\rho^{1/2}})|,
\end{align*}
which by summing over $j$ means
\begin{align*}
|\cup_{i,j,k} \bar{Q}_{i,j,k}| \lesssim |U(\F_{\rho^{1/2}}, Y_{\rho^{1/2}})| \leq (\rho^{1/2})^{\sigma - \eps_2}. 
\end{align*}
Therefore, $U(\F_{\rho^{1/2}}, Y_{\rho^{1/2}})$ can be covered by a set $\bar{\Q}$ of $\rho^{1/2}$-cubes of cardinality $\lesssim (\rho^{1/2})^{\sigma - \eps_2}/(\rho^{1/2})^n = \rho^{(\sigma - n)/2 - \eps_2/2}$, as desired.

Let $$\bar{\Q}' \coloneqq \{\bar{Q} \in \bar{\Q}: |\int_{\bar{Q}}\mu(\F_{\rho,3}, Y_{\rho,3})(x) \D x \geq \rho^{(n-\sigma)/2 + \eps_2} \},$$ then 
\begin{align*}
\int_{U(\bar{\Q} \setminus \bar{\Q}')} \mu(\F_{\rho,3}, Y_{\rho,3})(x) \D x 
&\leq 
\rho^{(\sigma - n)/2 - \eps_2/2} \rho^{(n-\sigma)/2 + \eps_2} \\
&\leq 
\rho^{5\eps_1} 
\leq 
\rho^{\eps_1}
\sum_{\F_{\rho,1}} |Y_{\rho,1}(F)| 
\leq 
(1/2) 
\sum_{\F_{\rho,3}} |Y_{\rho,3}(F)|.
\end{align*}
So if we define $Y_{\rho,4}(F) \coloneqq Y_{\rho,3}(F) \cap U(\bar{\Q}')$ for each $F \in \F_{\rho,3}$, then $(\F_{\rho,3}, Y_{\rho,4})$ is a $\gtrsim 1$ refinement of $(\F_{\rho,3}, Y_{\rho,3})$ (thus a $\gtrapprox 1$ refinement of $(\F_{\rho,2}, Y_{\rho,2})$). The problem is that $(\F_{\rho,3}, Y_{\rho,4})$ may no longer be formed by unions of $(\rho/10n)$-cubes (although $(\F_{\rho,2}, Y_{\rho,2})$ is). To fix this issue, let $\F_{\rho,5} \coloneqq \F_{\rho,3}$, and for each $F_\rho \in \F_{\rho,5}$, define $Y_{\rho,5}(F_\rho)$ to be the union of all $(\rho/10n)$-cubes intersecting $Y_{\rho,4}(F_\rho)$. It is easy to see that the resulting $(\F_{\rho,5}, Y_{\rho,5})$ is still a $\gtrapprox 1$ refinement of $(\F_{\rho,2}, Y_{\rho,2})$, but now formed by unions of $(\rho/10n)$-cubes. Most importantly, $\bar{\Q}'$ covers $(\F_{\rho,5}, Y_{\rho,5})$, and for each $\rho^{1/2}$-cube $\bar{Q} \in \bar{\Q}'$, we have $$\int_{\bar{Q}}\mu(\F_{\rho,5}, Y_{\rho,5})(x) \D x \gtrsim \int_{\bar{Q}}\mu(\F_{\rho,3}, Y_{\rho,3})(x) \D x 
\geq \rho^{(n-\sigma)/2 + \eps_2}.$$

Note that for each $x \in \R^n$, we have $\mu(\F_{\rho,5}, Y_{\rho,5})(x) \leq \mu(\F_{\rho,2}', Y_{\rho,2}') \leq \rho^{-\sigma-2\eps_2}$. Also, the two itemized properties previously described (for $(\F_{\rho,2}, Y_{\rho,2})$) continue to hold for $(\F_{\rho,5}, Y_{\rho,5})$. Therefore, if $\bar{Q}$ is a $\rho^{1/2}$-cube that intersects $U(\F_{\rho,5}, Y_{\rho,5})$, then by \eqref{Eq:rho_1/2_counting}, $U(\F_{\rho,5}, Y_{\rho,5}) \cap \bar{Q}$ is contained in a union of at most $\rho^{-\eps_2} \rho^{(\sigma-1)/2}$ ``local grains'', which are sets of the form $\bar{Q} \cap N_{\rho}(t+\HH(\bar{Q})^\perp)$, where $t \in \HH(\bar{Q})$. By pigeonholing inside each $\bar{Q}$, we may assume that $U(\F_{\rho,5}, Y_{\rho,5}) \cap \bar{Q} \subset \bigsqcup_{j=1}^J G_{\bar{Q}, j}$, where $J \leq \rho^{-\eps_2} \rho^{(\sigma-1)/2}$, $G_{\bar{Q}, j} \coloneqq \bar{Q} \cap N_{\rho}(t_j+\HH(\bar{Q})^\perp)$, and $\{t_j\}_{j=1}^J \subset \HH(\bar{Q})$ is $10n\rho$-separated. Let $\G(\bar{Q}) \coloneqq \{ G_{\bar{Q}, j} \}_{j=1}^J$ and $$\G(\bar{Q})' \coloneqq \{ G_{\bar{Q}, j} \in \G(\bar{Q}): \int_{G_{\bar{Q}, j}} \mu(\F_{\rho,5}, Y_{\rho,5})(x) \D x \geq \rho^{(n+1 - 2\sigma)/2 + 2\eps_2} \},$$ 
then 
\begin{align*}
\int_{U(\G(\bar{Q}) \setminus \G(\bar{Q})')} \mu(\F_{\rho,5}, Y_{\rho,5})(x) \D x 
&\leq \rho^{-\eps_2}\rho^{(\sigma-1)/2} \rho^{(n+1 - 2\sigma)/2 + 2\eps_2}\\
&= 
\rho^{\eps_2} \rho^{(n-\sigma)/2 + \eps_2} 
\leq (1/2)\int_{\bar{Q}}\mu(\F_{\rho,5}, Y_{\rho,5})(x) \D x.
\end{align*}
This means that by only retaining the mass of $(\F_{\rho,5}, Y_{\rho,5})$ over all $U(\G(\bar{Q})')$, we obtain a $\gtrsim 1$ refinement $(\F_{\rho,6}, Y_{\rho,6})$ such that for each $\rho^{1/2}$-cube $\bar{Q}$ intersecting $U(\F_{\rho,6}, Y_{\rho,6})$, we have $U(\F_{\rho,6}, Y_{\rho,6}) \cap \bar{Q} \subset U(\G({\bar{Q}})')$ and $\int_{\bar{Q}}\mu(\F_{\rho,6}, Y_{\rho,6})(x) \D x 
\gtrsim \rho^{(n-\sigma)/2 + \eps_2}$. Moreover, for each $G_{\bar{Q}, j} \in \G(\bar{Q})'$, we have 
\begin{align}\label{Eq:fill_out_grain}
\int_{G_{\bar{Q}, j}} \mu(\F_{\rho,6}, Y_{\rho,6})(x) \D x \geq \rho^{(n+1 - 2\sigma)/2 + 2\eps_2}.
\end{align}
All the refinements in this paragraph can be suitably done to make sure that $(\F_{\rho,6}, Y_{\rho,6})$ is formed by unions of $(\rho/10n)$-cubes. Finally, as $\mu(\F_{\rho,6}, Y_{\rho,6})(x) \leq \mu(\F_{\rho,5}, Y_{\rho,5})(x) \leq \rho^{-\sigma - 2\eps_2}$ for each $x\in\R^n$, one can further conclude that 
\begin{align*}
    |U(\F_{\rho,6}, Y_{\rho,6}) \cap \bar{Q}| \geq \rho^{\sigma + 2\eps_2} \int_{\bar{Q}} \mu(\F_{\rho,6}, Y_{\rho,6})(x) \D x \gtrsim \rho^{(\sigma+n)/2 + 3\eps_2}.
\end{align*}

Next, let $\bar{Q}$ be one of the $\rho^{1/2}$-cubes described above. By Fubini's theorem, there must be a point $x \in \bar{Q}$ so that the line $\mathscr{L} \coloneqq x + \HH(\bar{Q})$ satisfies 
\begin{align}\label{Eq:pivot_line}
    |\mathscr{L} \cap U(\F_{\rho,6}, Y_{\rho,6}) \cap \bar{Q}| \gtrsim \rho^{(\sigma+1)/2 + 3\eps_2}, 
\end{align}
and so $\mathscr{L}$ intersects $U(\F_{\rho,6}, Y_{\rho,6})$ in $\gtrsim \rho^{(\sigma-1)/2 + 3\eps_2}$ local grains from $\G(\bar{Q})'$, which we denote by $\G(\bar{Q})''$. Then by \eqref{Eq:fill_out_grain}, we have 
\begin{align*}
\int_{U(\G(\bar{Q})'')} \mu(\F_{\rho,6}, Y_{\rho,6})(x) \D x 
\gtrsim 
\rho^{(\sigma-1)/2 + 3\eps_2} \rho^{(n+1 - 2\sigma)/2 + 2\eps_2} = \rho^{(n-\sigma)/2 + 5\eps_2}. 
\end{align*}
On the other hand, 
\begin{align*}
\int_{U(\G(\bar{Q})')} \mu(\F_{\rho,6}, Y_{\rho,6})(x) \D x 
&\leq 
\rho^{-\sigma-2\eps_2} |U(\G(\bar{Q})')|\\ &\lesssim 
\rho^{-\sigma-2\eps_2} \cdot \rho^{-\eps_2}\rho^{(\sigma-1)/2} \cdot \rho^{(n+1)/2}\\
&= \rho^{(n-\sigma)/2-3\eps_2}. 
\end{align*}
(Note that $\rho^{(n+1)/2}$ is the volume of a $\rho \times \rho^{1/2} \times \rho^{1/2}\1$ grain.) Thus 
\begin{align*}
\int_{U(\G(\bar{Q})'')} \mu(\F_{\rho,6}, Y_{\rho,6})(x) \D x \gtrsim 
\rho^{8\eps_2}
\int_{U(\G(\bar{Q})')} \mu(\F_{\rho,6}, Y_{\rho,6})(x) \D x. 
\end{align*}
This means that by only retaining the mass of $(\F_{\rho,6}, Y_{\rho,6})$ over all $U(\G(\bar{Q})'')$, we obtain a $\gtrsim \rho^{8\eps_2}$ refinement $(\F_{\rho,7}, Y_{\rho,7})$ formed by unions of $(\rho/10n)$-cubes. 

Unwinding definitions, we see that if $\bar{Q}$ is a $\rho^{1/2}$-cube and if $J \subset \HH(\bar{Q})$ is a line segment, then 
\begin{align*}
\EE_\rho\left( J \cap \pi_{\HH(\bar{Q})}(U(\F_{\rho,7}, Y_{\rho,7}) \cap \bar{Q}) \right) \lesssim 
\EE_\rho\left( \pi_{\HH(\bar{Q})}( \mathscr{L}_J \cap U(\F_{\rho,7}, Y_{\rho,7})) \right),
\end{align*}
where $\mathscr{L}_J \coloneqq \mathscr{L} \cap \pi_{\HH(\bar{Q})}^{-1}(J)$, and $\mathscr{L}$ is the line from \eqref{Eq:pivot_line}. In particular, if $J$ has length $\tau$, then the line segment $\mathscr{L}_J$ is contained in a union of $O(1)$ $\tau$-cubes $\Tilde{Q}$, and hence 
\begin{align*}
\EE_\rho\left( J \cap \pi_{\HH(\bar{Q})}(U(\F_{\rho,7}, Y_{\rho,7}) \cap \bar{Q})\right) \lesssim 
\sum_{\Tilde{Q}} 
\EE_{\rho}\big( \pi_{\HH(\bar{Q})}(U(\F_{\rho, 7}, Y_{\rho, 7})  \cap \Tilde{Q}) \big),
\end{align*}
where the sum is taken over the $O(1)$ $\tau$-cubes $\Tilde{Q}$ described above. Since the contribution from each of these $\tau$-cubes is $\lesssim \delta^{-\eta} \rho^{-\eps_2}(\tau/\rho)^{1-\sigma}$ by \eqref{Eq:tau_counting}, we conclude that 
\begin{align}\label{Eq:grains_single_rho_tau'}
\EE_\rho\left( J \cap \pi_{\HH(\bar{Q})}(U(\F_{\rho,7}, Y_{\rho,7}) \cap \bar{Q})\right) 
\lesssim 
\delta^{-\eta} \rho^{-\eps_2}(\tau/\rho)^{1-\sigma}.
\end{align}

Finally, by tracing back all the refinement procedures, $(\F_{\rho,7}, Y_{\rho,7})$ is a $\gtrapprox \rho^{8\eps_2}$ refinement of $(\F_{\rho,1}, Y_{\rho,1})$ (formed by unions of $(\rho/10n)$-cubes), and so induces a $\gtrapprox \rho^{8\eps_2}$ refinement $(\F',Y')$ of $(\F_1,Y_1)$ (thus of $(\F,Y)$). 
Since $U(\F',Y') \subset U(\F_{\rho,7}, Y_{\rho,7})$ and $\delta \geq \rho^{1/\eps}$, by taking $\eta$ and $\eps_2$ sufficiently small depending on $\eps$, we get from \eqref{Eq:grains_single_rho_tau'} that \eqref{Eq:grains_single_rho_tau} holds true. The smallness of $\eta$ and $\eps_2$ would also guarantee that $(\F',Y')$ is $\eps$-extremal. And so the proof is completed.
\end{proof}

To prove Proposition~\ref{Prop:local_grains}, the remaining issues are to upgrade Lemma~\ref{Lem:grains_single_rho_tau} to a version that handles multiple scales $\rho$ and $\tau$ simultaneously. For this purpose, there is nothing special about the geometry of flakes, and the iteration steps are exactly the same as those in \cite[Lemma~4.11~and~4.12]{WZ26sticky}. To avoid unnecessary repetitions, we only state the corresponding results for flakes. 

Firstly, for a fixed $\rho$, Lemma~\ref{Lem:grains_single_rho} records what happens when we iterate Lemma~\ref{Lem:grains_single_rho_tau} 
for many different scales $\tau$.
\begin{lem}\label{Lem:grains_single_rho}
For all $\eps>0$, there exist $\eta, \delta_0 >0$ so that the following holds for all $\delta \in (0, \delta_0]$. Let $(\F,Y)$ be an $\eta$-extremal set of $\delta$-flakes with $\delta^{-\eta}$-Lipschitz plane map $\HH$. Let $\rho \in [\delta^{1-\eps}, \delta^\eps]$. For each $\rho^{1/2}$-cube $\bar{Q}$ intersecting $U(\F,Y)$, define $\HH(\bar{Q}) \coloneqq \HH(y_{\bar{Q}})$, where $y_{\bar{Q}}$ is any fixed point in $\bar{Q} \cap U(\F,Y)$. 

Then there exists an $\eps$-extremal refinement $(\F',Y')$, so that for each $\rho^{1/2}$-cube $\bar{Q}$ intersecting $U(\F,Y)$ and each line segment $J \subset \HH(\bar{Q})$ of length $\tau \geq \rho$, we have
\begin{align*}
\EE_\rho\left( J \cap \pi_{\HH(\bar{Q})}(U(\F',Y') \cap \bar{Q})\right) \leq \rho^{-\eps}(\tau/\rho)^{1-\sigma}.
\end{align*}
\end{lem}

Next, Lemma~\ref{Lem:grains_multi_scale} records what happens when we iterate Lemma~\ref{Lem:grains_single_rho} for many different scales $\rho$.

\begin{lem}\label{Lem:grains_multi_scale}
For all $\eps>0$, there exist $\eta, \delta_0 >0$ so that the following holds for all $\delta \in (0, \delta_0]$. Let $(\F,Y)$ be an $\eta$-extremal set of $\delta$-flakes with $\delta^{-\eta}$-Lipschitz plane map $\HH$.

Then there exists an $\eps$-extremal refinement $(\F',Y')$, so that for all $x \in U(\F,Y)$, all $\rho \geq \delta$, and all line segments $J \subset \HH(x)$ of length $\tau \geq \rho$, we have 
\begin{align*}
\EE_\rho\left( J \cap \pi_{\HH(x)}(U(\F',Y') \cap B(x, \rho^{1/2})) \right) \leq \delta^{-\eps}(\tau/\rho)^{1-\sigma}.
\end{align*}
In particular, $\pi_{\HH(x)}(U(\F',Y') \cap B(x, \rho^{1/2}))$ is a $(\rho, 1-\sigma, \delta^{-\eps})_1$-ADset.
\end{lem}

We are now ready to prove Proposition~\ref{Prop:local_grains}. 

\begin{proof}[Proof of Proposition~\ref{Prop:local_grains}]
Combining Lemma~\ref{Lem:exist_eps_Lips} with Lemma~\ref{Lem:grains_multi_scale} yields a set of $\delta$-flakes $(\F,Y)$ satisfying all the requirements, except that the plane map $\HH$ for it is only $\delta^{-\eta}$-Lipschitz (for some $0 < \eta \ll \eps$). This can be fixed via a mild rescaling similar to \cite[Lemma~3.5]{WZ26sticky}. The arguments are straightforward (though somewhat tedious), and so are omitted.
\end{proof}

\section{An unusual slicing argument}\label{Sec:slicing_set_up}

Recall that for the proof of Theorem~\ref{Thm:Main_Sticky}, we start by assuming $\sigma>0$ and aim for some contradiction. So far we have established the planiness and the local grains structure. To continue the proof, we will use certain slicing arguments, and here is an outline of the key steps:
\begin{enumerate}[label=(\rm{\roman*)}]
    \item Choose a suitable $X \in \Gr(n,n-3)$, and cut the $\delta$-flakes up by $(\delta,3)$-flakes orthogonal to $X$ in $\R^n$. We get a bunch of $\delta$-tubes in each slice ($(\delta,3)$-flake). 

    \item Use a projection method and $(n-3)$-plane Kakeya in $\R^{n-1}$ (induction hypothesis: Theorem~\ref{Thm:Main} when $n$ is replaced with $n-1$) to show that the $\delta$-tubes in each slice is essentially distinct.

    \item Apply the factoring lemma (Lemma~\ref{Lem:factoring}) in each slice to get a bunch of convex sets of dimensions $\delta\one \times a \times b \times 1$ ($\delta \leq a \leq b \leq 1$) that factors the $\delta$-tubes. The convex sets themselves are Katz--Tao, while the $\delta$-tubes in each convex set is Frostman. 

    \item If $b \approx \delta$, then the $\delta$-tubes in each slice are Katz--Tao, and so are essentially disjoint by $\R^3$ Kakeya for tubes (Theorem~\ref{Thm:R3_Kakeya}). As a result, the original $\delta$-flakes are essentially disjoint, contradicting the extremal property. Therefore, we may assume that $b \gg \delta$.
    
    \item 
    By $\R^3$ Kakeya for planks (\cite[Lemma~6.4]{GWZ26streamlined}, a by-product of $\R^3$ Kakeya for tubes), in each slice, the union of $\delta$-tubes must fill out all the convex sets from Step~(iii). This can be viewed as an initial step that grows grains.

    \item If $a \gg \delta$, then by combining the growth in all slices, each $\delta$-flake must fill out a much thicker $a$-flake containing it, which allows us to go efficiently from scale $\delta$ to scale $a$ (by self-similarity) and gain an exponent, contradicting the extremal property.

    \item Otherwise, we are in the case when $a \approx \delta$ and $b \gg \delta$. We know that in each slice, the $\delta$-tubes have grown into much wider $\delta\1 \times b \times 1$ grains, which we call \textit{medium grains}\footnote{We choose the terminology to distinguish these grains from \textit{local grains} and \textit{global grains} in \cite{WZ26sticky}, and to indicate that they are intermediate objects subject to further manipulations.}. However, since the direction of such anisotropic growth in $\R^n$ may not be uniform over all slices, we cannot directly grow each $\delta$-flake into a much wider $\delta \times b \times \1$ plank. To overcome this difficulty, we need to play with both medium grains and local grains, the compatibility of which provides us with some additional rigidity. More precisely, the medium grains can be used to regularize the orientations of local grains, i.e., the plane map is kind of concentrated when we walk along each medium grain.

    \item 
    The strategy in Step~(vii) is good but still not enough, as we need to regularize the behaviour for different medium grains. Fortunately, by choosing $(n-2)$ ``transversal'' $X \in \Gr(n,n-3)$ and applying all the previous discussions with respect to each of them, we can make sure that the plane map is kind of concentrated over each $\delta$-flake. 
    
    \item 
    By the concentration property of the plane map obtained in Step~(viii), we can combine the local grains along each $\delta$-flake to grow it into a much wider $\delta \times b' \times \1$ plank for some $b' \gg \delta$ (although $b'$ may be much smaller than $b$).

    \item By self-similarity, the $\delta$-flakes in each $b'$-flake (from the $b'$-cover) are Frostman, and so the corresponding $\delta \times b' \times \1$ planks are also Frostman, thus filling out the $b'$-flake via a local C\'ordoba-type $L^2$ argument. This allows us to go efficiently from scale $\delta$ to scale $b'$ (still by self-similarity) and gain an exponent, contradicting the extremal property. The whole proof is completed.
\end{enumerate}

In a nutshell, we use the known $\R^3$ Kakeya in many scenarios (projection when $n=4$, slicing, growing medium grains), and by combining the grains structure with the slicing method in an organic way, we are able to prove the $(n-2)$-plane sticky Kakeya theorem without resorting to advanced projection theory and additive combinatorics. More precisely, these ingredients (that distinguish $\R$ and $\C$) are implicitly used whenever we apply $\R^3$ Kakeya, but we need nothing more for $(n-2)$-plane Kakeya in $\R^n$ than for line Kakeya in $\R^3$. 

The above outline involves many informal reductions, e.g., often we can only expect a property to hold for ``most'' objects but not ``all'' of them. And our main task will be to flesh out the outline in a rigorous way by the language of shadings. We will finish Steps~(i)--(vi) in this section, and leave Steps~(vii)--(x) to the next section. {\color{black} The planiness and the local grains structure are only used in the last four steps (Section~\ref{Sec:medium_grains_end}), while the self-similarity property is only used in Step~(vi)~and~(x).}

We now begin this long proof. It is important to keep in mind that $(\F,Y)$ is fixed and is $\eta$-extremal. $\eta$ will be reserved as the extremal parameter throughout the proof and  {\color{black} it will be chosen to be  sufficiently small later on.} First of all, by Proposition~\ref{Prop:local_grains}, for any $\eta, \delta_0 > 0$, there exist $\delta \in (0,\delta_0]$ and an $\eta$-extremal set of $\delta$-flakes $(\F,Y)$, such that for all $\rho \in [\delta,1]$ and $x\in U(\F,Y)$, $\pi_{\HH(x)}(U(\F,Y) \cap B(x,\rho^{1/2}))$ is a $(\rho, 1-\sigma, \delta^{-\eta})_1$-ADset. Here $\HH$ is a $1$-Lipschitz plane map for $(\F,Y)$.

We will also introduce many small parameters, whose relationship can be roughly summarized as follows:
\begin{align}\label{Eq:parameter_relation}
\boxed{\eta \ll \eps_1 \ll \eps_1^{(2)} \ll \cdots \ll \eps_1^{(n-2)} \leq \eta_0 \ll \eps_3 \ll \eta_2 \ll \eta_1 \ll \eps_2 \ll \sigma.}
\end{align}
{\color{black} To help the reader follow the proof, we provide a formal proposition summarizing what we will do at the beginning of each step (except for the last Step~(x)).}

{\bf Step~(i).} The first thing is to choose good direction subspaces to perform slicing. 
More explicitly, we will prove the following result:
\begin{proposition}\label{Prop:step_1}
There exist $D > 0$, a $D$-flake $F_D^*$, an orthonormal basis $\{v_k\}_{k=1}^{n-2}$ of $\dir(F_D^*)$, and a $\gtrsim 1$ refinement $(\F',Y)$ of $(\F,Y)$, such that the following holds. 

For each $F \in \F'$, we have $F \subset F_D^*$, $|Y(F)| \gtrapprox \delta^\eta |F|$, and 
\begin{align}\label{Eq:choice_D}
\angle(\dir(F), \dir(F_D^*)) \leq (100n)^{-1}.
\end{align}

Besides, for each $1 \leq k \leq n-2$, if we let 
\begin{itemize}
    \item the $k$-th direction subspace be 
    \begin{align}\label{Def:k-th_subspace}
        X_k \coloneqq \Span(\{v_j\}_{j=1}^{n-2}\setminus\{v_k\}),
    \end{align}

    \item the lattice of mesh size $\delta$ in $X_k$ be 
    \begin{align}\label{Def:k-th_lattice}
    \Lambda_k \coloneqq \Big\{\sum_{\substack{1 \leq j \leq n-2\\ j \neq k}} m_j \delta v_j: \text{ each } m_j \in [-\delta^{-1}, \delta^{-1}] \cap \Z \Big\},
    \end{align} 

    \item the \textit{adapted $\rho$-neighborhood} $(\rho>0)$ of a set $E \subset \R^n$ be \begin{align}\label{Def:k-th_neighbor}
    \Bar{N}_\rho(E) \coloneqq \Big\{ x + \sum_{k=1}^{n-2} a_k v_k: x \in E, \text{ each } a_k \in [-\rho, \rho) \Big\},
    \end{align}

    \item and the slice of thickness $\delta$ with direction $X_k^\perp$ and offset $\omega \in \Lambda_k$ be     \begin{align}\label{Def:k-th_slices}
    \Pi_\omega^{(k)} \coloneqq \Bar{N}_{\delta/2}(X_k^\perp + \omega), 
    \end{align}
\end{itemize}  
then $F \cap \Pi_\omega^{(k)}$ is comparable to a $\delta$-tube\footnote{In fact, each tube is contained in and comparable to a $\sqrt{n}\delta$-tube. So strictly speaking, one may need to (harmlessly) replace $\delta$ with $\sqrt{n}\delta$ at some places later on. We ignore such technical issues for simplicity.} for all $\omega \in \Lambda_k$ and all $F \in \F'$. 
\end{proposition}

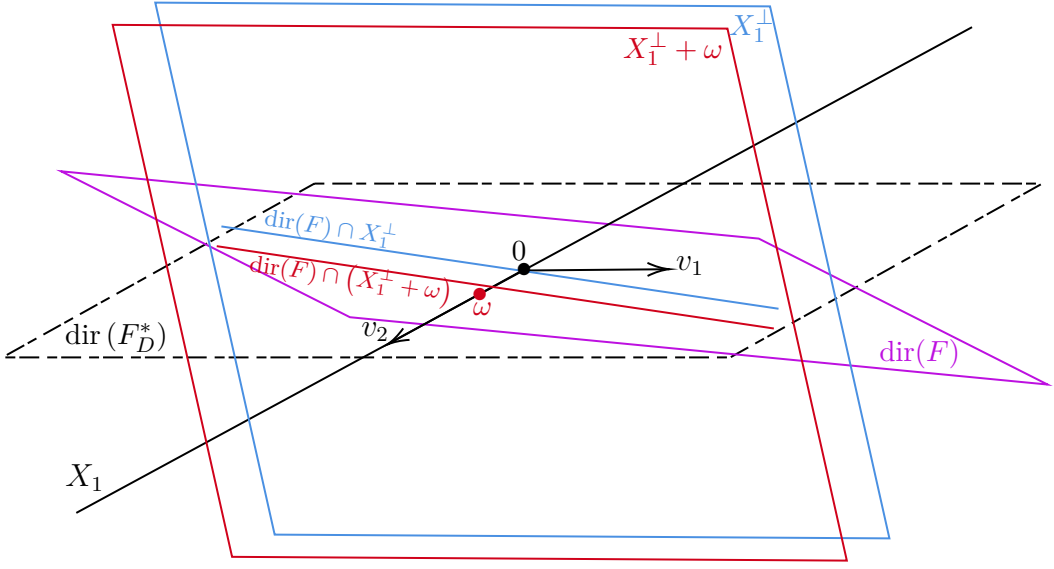
\begin{figure}[htbp]
  \centering
  
    \hspace*{0.5cm}

\tikzset{every picture/.style={line width=0.75pt}} %set default line width to 0.75pt        

\begin{tikzpicture}[x=0.75pt,y=0.75pt,yscale=-1,xscale=1]
%uncomment if require: \path (0,1303); %set diagram left start at 0, and has height of 1303

%Shape: Parallelogram 
\draw  [dash pattern={on 3.75pt off 3pt on 7.5pt off 1.5pt}] (281.5,143.32) -- (646.66,143.32) -- (490.16,230.51) -- (125,230.51) -- cycle ;
%Shape: Parallelogram  
\draw  [color={rgb, 255:red, 189; green, 16; blue, 224 }  ,draw opacity=1 ] (299.33,210.35) -- (649.33,244.12) -- (504.42,170.91) -- (154.42,137.15) -- cycle ;
%Straight Lines 
\draw    (385.83,186.92) -- (319.81,223.12) ;
\draw [shift={(318.06,224.09)}, rotate = 331.26] [color={rgb, 255:red, 0; green, 0; blue, 0 }  ][line width=0.75]    (10.93,-3.29) .. controls (6.95,-1.4) and (3.31,-0.3) .. (0,0) .. controls (3.31,0.3) and (6.95,1.4) .. (10.93,3.29)   ;
%Straight Lines 
\draw    (385.83,186.92) -- (457.55,186.29) ;
\draw [shift={(459.55,186.27)}, rotate = 179.5] [color={rgb, 255:red, 0; green, 0; blue, 0 }  ][line width=0.75]    (10.93,-3.29) .. controls (6.95,-1.4) and (3.31,-0.3) .. (0,0) .. controls (3.31,0.3) and (6.95,1.4) .. (10.93,3.29)   ;
%Straight Lines 
\draw [line width=0.75]    (611.43,64.85) -- (162.06,308.52) ;
%Shape: Parallelogram 
\draw  [color={rgb, 255:red, 74; green, 144; blue, 226 }  ,draw opacity=1 ] (261.57,319.41) -- (569.98,321.3) -- (510.09,54.42) -- (201.68,52.53) -- cycle ;
%Shape: Parallelogram [id:dp9043628715228554] 
\draw  [color={rgb, 255:red, 208; green, 2; blue, 27 }  ,draw opacity=1 ] (240.21,330.56) -- (548.62,332.45) -- (488.73,65.57) -- (180.32,63.69) -- cycle ;
%Straight Lines [id:da667325585856684] 
\draw [color={rgb, 255:red, 74; green, 144; blue, 226 }  ,draw opacity=1 ]   (234.86,164.8) -- (514.35,206.1) ;
%Straight Lines [id:da20184788869446157] 
\draw [color={rgb, 255:red, 208; green, 2; blue, 27 }  ,draw opacity=1 ]   (232.51,174.71) -- (512,216.01) ;

% Text Node
\draw (182.22,220.28) node  [font=\small]  {$\mathrm{dir}\left( F_{D}^{*}\right)$};
% Text Node
\draw (585.27,228.89) node  [font=\small,color={rgb, 255:red, 189; green, 16; blue, 224 }  ,opacity=1 ]  {$\mathrm{dir}( F)$};
% Text Node
\draw (166.38,290.96) node    {$X_{1}$};
% Text Node
\draw (470.01,184.41) node    {$v_{1}$};
% Text Node
\draw (312.08,217.45) node    {$v_{2}$};
% Text Node
\draw (501.01,64.29) node  [font=\small,color={rgb, 255:red, 74; green, 144; blue, 226 }  ,opacity=1 ]  {$X_{1}^{\perp }$};
% Text Node
\draw (384.18,177.29) node  [font=\small]  {$0$};
% Text Node
\draw (386.74,186.69) node    {$\smblkcircle $};
% Text Node
\draw (364.47,199.07) node    {$\textcolor[rgb]{0.82,0.01,0.11}{\smblkcircle }$};
% Text Node
\draw (461.44,75.68) node  [font=\small,color={rgb, 255:red, 208; green, 2; blue, 27 }  ,opacity=1 ]  {$X_{1}^{\perp } +\omega $};
% Text Node
\draw (365.97,206.07) node  [font=\small,color={rgb, 255:red, 208; green, 2; blue, 27 }  ,opacity=1 ]  {$\omega $};
% Text Node
\draw (290.08,165.42) node  [font=\scriptsize,color={rgb, 255:red, 74; green, 144; blue, 226 }  ,opacity=1 ,rotate=-8.45]  {$\mathrm{dir}( F) \cap X_{1}^{\perp }$};
% Text Node
\draw (300,191.5) node  [font=\scriptsize,color={rgb, 255:red, 208; green, 2; blue, 27 }  ,opacity=1 ,rotate=-8.37]  {$\mathrm{dir}( F) \cap \left( X_{1}^{\perp } +\omega \right)$};

\end{tikzpicture}
\caption{A rough sketch of how we do slicing when $n=4$. For simplicity, we focus on the $k=1$ case and only draw one $\dir(F)$. Here $\{v_1, v_2\}$ is an orthonormal basis of $\dir(F_D^*)$, $X_1 = \langle v_2 \rangle$, and $\omega \in \Lambda_1 = \{m\delta v_2 : \, m \in [-\delta^{-1}, \delta^{-1}] \cap \Z \}$. The reader should keep in mind that both $X_1^\perp$ and $X_1^\perp + \omega$ have dimension $3$, although for technical reasons we can only draw them as $2$-dimensional planes.}
\label{fig:slicing}
\end{figure}

Before everything, let us first digest what Proposition~\ref{Prop:step_1} tells us. Intuitively, for each $F \in \F'$, since $X_k \subset \dir(F_D^*)$ is quite close to $\dir(F)$, we should have transversality between $\dir(F)$ and each $X_k^\perp + \omega$. See Figure~\ref{fig:slicing} for an illustration of the special case when $n=4$. This gives the transversality between $F$ and each $\Pi_\omega^{(k)}$, which is compatible with the claim that $F \cap \Pi_\omega^{(k)}$ is a $\delta$-tube. 

\begin{proof}[Proof of Proposition~\ref{Prop:step_1}]
Firstly, let $D > 0$ be a small constant depending only on $n$, so that for any $\delta$-flake $F$ contained in a $D$-flake $F_D$, we have 
$$\angle(\dir(F), \dir(F_D)) \leq (100n)^{-1}.$$ 
This can be done because Lemma~\ref{Lem:flake_contain} implies $\angle(\dir(F), \dir(F_D)) \lesssim D$ whose constant depends only on $n$.

Since $\lambda(\F, Y) \geq \delta^\eta$, by discarding those flakes in $\F$ with low density of shadings, we may assume that $|Y(F)| \gtrapprox \delta^\eta|F|$ for each $F \in \F$.

By Remark~\ref{Rem:major_metric}, the major angle yields a metric on the finite-dimensional manifold $\Gr(n, n-2)$, which allows us to decompose $\Gr(n, n-2)$ into a finitely overlapping union of $O(1)$ small caps of diameter $D/10$. Therefore, we can cover $(\F,Y)$ by a set of $O(1)$ $D$-flakes. By pigeonholing, we can find one of them, say $F_D^*$, so that $\sum_{\F[F_D^*]} |Y(F)| \gtrsim \sum_{\F} |Y(F)|$. Let $\{v_k\}_{k=1}^{n-2}$ be any (fixed) orthonormal basis of $\dir(F_D^*)$. Let $\F' \coloneqq \F[F_D^*]$. It remains to check that $D$, $F_D^*$, $\{v_k\}_{k=1}^{n-2}$, and $(\F', Y)$ satisfy the desired properties.

By the construction of $\F'$ and our choice of $D$, for each $F \in \F'$, we clearly have $F \subset F_D^*$, $|Y(F)| \gtrapprox \delta^\eta |F|$, and \eqref{Eq:choice_D}.

For each $k$, let $X_k$ be defined as in \eqref{Def:k-th_subspace}. For each $F \in \F'$, \eqref{Eq:choice_D} together with Remark~\ref{Rem:major_equiv} gives that $X_k \cap B_1 \subset \dir(F_D^*) \cap B_1 \subset N_{1/10}(\dir(F))$. Therefore, by Lemma~\ref{Lem:slice_transversality} with $V = \dir(F)$ and $W = X_k$, we have $\dim(\dir(F) \cap X_k^\perp) = 1$ and $N_\delta(\dir(F)) \cap N_\delta(X_k^\perp) \subset N_{2\delta}(\dir(F)\cap X_k^\perp)$. This implies that $F \cap \Pi_\omega^{(k)}$ is comparable to a $\delta$-tube, as desired. 
\end{proof}

Let $D$, $F_D^*$, $\{v_k\}_{k=1}^{n-2}$, and $(\F', Y)$ be the output of Proposition~\ref{Prop:step_1}. Let $X_k$, $\Lambda_k$, $\Bar{N}_\rho(\cdot)$, and $\{\Pi_\omega^{(k)}\}_{\omega\in\Lambda_k}$ be defined as in \eqref{Def:k-th_subspace}, \eqref{Def:k-th_lattice}, \eqref{Def:k-th_slices}, and \eqref{Def:k-th_neighbor}, respectively. 

In particular, inside each $\Pi_\omega^{(k)}$, we get a bunch of $\delta$-tubes 
\begin{align*}
    \T^{\omega, (k)} \coloneqq \{ F \cap \Pi_\omega^{(k)}: F \in \F'\},
\end{align*}
whose shading $Y^{\omega, (k)}$ is defined by 
\begin{align*}
    Y^{\omega, (k)}(F \cap \Pi_\omega^{(k)}) \coloneqq Y(F) \cap \Pi_\omega^{(k)}.
\end{align*} 
Thus each $F \in \F'$ is decomposed into $\sim \delta^{-(n-3)}$ parallel $\delta$-tubes pointing in $\dir(F) \cap X_k^\perp$, which we denote by \begin{align*}
    \T^{(k)}(F) \coloneqq \{ F \cap \Pi_\omega^{(k)}: \omega \in \Lambda_k\}.
\end{align*}
Moreover, note that 
$$U(\F', Y) = \bigsqcup_{\omega\in\Lambda_k} U(\T^{\omega, (k)}, Y^{\omega, (k)}),$$
and 
$$\sum_{F \in \F'} |Y(F)| = \sum_{\omega\in\Lambda_k} \sum_{T \in \T^{\omega, (k)}} |Y^{\omega, (k)}(T)|.$$ 
Thus the study of $\delta$-flakes is reduced to that of $\delta$-tubes in each slice, which is essentially a $3$-dimensional problem. 

This finishes Step~(i).

In Step~(ii)--(vii), let us only focus on $X_1$, and omit all the superscripts ``$(1)$'' for simplicity. For Step~(viii) in the next section, we will further consider all the other $X_k$ and need to use the full notation.

{\bf Step~(ii).} This is the only place in the whole proof that we use the induction hypothesis (Theorem~\ref{Thm:Main} holds when $n$ is replaced with $n-1$).

The caveat in Step (i) is that $\T^\omega$ may not be essentially distinct. This could be a problem (at least in some slices) if we want to show the largeness of $U(\T^\omega, Y^\omega)$. In general, we need to pass to a refinement to forbid such a bad scenario. 
More explicitly, we will prove the following result:
\begin{proposition}\label{Prop:step_2}
Suppose that $0 < \eta \ll \eps_1$, then there exists $\delta_0 > 0$ so that the following holds for all $\delta \in (0,\delta_0]$. 

There exist $\lambda \gtrapprox \delta^\eta$ and a $\gtrapprox \delta^{2\eps_1 + \eta}$ refinement $(\F_1,Y_1)$ of $(\F, Y)$, such that 
\begin{itemize}    
    \item For each $\omega \in \Lambda_1$, \begin{align}\label{Def:T^w,1}
        \T^{\omega,1} \coloneqq \{F \cap \Pi_\omega: F \in \F_1, Y_1(F) \cap \Pi_\omega \neq \varnothing \}
    \end{align} 
    is a set of essentially distinct $\delta$-tubes. 

    \item For each $\omega \in \Lambda_1$, the shading $Y^{\omega,1}$ on $\T^{\omega,1}$ defined by 
    \begin{align}\label{Def:Y^w,1}
        Y^{\omega,1}(F \cap \Pi_\omega) \coloneqq Y_1(F) \cap \Pi_\omega
    \end{align}
    has the property that $|Y^{\omega,1}(T)| \sim \lambda |T|$ for all $T \in \T^{\omega,1}$.
\end{itemize}  
\end{proposition}
\begin{proof}
Recall that $\dir(T) = \dir(F) \cap X_1^\perp$ for all $T \in \T(F)$. This leads us to consider the well-defined map 
\begin{align*}
    \Phi: \F' &\rightarrow  \Gr(X_1^\perp, 1)\\
    F & \mapsto \dir(F) \cap X_1^\perp.
\end{align*}
Roughly speaking, for $F_1, F_2 \in \F'$, if $\Phi(F_1)$ and $\Phi(F_2)$ are $\delta$-separated, then any $T_1 \in \T(F_1)$ and $T_2 \in \T(F_2)$ must be essentially distinct, as they have essentially distinct directions. Therefore, we only need to handle those $F$ with essentially the same image under $\Phi$.

Choose a maximal set of $\delta$-separated elements $\{\ell_i\}_{i=1}^M$ on $\Gr(X_1^\perp, 1)$ (major angle yields the metric), then $\Gr(X_1^\perp, 1) \subset \cup_{i=1}^M B(\ell_i, 2\delta)$ (ball in the metric space), and so $\F = \cup_{i=1}^M \Phi^{-1}(B(\ell_i, 2\delta))$. Divide $\{\ell_i\}_{i=1}^M$ into $J \lesssim 1$ groups $\{\ell_{i,j}\}_{i=1}^{M_j}$ ($j = 1,\dots, J$) so that the points in each group are $10\delta$-separated. By pigeonholing, there exists some $j_0$ such that 
$$\F_1 \coloneqq \bigsqcup_{i=1}^{M_{j_0}} \Phi^{-1}(B(\ell_{i,j_0}, 2\delta))$$ 
satisfies $\sum_{\F_1}|Y(F)| \gtrsim \sum_{\F} |Y(F)|$. For each $i$, we will refine the shading $Y$ on $$\F_{1,i} \coloneqq \Phi^{-1}(B(\ell_{i,j_0}, 2\delta))$$ 
to make sure that the resulting $\delta$-tubes with nonempty shadings are essentially distinct. Without loss of generality, we may assume that for each $1\leq i \leq M_{j_0}$, there exists $F \in \F_{1,i}$ such that $Y(F) \neq \varnothing$, which ensures $|Y(F)| \gtrapprox \delta^\eta |F|$ as $(\F', Y)$ comes from Proposition~\ref{Prop:step_1}. Combining the refinements for all $(\F_{1,i}, Y)$ will yield the desired refinement for $(\F, Y)$. 

Fix an arbitrary $i$. Let $\pi_{i,j_0}$ be the projection along the direction $\ell_{i,j_0}$ to the $(n-1)$-plane $\ell_{i,j_0}^\perp$. For each $F \in \F_{1,i}$, since $\Phi(F) \in B(\ell_{i,j_0}, 2\delta)$, we know that each $T \in \T(F)$ makes an angle $\leq 2\delta$ with $\ell_{i,j_0}$, and so $\pi_{i,j_0}(F)$ is contained in (and comparable to) a $(4\delta,n-3)$-flake $\Tilde{F}(F)$ in $\ell_{i,j_0}^\perp$. Let 
$$\Tilde{\F}_{1, i} \coloneqq \{\Tilde{F}(F): F \in \F_{1,i}\}.$$ 
To define a reasonable shading on $\Tilde{\F}_{1, i}$ that is amenable to various manipulations, we need to uniformize the contribution of each fibre. Consider the (measurable) function $f(z) \coloneqq \int_{\pi_{i,j_0}^{-1}(z)} \mu(\F_{1,i}, Y)(x) \D x$ on $\ell_{i,j_0}^\perp$, where the integration is $1$-dimensional. Then by Fubini's theorem, we have 
\begin{align*}
    \int_{\ell_{i,j_0}^\perp} f(z) \D z = \int_{\R^n} \mu(\F_{1,i}, Y)(x) \D x.
\end{align*}
By first discarding the mass over $\{z: f(z) \leq \delta^{100} \}$ and then apply dyadic pigeonholing with respect to the value of $f$, we get a set $A \subset \ell_{i,j_0}^\perp$ so that $f(z)$ is comparable for all $z \in A$, and 
\begin{align*}
    \int_{A} f(z) \D z \gtrapprox \int_{\R^n} \mu(\F_{1,i}, Y)(x) \D x.
\end{align*}
For each $F \in \F_{1,i}$, define $$Y_{1,i}(F) \coloneqq Y(F) \cap \pi_{i,j_0}^{-1}(A),$$ then $(\F_{1,i}, Y_{1,i})$ is a $\gtrapprox 1$ refinement of $(\F_{1,i}, Y)$. Note that $\lambda(\F_{1,i}, Y_{1,i}) \gtrapprox \lambda(\F_{1,i}, Y) \gtrapprox \delta^\eta$. For each $\Tilde{F} = \Tilde{F}(F) \in \Tilde{\F}_{1,i}$, define 
$$\Tilde{Y}_{1,i}(\Tilde{F}) \coloneqq \pi_{i,j_0}(Y_{1,i}(F)).$$ 
Clearly $A = U(\Tilde{\F}_{1,i}, \Tilde{Y}_{1,i})$. Since the length of each fibre is $\lesssim 1$, we have 
\begin{align*}
|\Tilde{Y}_{1,i}(\Tilde{F})| 
&= |\pi_{i,j_0}(Y_{1,i}(F))| \\
&= \int_{\pi_{i,j_0}(Y_{1,i}(F))} 1 \D z \\
&\gtrsim 
\int_{\pi_{i,j_0}(Y_{1,i}(F))} |\pi_{i,j_0}^{-1}(z) \cap Y_{1,i}(F)| \D z \\
&= |Y_{1,i}(F)|, 
\end{align*}
which implies that $$\lambda(\Tilde{\F}_{1,i}, \Tilde{Y}_{1,i}) = \frac{\sum_{\Tilde{\F}_{1,i}} |\Tilde{Y}_{1,i}(\Tilde{F})|}{\sum_{\Tilde{\F}_{1,i}}|\Tilde{F}|}
\gtrsim 
\frac{\sum_{\F_{1,i}}|Y_{1,i}(F)|}{\sum_{\F_{1,i}}|F|} = \lambda(\F_{1,i}, Y_{1,i}) 
\gtrapprox \delta^\eta.$$

By Definition~\ref{Def:strong_sticky}, $\F$ is $\delta^{-\eta}$-sticky implies that $\Delta_{max}(\F) \leq \delta^{-\eta}$. By Proposition~\ref{Prop:inherit}(B), this means $\Delta_{max}(\F_{1,i}) \leq \delta^{-\eta}$. Thus for any convex set $K \subset \ell_{i,j_0}^\perp$, we have 
\begin{align*}
    |\Tilde{\F}_{1,i}[K]| &\leq |\{ F \in \F_{1,i}: \pi_{i,j_0}(F) \subset K \}| \\
    &= |\F_{1,i}[\pi_{i,j_0}^{-1}(K) \cap B_1]| \\
    &\leq \Delta_{max}(\F_{1,i}) |\pi_{i,j_0}^{-1}(K) \cap B_1|/|F|\\
    &\lesssim \delta^{-\eta} |K|/|F|. 
\end{align*}
This estimate immediately yields $$\Delta_{max}(\Tilde{\F}_{1,i}) = \sup_{K \text{ convex}} \frac{|\Tilde{F}| |\Tilde{\F}_{1,i}[K]|}{|K|} \lesssim \frac{\delta^{-\eta} |\Tilde{F}||K|}{|F| |K|} 
\lesssim 
\delta^{-\eta}.$$ 
Therefore, for any $\eps_1>0$, {\color{black} by choosing $\eta \ll \eps_1$} and $\delta_0$ sufficiently small and applying the induction hypothesis (Theorem~\ref{Thm:Main} when $n$ is replaced with $n-1$) to $(\Tilde{\F}_{1,i}, \Tilde{Y}_{1,i})$, we have $|A| = |U(\Tilde{\F}_{1,i}, \Tilde{Y}_{1,i})| \geq \delta^{\eps_1} |\Tilde{\F}_{1,i}| |\Tilde{F}|$. 

Recall that 
\begin{align*}
\int_{A} \mu(\Tilde{\F}_{1,i}, \Tilde{Y}_{1,i})(z) \D z = \int_{\ell_{i,j_0}^\perp} \mu(\Tilde{\F}_{1,i}, \Tilde{Y}_{1,i})(z) \D z = \sum_{\Tilde{\F}_{1,i}} |\Tilde{Y}_{1,i}(\Tilde{F})| \leq |\Tilde{\F}_{1,i}| |\Tilde{F}|,
\end{align*}
and so the high-multiplicity set 
$$A' \coloneqq \{z\in A: \mu(\Tilde{\F}_{1,i}, \Tilde{Y}_{1,i})(z) \geq \delta^{-2\eps_1} \}$$ 
satisfies $\delta^{-2\eps_1}|A'| \leq |\Tilde{\F}_{1,i}| |\Tilde{F}|$. In particular, this means $|A'| \leq \delta^{2\eps_1} |\Tilde{\F}_{1,i}| |\Tilde{F}| \leq \delta^{\eps_1} |A|$. For each $F \in \F_{1,i}$, define $$Y_{1,i}'(F) \coloneqq Y_{1,i}(F) \cap \pi_{i,j_0}^{-1}(A\setminus A').$$ Since $f(z)$ is comparable over $A$, we know that $(\F_{1,i}, Y_{1,i}')$ is a $\gtrsim 1$ refinement of $(\F_{1,i}, Y_{1,i})$. Moreover, by our construction, for any $x \in U(\F_{1,i}, Y_{1,i}')$, we have $$\mu(\F_{1,i}, Y_{1,i}')(x) \leq \delta^{-2\eps_1}.$$ 

Note that $\lambda(\F_{1,i}, Y_{1,i}') \gtrsim \lambda(\F_{1,i}, Y_{1,i}) \gtrapprox \delta^\eta$, where 
\begin{align*}
\lambda(\F_{1,i}, Y_{1,i}') = \frac{\sum_{\F_{1,i}}|Y_{1,i}'(F)|}{\sum_{\F_{1,i}}|F|} = \frac{\sum_{\F_{1,i}}\sum_{T \in \T(F)} |Y_{1,i}'(F) \cap T|}{\sum_{\F_{1,i}}\sum_{T \in \T(F)}|T|}. 
\end{align*}
So by pigeonholing, there exists some constant $\lambda_i$ and $\T(F)' \subset \T(F)$ for each $F \in \F_{1,i}$, such that $|Y_{1,i}'(F) \cap T|/|T| \sim \lambda_i \gtrapprox \delta^\eta$ for any $F \in \F_{1,i}$ and $T \in \T(F)'$, and 
\begin{align*}
\sum_{\F_{1,i}}\sum_{T \in \T(F)'} |Y_{1,i}'(F) \cap T| \gtrapprox 
\sum_{\F_{1,i}}\sum_{T \in \T(F)} |Y_{1,i}'(F) \cap T|.
\end{align*}
In other words, if we define 
$$Y_{1,i}''(F) \coloneqq \bigsqcup_{\T(F)'} Y_{1,i}'(F) \cap T$$ for each $F \in \F_{1,i}$, then $(\F_{1,i}, Y_{1,i}'')$ is a $\gtrapprox 1$ refinement of $(\F_{1,i}, Y_{1,i}')$. Moreover, for each $\omega \in \Lambda_1$ and $F \in \F_{1,i}$, either $Y_{1,i}''(F) \cap \Pi_\omega = \varnothing$, or $|Y_{1,i}''(F) \cap \Pi_\omega|/ |F \cap \Pi_\omega| \sim \lambda_i \gtrapprox \delta^\eta$. 
This naturally induces a subset 
$$\T_i^\omega \coloneqq \{F \cap \Pi_\omega: F \in \F_{1,i}, Y_{1,i}''(F) \cap \Pi_\omega \neq \varnothing\} \subset \T^\omega$$ 
and a shading $Y_{i}^\omega$ on it defined by $$Y_{i}^\omega(F \cap \Pi_\omega) \coloneqq Y_{1,i}''(F) \cap \Pi_\omega.$$ 
Besides, for any $x \in \Pi_\omega$, we have 
$$\mu(\T_i^\omega, Y_i^\omega)(x) = \mu(\F_{1,i}, Y_{1,i}'')(x) \leq \mu(\F_{1,i}, Y_{1,i}')(x) \leq \delta^{-2\eps_1}.$$

Now the key observation is that for each fixed $\omega$, although $\T_{i}^\omega$ may not be essentially distinct in itself, it has nicely bounded repetition: For any $10\delta$-tube $\Bar{T}$ in $\R^n$, there cannot be too many $\T_i^\omega$ in $\Bar{T}$. More precisely, we have 
\begin{align*}
|\T_i^\omega[\Bar{T}]| \cdot \lambda_i |T| 
\sim 
\sum_{T \in \T_i^\omega[\Bar{T}]} |Y_i^\omega(T)| 
= 
\int_{\Bar{T} \cap \Pi_\omega} \mu(\T_i^\omega[\Bar{T}], Y_i^\omega)(x) \D x 
\leq 
\delta^{-2\eps_1} |\Bar{T} \cap \Pi_\omega|,
\end{align*}
which implies 
\begin{align*}
|\T_i^\omega[\Bar{T}]| \lesssim 
\lambda_i^{-1}\delta^{-2\eps_1} |\Bar{T} \cap \Pi_\omega|/|T| \lessapprox 
\delta^{-2\eps_1-\eta}. 
\end{align*}
Therefore, by first passing to a (dominating) $100\delta$-separated cover composed of $10\delta$-tubes and then choosing (at most) one element of $\T_i^\omega$ in each of them, we get an essentially distinct subset $\T_i^{\omega,1} \subset \T_i^\omega$ with $$|\T_i^{\omega,1}| \gtrapprox \delta^{2\eps_1+\eta} |\T_i^\omega|.$$ 
Since $|Y_i^\omega(T)| \sim \lambda_i |T|$ for each $T \in \T_i^\omega$, we know that $(\T_i^{\omega,1}, Y_i^\omega)$ is a $\gtrapprox \delta^{2\eps_1+\eta}$ refinement of $({\T_i^\omega}, Y_i^\omega)$. Let $$\T_i \coloneqq \bigsqcup_\omega \T_i^{\omega,1}$$ 
be the resulting $\delta$-tubes from all slices. If we define 
$$Y_{1,i}'''(F) \coloneqq \bigsqcup_{\T(F)' \cap \T_i} Y_{1,i}'(F) \cap T$$ 
for each $F \in \F_{1,i}$, then 
\begin{align*}
\sum_{\F_{1,i}}|Y_{1,i}'''(F)| = \sum_\omega \sum_{\T_i^{\omega,1}} |Y_i^\omega(T)| \gtrapprox 
\delta^{2\eps_1+\eta} 
\sum_\omega \sum_{{\T_i^\omega}} |Y_i^\omega(T)| 
= 
\delta^{2\eps_1+\eta}
\sum_{\F_{1,i}}|Y_{1,i}''(F)|,
\end{align*}
and so $(\F_{1,i}, Y_{1,i}''')$ is a $\gtrapprox \delta^{2\eps_1+\eta}$ refinement of $(\F_{1,i}, Y_{1,i}'')$ (thus of $(\F_{1,i}, Y)$). 

Finally, for each $\omega \in \Lambda_1$, let $$\T^{\omega,1} \coloneqq \bigsqcup_{i=1}^{M_{j_0}} \T_i^{\omega,1},$$ then it is easy to verify that $\T^{\omega,1}$ is essentially distinct. Define a shading $Y^{\omega,1}$ on it by 
$$Y^{\omega,1}(T) = Y_i^{\omega}(T) \quad \text{for any } T \in \T_i^{\omega,1}.$$ 
Also, define a shading $Y_1$ on $\F_1 = \bigsqcup_{i=1}^{M_{j_0}} \F_{1,i}$ by $$Y_1(F) \coloneqq Y_{1,i}'''(F) \quad \text{for any } F \in \F_{1,i}.$$ 
Then $(\F_1, Y_1)$ is a $\gtrapprox \delta^{2\eps_1+\eta}$ refinement of $(\F_1, Y)$ (thus of $(\F,Y)$), and is ``compatible'' with $\{ (\T^{\omega,1}, Y^{\omega,1}) \}_\omega$ in the sense that \eqref{Def:T^w,1} and \eqref{Def:Y^w,1} are satisfied. Moreover, for any $T \in \T^{\omega,1}$, we have $|Y^{\omega,1}(T)| \gtrapprox \delta^\eta |T|$. So by a harmless pigeonholing step, we may further assume that for some constant $\lambda$, we have 
$$|Y^{\omega,1}(T)|/|T| \sim \lambda \gtrapprox \delta^\eta$$ 
for all $\omega \in \Lambda_1$ and all $T \in \T^{\omega,1}$. It is easy to verify that $\lambda$ and $(\F_1, Y_1)$ satisfy all the desired properties.
\end{proof}

For $0 < \eta \ll \eps_1$, let $\lambda$ and $(\F_1, Y_1)$ be the output of Proposition~\ref{Prop:step_2}, and let $\{(\T^{\omega,1}, Y^{\omega,1})\}_{\omega\in \Lambda_1}$ be defined as in \eqref{Def:T^w,1} and \eqref{Def:Y^w,1}.

This finishes Step~(ii).

{\bf Step~(iii).} 
For each $\omega \in \Lambda_1$, we do not have any a priori structural information about $\T^{\omega,1}$ (except that they are essentially distinct), but we can try to apply the factoring lemma to find some structure. For this purpose, a technical issue is that although the distribution of $\T^{\omega,1}$ is intuitively compressed into a $3$-dimensional space, the shadings may not be evenly distributed along $X_1$, the direction of $\delta$-thickness. To fix this issue, we need a projection argument. Define $$\Tilde{\T}^{\omega,1} \coloneqq \{ \pi_{X_1^\perp}(T): T \in \T^{\omega,1} \}.$$ More explicitly, we will prove the following result:
\begin{proposition}\label{Prop:step_3}
There exist $a, b \in [\delta,1]$, $N \in \N^+$, $\Gamma \subset \Lambda_1$, as well as $\T^{\omega, 2} \subset \T^{\omega, 1}$ and a shading $\Tilde{Y}^{\omega,1}$ on $\tilde{\T}^{\omega,1}$ for each $\omega \in \Gamma$, such that 
\begin{itemize}
    \item $N \gtrapprox \delta^{-2+2\eps_1+3\eta}$, $|\Gamma| \gtrapprox \delta^{-(n-3) + 2\eps_1 + 4\eta}$, and \begin{align}\label{Eq:total_mass_0}
    \sum_{\omega \in \Gamma} \sum_{\T^{\omega,2}}|Y^{\omega,1}(T)|  \gtrapprox \delta^{2\eps_1+3\eta}.
    \end{align}
    
    \item For each $\omega \in \Gamma$, we have $|\T^{\omega, 2}| \sim N$, and $\Tilde{\T}^{\omega,2} \coloneqq \{ \pi_{X_1^\perp}(T): T \in \T^{\omega,2} \} \subset \Tilde{\T}^{\omega,1}$ is factored by a set of convex sets $\tilde{\W}_\omega$ of dimensions $a \times b \times 1$. In particular, for each $\omega \in \Gamma$, there is a partition 
    \begin{align*}    \Tilde{\T}^{\omega,2} 
    =     \bigsqcup_{\tilde{W} \in \Tilde{\W}_\omega} \Tilde{\T}_{\tilde{W}}^{\omega,2}
    \end{align*}   such that $\Tilde{\W}_\omega$ is Katz--Tao, and $\Tilde{\T}_{\tilde{W}}^{\omega,2}$ is Frostman in $\tilde{W}$ for each $\tilde{W} \in \Tilde{\W}_\omega$.

    \item For each $\omega \in \Gamma$ and $\tilde{T} \in \Tilde{\T}^{\omega,1}$, we have $|\Tilde{Y}^{\omega,1}(\Tilde{T})| \geq \delta^{2\eta} |\tilde{T}|$, and \begin{align}\label{Eq:regular_fibre}
    |\pi_{X_1^\perp}^{-1}(z) \cap Y^{\omega,1}(T)| \geq \delta^{(n-3)+2\eta}
    \end{align}
    for any $z \in \Tilde{Y}^{\omega,1}(\Tilde{T})$.
\end{itemize}
\end{proposition}

\begin{proof}
Since $|Y^{\omega, 1}(T)| \sim \lambda |T| \gtrapprox \delta^\eta |T|$ for $T \in \T^{\omega, 1}$, by pigeonholing with respect to fibres of $\pi_{X_1^\perp}$, we immediately know that for each $\tilde{T} \in \Tilde{\T}^{\omega,1}$, there is a subset $\Tilde{Y}^{\omega,1}(\Tilde{T}) \subset \pi_{X_1^\perp}(Y^{\omega,1}(T))$ such that the last bullet point holds.

For each $\omega \in \Lambda_1$, by applying Lemma~\ref{Lem:factoring} to the set of $\delta$-tubes $\Tilde{\T}^{\omega,1}$ in $X_1^\perp$, we get a subset $\Tilde{\T}^{\omega,2} \subset \Tilde{\T}^{\omega,1}$ with $|\Tilde{\T}^{\omega,2}| \gtrapprox |\Tilde{\T}^{\omega,1}|$; a family $\Tilde{\W}_\omega$ of convex sets in $X_1^\perp$, each of which has dimensions $a_\omega \times b_\omega \times 1$ (with $\delta \leq a_\omega \leq b_\omega \leq 1$); and a partition 
\begin{align*}
\Tilde{\T}^{\omega,2} 
= 
\bigsqcup_{\tilde{W} \in \Tilde{\W}_\omega} \Tilde{\T}_{\tilde{W}}^{\omega,2},
\end{align*}
where each set in $\Tilde{\T}_{\tilde{W}}^{\omega,2}$ is contained in $\tilde{W}$. This partition has the following properties. 
\begin{enumerate}[label=\rm{(\roman*)}]
    \item $\Tilde{\W}_\omega$ is Katz--Tao, i.e., $\Delta_{max}(\Tilde{\W}_\omega) \sim 1$.

    \item $\Tilde{\T}_{\tilde{W}}^{\omega,2}$ is Frostman in $\tilde{W}$, i.e., $\Delta(\Tilde{\T}_{\tilde{W}}^{\omega,2}, \tilde{W}) \sim \Delta_{max}(\Tilde{\T}_{\tilde{W}}^{\omega,2})$ for each $\tilde{W} \in \Tilde{\W}_\omega$. Furthermore, $\Delta(\Tilde{\T}_{\tilde{W}}^{\omega,2}, \tilde{W}) \sim \Delta_{max}(\Tilde{\T}^{\omega,2})$.
\end{enumerate} 
Define 
$$\T^{\omega,2} \coloneqq \{ T \in \T^{\omega,1}: \pi_{X_1^\perp}(T) \in \Tilde{\T}^{\omega,2} \}.$$ 
Then for each $T \in \T^{\omega,2}$, we may simply use $\tilde{T}$ to denote $\pi_{X_1^\perp}(T)$. Note that $|{\T}^{\omega,1}| = |\Tilde{\T}^{\omega,1}|$ and $|{\T}^{\omega,2}| = |\Tilde{\T}^{\omega,2}|$.

By pigeonholing with respect to $a_\omega$, $b_\omega$, and $|{\T}^{\omega,2}|$ over $\omega\in\Lambda_1$, we get a subset $\Gamma \subset \Lambda_1$ such that for some constants $a$, $b$ and $N$ (with $\delta \leq a \leq b \leq 1$), we have $a_\omega \sim a$, $b_\omega \sim b$, $|{\T}^{\omega,2}| \sim N$ for all $\omega \in \Gamma$; also, 
\begin{align}\label{Eq:total_mass'}
\sum_{\omega \in \Gamma} \sum_{\T^{\omega,2}}|Y^{\omega,1}(T)| \gtrapprox 
\sum_{\omega} \sum_{\T^{\omega,2}}|Y^{\omega,1}(T)| \gtrapprox 
\sum_{\omega} \sum_{\T^{\omega,1}}|Y^{\omega,1}(T)|,
\end{align}
where the first inequality is what we require $\Gamma$ to satisfy, and the second inequality is by $|{\T}^{\omega,2}| \gtrapprox |{\T}^{\omega,1}|$ and $|Y^{\omega,1}(T)| \sim \lambda|T|$ for all $T$. Thus the second bullet point holds.

Since $(\F_1, Y_1)$ is a $\gtrapprox \delta^{2\eps_1 + \eta}$ refinement of $(\F, Y)$, and $\sum_{\F}|Y(F)| = \lambda(\F, Y) |F| |\F| \geq \delta^{2\eta}$, we also have 
\begin{align}\label{Eq:total_mass}
    \sum_\omega \sum_{\T^{\omega,1}}|Y^{\omega,1}(T)| = \sum_{\F_1}|Y_1(F)| \gtrapprox \delta^{2\eps_1+\eta}  \sum_{\F}|Y(F)| \geq \delta^{2\eps_1+3\eta}.
    \end{align}
Combining \eqref{Eq:total_mass'} and \eqref{Eq:total_mass}, we get \eqref{Eq:total_mass_0}.

Finally, observe that \eqref{Eq:total_mass'} can be rewritten as $|\Gamma| N \lambda |T| \gtrapprox \delta^{2\eps_1+3\eta}$, where $\lambda \leq 1$ and $|T|\sim \delta^{n-1}$. As $|\Gamma| \lesssim \delta^{-(n-3)}$, we must have $N \gtrapprox \delta^{-2+2\eps_1+3\eta}$. On the other hand, as $N \leq |\T^{\omega,1}| \leq |\T^\omega| \leq |\F'| \leq |\F| \leq \delta^{-2-\eta}$ ($\F$ is $\eta$-extremal), we also have $|\Gamma| \gtrapprox \delta^{-(n-3) + 2\eps_1 + 4\eta}$. So the first bullet point holds.
\end{proof}

For $0 < \eta \ll \eps_1$, let $a, b, N, \Gamma, \{\T^{\omega,2}\}_{\omega\in\Gamma}, \{\tilde{Y}^{\omega,1}\}_{\omega\in\Gamma}$ be the output of Proposition~\ref{Prop:step_3}, and let $\{\tilde{\T}^{\omega,2}\}_{\omega\in\Gamma}, \{\tilde{\W}_\omega\}_{\omega\in\Gamma}, \{\Tilde{\T}_{\tilde{W}}^{\omega,2}\}_{\omega\in\Gamma, \tilde{W} \in \tilde{\W}_\omega}$ be determined as in the second bullet point. Roughly speaking, the first bullet point means that in ``most'' of the slices, there are ``around'' $\delta^{-2}$ tubes.

This finishes Step~(iii).

{\bf Step~(iv).} Now we start to analyze the outputs obtained from the factoring. The first observation is that $b$ cannot be too close to $\delta$. More explicitly, we will prove the following result:
\begin{proposition}\label{Prop:step_4}
Suppose $0 < \eta \ll \eps_1, \eta_1 \ll \sigma$, then there exists $\delta_0 > 0$ so that $b \in [\delta^{1-\eta_1}, 1]$ holds for all $\delta \in (0, \delta_0]$.
\end{proposition}
\begin{proof}
{\color{black} Let $\eta_1, \eps_2 > 0$ be small parameters to be determined later on.} 
For each $\omega \in \Gamma$ and $\tilde{W} \in \Tilde{\W}_\omega$, we know that $\tilde{W}$ has dimensions $a \times b \times 1$, and $\Tilde{\W}_\omega$ is Katz--Tao. Suppose $b \leq \delta^{1-\eta_1}$, then intuitively $\Tilde{\T}^{\omega,2}$ is not far from $\Tilde{\W}_\omega$, which should also be Katz--Tao with small error. To make this point rigorous, recall that $\Tilde{\T}^{\omega,2}$ is essentially distinct, and so for each $\tilde{W} \in \Tilde{\W}_\omega$, we have $$|\Tilde{\T}^{\omega,2}[\tilde{W}]| \lesssim (|\tilde{W}|/|\tilde{T}|)^2 \lesssim (b/\delta)^4 \leq \delta^{-4\eta_1}.$$ This implies that for any convex set $\tilde{K}$ of dimensions $a' \times b' \times 1$ (with $\delta \leq a' \leq b' \leq 1$) in $X_1^\perp$, we have 
\begin{align*}
|\Tilde{\T}^{\omega,2}[\tilde{K}]| 
&\leq |\Tilde{\W}_\omega[N_b(\tilde{K})]| 
\sup_{\tilde{W} \in \Tilde{\W}_\omega} 
|\Tilde{\T}^{\omega,2}[\tilde{W}]| \\
&\lesssim 
\Delta_{max}(\Tilde{\W}_\omega) \frac{|N_b(\tilde{K})|}{|\tilde{W}|} \delta^{-4\eta_1} \\
&\lesssim 
\frac{(b/\delta)^2|\tilde{K}|}{\delta^2} \delta^{-4\eta_1} \\
&\lesssim 
\delta^{-6\eta_1} \frac{|\tilde{K}|}{|\tilde{T}|},
\end{align*}
and so $\Delta_{max}(\Tilde{\T}^{\omega,2}) \lesssim \delta^{-6\eta_1}$. Since $\lambda(\Tilde{\T}^{\omega,2}, \tilde{Y}^{\omega,1}) \geq \delta^{2\eta}$, by {\color{black} choosing $\max\{\eta, \eta_1\} \ll \eps_2$} and applying Theorem~\ref{Thm:R3_Kakeya} with $\eps_2$ in place of $\eps$, we have $$|U(\Tilde{\T}^{\omega,2}, \tilde{Y}^{\omega,1})| \geq \delta^{\eps_2} |\Tilde{\T}^{\omega,2}||\tilde{T}| \sim \delta^{\eps_2} N \delta^2 \gtrapprox \delta^{2\eps_1 + \eps_2 + 3\eta}.$$ 
By \eqref{Eq:regular_fibre}, this implies that 
\begin{align*}
|U(\T^{\omega,1}, Y^{\omega,1})| &\geq |U(\T^{\omega,2}, Y^{\omega,1})|\\
& \geq \delta^{(n-3)+2\eta} |U(\Tilde{\T}^{\omega,2}, \tilde{Y}^{\omega,1})|\\
& \gtrapprox  
\delta^{(n-3)+2\eps_1+\eps_2+5\eta}.
\end{align*}
Therefore, by summing over all $\omega \in \Gamma$, we get 
\begin{align*}
|U(\F_1,Y_1)| = \sum_{\omega} |U(\T^{\omega,1}, Y^{\omega,1})| \geq \sum_{\omega \in \Gamma} |U(\T^{\omega,1}, Y^{\omega,1})| \gtrapprox \delta^{(n-3)+2\eps_1+\eps_2 + 5\eta}|\Gamma| \gtrapprox \delta^{4\eps_1 + \eps_2 + 9\eta}.
\end{align*}
On the other hand, since $(\F, Y)$ is $\eta$-extremal, we have 
\begin{align*}
|U(\F_1,Y_1)| \leq |U(\F,Y)| \leq \delta^{\sigma - \eta}.
\end{align*}
As $\sigma >0$ is fixed, by first {\color{black} choosing $\max\{ \eta,\eps_1,\eps_2 \} \ll \sigma$} and then letting $\delta_0$ be small enough, we obtain a contradiction between the upper and lower bounds. Therefore, we must have $b \geq \delta^{1-\eta_1}$.
\end{proof}

In other words, for $0 < \eta \ll \eps_1, \eta_1 \ll \sigma$, by Proposition~\ref{Prop:step_4}, we can always assume $b \geq \delta^{1-\eta_1}$, where $b$ is from Proposition~\ref{Prop:step_3}.

This finishes Step~(iv).

{\bf Step~(v).} To proceed, we grow medium grains by using Kakeya for planks in $\R^3$. For technical reasons, we also want certain regularity of the shadings of these grains. More explicitly, we will prove the following result: 
\begin{proposition}\label{Prop:step_5}
Suppose $0 < \eta \ll \min\{\eps_1, \eta_1\}$ and $\max\{\eps_1, \eta_1\} \ll \sigma$, then there exists $\delta_0 > 0$ so that the following holds for all $\delta \in (0,\delta_0]$. 

There exist a refinement $(\T^{\omega, 3}, Y^{\omega, 3})$ of $(\T^{\omega, 2}, Y^{\omega, 1})$ and a set of convex sets $\W_\omega$ of dimensions $\delta\one \times a \times b \times 1$ in $\Pi_\omega$ for each $\omega \in \Gamma$, such that 
\begin{itemize}
    \item For each $\omega \in \Gamma$ and $T \in \T^{\omega,3}$, we have \begin{align}\label{Eq:single_tube_density}
    |Y^{\omega,3}(T)| \gtrapprox  \delta^{4\eta}|T|,
    \end{align}
    and there exists some $W \in \W_\omega$ with $T \subset W$.

    \item The total mass satisfies \begin{align}\label{Eq:total_mass''}
    \sum_{\omega\in\Gamma} \sum_{\T^{\omega,3}} |Y^{\omega,3}(T)| 
    \gtrapprox \delta^{2\eps_1 + 7\eta}.
    \end{align}

    \item For each $\omega \in \Gamma$ and $W \in \W_\omega$, we have not only \begin{align}\label{Eq:grains_in_slices}
    |U(\T^{\omega,2}, Y^{\omega,1}) \cap W| \geq \delta^{\eps_1 + 2\eta} |W|,
    \end{align}
    but also \begin{align}\label{Eq:regular_grains_in_slices}
    |U(\T^{\omega,2}, Y^{\omega,1}) \cap W \cap B(x,r)| 
    \gtrapprox 
    \delta^{\eps_1 +2\eta} |W \cap B(x,r)|
    \end{align}
    for each $x \in U(\T^{\omega,3}[W], Y^{\omega,3})$ and each $r \in [\delta,1]$. 
\end{itemize}
\end{proposition}
\begin{proof}
For each $\omega \in \Gamma$ and $\tilde{W} \in \Tilde{\W}_\omega$, we know that $\tilde{W}$ has dimensions $a \times b \times 1$, and $\Tilde{\T}_{\tilde{W}}^{\omega,2}$ is Frostman in $\tilde{W}$. Let $(\Tilde{\T}_{\tilde{W}}^{\omega,2}, \Tilde{Y}_0^{\omega,1})$ be the constant multiplicity refinement of $(\Tilde{\T}_{\tilde{W}}^{\omega,2}, \Tilde{Y}^{\omega,1})$. Under a unit rescaling $\phi_{\tilde{W}}$ sending $\tilde{W}$ to a $X_1^\perp \cap B_1$, $(\Tilde{\T}_{\tilde{W}}^{\omega,2}, \Tilde{Y}_0^{\omega,1})$ becomes a set of $\delta/b \times \delta/a \times 1$ planks $\phi_{\tilde{W}} (\Tilde{\T}_{\tilde{W}}^{\omega,2}, \Tilde{Y}_0^{\omega,1})$ that is Frostman and has averaged density $\gtrapprox \delta^{2\eta}$. Recall that in proving the $\R^3$ Kakeya for tubes, Wang and Zahl \cite[Section~5]{WZ25volume} also obtained corresponding results for planks in $\R^3$ as a by-product. In other words, we know that the assertion $K_F(\beta)$ in \cite[Definition~3.5]{GWZ26streamlined} holds for all $\beta > 0$, and so by \cite[Lemma~6.4]{GWZ26streamlined}, we have the following result:
\begin{lem}\label{Lem:R3_plank_Frostman}
For every $\eps>0$, there exist $\eta,a_0>0$ such that
the following holds whenever $0<a\leq a_0$ and $a\leq b\leq 1$.
Let $\LL$ be a set of $a\times b\times 1$ planks in
$B_1\subset\mathbb R^3$ that is Frostman, and let $Y$ be a shading
on $\LL$ with $\lambda(\LL,Y)\geq a^\eta$.
Then $|U(\LL,Y)|\geq a^\eps$.
\end{lem}

\begin{rem}\label{Rem:a_0_b_0}
\label{rem:small-a}
\cite[Lemma~6.4]{GWZ26streamlined} initially gives Lemma \ref{Lem:R3_plank_Frostman} in the range $0<a\leq b\leq b_*$ for suitable $\eta_*,b_*>0$ depending on $\eps$. However, only $a$ needs to be sufficiently small. The point is that when $b$ is bounded below, an $a\times b\times 1$ plank has only one thin direction. This regime is controlled by the standard C\'ordoba-type $L^2$ estimates for slabs, which estimate overlap using pairwise slab intersections (see Lemma~\ref{Lem:slab_2}).

Indeed, choose $\eta\le \min\left\{\eta_*,\frac{\eps}{2C_3}\right\}$, where $C_3$ is the constant in Lemma~\ref{Lem:slab_2}. If $b\leq b_*$, then \cite[Lemma~6.4]{GWZ26streamlined} applies directly. If $b>b_*$, choose $a_0$ sufficiently small so that $a_0^\eta\leq b_*$. After decreasing $a_0$ further, the Frostman hypothesis also gives
$C_F(\LL,B_1)\leq a^{-\eta}$ for $a\leq a_0$.
We then have $a^\eta\leq b$, so Lemma~\ref{Lem:slab_2}, with
$(\delta,b_1,b_2)=(a,b,1)$ gives $|U(\LL,Y)| \gtrsim a^{C_3\eta} \geq a^{\eps/2}$. We can decrease $a_0$ to absorb the implicit constant and get $|U(\LL,Y)|\geq a^\eps$.

In the application immediately below, the rescaled planks have
dimensions $(\delta/b)\times(\delta/a)\times 1$.
Thus the parameter called $a$ in the preceding lemma
corresponds to $\delta/b$ here, while its parameter $b$
corresponds to $\delta/a$.
Since $\delta\leq a\leq b$, we have $\frac{\delta}{b}\leq\frac{\delta}{a}\leq 1.$ Moreover, $b\geq\delta^{1-\eta_1}$ implies
$\delta/b\leq\delta^{\eta_1}\leq a_0$ for sufficiently small
$\delta$. This verifies the required smallness of the shortest
rescaled side.
\end{rem}

To apply Lemma~\ref{Lem:R3_plank_Frostman} to the planks $\phi_{\tilde{W}} (\Tilde{\T}_{\tilde{W}}^{\omega,2}, \Tilde{Y}_0^{\omega,1})$, note that since $b \geq \delta^{1-\eta_1}$, we have $\lambda(\phi_{\tilde{W}} (\Tilde{\T}_{\tilde{W}}^{\omega,2}, \Tilde{Y}_0^{\omega,1})) \gtrapprox \delta^{2\eta} \geq (\delta/b)^{2\eta/\eta_1}$. 
By {\color{black} choosing $\eta \ll \min\{ \eta_1, \eps_1 \}$}, Lemma~\ref{Lem:R3_plank_Frostman} yields $|U(\phi_{\tilde{W}} (\Tilde{\T}_{\tilde{W}}^{\omega,2}, \Tilde{Y}_0^{\omega,1}))| \geq (\delta/b)^{\eps_1}$. Undoing the rescaling, we get $$|U(\Tilde{\T}_{\tilde{W}}^{\omega,2}, \Tilde{Y}_0^{\omega,1})| \geq (\delta/b)^{\eps_1}|\Tilde{W}| \geq \delta^{\eps_1}|\Tilde{W}|.$$
For technical reasons, we also need certain regularity of $U(\Tilde{\T}_{\tilde{W}}^{\omega,2}, \Tilde{Y}_0^{\omega,1})$ in $\tilde{W}$. We record the following definition:
\begin{definition}[{\cite[Definition~2.2]{WZ25assouad}}]\label{Def:regular_shading}
Let $\Omega \subset \R^n$, let $Y(\Omega) \subset \Omega$ be a shading, and let $\delta >0$. We say that the shading $Y(\Omega)$ is \textit{regular at scales $\geq \delta$} if for each $x \in Y(\Omega)$ and each $r \in [\delta,1]$, we have 
\begin{align*}
|Y(\Omega) \cap B(x,r)| \geq (100^n\log(1/\delta))^{-1}|Y(\Omega)| \Big( \frac{|\Omega \cap B(x,r)|}{|\Omega|} \Big).
\end{align*}
If the quantity $\delta$ is apparent from context, then we will omit it and say that $Y(\Omega)$ is
\textit{regular}.
\end{definition}

The next lemma says after a harmless refinement, every shading has a regular subshading. 
\begin{lem}[{\cite[Lemma~2.3]{WZ25assouad}}]\label{Lem:regular_subshading}
Let $\delta>0$, let $\Omega \subset B_1 \subset \R^n$, and let $Y(\Omega) \subset \Omega$ be a shading. Then there is a subshading $Y'(\Omega) \subset Y(\Omega)$ that is regular at scales $\geq \delta$ and satisfies $|Y'(\Omega)| \geq (1/2)|Y(\Omega)|$.
\end{lem}

Applying Lemma~\ref{Lem:regular_subshading} to $\Omega = \tilde{W}$ and $Y(\Omega) = U(\Tilde{\T}_{\tilde{W}}^{\omega,2}, \Tilde{Y}_0^{\omega,1})$, we get a subset $\tilde{A} \subset U(\Tilde{\T}_{\tilde{W}}^{\omega,2}, \Tilde{Y}_0^{\omega,1})$ that is regular and satisfies $|\tilde{A}| \geq (1/2)|U(\Tilde{\T}_{\tilde{W}}^{\omega,2}, \Tilde{Y}_0^{\omega,1})|$. For each $\tilde{T} \in \Tilde{\T}_{\tilde{W}}^{\omega,2}$, define $$\Tilde{Y}_{00}^{\omega,1}(\tilde{T}) \coloneqq \Tilde{Y}_0^{\omega,1}(\tilde{T}) \cap \tilde{A}.$$ 
Then since $(\Tilde{\T}_{\tilde{W}}^{\omega,2}, \Tilde{Y}_0^{\omega,1})$ has constant multiplicity and $$|U(\Tilde{\T}_{\tilde{W}}^{\omega,2}, \Tilde{Y}_{00}^{\omega,1})| = |\tilde{A}| \geq (1/2)|U(\Tilde{\T}_{\tilde{W}}^{\omega,2}, \Tilde{Y}_0^{\omega,1})|,$$ 
we know that $(\Tilde{\T}_{\tilde{W}}^{\omega,2}, \Tilde{Y}_{00}^{\omega,1})$ is a $\gtrsim 1$ refinement of $(\Tilde{\T}_{\tilde{W}}^{\omega,2}, \Tilde{Y}_0^{\omega,1})$. By pigeonholing with respect to the density of shadings, we get a subset $\Tilde{\T}_{\tilde{W}}^{\omega,3} \subset \Tilde{\T}_{\tilde{W}}^{\omega,2}$ such that $(\Tilde{\T}_{\tilde{W}}^{\omega,3}, \Tilde{Y}_{00}^{\omega,1})$ is a $\gtrsim 1$ refinement of $(\Tilde{\T}_{\tilde{W}}^{\omega,2}, \Tilde{Y}_{00}^{\omega,1})$, and for any $\tilde{T} \in \Tilde{\T}_{\tilde{W}}^{\omega,3}$, we have $|\Tilde{Y}_{00}^{\omega,1}(\tilde{T})|/|\tilde{T}| \gtrapprox \delta^{2\eta}$. Note that $(\Tilde{\T}_{\tilde{W}}^{\omega,3}, \Tilde{Y}_{00}^{\omega,1})$ is a $\gtrapprox 1$ refinement of $(\Tilde{\T}_{\tilde{W}}^{\omega,2}, \Tilde{Y}^{\omega,1})$.

For each $\omega \in \Gamma$, define 
$$\W_\omega \coloneqq \{ \pi_{X_1^\perp}^{-1}(\Tilde{W}) \cap \Pi_\omega: \Tilde{W} \in \Tilde{\W}_\omega \}.$$ 
Then for each $W \in \W_\omega$, we may simply use $\tilde{W}$ to denote $\pi_{X_1^\perp}(W)$. Under these notations, define $$\T^{\omega,3} \coloneqq \bigsqcup_{W\in\W_\omega} \T_W^{\omega,3},$$ where  
$$\T_W^{\omega,3} \coloneqq \{ T \in \T^{\omega,2}: \pi_{X_1^\perp}(T) \in \Tilde{\T}_{\tilde{W}}^{\omega,3} \}.$$ 
Also, define a shading $Y^{\omega,3}$ on $\T^{\omega,3}$ by $$Y^{\omega,3}(T) \coloneqq Y^{\omega,1}(T) \cap \pi_{X_1^\perp}^{-1}(\tilde{Y}_{00}^{\omega,1}(\tilde{T}))$$ for each $T \in \T^{\omega,3}$. It remains to check the three bullet points in the statement.

Firstly, by \eqref{Eq:regular_fibre}, for each $\omega \in \Gamma$ and $T \in \T^{\omega,3}$, we have 
\begin{align*}
|Y^{\omega,3}(T)| \geq \delta^{(n-3)+2\eta} |\tilde{Y}_{00}^{\omega,1}(\tilde{T})| \gtrapprox \delta^{(n-3)+4\eta}|\tilde{T}| \sim \delta^{4\eta}|T|.
\end{align*}
And by our construction of $\T^{\omega,3}$, there must exist some $W \in \W_\omega$ such that $T \in \T_W^{\omega,3}$, which means $\pi_{X_1^\perp}(T) \in \tilde{\T}_{\tilde{W}}^{\omega,3}$. By our construction of $\tilde{\T}_{\tilde{W}}^{\omega,3}$, we have $\pi_{X_1^\perp}(T) \subset \tilde{W}$, and so 
$$ T \subset \pi_{X_1^\perp}^{-1}(\pi_{X_1^\perp}(T)) \cap \Pi_\omega \subset \pi_{X_1^\perp}^{-1}(\tilde{W}) \cap \Pi_\omega = W. $$ 
Thus the first bullet point holds. 

Secondly, by \eqref{Eq:regular_fibre}, we also have 
\begin{align*}
\sum_{\T_W^{\omega,3}} |Y^{\omega,3}(T)| 
&\geq 
\delta^{(n-3)+2\eta} \sum_{\Tilde{\T}_{\tilde{W}}^{\omega,3}} |\Tilde{Y}_{00}^{\omega,1}(\tilde{T})| \\
&\gtrapprox 
\delta^{(n-3)+2\eta} \sum_{\Tilde{\T}_{\tilde{W}}^{\omega,2}} |\Tilde{Y}^{\omega,1}(\tilde{T})| \\
&\geq 
\delta^{(n-3)+4\eta} \sum_{\Tilde{\T}_{\tilde{W}}^{\omega,2}} |\tilde{T}| 
\sim 
\delta^{4\eta} \sum_{\T_W^{\omega,2}} |T| 
\geq 
\delta^{4\eta} \sum_{\T_W^{\omega,2}} |Y^{\omega,1}(T)|.
\end{align*}
Summing this over all $\omega \in \Gamma$ and $W\in \W_\omega$, we get 
\begin{align*}
\sum_{\omega\in\Gamma} 
\sum_{\T^{\omega,3}} |Y^{\omega,3}(T)| 
\gtrapprox 
\delta^{4\eta} 
\sum_{\omega\in\Gamma} 
\sum_{\T^{\omega,2}} |Y^{\omega,1}(T)| \gtrapprox \delta^{2\eps_1 + 7\eta},
\end{align*}
where we used \eqref{Eq:total_mass_0} in the last inequality. Thus the second bullet point holds.

Finally, by \eqref{Eq:regular_fibre}, for each $\omega \in \Gamma$ and $W \in \W_\omega$, we have 
\begin{equation*}
\begin{split}
|U(\T^{\omega,2}, Y^{\omega,1}) \cap W|
&\geq |U(\T^{\omega,2}[W], Y^{\omega,1})| \\
&\geq \delta^{(n-3)+2\eta} |U(\Tilde{\T}_{\tilde{W}}^{\omega,2}, \Tilde{Y}^{\omega,1})| \\
&\geq \delta^{\eps_1 + 2\eta} (\delta^{n-3}|\Tilde{W}|) 
= \delta^{\eps_1 + 2\eta} |W|.
\end{split}
\end{equation*}
Moreover, by the regularity of $\tilde{A} = U(\Tilde{\T}_{\tilde{W}}^{\omega,2}, \Tilde{Y}_{00}^{\omega,1})$ and \eqref{Eq:regular_fibre}, we know that for each $x \in U(\T^{\omega,3}[W], Y^{\omega,3})$ and each $r \in [\delta,1]$, we have 
\begin{equation*}
\begin{split}
|U(\T^{\omega,2}, Y^{\omega,1}) \cap W \cap B^n(x,r)| 
&\geq 
|U(\T^{\omega,2}[W], Y^{\omega,1}) \cap B^n(x,r)| \\
&\geq 
\delta^{(n-3)+2\eta} 
|U(\Tilde{\T}_{\tilde{W}}^{\omega,2}, \Tilde{Y}^{\omega,1}) \cap B^3(\tilde{x},r)| \\
&\geq \delta^{(n-3)+2\eta} 
|\tilde{A} \cap B^3(\tilde{x},r)|\\
&\gtrapprox 
\delta^{(n-3)+2\eta}
|\tilde{A}| |\tilde{W} \cap B^3(\tilde{x},r)|/|\tilde{W}|\\
&\gtrsim 
\delta^{2\eta}
|U(\Tilde{\T}_{\tilde{W}}^{\omega,2}, \Tilde{Y}_0^{\omega,1})||W \cap B^n(x,r)|/|\tilde{W}|\\
&\geq 
\delta^{\eps_1 +2\eta} |W \cap B^n(x,r)|,
\end{split}
\end{equation*}
where $\tilde{x} = \pi_{X_1^\perp}(x) \in \pi_{X_1^\perp}(U(\T^{\omega,3}, Y^{\omega,3}))  \subset \tilde{A}$. Thus the last bullet point holds.
\end{proof}

For $0 < \eta \ll \min\{\eps_1, \eta_1\}$ and $\max\{\eps_1, \eta_1\} \ll \sigma$, let $\{(\T^{\omega,3}, Y^{\omega,3})\}_{\omega\in\Gamma}$ and $\{\W_\omega\}_{\omega\in\Gamma}$ be the output of Proposition~\ref{Prop:step_5}. Intuitively, each $T \in \T^{\omega,3}$ grows into a significantly larger $\delta\one \times a \times b \times 1$ (recall that $b \geq \delta^{1-\eta_1}$) grain $W \in \W_\omega$ (with $T \subset W$), whose shading is regular. Note that $W \cap B(x,r)$ in \eqref{Eq:regular_grains_in_slices} is a rectangular prism in $\Pi_\omega$ of size $\delta\one \times \min\{a, r\} \times \min\{b, r\} \times r$.

This finishes Step~(v).

{\bf Step~(vi).} 
Now let us show that for the factoring, $a$ cannot be too large than $\delta$. More explicitly, we will prove the following result:
\begin{proposition}\label{Prop:step_6}
Suppose $0 < \eta \ll \eps_1 \ll \eta_i \ll \sigma$ where $i=1,2$, then there exists $\delta_0 > 0$ so that $a \in [\delta, \delta^{1-\eta_2}]$ and $b \in [\delta^{1-\eta_1}, 1]$ both hold for all $\delta \in (0,\delta_0]$. 
\end{proposition}
For the proof, we need to understand how information from slices translates into information about the flakes. It would be convenient to first introduce some notation. For each $F \in \F_1$, define a shading $Y_2$ by 
$$Y_2(F) \coloneqq \bigsqcup_{\T(F) \cap (\sqcup_{\omega\in \Gamma} \T^{\omega,3})} Y^{\omega,3}(T).$$ 
Note that \eqref{Eq:total_mass''} implies that 
\begin{align}\label{Eq:Y_2_density}
\sum_{\F_1} |Y_2(F)| = 
\sum_{\omega \in \Gamma} \sum_{\T^{\omega,3}}|Y^{\omega,3}(T)|
\gtrapprox 
\delta^{2\eps_1+7\eta}.
\end{align} 

\begin{proof}[Proof of Proposition~\ref{Prop:step_6}]
{\color{black} Let $\eta_2 > 0$ be a small parameter to be determined later on.} 
For each $\omega \in \Gamma$, we have already grown tubes into a bunch of $\delta\one \times a \times b \times 1$ ($\delta \leq a \leq b \leq 1$) grains $\W_\omega$ in $\Pi_\omega$. By Proposition~\ref{Prop:step_4}, we already know that $b \in [\delta^{1-\eta_1}, 1]$, so it remains to focus on $a$. 

For technical reasons, let us first handle a trivial case. If $a \geq \delta^{\eta_2}$, then by merely summing up the contributions from all slices, we have 
\begin{align*}
|U(\F_1, Y_1)| 
&\geq 
\sum_{\omega \in \Gamma} |U(\T^{\omega,2}, Y^{\omega,1})| \\
&\geq 
|\Gamma| \delta^{\eps_1+2\eta} |W| \\
&\gtrapprox 
\delta^{-(n-3) + 2\eps_1 + 4\eta} \delta^{\eps_1+2\eta} \delta^{(n-3)+2\eta_2} \\
&= 
\delta^{3\eps_1 + 6\eta + 2\eta_2}, 
\end{align*}
where we used the fact that $\W_\omega \neq \varnothing$ for all $\omega \in \Gamma$ and \eqref{Eq:grains_in_slices}. On the other hand, since $(\F, Y)$ is $\eta$-extremal, we have 
\begin{align*}
|U(\F_1,Y_1)| \leq |U(\F,Y)| \leq \delta^{\sigma - \eta}.
\end{align*}
As $\sigma >0$ is fixed, by first {\color{black} choosing $\max\{ \eta, \eta_2, \eps_1\} \ll \sigma $} and then letting $\delta_0$ be small enough, we obtain a contradiction between the upper and lower bounds. Therefore, we may always assume that $a \leq \delta^{\eta_2}$.

Suppose $\delta^{1-\eta_2} \leq a \leq \delta^{\eta_2}$, then intuitively a typical $\delta$-flake still grows significantly by $a$ in its two normal directions (as it does in each slice), and so the union must fill out an $a$-neighborhood of it. To make this point rigorous, we need to do some pigeonholing on $(\F_1, Y_1)$.  
More precisely, since we have \eqref{Eq:Y_2_density} and $|\F_1| \leq |\F| \leq \delta^{-2-\eta}$, there exists a subset $\F_1' \subset \F_1$ such that $|Y_2(F)| \geq \delta^{2 + 2\eps_1 + 9\eta}$ for each $F \in \F_1'$, and 
\begin{align*}
\sum_{\F_1'} |Y_2(F)| \gtrapprox \delta^{2\eps_1+7\eta}.
\end{align*}
Due to $|Y_2(F)| \leq |F| \sim \delta^2$, we have $|\F_1'| \gtrapprox \delta^{-2+2\eps_1+7\eta} \geq \delta^{2\eps_1+8\eta}|\F|$. Thus by Lemma~\ref{Lem:split_sticky} with $a$ in place of $\rho$ and $\eta_2$ in place of $\eps$, we know that $\F_1'$ contains a uniform subset $\F_2 \subset \F_1'$ with $|\F_2| \gtrapprox |\F_1'|$ that is $\delta^{-(2\eps_1 + 12\eta)}$-sticky. Moreover, $\F_2$ is uniform at scale $a$, and $(\F_2)_a$ is an $a^{-(2\eps_1 + 13\eta)/\eta_2}$-sticky family of $a$-flakes, which is also spaced by the proof of Lemma~\ref{Lem:split_sticky}. So if we thicken each $a$-flake $F_a \in (\F_2)_a$ to be a $2a$-flake $F_{2a}$ with respect to $\dir(F_a)$, then the resulting set of $2a$-flakes $(\F_2)_{2a}$ (slightly abusing notation) is $(2a)^{-(2\eps_1 + 14\eta)/\eta_2}$-sticky.

Now for each $F_{2a} \in (\F_2)_{2a}$, take one $F \in \F_2[F_a]$. Since $|Y_2(F)| \geq \delta^{2 + 2\eps_1 + 9\eta}$ and $|Y^{\omega,3}(T)| \leq |T| \lesssim \delta^{n-1}$, we know that 
\begin{align}\label{Eq:layer_number}
\big|\T(F) \cap \bigsqcup_{\omega\in \Gamma} \T^{\omega,3}\big| \gtrsim \delta^{-(n-3) + 2\eps_1 + 9\eta}.
\end{align}
On the other hand, by \eqref{Eq:regular_grains_in_slices} with $r = a$, for each $\omega \in \Gamma$ and $T \in \T^{\omega,3}$, there is a $W \in \W_\omega$ with $T \subset W$, and we have $$|U(\T^{\omega,2}, Y^{\omega,1}) \cap W \cap B(x,a)| \gtrapprox \delta^{\eps_1 + 2\eta}|W \cap B(x,a)|$$ 
for each $x \in Y^{\omega,3}(T)$. Recall from \eqref{Eq:single_tube_density} that $|Y^{\omega,3}(T)| \gtrapprox \delta^{4\eta}|T|$, and so there exists a set of $2a$-separated points $\{x_k\} \subset Y^{\omega,3}(T)$ with cardinality $\gtrapprox \delta^{4\eta}/a$. By summing up the contributions from all $a$-balls, for each $\omega \in \Gamma$ with $\T(F) \cap \T^{\omega,3} \neq \varnothing$, we have 
\begin{align*}
|U(\T^{\omega,2}, Y^{\omega,1}) \cap W \cap F_{2a}| 
&\geq 
\big|U(\T^{\omega,2}, Y^{\omega,1}) \cap W \cap \bigcup_k B(x_k,a)\big| \\
&=
\sum_k |U(\T^{\omega,2}, Y^{\omega,1}) \cap W \cap B(x_k,a)| \\
&\gtrapprox 
(\delta^{4\eta}/a) \cdot \delta^{\eps_1+2\eta} |W \cap B(x,a)| \\
&\sim 
\delta^{\eps_1+6\eta} (\delta^{n-3} a^3)/a\\
&= \delta^{(n-3)+\eps_1+6\eta} a^2.
\end{align*}
Thus by \eqref{Eq:layer_number} and further summing the above estimates over all $\omega \in \Gamma$ with $\T(F) \cap \T^{\omega,3} \neq \varnothing$, we obtain 
\begin{align*}
|U(\F_1, Y_1) \cap F_{2a}| 
&\geq 
\sum_{\omega}|U(\T^{\omega,2}, Y^{\omega,1}) \cap F_{2a}| \\
&\gtrapprox 
\delta^{-(n-3) + 2\eps_1 + 9\eta} \cdot \delta^{(n-3) + \eps_1 + 6\eta} a^2 \\
&\sim \delta^{3\eps_1 + 15\eta} |F_{2a}|.
\end{align*}
If we define a shading $\Bar{Y}_{2a}$ on $(\F_2)_{2a}$ by $\Bar{Y}_{2a}(F_{2a}) \coloneqq U(\F_1, Y_1) \cap F_{2a}$ for each $F_{2a} \in (\F_2)_{2a}$, then $\lambda((\F_2)_{2a}, \Bar{Y}_{2a}) \gtrapprox \delta^{3\eps_1 + 15\eta} \gtrsim (2a)^{(3\eps_1 + 15\eta)/\eta_2}$. Therefore, $((\F_2)_{2a}, \Bar{Y}_{2a})$ satisfies Items~\ref{item_a}~and~\ref{item_b} in Definition~\ref{Def:sigma} with $\delta = 2a$ and $\eps = (3\eps_1 + 16\eta)/\eta_2$. 

{\color{black} Let $\eps_3 > 0$ be a small parameter to be determined later on.} By Remark~\ref{Rem:extreme_mult}, if we {\color{black} choose $\max\{\eta,\eps_1\} \ll \min\{\eta_2,\eps_3\}$}, then $$|U((\F_2)_{2a}, \Bar{Y}_{2a})| \geq (2a)^{\sigma + \eps_3} \geq \delta^{(1-\eta_2)(\sigma+\eps_3)} = \delta^{\sigma - (\eta_2\sigma + \eta_2 \eps_3 - \eps_3)}.$$
{\color{black} By choosing $\eps_3 \leq \eta_2 \sigma/2$,} the lower bound becomes $|U((\F_2)_{2a}, \Bar{Y}_{2a})| \geq \delta^{\sigma - \eta_2\sigma/2}$. 
On the other hand, since $U((\F_2)_{2a}, \Bar{Y}_{2a}) \subset U(\F_1, Y_1) \subset U(\F, Y)$ and $(\F, Y)$ is $\eta$-extremal, we have the upper bound $|U(\F_1,Y_1)| \leq |U(\F,Y)| \leq \delta^{\sigma - \eta}$. {\color{black} By choosing $\eta \leq \eta_2\sigma/4$} and then letting $\delta_0$ be small enough, we obtain a contradiction between the upper and lower bounds. Therefore, we must have $a \leq \delta^{1-\eta_2}$. 
\end{proof}

In other words, for $0 < \eta \ll \eps_1 \ll \eta_i \ll \sigma$ where $i=1,2$, by Proposition~\ref{Prop:step_6}, we can always assume $a \leq \delta^{1-\eta_2}$ and $b \geq \delta^{1-\eta_1}$, where $a, b$ are from Proposition~\ref{Prop:step_3}. {\color{black} Recall that $\eta$ is the extremal parameter for $(\F,Y)$ at the beginning of the proof (before Step~(i)), and $\eps_1$ first appears when we apply the induction hypothesis during the proof of Proposition~\ref{Prop:step_2} (in Step~(ii)).}

This finishes Step~(vi).

\section{Medium grains and the end of the proof}\label{Sec:medium_grains_end}

In this section, we will continue the proof in Section~\ref{Sec:slicing_set_up} and finish Steps~(vii)--(x) therein. {\color{black} Recall from Step~(v) that each $T \in \T^{\omega, 3}$ grows into a $\delta\one \times a \times b \times 1$ grain $W \in \W_\omega$ (with $T \subset W$)}, and by Step~(iv) and Step~(vi), we may assume that $b \geq \delta^{1-\eta_1}$ and $a \leq \delta^{1-\eta_2}$. {\color{black} Now we choose $\eta_2 \ll \eta_1$.} This essentially means that $a \approx \delta$ and $b \gg \delta$, and the grains roughly have dimensions $\delta\1 \times b \times 1$. From now on, we will call them \textit{medium grains}, which extends the $\delta$-tube $T$ in some normal direction of a $\delta$-flake $F$ containing $T$. As outlined at the beginning of Section~\ref{Sec:slicing_set_up}, we need to analyze the compatibility of medium grains and local grains to regularize the plane map, which is vital for further grains-growing procedures. Now we turn to the details.

{\bf Step~(vii).} Let us first show that the plane map $\HH$ cannot change too much along the tubes corresponding to $X_1$ (up to a refinement). More explicitly, we will prove the following result: 
\begin{proposition}\label{Prop:step_7}
Suppose $0 < \eta \ll \eps_1 \ll \eta_2 \ll \eta_1 \ll \sigma$, then there exists $\delta_0 > 0$ so that the following holds for all $\delta \in (0,\delta_0]$. 

For each $\omega \in \Gamma$ and $T \in \T^{\omega,3}$, there exists some $\HH_T \in \Gr(n,1)$ depending only on $T$, such that $\angle(\HH(x), \HH_T) \lesssim \delta^{\eta_1/2}$ 
for all $x \in Y^{\omega,3}(T)$.
\end{proposition}

\begin{proof}
Fix any $\omega \in \Gamma$ and $T \in \T^{\omega,3}$. Then there exists a $W \in \W_\omega$ with $T \subset W$, and an $F \in \F_1$ with $T \in \T(F)$. Let $s \coloneqq \min\{b, \delta^{1/2}\}$, then $s \geq \delta^{1-\eta_1} \geq a$ (as $\eta_2 \leq \eta_1 \leq 1/2$). By \eqref{Eq:regular_grains_in_slices} with $s$ in place of $r$, for each $x \in Y^{\omega,3}(T)$, we have 
\begin{align*}
|U(\T^{\omega,2}, Y^{\omega,1}) \cap W \cap B(x,s)| 
\gtrapprox 
\delta^{\eps_1 + 2\eta} 
|W \cap B(x,s)|,
\end{align*}
where $W \cap B(x,s)$ is a rectangular box in $\Pi_\omega$ of size $\delta\one \times a \times s \times s$. 

Let $P_{W,x} \subset B(x,s)$ be the $a\1 \times s \times s$ plate segment obtained by thickening the sides of $W \cap B(x,s)$ along $X_1$ from $\delta$ to $a$. We may assume that $P_{W,x}$ is centered at $x$. Since $s \gg a$, by a rescaled version of Definition~\ref{Def:flake_direction}, it makes sense to talk about its direction subspace $\dir(P_{W,x})$, which is well defined up to error $a/s$ by (a rescaled version of) Lemma~\ref{Lem:flake_contain}. We may assume that $\dir(T) \subset \dir(P_{W,x}) \subset X_1^\perp$. Also, since each $P_{W,x}$ is parallel to $W$, we can make sure that $\dir(P_{W,x})$ is the same for all $x \in Y^{\omega,3}(T)$. In particular, each $P_{W,x}$ can be tiled by $\sim s/a$ parallel $a\1 \times a \times s$ tube segments pointing in $\dir(T)=\dir(F) \cap X_1^\perp$, as well as in some $\ell_{W,x} \in \Gr(n,1)$ perpendicular to $\dir(T)$, where $\ell_{W,x}$ are the same for all $x\in Y^{\omega,3}(T)$. Let $\ell_W$ be the common value. Finally, note that 
\begin{equation}
\begin{split}
|U(\F, Y) \cap P_{W,x}| &\geq |U(\T^{\omega,2}, Y^{\omega,1}) \cap P_{W,x}| \\
&\geq 
|U(\T^{\omega,2}, Y^{\omega,1}) \cap W \cap B(x,s)| \\
&\gtrapprox 
\delta^{\eps_1 + 2\eta} 
|W \cap B(x,s)| \\
&\gtrsim 
\delta^{\eps_1 + 2\eta} (\delta/a)^{n-3}|P_{W,x}| \\
&\geq 
\delta^{\eps_1 + (n-3)\eta_2 + 2\eta} |P_{W,x}|
\end{split}\label{dense_flake_seg}
\end{equation}

On the other hand, the graininess of $(\F,Y)$ at scale $\rho = a \in [\delta,1]$ implies that for all $x \in U(\F,Y)$, $\pi_{\HH(x)}(U(\F,Y) \cap B(x,a^{1/2}))$ is an $(a, 1-\sigma, \delta^{-\eta})_1$-ADset. Since $s \leq \delta^{1/2} \leq a^{1/2}$, we have $U(\F, Y) \cap P_{W,x} \subset U(\F,Y) \cap B(x,a^{1/2})$, and so $\pi_{\HH(x)}(U(\F,Y) \cap P_{W,x})$ is an $(a, 1-\sigma, \delta^{-\eta})_1$-ADset. 
Since $x \in Y^{\omega,3}(T) \subset Y^{\omega,1}(T) \subset Y_1(F) \subset Y(F)$, we know from Definition~\ref{Def:plane_map} that 
\begin{align}\label{Eq:2-plane_HH}
\dir(F) \cap B_1 \subset N_{2\delta}(\HH(x)^\perp). 
\end{align}

\begin{rem}
Intuitively, $\HH(x)$ is close to $\dir(F)^\perp \in \Gr(n,2)$ with error $\sim\delta$. So there is essentially one degree of freedom (corresponding to the unit circle in $\dir(F)^\perp$). We will show that $\ell_W$ and $\HH(x)$ are approximately orthogonal (the cosine of the angle is $\leq \delta^{\eta_1/2}$), which combined with the fact that $\ell_W$ is transverse to $\dir(F)$ would yield an additional effective restriction on $\HH(x)$ (see \eqref{Eq:NN_HH}). See Figure~\ref{fig:medium_grains} for an illustration of the special case when $n=4$. For this purpose, we need to exploit the local grains structure.

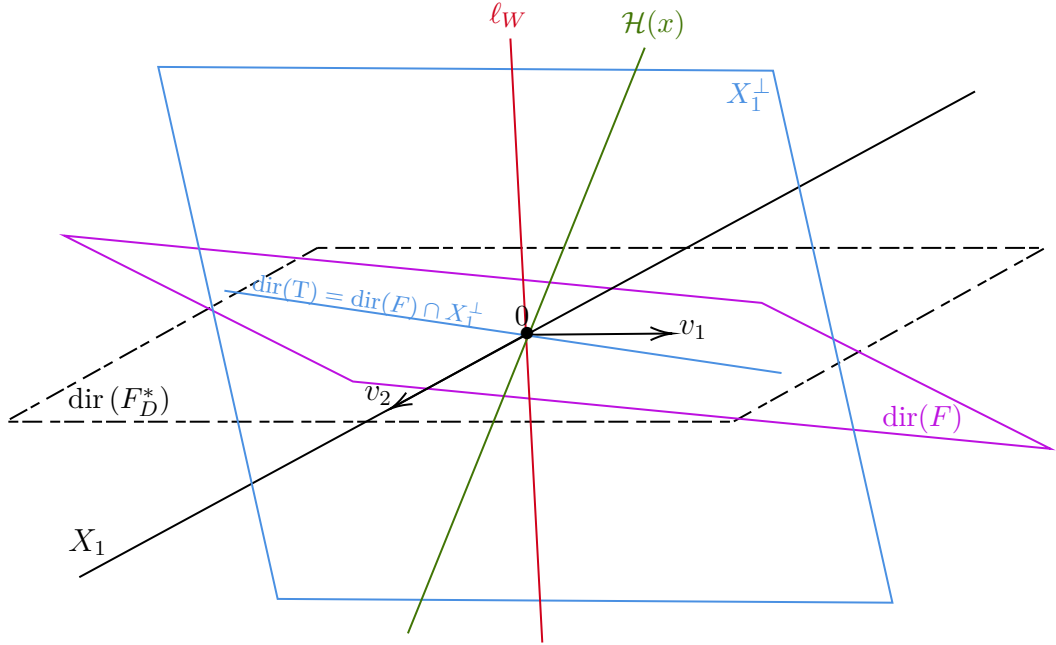
\begin{figure}[htbp]
  \centering
  
    \hspace*{0.5cm}

\tikzset{every picture/.style={line width=0.75pt}} %set default line width to 0.75pt        

\begin{tikzpicture}[x=0.75pt,y=0.75pt,yscale=-1,xscale=1]
%uncomment if require: \path (0,1645); %set diagram left start at 0, and has height of 1645

%Shape: Parallelogram 
\draw  [dash pattern={on 3.75pt off 3pt on 7.5pt off 1.5pt}] (309.5,758.48) -- (674.66,758.48) -- (518.16,845.67) -- (153,845.67) -- cycle ;
%Shape: Parallelogram 
\draw  [color={rgb, 255:red, 189; green, 16; blue, 224 }  ,draw opacity=1 ] (327.33,825.51) -- (677.33,859.27) -- (532.42,786.07) -- (182.42,752.31) -- cycle ;
%Straight Lines 
\draw    (413.83,802.07) -- (347.81,838.28) ;
\draw [shift={(346.06,839.24)}, rotate = 331.26] [color={rgb, 255:red, 0; green, 0; blue, 0 }  ][line width=0.75]    (10.93,-3.29) .. controls (6.95,-1.4) and (3.31,-0.3) .. (0,0) .. controls (3.31,0.3) and (6.95,1.4) .. (10.93,3.29)   ;
%Straight Lines 
\draw    (413.83,802.07) -- (485.55,801.45) ;
\draw [shift={(487.55,801.43)}, rotate = 179.5] [color={rgb, 255:red, 0; green, 0; blue, 0 }  ][line width=0.75]    (10.93,-3.29) .. controls (6.95,-1.4) and (3.31,-0.3) .. (0,0) .. controls (3.31,0.3) and (6.95,1.4) .. (10.93,3.29)   ;
%Straight Lines 
\draw [line width=0.75]    (639.43,680.01) -- (190.06,923.68) ;
%Shape: Parallelogram 
\draw  [color={rgb, 255:red, 74; green, 144; blue, 226 }  ,draw opacity=1 ] (289.57,934.57) -- (597.98,936.46) -- (538.09,669.58) -- (229.68,667.69) -- cycle ;
%Straight Lines 
\draw [color={rgb, 255:red, 74; green, 144; blue, 226 }  ,draw opacity=1 ]   (262.86,779.96) -- (542.35,821.26) ;
%Straight Lines 
\draw [color={rgb, 255:red, 208; green, 2; blue, 27 }  ,draw opacity=1 ]   (406.33,653.45) -- (422.33,956.45) ;
%Straight Lines 
\draw [color={rgb, 255:red, 65; green, 117; blue, 5 }  ,draw opacity=1 ]   (473.75,658.04) -- (354.92,951.85) ;

% Text Node
\draw (210.22,835.44) node  [font=\small]  {$\mathrm{dir}\left( F_{D}^{*}\right)$};
% Text Node
\draw (613.27,844.05) node  [font=\small,color={rgb, 255:red, 189; green, 16; blue, 224 }  ,opacity=1 ]  {$\mathrm{dir}( F)$};
% Text Node
\draw (194.38,906.12) node    {$X_{1}$};
% Text Node
\draw (498.01,799.57) node    {$v_{1}$};
% Text Node
\draw (340.08,832.61) node    {$v_{2}$};
% Text Node
\draw (526.01,681.45) node  [font=\small,color={rgb, 255:red, 74; green, 144; blue, 226 }  ,opacity=1 ]  {$X_{1}^{\perp }$};
% Text Node
\draw (412.18,792.44) node  [font=\small]  {$0$};
% Text Node
\draw (414.74,801.85) node    {$\smblkcircle $};
% Text Node
\draw (335.08,784.58) node  [font=\scriptsize,color={rgb, 255:red, 74; green, 144; blue, 226 }  ,opacity=1 ,rotate=-8.45]  {$\mathrm{dir( T) =dir}( F) \cap X_{1}^{\perp }$};
% Text Node
\draw (405.34,641.52) node  [font=\small,color={rgb, 255:red, 208; green, 2; blue, 27 }  ,opacity=1 ]  {$\ell _{W}$};
% Text Node
\draw (478.08,646.63) node  [font=\small,color={rgb, 255:red, 65; green, 117; blue, 5 }  ,opacity=1 ]  {$\mathcal{H}( x)$};

\end{tikzpicture}
\caption{This can be seen as a continuation of Figure~\ref{fig:slicing}. The point here is that we want to ``determine'' $\HH(x)$ by using $\ell_W$ and $\dir(F)$.}
\label{fig:medium_grains}
\end{figure}
\end{rem}

We need a mild generalization of \cite[Lemma~6.2]{WZ26sticky}, which says that if $T_{seg}$ is
a segment of a tube $T$, and if the image of a large subset of $T_{seg}$ under $\pi_\HH$ is an ADset, then $\HH$ must be nearly orthogonal to $\dir(T)$. The
precise statement is as follows.

\begin{lem}\label{Lem:tube_seg_AD}
Let \( \varepsilon, \eta > 0 \), \( 0 \leq \alpha < 1 \), and \( 0 < \delta \leq a \leq s \leq 1 \). Let \( T \) be an \(a\)-tube, let \( T_{seg} \subset T \) be a tube segment of length \( s \), and let \( U \subset T_{seg} \) with
\begin{align}\label{Eq:dense_tube_seg}
|U| \geq \delta^\varepsilon |T_{seg}|. 
\end{align}
Let \( \ell \in \Gr(n,1) \), and suppose \(\pi_\ell(U)\) is an \( (a, \alpha, \delta^{-\eta})_1 \)-ADset. Then
\begin{align}\label{Eq:dense_tube_seg'}
\dir(T) \cap B_1 \subset N_{\tau}(\ell^\perp).
\end{align}
for some $\tau \lesssim \delta^{-\frac{\eps + \eta}{1-\alpha}} a/s$.
\end{lem}

\begin{proof}
Let \( \tau\) be the smallest number such that \eqref{Eq:dense_tube_seg'} holds. It suffices to bound this $\tau$. We will suppose that \( \tau \geq a/s \), since otherwise there is nothing to prove. Divide \( T_{seg} \) into subsegments of length \( a/\tau \) (by hypothesis \( a/\tau \leq s \), so there is at least one such segment), and let \( \mathbb{J} \) be this set of segments. For each \( J \in \mathbb{J} \), \( \pi_\ell(J) \) is a line segment of length \( \sim a \), and each point in $\ell$ is contained in $\lesssim 1$ of the sets \( \{ \pi_\ell(J): J \in \mathbb{J}\} \). By \eqref{Eq:dense_tube_seg}, we have
\[
\#\{J \in \mathbb{J}: J \cap U \neq \varnothing\} \geq \delta^\varepsilon \tau s/a,
\]
and thus 
\[
\mathcal{E}_a(\pi_\ell(U)) \gtrsim \delta^\varepsilon \tau s/a.
\]
On the other hand, \( \pi_\ell(U) \subset \pi_\ell(T_{seg}) \) is contained in a line segment \( I \subset \ell \) of length \( \sim \tau s\), and hence
\[
\mathcal{E}_a(\pi_\ell(U)) \leq \mathcal{E}_a(\pi_\ell(U) \cap I) \leq \delta^{-\eta}(|I|/a)^\alpha \lesssim \delta^{-\eta}(\tau s/a)^\alpha.
\]
Comparing the lower and upper bounds of $\mathcal{E}_a(\pi_\ell(U))$, we conclude that \( \tau \lesssim \delta^{-\frac{\eps + \eta}{1-\alpha}} a/s \), as desired.
\end{proof}

By \eqref{dense_flake_seg} and pigeonholing, there exists an $a\1 \times a \times s$ tube segment $T_{W,x} \subset P_{W,x}$ pointing in $\ell_{W}$, such that \begin{align*}
|U(\F,Y) \cap T_{W,x}| \gtrapprox 
\delta^{\eps_1 + (n-3)\eta_2 + 2\eta} |T_{W,x}|.
\end{align*}
Also, $\pi_{\HH(x)}(U(\F,Y) \cap T_{W,x}) \subset \pi_{\HH(x)}(U(\F,Y) \cap P_{W,x})$ is an $(a, 1-\sigma, \delta^{-\eta})_1$-ADset. So by applying Lemma~\ref{Lem:tube_seg_AD} with $T_{seg} = T_{W,x}$, $U = U(\F,Y) \cap T_{W,x}$, $\ell = \HH(x)$, $\eps = \eps_1 + (n-3) \eta_2 + 3\eta$, and $\alpha = 1 - \sigma$, we obtain
\begin{align*}
    \ell_W \cap B_1 \subset N_\tau(\HH(x)^\perp)
\end{align*}
for some 
$$\tau 
\lesssim 
\delta^{-[\eps_1 + (n-3)\eta_2 + 4\eta]/\sigma} a/s 
\leq 
\delta^{\eta_1 - \eta_2 - [\eps_1 + (n-3)\eta_2 + 4\eta]/\sigma}.$$ 
If we {\color{black} choose $\max\{ \eta, \eta_2, \eps_1\} \ll \min\{ \eta_1, \sigma \}$}, then this becomes 
\begin{align}\label{Eq:NN_HH}
    \ell_W \cap B_1 \subset N_{\delta^{\eta_1/2}}(\HH(x)^\perp).
\end{align}

Let $V_0 \coloneqq \dir(F) + \ell_W$. We want to apply Lemma~\ref{Lem:geo_fact_1} with
\begin{align*}
    V_1 = \dir(F), & \quad V_1' = \pi_{\HH(x)^\perp}(\dir(F)), \\
    V_2 = \ell_W, & \quad V_2' = \pi_{\HH(x)^\perp}(\ell_W).
\end{align*} 
Let us check the assumptions. By \eqref{Eq:2-plane_HH} and Lemma~\ref{Lem:faithful_proj} with $\rho = 2\delta$, we get $\dim V_1' = n-2$ and $\angle(V_1, V_1') \leq 2\delta \leq \delta^{\eta_1/2}$. Similarly, by \eqref{Eq:NN_HH} and Lemma~\ref{Lem:faithful_proj} with $\rho = \delta^{\eta_1/2}$, we get $\dim V_2' = 1$ and $\angle(V_2, V_2') \leq \delta^{\eta_1/2}$. Finally, note that $X_1 \oplus \dir(T) \subset \ell_W^\perp$, and by Remark~\ref{Rem:major_equiv} and \eqref{Eq:choice_D}, we have 
$$\angle(X_1 \oplus \dir(T), \dir(F)) \leq \angle(\dir(F_D^*), \dir(F)) \leq (100n)^{-1}.$$ 
These mean that $V_1 \cap B_1 \subset N_{1/20}(X_1 \oplus \dir(T)) \subset N_{1/20}(V_2^\perp)$. Therefore, we can apply Lemma~\ref{Lem:geo_fact_1} with $\rho = \delta^{\eta_1/2}$ to conclude that $V_1 + V_2 = V_0$ and $V_1' + V_2' \subset \HH(x)^\perp$ are both direct sums, which implies $V_0 = V_1 \oplus V_2 \in \Gr(n,n-1)$ and $\HH(x)^\perp = V_1' \oplus V_2'$ by counting dimensions. Moreover, $\angle(V_0, \HH(x)^\perp) = \angle(V_1 \oplus V_2, V_1' \oplus V_2') \leq 3\delta^{\eta_1/2}$.

Let $\HH_W \coloneqq V_0^\perp \in \Gr(n,1)$. Then by Item~(iii) of Theorem~\ref{Thm:Principal_angle}, $\angle(\HH_W, \HH(x)) = \angle(\HH_W^\perp, \HH(x)^\perp) \leq 3\delta^{\eta_1/2}$. The key fact is that this uniformly holds over all $x\in Y^{\omega,3}(T)$. If we fix a choice of $W$ for each $T$, then we may also write $\HH_W$ as $\HH_T$ to emphasize such dependence. In a nutshell, for each $T \in \T^{\omega,3}$, we have proved that 
$$\angle(\HH(x), \HH_T) \lesssim \delta^{\eta_1/2}$$ 
for all $x \in Y^{\omega,3}(T)$. 
\end{proof}

Proposition~\ref{Prop:step_7} shows that $\HH$ is stable along each $T \in \T^{\omega,3}$ (up to a refinement).

This finishes Step~(vii).

{\bf Step~(viii).} 
As forecast in the outline, by performing slicing arguments to all direction subspaces $X_k$, the plane map $\HH$ is also forced to not change too much along each $\delta$-flake (up to a refinement). More explicitly, we will prove the following result:
\begin{proposition}\label{Prop:step_8}
Suppose $0 < \eta \ll \eta_0 \ll \eta_1 \ll \sigma$, then there exists $\delta_0 > 0$ so that the following holds for all $\delta \in (0,\delta_0]$. 

There exists a refinement $(\F, Y_3)$ of $(\F, Y)$ such that 
\begin{itemize}
    \item $\sum_{\F} |Y_3(F)| \gtrapprox  \delta^{4n\eta_0}$. In particular, $(\F, Y_3)$ is $4n\eta_0$-extremal.

    \item For any $F \in \F$ with $Y_3(F) \neq \varnothing$, we have 
    $$\angle(\HH(x), \HH_{F}) \lesssim \delta^{\eta_1/2}$$
    for each $x \in Y_3(F)$, where $\HH_{F} \in \Gr(\dir(F)^\perp,1)$ only depends on $F$.
\end{itemize}
\end{proposition}

\begin{proof}
For clarity of notation later on, we rewrite $Y_2$ as $Y^{(1)}$, $\Gamma$ as $\Gamma^{(1)}$, $\eta$ as $\eta^{(1)}$, $\eps_1$ as $\eps_1^{(1)}$, and $(\T^{\omega,3}, Y^{\omega,3})$ as $(\T_*^{\omega,(1)}, Y_*^{\omega,(1)})$. Recall from Step~(vi) that for each $F \in \F_1$, we have 
$$Y^{(1)}(F) = \bigsqcup_{\T^{(1)}(F) \cap (\sqcup_{\omega\in \Gamma^{(1)}} \T_*^{\omega, (1)})} Y_*^{\omega, (1)}(T),$$ 
and \eqref{Eq:Y_2_density} becomes 
\begin{align*}
\sum_{\F_1} |Y^{(1)}(F)| 
\gtrapprox 
\delta^{2\eps_1^{(1)}+7\eta^{(1)}}.
\end{align*}
Also, by \eqref{Eq:single_tube_density}, $|Y_*^{\omega, (1)}(T)| \gtrapprox \delta^{4\eta^{(1)}}|T|$ for each $T \in \bigsqcup_{\omega \in \Gamma^{(1)}} \T_*^{\omega,(1)}$. 
For technical clarity, let us trivially extend the definition of $Y^{(1)}$ to the whole $\F$ by requiring $Y^{(1)}(F) \coloneqq \varnothing$ for each $F \in \F\setminus\F_1$. Then 
\begin{align*}
\lambda(\F,Y^{(1)}) = \frac{\sum_{\F} |Y^{(1)}(F)|}{\sum_{\F}|F|} = \frac{\sum_{\F_1} |Y^{(1)}(F)|}{\sum_{\F}|F|} 
\geq 
\delta^{2\eps_1^{(1)} + 9\eta^{(1)}}.
\end{align*}
And so $(\F,Y^{(1)})$ is $(2\eps_1^{(1)} + 9\eta^{(1)})$-extremal, with the same plane map $\HH$ as $(\F,Y)$.

Now we consider other $X_k$. {\color{black} We divide our arguments into three smaller steps, each of which involves an iterative algorithm. In Step~(A), we iteratively refine $(\F, Y^{(1)})$ for each $X_k$, which will yield a preparatory refinement $(\F, Y^{(n-2)})$. In Step~(B), we iteratively construct the shading in each $F$ to get the desirable refinement $(\F, Y_3)$, after which we verify the first bullet point. In Step~(C), we iteratively construct an admissible route connecting any pair of points in each $Y_3(F)$, by which we verify the second bullet point.}

{\bf Step~(A).} We run the following iterative algorithm. Let $k$ be the loop index, with $2$ as its initial value. 

Repeat Step~(ii)--(vii) with $(\F, Y^{(k-1)})$ in place of $(\F,Y)$, 
$$\eta^{(k)} \coloneqq \sum_{j=1}^{k-1} c_j^{(k)} \eps_1^{(j)} + c_0^{(k)} \eta^{(k-1)}$$ 
in place of $\eta$ ($c_1^{(2)}=2$, $c_0^{(2)}=9$), and $X_k$ in place of $X_1$. Recall that for each $\omega\in \Lambda_k$, the slicing procedure in Step~(i) gives $(\T^{\omega, (k)},Y^{\omega, (k)})$, and $(\F, Y^{(k-1)})$ naturally induces a refinement $(\T^{\omega, (k)},Y_0^{\omega, (k)})$ of it by restricting the shading of flakes to the slice $\Pi_\omega^{(k)}$. 

Throughout this process, we will introduce an additional collection of small parameters, whose relationship can be roughly summarized as follows:
\begin{align*}
\boxed{ \sum_{j=1}^{k-1} c_j^{(k)} \eps_1^{(j)} + c_0^{(k)} \eta^{(k-1)} = \eta^{(k)} \ll \eps_1^{(k)} \ll \eps_3^{(k)} \ll \eta_2^{(k)} \ll \eta_1^{(k)} \ll \eps_2^{(k)} \ll \sigma.}
\end{align*}
Here and henceforth we put superscripts ``$(k)$'' to emphasize that it is the $k$-th time running Step~(ii)--(vii), and to indicate the parallel structure between objects, e.g., $\eta_1^{(k)}$ plays the same role as that of $\eta_1^{(1)} = \eta_1$. Note that by choosing $$\max\{ \eps_1^{(1)}, \ldots, \eps_1^{(k-1)} \} \leq \eta^{(k)} \ll \eps_1^{(k)} \ll \eps_3,$$ 
we may ``take'' $\eps_3^{(k)} = \eps_3$, $\eta_2^{(k)} = \eta_2$, $\eta_1^{(k)} = \eta_1$, and $\eps_2^{(k)} = \eps_2$.\footnote{Rigorously speaking, we may need to modify the exact values of small parameters, but the important thing is that their ``$\ll$''-relationship is preserved. We write in this way for simplicity. The reader is invited to think about the logic behind the choices.} As a result, in each slice $\Pi_{\omega}^{(k)}$ with $\omega \in \Gamma^{(k)} \subset \Lambda_k$, there is a refinement $(\T_*^{\omega, (k)}, Y_*^{\omega, (k)})$ of $(\T^{\omega, (k)},Y_0^{\omega, (k)})$, such that for each $T \in \T_*^{\omega, (k)}$, we have 
$$|Y_*^{\omega, (k)}(T)| \gtrapprox \delta^{4\eta^{(k)}}|T| = \delta^{
4\left(\sum_{j=1}^{k-1} c_j^{(k)} \eps_1^{(j)} + c_0^{(k)} \eta^{(k-1)} \right)
}|T|,$$
and 
$$
\angle(\HH(x), \HH_T) \lesssim \delta^{\eta_1/2}
$$
for all $x \in Y_*^{\omega, (k)}(T)$, where $\HH_T \in \Gr(n,1)$ depends only on $T$. Moreover, for each $F \in \F$, if we define 
$$Y^{(k)}(F) \coloneqq \bigsqcup_{\T^{(k)}(F) \cap (\sqcup_{\omega \in \Gamma^{(k)}} \T_*^{\omega, (k)})} Y_*^{\omega, (k)}(T),$$
then $(\F,Y^{(k)})$ is a refinement of $(\F,Y^{(1)})$, and 
\begin{align*}
\sum_{\F}|Y^{(k)}(F)| \gtrapprox \delta^{2\eps_1^{(k)} + 7\eta^{(k)}} = \delta^{ 2\eps_1^{(k)} + 7 \left( \sum_{j=1}^{k-1} c_j^{(k)} \eps_1^{(j)} + c_0^{(k)} \eta^{(k-1)}
\right) }.
\end{align*} 
Now if we take $c_k^{(k+1)} \coloneqq 2$, $c_j^{(k+1)} \coloneqq 7c_j^{(k)}$, and $c_0^{(k+1)} \coloneqq 7 c_0^{(k)} + 2$, then $(\F, Y^{(k)})$ is $\eta^{(k+1)}$-extremal for 
$$\eta^{(k+1)} \coloneqq \sum_{j=1}^{k} c_j^{(k+1)} \eps_1^{(j)} + c_0^{(k+1)} \eta^{(k)}.$$ 
If $k=n-2$, then we stop. Otherwise, we go back to the start of the loop and proceed with $k+1$ in place of $k$.

Let $\eta_0 \coloneqq \eta^{(n-1)}$. Then the output of the above algorithm is an $\eta_0$-extremal refinement $(\F, Y^{(n-2)})$ of $(\F, Y^{(1)})$, where 
$$\eta_0 = \sum_{j=1}^{n-2} c_j^{(n-1)} \eps_1^{(j)} + c_0^{(n-1)} \eta^{(n-2)}.$$ 
By induction, one can easily see that $\eta_0 \geq \eta^{(k)}$ for any $2 \leq k \leq n-2$, and $\eta_0$ is actually a (determinate) linear combination of $\{\eps_1^{(j)}\}_{j=1}^{n-2}$ and $\eta^{(1)}$ with positive integer coefficients. We can also see the previous claimed relation \eqref{Eq:parameter_relation} between various small parameters. 

\begin{rem}
Intuitively, for each $F \in \F$, the shading $Y^{(n-2)}(F)$ is simultaneously well-behaved under $(n-2)$ transverse (to be checked!) tilings of $\delta$-tubes, in the sense that $\HH(x)$ deviates by $\lesssim \delta^{\eta_1/2}$ when we walk along $T \in \T^{(k)}(F) \cap \T_*^{\omega, (k)}$ ($1\leq k \leq n-2$), as long as $x$ remains confined to the set $Y^{(n-2)}(F)$. One may suspect that this already guarantees the existence of some $\HH_F \in \Gr(n,1)$, so that $\angle(\HH(x), \HH_F) \lesssim \delta^{\eta_{1/2}}$ for all $x \in Y^{(n-2)}(F)$, as any two points $x, y \in Y^{(n-2)}(F)$ can reach each other by walking along a broken line aligned with the tilings (like a contagious structure). However, there is a caveat that the corner points of this broken line may not lie entirely in $Y^{(n-2)}(F)$. To settle this problem, we should employ certain backward iteration to construct another suitable refinement $Y_3(F) \subset Y^{(1)}(F)$, which may not have any immediate inclusion relation with $Y^{(n-2)}(F)$. This is why we still need Step~(B)~and~(C).
\end{rem}

{\bf Step~(B).} To construct the desired $(\F, Y_3)$, we will need the transversality of $X_k$. More precisely, since $\angle(\dir(F), \dir(F_D^*)) \leq (100n)^{-1}$ by \eqref{Eq:choice_D}, for each $k$, we can apply Lemma~\ref{Lem:geo_fact_3} with $V = \dir(F)$, $W = \dir(F_D^*)^\perp$, $\ell = \langle v_k \rangle$, and $\rho = (100n)^{-1}$. This tells us that $\langle v_k \rangle$ is obliquely projected onto a line $\ell_k' \in \Gr(\dir(F),1)$ along $\dir(F_D^*)^\perp$, and $\angle(\langle v_k \rangle, \ell_k') \leq (50n)^{-1}$. Note that $\ell_k'$ is exactly $\dir(F) \cap X_k^\perp$, the direction of tubes in $\T^{(k)}(F)$. Since $\{v_k\}_{k=1}^{n-2}$ is orthonormal, we know that $\{\ell_k'\}_{k=1}^{n-2}$ is transverse in the common sense. See Figure~\ref{fig:transversality} for an illustration of 
the special case when $n = 4$. 

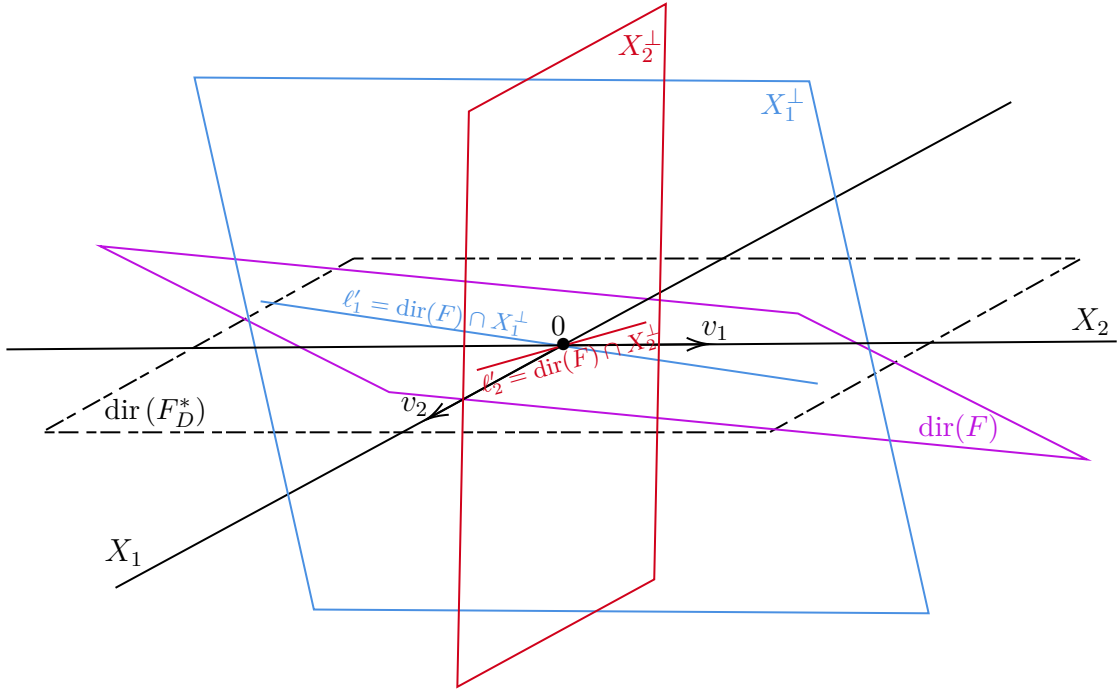
\begin{figure}[htbp]
  \centering
  
    \hspace*{0.5cm}

\tikzset{every picture/.style={line width=0.75pt}} %set default line width to 0.75pt        

\begin{tikzpicture}[x=0.75pt,y=0.75pt,yscale=-1,xscale=1]
%uncomment if require: \path (0,1645); %set diagram left start at 0, and has height of 1645

%Shape: Parallelogram 
\draw  [dash pattern={on 3.75pt off 3pt on 7.5pt off 1.5pt}] (345.5,1156.48) -- (710.66,1156.48) -- (554.16,1243.67) -- (189,1243.67) -- cycle ;
%Shape: Parallelogram 
\draw  [color={rgb, 255:red, 189; green, 16; blue, 224 }  ,draw opacity=1 ] (363.33,1223.51) -- (713.33,1257.27) -- (568.42,1184.07) -- (218.42,1150.31) -- cycle ;
%Straight Lines 
\draw    (449.83,1200.07) -- (383.81,1236.28) ;
\draw [shift={(382.06,1237.24)}, rotate = 331.26] [color={rgb, 255:red, 0; green, 0; blue, 0 }  ][line width=0.75]    (10.93,-3.29) .. controls (6.95,-1.4) and (3.31,-0.3) .. (0,0) .. controls (3.31,0.3) and (6.95,1.4) .. (10.93,3.29)   ;
%Straight Lines 
\draw    (449.83,1200.07) -- (521.55,1199.45) ;
\draw [shift={(523.55,1199.43)}, rotate = 179.5] [color={rgb, 255:red, 0; green, 0; blue, 0 }  ][line width=0.75]    (10.93,-3.29) .. controls (6.95,-1.4) and (3.31,-0.3) .. (0,0) .. controls (3.31,0.3) and (6.95,1.4) .. (10.93,3.29)   ;
%Straight Lines 
\draw [line width=0.75]    (675.43,1078.01) -- (226.06,1321.68) ;
%Shape: Parallelogram 
\draw  [color={rgb, 255:red, 74; green, 144; blue, 226 }  ,draw opacity=1 ] (325.57,1332.57) -- (633.98,1334.46) -- (574.09,1067.58) -- (265.68,1065.69) -- cycle ;
%Straight Lines 
\draw [color={rgb, 255:red, 74; green, 144; blue, 226 }  ,draw opacity=1 ]   (298.86,1177.96) -- (578.35,1219.26) ;
%Straight Lines 
\draw    (171.33,1202.07) -- (728.33,1198.07) ;
%Shape: Parallelogram 
\draw  [color={rgb, 255:red, 208; green, 2; blue, 27 }  ,draw opacity=1 ] (502.12,1028.78) -- (496.32,1317.46) -- (397.54,1371.37) -- (403.33,1082.69) -- cycle ;
%Straight Lines 
\draw [color={rgb, 255:red, 208; green, 2; blue, 27 }  ,draw opacity=1 ]   (492.33,1188.45) -- (407.33,1212.45) ;

% Text Node
\draw (246.22,1233.44) node  [font=\small]  {$\mathrm{dir}\left( F_{D}^{*}\right)$};
% Text Node
\draw (649.27,1242.05) node  [font=\small,color={rgb, 255:red, 189; green, 16; blue, 224 }  ,opacity=1 ]  {$\mathrm{dir}( F)$};
% Text Node
\draw (230.38,1304.12) node    {$X_{1}$};
% Text Node
\draw (527.01,1192.57) node    {$v_{1}$};
% Text Node
\draw (376.08,1230.61) node    {$v_{2}$};
% Text Node
\draw (561.01,1079.45) node  [font=\small,color={rgb, 255:red, 74; green, 144; blue, 226 }  ,opacity=1 ]  {$X_{1}^{\perp }$};
% Text Node
\draw (448.18,1190.44) node  [font=\small]  {$0$};
% Text Node
\draw (450.74,1199.85) node    {$\smblkcircle $};
% Text Node
\draw (387.08,1183.58) node  [font=\scriptsize,color={rgb, 255:red, 74; green, 144; blue, 226 }  ,opacity=1 ,rotate=-8.45]  {$\ell _{1} '=\mathrm{dir}( F) \cap X_{1}^{\perp }$};
% Text Node
\draw (715.38,1188.12) node    {$X_{2}$};
% Text Node
\draw (489.43,1049.94) node  [font=\small,color={rgb, 255:red, 208; green, 2; blue, 27 }  ,opacity=1 ]  {$X_{2}^{\perp }$};
% Text Node
\draw (455.08,1207.58) node  [font=\scriptsize,color={rgb, 255:red, 208; green, 2; blue, 27 }  ,opacity=1 ,rotate=-344.2]  {$\ell _{2} '=\mathrm{dir}( F) \cap X_{2}^{\perp }$};

\end{tikzpicture}
\caption{This can be seen as a continuation of Figure~\ref{fig:slicing}. Here $\langle v_1 \rangle$ and $\langle v_2 \rangle$ are  obliquely projected to $\ell_1'$ and $\ell_2'$, respectively.}
\label{fig:transversality}
\end{figure}

\begin{rem}
Intuitively, the whole collection of slices $\bigsqcup_{k=1}^{n-2} \{\Pi_\omega^{(k)}: \omega \in \Lambda_k \}$ perfectly divides $\dir(F_D^*)$ into a lattice $\Lambda$ of mesh size $\delta$, whose fundamental domain is the $\delta$-cube centered at $0$ and aligned with $\{v_j\}_{j=1}^{n-2}$. It also divides each $\dir(F)$ into a tilted lattice $\Lambda_F$, which is aligned with the whole collection of tubes 
$\bigsqcup_{k=1}^{n-2} \T^{(k)}(F)$. Note that $\Lambda_F$ is the oblique projection of $\Lambda$ onto $\dir(F)$ along $\dir(F_D^*)^\perp$. Since $\angle(\dir(F), \dir(F_D^*)) \leq (100n)^{-1}$, we should also expect $\Lambda_F$ to be close to $\Lambda$, and so has its fundamental domain comparable to a $\delta$-cube.
\end{rem}

Now for each $F \in \F$ with $Y^{(n-2)}(F) \neq \varnothing$, there exists some $\omega_0 \in \Gamma^{(n-2)}$ and $T_0 \in \T^{(n-2)}(F) \cap \T_*^{\omega_0, (n-2)}$ with $|Y_*^{\omega_0, (n-2)}(T_0)| \gtrapprox \delta^{4\eta^{(n-2)}}|T_0|$. Fix one such $T_0$ and let $E_F^{(n-2)} \coloneqq Y_*^{\omega_0, (n-2)}(T_0)$. 

We run the following iterative algorithm. Let $k$ be the loop index, with $n-2$ as its initial value. 

By transversality of $\T^{(k)}(F)$ and $\T^{(k-1)}(F)$, we know from $E_F^{(k)} \subset Y^{(k-1)}(F)$ that the set  
\begin{align*}
\Gamma_{F}^{(k-1)} \coloneqq \{ \omega \in \Gamma^{(k-1)}: (Y^{(k-1)}(F) \cap \Pi_{\omega}^{(k-1)}) \cap E_F^{(k)} \neq \varnothing \}
\end{align*}
has cardinality $\gtrsim \delta^{-n}|E_F^{(k)}|$. Note that if $\omega \in \Gamma_F^{(k-1)}$, then $Y^{(k-1)}(F) \cap \Pi_{\omega}^{(k-1)} \neq \varnothing$, which means $F \cap \Pi_\omega^{(k-1)} \in \T_*^{\omega,(k-1)}$ and $Y^{(k-1)}(F) \cap \Pi_\omega^{(k-1)} = Y_*^{\omega,(k-1)}(F \cap \Pi_\omega^{(k-1)})$ by the construction of $Y^{(k-1)}$. Define 
\begin{align*}
E_F^{(k-1)} \coloneqq 
\bigsqcup_{\omega \in \Gamma_{F}^{(k-1)}} Y_*^{\omega,(k-1)}(F \cap \Pi_\omega^{(k-1)}).
\end{align*}
Then $E_F^{(k-1)} \subset Y^{(k-1)}(F)$, and 
\begin{align*}
|E_F^{(k-1)}| \gtrapprox \delta^{-n}|E_F^{(k)}| \cdot \delta^{4\eta^{(k-1)}} |T| 
\sim 
\delta^{-1 + 4\eta^{(k-1)}} |E_F^{(k)}|.
\end{align*}
If $k=2$, then we stop. Otherwise, we go back to the start of the loop and proceed with $k - 1$ in place of $k$.

The output of the above algorithm is a set $E_F^{(1)} \subset F$ satisfying 
\begin{align*}
|E_F^{(1)}| \gtrapprox |E_F^{(n-2)}| \cdot \prod_{j=2}^{n-2} \delta^{-1+4\eta^{(j-1)}} 
\gtrapprox 
\delta|T_0| \cdot \prod_{j=2}^{n-1} \delta^{-1+4\eta^{(j-1)}} \sim \delta^{4\sum_{j=1}^{n-2} \eta^{(j)}} \delta^2 \geq \delta^{4(n-2)\eta_0} |F|.
\end{align*}
Thus for each $F \in \F$ with $Y^{(n-2)}(F) \neq \varnothing$, we can define $Y_3(F) \coloneqq E_F^{(1)}$. What we just showed is $|Y_3(F)| \gtrapprox \delta^{4(n-2) \eta_0} |F|$. For technical clarity, let us trivially extend the definition of $Y_3$ to the whole $\F$ by setting $Y_3(F) \coloneqq \varnothing$ for each $F$ with $Y^{(n-2)}(F) = \varnothing$. 

Recall that $(\F, Y^{(n-2)})$ is $\eta_0$-extremal, which implies $$\sum_\F |Y^{(n-2)}(F)| = \lambda(\F, Y^{(n-2)}) |\F||F| \geq \delta^{2\eta_0}.$$
So we must have 
$$|\{F \in \F: Y^{(n-2)}(F) \neq \varnothing\}| \geq \delta^{-2 + 2\eta_0}.$$ 
Therefore, we can conclude that 
\begin{align*}
\sum_{\F} |Y_3(F)| 
\gtrapprox \delta^{-2 + 2\eta_0} \cdot 
\delta^{4(n-2)\eta_0} |F| = \delta^{(4n-6)\eta_0}.
\end{align*}
In particular, $(\F,Y_3)$ is $4n\eta_0$-extremal. Thus the first bullet point holds.

{\bf Step~(C).} It remains to check the desired concentration property of the plane map $\HH$ with respect to $(\F, Y_3)$. By the construction of $Y_3$, if $Y_3(F) \neq \varnothing$\footnote{If $Y_3(F) = \varnothing$, then there is nothing to do.}, then $Y^{(n-2)}(F) \neq \varnothing$, and $Y_3(F)$ grows out of some $T_0$ (depending on $F$) with $Y_*^{\omega_0, (n-2)}(T_0) = E_F^{(n-2)}$. Fix any point $x_0 \in Y_*^{\omega_0, (n-2)}(T_0) \subset Y^{(n-2)}(F)$. For each $x \in Y_3(F) = E_F^{(1)}$, we can use the following iterative algorithm to find an admissible broken line connecting $x$ to $x_0$. Let $x^{(1)} \coloneqq x$. Let $k$ be the loop index, with $1$ as its initial value.

By the definition of $E_F^{(k)}$, there is some $\omega \in \Gamma_F^{(k)}$ such that $x^{(k)} \in Y^{(k)}(F) \cap \Pi_\omega^{(k)} = Y_*^{\omega,(k)}(F \cap \Pi_\omega^{(k)})$, and $Y_*^{\omega,(k)}(F \cap \Pi_\omega^{(k)}) \cap E_F^{(k+1)} \neq \varnothing$. Let $T^{(k)} \coloneqq F \cap \Pi_\omega^{(k)}$. Pick some $x^{(k+1)} \in Y_*^{\omega,(k)}(T^{(k)}) \cap E_F^{(k+1)}$. By the stability of $\HH$ along $T^{(k)} \in \T^{(k)}(F)$ (Step~(vii)), we have 
\begin{align*}
\angle(\HH(x^{(k)}), \HH(x^{(k+1)})) \leq \angle(\HH(x^{(k)}), \HH_{T^{(k)}}) + \angle(\HH_{T^{(k)}}, \HH(x^{(k+1)})) \lesssim \delta^{\eta_1/2}.
\end{align*}
If $k=n-3$, then we stop. Otherwise, we go back to the start of the loop and proceed with $k + 1$ in place of $k$. 

The output of the above algorithm is a sequence 
$$x = x^{(1)} \xrightarrow{T^{(1)}} x^{(2)} \xrightarrow{T^{(2)}} \cdots \xrightarrow{T^{(n-3)}} x^{(n-2)}.$$
By the triangle inequality, we have 
\begin{align*}
\angle(\HH(x), \HH(x^{(n-2)})) \leq \sum_{k=1}^{n-3} \angle(\HH(x^{(k)}), \HH(x^{(k+1)})) \lesssim \delta^{\eta_1/2}.
\end{align*}
Also, note that $x^{(n-2)} \in Y_*^{\omega,(n-3)}(T^{(n-3)}) \cap E_F^{(n-2)} \subset E_F^{(n-2)} \subset T_0$, and so by the stability of $\HH$ along $T_0 \in \T^{(n-2)}(F)$ (Step~(vii)), we have 
\begin{align*}
\angle(\HH(x^{(n-2)}), \HH(x_0)) \leq \angle(\HH(x^{(n-2)}), \HH_{T_0}) + \angle(\HH_{T_0}, \HH(x_0))  \lesssim \delta^{\eta_1/2}. 
\end{align*}

Finally, define 
$$\HH_F \coloneqq \pi_{\dir(F)^\perp}(\HH(x_0)) \in \Gr(\dir(F)^\perp,1),$$ which depends only on $F$. Since $\dir(F) \cap B_1 \subset N_{2\delta}(\HH(x_0)^\perp)$, by applying Lemma~\ref{Lem:faithful_proj} with $V = \HH(x_0)^\perp$, $W = \dir(F)$, and $\rho = 2\delta$, we get $\HH(x_0) \cap B_1 \subset N_{2\delta}(\dir(F)^\perp)$, which by Lemma~\ref{Lem:faithful_proj} again with $V = \dir(F)^\perp$, $W = \HH(x_0)$, and $\rho = 2\delta$ implies $$\angle(\HH_F, \HH(x_0)) \leq 2\delta \leq \delta^{\eta_1/2}.$$ 
So by the triangle inequality, we get $$\angle(\HH(x), \HH_F) \leq \angle(\HH(x), \HH(x^{(n-2)})) + \angle(\HH(x^{(n-2)}), \HH(x_0)) + \angle(\HH(x_0), \HH_F) \lesssim \delta^{\eta_1/2},$$which uniformly holds for all $x \in Y_3(F)$. Thus the second bullet point holds.
\end{proof}

For $0 < \eta \ll \eta_0 \ll \eta_1 \ll \sigma$, let $(\F,Y_3)$ be the output of Proposition~\ref{Prop:step_8}. 

This finishes Step~(viii). 

From now on, we can forget about the slicing structures, because Steps~(ix)~and~(x) depend only on the concentration property of $\HH$ we get in this step, which is preserved under any refinement of $(\F,Y_3)$.

{\bf Step~(ix).} 
Now we are ready to grow the flakes into wider planks (up to a refinement). More precisely, we will prove the following result: 
\begin{proposition}\label{Prop:step_9}
Suppose $0 < \eta \ll\eps_3 \ll \eta_1 \ll \sigma$, then there exists $\delta_0 > 0$ so that the following holds for all $\delta \in (0,\delta_0]$. 

There exists $\F_3 \subset \F$ with $|\F_3| \geq \delta^{\eps_3}|\F|$, such that each $F \in \F_3$ is contained in a plank $L_F$ of dimensions $\delta \times \delta^{1-\eta_1/2} \times \1$ with $|U(\F, Y) \cap L_F| \gtrapprox \delta^{5\eps_3}|L_F|$.
\end{proposition}

\begin{proof}
The idea is to piece together the local grains along a typical $\delta$-flake, whose orientations are concentrated in a $\lesssim \delta^{\eta_1/2}$ range by Step~(viii). To make this rigorous, we need some additional pigeonholing steps to make the shading in each local grain regular. 

Let $(\F,Y_3')$ be the constant multiplicity refinement of $(\F,Y_3)$. Note that $(\F,Y_3')$ is $4n\eta_0$-extremal, so by the proof of Lemma~\ref{Lem:grains_number_given_scale} and by {\color{black} taking $\eta_0 \ll \eps_3$} (which can be ensured by taking 
$$\eta \ll \eps_1 \ll \eps_1^{(2)} \ll \cdots \ll \eps_1^{(n-2)} \ll \eps_3$$ 
as in \eqref{Eq:parameter_relation}), we get a $\gtrapprox 1$ refinement $(\F,Y_4)$ of $(\F,Y_3')$, and an $\eps_3$-extremal cover $(\F_{\delta^{1/2}}, Y_{\delta^{1/2}})$ of it, which is formed by unions of $(\delta^{1/2}/10n)$-cubes. By similar arguments as in the proof of Lemma~\ref{Lem:grains_single_rho_tau}, $U(\F_{\delta^{1/2}}, Y_{\delta^{1/2}})$, and thus $U(\F,Y_4)$, is contained in a union of $(\delta^{1/2}/10n)$-cubes of cardinality $\lesssim \delta^{(\sigma-n)/2 - \eps_3/2}$. Let $\Bar{\Q}$ be this set of cubes. By a further harmless $\gtrsim 1$ refinement, we may assume that $\Bar{\Q}$ is $10n\delta^{1/2}$-separated. Let $$\bar{\Q}' \coloneqq \{\bar{Q} \in \bar{\Q}: \int_{\bar{Q}}\mu(\F, Y_4)(x) \D x \geq \delta^{(n-\sigma)/2 + \eps_3} \},$$ then 
\begin{align*}
\int_{U(\bar{\Q} \setminus \bar{\Q}')} \mu(\F, Y_4)(x) \D x 
&\leq 
\delta^{(\sigma - n)/2 - \eps_3/2} \delta^{(n-\sigma)/2 + \eps_3} \\
&= 
\delta^{\eps_3/2} 
\leq 
\delta^{4n\eta_0 + \eta} 
\lessapprox 
\delta^{\eta} \sum_{\F}|Y_3(F)|
\lessapprox
\delta^{\eta} \sum_{\F}|Y_4(F)|, 
\end{align*}
as $\eta$ and $\eta_0$ are very small compared to $\eps_3$. 
So if we define $Y_5(F) \coloneqq Y_4(F) \cap U(\bar{\Q}')$ for each $F \in \F$, then $(\F, Y_5)$ is a $\gtrsim 1$ refinement of $(\F, Y_4)$. 

We will follow the proof of Lemma~\ref{Lem:grains_single_rho_tau} to find local grains with dense shadings. For each $\bar{Q} \in \bar{\Q}$, fix a point $x_{\bar{Q}} \in \bar{Q} \cap U(\F,Y_5) \neq \varnothing$. Then 
$$\pi_{\HH(x_{\bar{Q}})}(U(\F,Y_5) \cap \bar{Q}) \subset \pi_{\HH(x_{\bar{Q}})}(U(\F,Y) \cap B(x_{\bar{Q}},\delta^{1/2}))$$ 
is a $(\delta, 1-\sigma, \delta^{-\eta})_1$-ADset by the planiness and the local grains structure of the original $(\F,Y)$. In particular, $U(\F,Y_5) \cap \bar{Q}$ is contained in at most $\delta^{-\eta}\delta^{(\sigma-1)/2}$ local grains of the form $\bar{Q} \cap N_{\delta/2}(t + \HH(x_{\bar{Q}})^\perp)$, where $t \in \HH(x_{\bar{Q}})$ and $\norm{t} \in \delta\Z$. In other words, we have 
$$U(\F, Y_5) \cap \bar{Q} \subset \bigsqcup_{j=1}^{J_{\bar{Q}}} G_{\bar{Q}, j},$$ 
where $J_{\bar{Q}} \lesssim \delta^{-\eta} \delta^{(\sigma-1)/2}$ and 
$$G_{\bar{Q}, j} \coloneqq \bar{Q} \cap N_{\delta/2}(t_j + \HH(x_{\bar{Q}})^\perp)$$ 
for some $t_j \in \HH(x_{\bar{Q}})$ with $\norm{t_j} = \delta\Z$. Let $\G(\bar{Q}) \coloneqq \{ G_{\bar{Q}, j} \}_{j=1}^{J_{\bar{Q}}}$ and $$\G(\bar{Q})' \coloneqq \{ G_{\bar{Q}, j} \in \G(\bar{Q}): \int_{G_{\bar{Q}, j}} \mu(\F, Y_5)(x) \D x \geq \delta^{(n+1 - 2\sigma)/2 + 2\eps_3} \},$$ 
then 
\begin{align*}
\int_{U(\G(\bar{Q}) \setminus \G(\bar{Q})')} \mu(\F, Y_5)(x) \D x  
&\lesssim \delta^{-\eta}\delta^{(\sigma-1)/2} \delta^{(n+1 - 2\sigma)/2 + 2\eps_3}\\
&= 
\delta^{\eps_3-\eta} \delta^{(n-\sigma)/2 + \eps_3} 
\leq \delta^\eta \int_{\bar{Q}}\mu(\F, Y_5)(x) \D x.
\end{align*}
This means that by only retaining the mass of $(\F, Y_5)$ over all $U(\G(\bar{Q})')$, we obtain a $\gtrsim 1$ refinement $(\F, Y_6)$ such that for each $(\delta^{1/2}/10n)$-cube $\bar{Q}$ intersecting $U(\F, Y_6)$, we have $U(\F, Y_6) \cap \bar{Q} \subset U(\G({\bar{Q}})')$ and 
$$\int_{\bar{Q}}\mu(\F, Y_6)(x) \D x 
\gtrsim \delta^{(n-\sigma)/2 + \eps_3}.$$
Moreover, for each $G_{\bar{Q}, j} \in \G(\bar{Q})'$, we have 
\begin{align}\label{Eq:fill_out_grain'}
\int_{G_{\bar{Q}, j}} \mu(\F, Y_6)(x) \D x \geq \delta^{(n+1 - 2\sigma)/2 + 2\eps_3}.
\end{align}

Therefore, if we define 
$$\V_{\bar{Q},j} \coloneqq \{F \cap G_{\bar{Q},j}: F \in \F\}$$ 
and 
$$Y_{\bar{Q},j}(F \cap G_{\bar{Q},j}) \coloneqq Y_6(F) \cap G_{\bar{Q},j},$$
then $(\V_{\bar{Q},j}, Y_{\bar{Q},j})$ is supported on $G_{\bar{Q},j}$, and \eqref{Eq:fill_out_grain'} tells that 
\begin{align*}
    \int_{G_{\bar{Q},j}} \mu(\V_{\bar{Q},j}, Y_{\bar{Q},j})(x) \D x \geq 
    \delta^{(n+1 - 2\sigma)/2 + 2\eps_3}. 
\end{align*}
On the other hand, by Remark~\ref{Rem:extreme_mult} and by {\color{black} taking $\eta_0 \ll \eps_3$}, we have $\mu(\F,Y_3') \leq \delta^{-\sigma-2\eps_3}$. Let $(\V_{\bar{Q},j}, Y_{\bar{Q},j}')$ be the constant multiplicity refinement of $(\V_{\bar{Q},j}, Y_{\bar{Q},j})$, then $$\mu(\V_{\bar{Q},j}, Y_{\bar{Q},j}')(x) \leq \mu(\V_{\bar{Q},j}, Y_{\bar{Q},j})(x) = \mu(\F,Y_6)(x) \leq \mu(\F, Y_3')(x) = \mu(\F,Y_3') \leq \delta^{-\sigma-2\eps_3}$$ 
for each $x \in G_{\bar{Q},j}$ (recall that $(\F,Y_3')$ has constant multiplicity). Apply Lemma~\ref{Lem:regular_subshading} to each $\Omega = G_{\bar{Q},j}$ and $Y(\Omega) = U(\V_{\bar{Q},j}, Y_{\bar{Q},j}')$, we get a subset $F_{\bar{Q},j} \subset U(\V_{\bar{Q},j}, Y_{\bar{Q},j}')$ that is regular and satisfies $|F_{\bar{Q},j}| \geq (1/2)|U(\V_{\bar{Q},j}, Y_{\bar{Q},j}')|$. Define $$Y_{\bar{Q},j}''(V) \coloneqq Y_{\bar{Q},j}'(V) \cap F_{\bar{Q},j}.$$ 
Then $(\V_{\bar{Q},j}, Y_{\bar{Q},j}'')$ is a $\gtrapprox 1$ refinement of $(\V_{\bar{Q},j}, Y_{\bar{Q},j})$. Note that 
\begin{align*}
|F_{\bar{Q},j}| = |U(\V_{\bar{Q},j}, Y_{\bar{Q},j}'')| \gtrapprox \frac{\sum_{\V_{\bar{Q},j}}|Y_{\bar{Q},j}(V)|}{\mu(\V_{\bar{Q},j}, Y_{\bar{Q},j}')} \geq 
\frac{\delta^{(n+1 - 2\sigma)/2 + 2\eps_3}}{\delta^{-\sigma - 2\eps_3}} = \delta^{(n+1)/2 + 4\eps_3}.
\end{align*}

Now if we define 
$$Y_7(F) \coloneqq \bigsqcup_{\bar{Q} \in \bar{\Q}'} \bigsqcup_{j=1}^{J_{\bar{Q}}} Y_{\bar{Q},j}''(F \cap G_{\bar{Q},j})$$
for each $F \in \F$, then $(\F,Y_7)$ is a $\gtrapprox 1$ refinement of $(\F,Y_6)$ (and thus of $(\F,Y_3)$). Moreover, 
by regularity of $F_{\bar{Q},j}$, for each $F \in \F$, if $Y_7(F) \cap G_{\bar{Q},j} \neq \varnothing$, then for any $x \in Y_7(F) \cap G_{\bar{Q},j}$ and $r \in [\delta,1]$, we have 
\begin{equation}
\begin{split}
|U(\F,Y_7) \cap G_{\bar{Q},j} \cap B(x,r)| &= |F_{\bar{Q},j} \cap B(x,r)| \\
&\gtrapprox |F_{\bar{Q},j}| \Big( \frac{|G_{\bar{Q},j} \cap B(x,r)|}{|G_{\bar{Q},j}|} \Big) \\
&\gtrapprox 
\delta^{(n+1)/2 + 4\eps_3} \Big( \frac{|G_{\bar{Q},j} \cap B(x,r)|}{\delta^{(n+1)/2}} \Big) \\
&= 
\delta^{4\eps_3} |G_{\bar{Q},j} \cap B(x,r)|. 
\end{split}\label{Eq:regular_local_grain}
\end{equation} 
Since $\sum_{\F}|Y_7(F)| \gtrapprox \sum_{\F}|Y_3(F)| \gtrapprox \delta^{4n\eta_0}$ and $|\F| \leq \delta^{-2-\eta}$, there is a subset $\F_3 \subset \F$, such that $(\F_3,Y_7)$ is a $\gtrsim 1$ refinement of $(\F,Y_7)$, 
and for each $F \in \F_3$, we have 
$|Y_7(F)| \geq \delta^{4n\eta_0 + 2\eta}|F|$. We also have $|\F_3| \gtrapprox \delta^{-2 + 4n\eta_0} \geq \delta^{4n\eta_0 + \eta} |\F|$. As $\eta, \eta_0 \ll \eps_3$, we conclude that $|Y_7(F)| \geq \delta^{\eps_3}|F|$ for each $F \in \F_3$, and $|\F_3| \geq \delta^{\eps_3}|\F|$.

It remains to check that each $F \in \F_3$ grows into a much wider plank. To see this, note that each $x \in Y_7(F)$ must lie in some $(\delta^{1/2}/10n)$-cube $\bar{Q} \in \Bar{\Q}'$, and $\angle(\HH(x), \HH_F) \lesssim \delta^{\eta_1/2}$ by Step~(viii). Since $\norm{x - x_{\bar{Q}}} \leq \delta^{1/2}$ and $\HH$ is $1$-Lipschitz over $U(\F,Y)$, we have $\angle(\HH(x), \HH(x_{\bar{Q}})) \leq \delta^{1/2}$, and so 
$$\angle(\HH(x_{\bar{Q}}), \HH_F) \lesssim \delta^{\eta_1/2}.$$ 
Pick $G_{\bar{Q},j} \in \G(\bar{Q})'$ such that $x \in Y_7(F) \cap G_{\bar{Q},j}$, then \eqref{Eq:regular_local_grain} with $r = \delta^{1 - \eta_1/2}$ implies 
\begin{align*}
|U(\F,Y_7) \cap S_{F,x}| \gtrapprox 
\delta^{4\eps_3} |S_{F,x}|, 
\end{align*}
where 
$$S_{F,x} \coloneqq G_{\bar{Q},j} \cap B(x,\delta^{1 - \eta_1/2})$$ 
is a slab segment of dimensions $\delta \times \delta^{1 - \eta_1/2} \times \delta^{1 - \eta_1/2}\1$. Since the normal direction of the central plane of $S_{F,x}$ is $\HH(x_{\bar{Q}})$ and $\angle(\HH(x_{\bar{Q}}), \HH_F) \lesssim \delta^{\eta_1/2}$, there is a slab segment $\bar{S}_{F,x}$ containing and comparable to $S_{F,x}$, so that the normal direction of $\bar{S}_{F,x}$ is exactly $\HH_F$. We may assume that $\bar{S}_{F,x}$ is centered at $x$. Note that 
\begin{align}\label{Eq:dense_slab_seg}
|U(\F,Y_7) \cap \bar{S}_{F,x}| \gtrapprox 
\delta^{4\eps_3} |\bar{S}_{F,x}| 
\end{align} 
still holds. As $|Y_7(F)| \geq \delta^{\eps_3}|F|$, there exists a $\geq C\delta^{1-\eta_1/2}$-separated set of points $\{x_k\}_k \subset Y_7(F)$ with cardinality $\gtrsim \delta^{-(n-2)(1-\eta_1/2) + \eps_3}$. Here $C \lesssim 1$ is chosen so that $\{\bar{S}_{F,x_k}\}_k$ are disjoint. By taking the union of all $\{\bar{S}_{F,x_k}\}_k$, we can grow $F$ into a much wider plank $L_F$. More precisely, let 
$$\ell_F \coloneqq \HH_F^\perp \cap \dir(F)^\perp \in \Gr(\dir(F)^\perp,1),$$ then there is a plank $L_F$ that satisfies the following properties:
\begin{itemize}
    \item $L_F$ contains $F$ and has central hyperplane parallel to $\HH_F^\perp = \dir(F) \oplus \ell_F$.
    \item $L_F$ contains all $\{\bar{S}_{F,x_k}\}_k$. 
    \item $L_F$ has dimensions $\delta \times \delta^{1-\eta_1/2} \times \1$, where $\delta$ is in $\HH_F$, $\delta^{1-\eta_1/2}$ is in $\ell_F$, and $\1$ is in $\dir(F)$.
\end{itemize}
The second bullet point above together with \eqref{Eq:dense_slab_seg} implies that 
\begin{align*}
|U(\F,Y_7) \cap L_F| &\geq \big|U(\F,Y_7) \cap \bigsqcup_k \bar{S}_{F,x_k}\big| \\
&= 
\sum_k |U(\F,Y_7) \cap \bar{S}_{F,x_k}| \\
&\gtrapprox 
\delta^{-(n-2)(1-\eta_1/2) + \eps_3} \cdot \delta^{4\eps_3} \delta^{(n-1)(1-\eta_1/2)+1} \\
&\sim  
\delta^{5\eps_3} |L_F|.
\end{align*} 
Thus for each $F \in \F_3$, we have $|U(\F, Y) \cap L_F| \gtrapprox \delta^{5\eps_3}|L_F|$, as desired.
\end{proof}

For $0 < \eta \ll\eps_3 \ll \eta_1 \ll \sigma$, let $\F_3$ and $\{L_F\}_{F \in \F_3}$ be the output of Proposition~\ref{Prop:step_9}.

Define $\LL \coloneqq \{L_F: F \in \F_3\}$ and $Y_\LL(L_F) \coloneqq U(\F, Y) \cap L_F$. Then $|Y_\LL(L_F)| \gtrapprox \delta^{5\eps_3}|L_F|$ 
for each $L_F \in \LL$. 

This finishes Step~(ix) with $b' \coloneqq \delta^{1-\eta_1/2}$.

{\bf Step~(x).}
In this final step, we will obtain a contradiction by further growing each $L_F$ into a $b'$-flake and using a local Córdoba-type $L^2$ argument. This contradiction will immediately lead to the validity of Theorem~\ref{Thm:Main_Sticky}.

Recall that $|\F_3| \geq \delta^{\eps_3} |\F|$. By Lemma~\ref{Lem:split_sticky} with $b'$ in place of $\rho$ and $\eta_1/2$ in place of $\eps$, we know that $\F_3$ contains a uniform subset $\F_4 \subset \F_3$ with $|\F_4| \gtrapprox |\F_3|$ that is $\delta^{-(\eps_3 + 4\eta)}$-sticky. Moreover, we have \begin{itemize}
    \item 
    $\F_4$ is uniform at scale $b'$, and $(\F_4)_{b'}$ is a ${b'}^{-(2\eps_3 + 10\eta)/\eta_1}$-sticky set of $b'$-flakes. Since $(\F_4)_{b'}$ is spaced by the proof of Lemma~\ref{Lem:split_sticky}, if we thicken each $b'$-flake $F_{b'} \in (\F_4)_{b'}$ to a $Cb'$-flake $F_{Cb'}$ with respect to its central plane (for $C \lesssim 1$), then the resulting set of $Cb'$-flakes $(\F_4)_{Cb'}$ (slightly abusing notation) is $(Cb')^{-(2\eps_3 + 11\eta)/\eta_1}$-sticky.

    \item  
    For each $F_{b'} \in (\F_4)_{b'}$, let $\widehat{\F}(F_{b'})$ be the image of the flakes in $\F_4[F_{b'}]$ under the rescaling taking $F_{b'}$ to $B_1$, then $\widehat{\F}(F_{b'})$ is a $(\delta/b')^{-(2\eps_3 + 8\eta)/\eta_1}$-sticky set of $(\delta/b')$-flakes. Since $(\F_4)_{b'}$ is spaced by the proof of Lemma~\ref{Lem:split_sticky}, if we rescale $\F_4[F_{Cb'}]$ by taking $F_{Cb'}$ to $B_1$ instead, then the resulting set $\widehat{\F}(F_{Cb'})$ of $(\delta/b')$-flakes, is $(\delta/b')^{-(2\eps_3 + 9\eta)/\eta_1}$-sticky.
\end{itemize} 
We technically introduce the constant $C \lesssim 1$ above to make sure that if $F \in \F_4$ is contained in some $F_{b'}$, then $L_F$ is fully contained in $F_{Cb'}$.

For each $F_{Cb'} \in (\F_4)_{Cb'}$, let $\LL(F_{Cb'}) \coloneqq \{L_F: F \in \F_4[F_{b'}]\}$.  Under the unit rescaling relative to $F_{Cb'}$ as in the second bullet point above, the set of $\delta \times b' \times \1$ planks $(\LL(F_{Cb'}), Y_\LL)$ becomes a set of $(\delta/b') \times 1 \times \1$-slabs $(\widehat{\LL}(F_{Cb'}), \widehat{Y}_{\widehat{\LL}(F_{Cb'})})$. By Step~(ix), we have $|\widehat{Y}_{\widehat{\LL}(F_{Cb'})}(\widehat{L})| \gtrapprox \delta^{5\eps_3}|\widehat{L}|$ for each $\widehat{L} \in \widehat{\LL}(F_{Cb'})$. Since $\widehat{\F}(F_{Cb'})$ is $(\delta/b')^{-(2\eps_3 + 9\eta)/\eta_1}$-sticky, we know from Remark~\ref{Rem:strong_sticky} that $C_F(\widehat{\F}(F_{Cb'})) \leq (\delta/b')^{-(4\eps_3 + 18\eta)/\eta_1}$. 
On the other hand, note that each $\widehat{F} \in \widehat{\F}(F_{Cb'})$ is contained in some $\widehat{L} \in \widehat{\LL}(F_{Cb'})$ (as it is before rescaling), so we can apply 
Proposition~\ref{Prop:inherit}(A) to get a subset $\widehat{\LL}(F_{Cb'})' \subset \widehat{\LL}(F_{Cb'})$ with 
$$C_F(\widehat{\LL}(F_{Cb'})') \lessapprox C_F(\widehat{\F}(F_{Cb'})) \leq (\delta/b')^{-(4\eps_3 + 18\eta)/\eta_1}.$$ 
Also observe that 
\begin{align*}
\lambda(\widehat{\LL}(F_{Cb'})', \widehat{Y}_{\widehat{\LL}(F_{Cb'})}) \gtrapprox \delta^{5\eps_3}.
\end{align*} 
Therefore, by applying Lemma~\ref{Lem:slab_0} (which does not require $\widehat{\LL}(F_{Cb'})'$ to be essentially distinct) to the set of $(\delta/b')$-slabs $(\widehat{\LL}(F_{Cb'})', \widehat{Y}_{\widehat{\LL}(F_{Cb'})})$, we obtain 
\begin{align*}
|U(\widehat{\LL}(F_{Cb'}), \widehat{Y}_{\widehat{\LL}(F_{Cb'})})| &\geq |U(\widehat{\LL}(F_{Cb'})', \widehat{Y}_{\widehat{\LL}(F_{Cb'})})| \\ 
&\gtrapprox 
C_F(\widehat{\LL}(F_{Cb'})')^{-1} \lambda(\widehat{\LL}(F_{Cb'}), \widehat{Y}_{\widehat{\LL}(F_{Cb'})})^2 \\
&\gtrapprox 
(\delta/b')^{(4\eps_3 + 18\eta)/\eta_1}\delta^{10\eps_3} \\
&= 
(\delta/b')^{(24\eps_3 + 18\eta)/\eta_1} \\
&\geq 
(Cb')^{(24\eps_3 + 18\eta)/\eta_1}, 
\end{align*}
where we also used $b' = \delta^{1-\eta_1/2}$ and $\delta/b' \geq Cb'$. 
By rescaling things back, this implies 
\begin{align*}
|U(\LL(F_{Cb'}), Y_\LL)| \gtrapprox (Cb')^{(24\eps_3 + 18\eta)/\eta_1}|F_{Cb'}|. 
\end{align*}

If we define a shading $\bar{Y}_{Cb'}$ on $(\F_4)_{Cb'}$ by $\bar{Y}_{Cb'}(F_{Cb'}) \coloneqq U(\LL(F_{Cb'}), Y_\LL) \cap F_{Cb'}$ for each $F_{Cb'} \in (\F_4)_{Cb'}$, then $\lambda((\F_4)_{Cb'}, \bar{Y}_{Cb'}) \gtrapprox (Cb')^{(24\eps_3 + 18\eta)/\eta_1}$. This together with the first bullet point means that $((\F_4)_{Cb'}, \bar{Y}_{Cb'})$ satisfies Items~\ref{item_a}~and~\ref{item_b} in Definition~\ref{Def:sigma} with $\delta = Cb'$ and $\eps = (24\eps_3 + 18\eta)/\eta_1$. By Remark~\ref{Rem:extreme_mult}, if we {\color{black} choose $\max\{\eta, \eps_3\} \ll \min\{\eta_1, \eta_2\}$}, then \begin{align*}
|U((\F_4)_{Cb'}, \bar{Y}_{Cb'})| \geq (Cb')^{\sigma + \eta_2} \geq \delta^{(1-\eta_1/2)(\sigma+\eta_2)} = \delta^{\sigma - (\eta_1\sigma/2 + \eta_1\eta_2/2 - \eta_2)}.
\end{align*}
{\color{black} By choosing $\eta_2 \leq \eta_1 \sigma/4$,} the lower bound becomes $|U((\F_4)_{Cb'}, \bar{Y}_{Cb'})| \geq \delta^{\sigma - \eta_1\sigma/4}$. 
On the other hand, since $U((\F_4)_{Cb'}, \Bar{Y}_{Cb'}) \subset U(\F, Y)$ and $(\F, Y)$ is $\eta$-extremal, we have the upper bound $|U(\F,Y)| \leq \delta^{\sigma - \eta}$. {\color{black} By choosing $\eta \leq \eta_1\sigma/8$} and then letting $\delta_0$ be small enough, we obtain a contradiction between the upper and lower bounds. This finishes Step~(x). 

In summary, under the initial assumption that $\sigma > 0$, we always get a contradiction by carrying out all Steps~(i)--(x). In other words, we must have $\sigma = 0$ at the very beginning. {\color{black} This completes the proof of Theorem~\ref{Thm:Main_Sticky}.}

\section{Proof of Theorem~\ref{Thm:2-plane_sticky}}\label{Sec:finish_sticky}

{\color{black} In this section, we will first deduce Part (B) of Theorem~\ref{Thm:2-plane_sticky} from Part (A), and then prove Part (A) by using Theorem~\ref{Thm:Main_Sticky}.} 

\subsection{Uniformity} 
Before everything, let us talk more about uniformity. In Section~\ref{subsec:uniform}, we have seen that any set of $\delta$-flakes has a uniform refinement by multiscale pigeonholing. If we modify the pigeonholing procedure so that at each step we pigeonhole not only $|\F[F_\rho]|$ but also $\Delta_{max}(\F[F_\rho])$, then we can also make $\Delta_{max}(\F[F_\rho])$ constant up to a factor $\sim 1$. In particular, this implies that $C_F(\F[F_\rho])$ is constant up to a factor $\sim 1$. In fact, we can do something more:
\begin{proposition}\label{Prop:CWA_uniform}
Every set $\F$ of {\color{black} essentially distinct} $\delta$-flakes in $B_1 \subset \R^n$ has a subset $\F'$ of cardinality $|\F'| \gtrapprox |\F|$ so that for any two scales $\delta \leq \tau \leq \theta \leq 1$, we have $|\F'_\tau[F'_\theta]|$, $C_F(\F'_\tau[F'_\theta])$, and $\Delta_{max}(\F'_\tau[F'_\theta])$ all constant up to a factor $\delta^{-C(n)(\log\log(1/\delta))^{-1}} \approx 1$ as $F'_\theta$ varies in $\F'_\theta$.
\end{proposition}
\begin{proof}
    Let $M$ and $\rho_k$ be as in Definition~\ref{Def:uniform_flakes}, then since $(\log(1/\delta))^{CM^2} \lessapprox 1$, by pigeonholing $\F_{\rho_j}[F_{\rho_k}]$ for all integer pairs $\{(k,j): 0 \leq k \leq j \leq M\}$ in lexicographically descending order, we immediately get a $\gtrapprox 1$ refinement $\F'$ that satisfies the conclusions at least for all $\tau=\rho_j$ and $\theta = \rho_k$ with $j \geq k$. Now we check the conclusions for general $\tau$ and $\theta$. For this purpose, let $\tau \in [\rho_{j}, \rho_{j-1}]$ and $\theta \in [\rho_{k}, \rho_{k-1}]$ ($1\leq j,k \leq M$). Fix any $F'_{\rho_{j}} \subset F'_{\tau} \subset F'_{\rho_{j-1}}$ and $F'_{\rho_{k}} \subset F'_{\theta} \subset F'_{\rho_{k-1}}$. Note that
    \begin{align*}  |\F'_{\rho_{j-1}}[F'_{\rho_{k}}]| &\leq |\F'_\tau[F'_\theta]| \leq |\F'_{\rho_{j}}[F'_{\rho_{k-1}}]|
        \lessapprox       |\F'_{\rho_{j-1}}[F'_{\rho_{k}}]|
        \\
        \Delta_{max}(\F'_{\rho_{j-1}}[F'_{\rho_{k}}]) &
        \leq \Delta_{max}(\F'_\tau[F'_\theta]) \leq \Delta_{max}(\F'_{\rho_{j}}[F'_{\rho_{k-1}}]) \lessapprox \Delta_{max}(\F'_{\rho_{j-1}}[F'_{\rho_{k}}]).
    \end{align*}
    Finally, by combining these with Definition~\ref{Def:Frostman}, we immediately get $C_F(\F'_\tau[F'_\theta]) \approx C_F(\F'_{\rho_{j-1}}[F'_{\rho_{k}}])$. And one can easily check that all factors are bounded by $\delta^{-C(n)M^{-1}}$. So the proof is completed.
\end{proof}

Therefore, we can modify Definition~\ref{Def:uniform_flakes} (and thus Definition~\ref{Def:uniform_shading}) in the following sense: Whenever we say $\F$ is \textit{uniform} (at some scales), we additionally assume multiscale uniformity on convex Wolff axiom constants (between these scales). In view of Proposition~\ref{Prop:CWA_uniform}, this convention is entirely harmless in practice. It also allows us to write expressions such as $|\F_\tau[F_\theta]|$, $\Delta_{max}(\F_\tau[F_\theta])$, $C_F(\F_\tau[F_\theta])$ without ambiguity (up to a $\approx 1$ factor), as long as $\F$ is uniform. {\color{black} And in this case, we will often omit words like ``for all $F_\theta \in \F_\theta$'' for simplicity.}

An important fact is that the stickiness property is stable under uniformization:
\begin{lem}\label{Lem:sticky_uniform}
Let $\eta>0$, and let $(\F, Y)$ be a set of {\color{black} essentially distinct} $\delta$-flakes in $B_1 \subset \R^n$, with $\lambda(\F,Y) \geq \delta^\eta$. 
\begin{enumerate}[label={\rm (\Alph*)}]
\item
If $\F$ is $\delta^{-\eta}$ Frostman at every scale, then there exist $\lambda \gtrapprox \delta^\eta$ and a uniform subset $\F' \subset \F$, such that $\lambda(\F', Y) \gtrapprox \lambda(\F, Y)$, $|Y(F)|/|F| \sim \lambda$ for all $F \in \F'$, and \( \F' \) is $\lessapprox \delta^{-2\eta}$ Frostman at every scale.

\item 
If \( \F \) is $\delta^{-\eta}$ Katz--Tao at every scale, then there exist $\lambda \gtrapprox \delta^\eta$ and a uniform subset $\F' \subset \F$, such that $\lambda(\F', Y) \gtrapprox \lambda(\F, Y)$, $|Y(F)|/|F| \sim \lambda$ for all $F \in \F'$, and \( \F' \) is $\delta^{-\eta}$ Katz--Tao at every scale.
\end{enumerate}
\end{lem}
\begin{proof}
The construction of $\F'$ is the same for both parts. 
Let $(\F'',Y)$ be a $\gtrapprox 1$ refinement of $(\F,Y)$ so that there is $\lambda > 0$ such that $|Y(F)|/|F| \sim \lambda$ for all $F \in \F''$. Then $|\F''| \gtrapprox \delta^{\eta} |\F|$ and $\lambda \sim \lambda(\F'',Y) \gtrapprox \lambda(\F,Y) \geq \delta^\eta$. By Proposition~\ref{Prop:CWA_uniform}, there is a subset $\F' \subset \F''$ with $|\F'| \gtrapprox |\F''|$ so that $\F'$ is uniform. 
Note that $(\F',Y)$ is still a $\gtrapprox 1$ refinement of $(\F,Y)$, and $|\F'| \gtrapprox \delta^\eta |\F|$, $\lambda(\F',Y) \sim \lambda(\F'',Y) \gtrapprox \delta^\eta$. It remains to check the required stickiness property for each part.

\begin{enumerate}[label={\rm (\Alph*)}]
\item For each $\delta \leq \rho \leq 1$, by pigeonholing, there is an $F_\rho$ such that $|\F'[F_\rho]| \gtrapprox \delta^\eta |\F[F_\rho]|$, which means $C_F(\F'[F_\rho]) \lessapprox \delta^{-\eta} C_F(\F[F_\rho]) \leq \delta^{-2\eta}$. Then by uniformity of $\F'$, this holds for all $F_\rho \in \F'_\rho$. This proves Part~(A). 

\item For each $\delta \leq \rho \leq 1$, since $\F'_\rho \subset \F_\rho$, we have $\Delta_{max}(\F'_\rho) \leq \Delta_{max}(\F_\rho) \leq \delta^{-\eta}$. This proves Part~(B).
\end{enumerate}
So the proof is completed.
\end{proof}

\subsection{Theorem~\ref{Thm:2-plane_sticky}(A) implies Theorem~\ref{Thm:2-plane_sticky}(B)} 

For this purpose, we need to first establish some facts about multiscale analysis.

\begin{lem}\label{Lem:KT_multi}
Let \( \F \) be a uniform set of \( \delta \)-flakes in \( B_1 \subset \mathbb{R}^n \), and \( \delta = \rho_M \leq \rho_{M-1} \leq \dots \leq \rho_1 \leq \rho_0 = 1 \) is a sequence of scales. Then
\[
\Delta_{max}(\F) \lessapprox_M \prod_{m=1}^M \Delta_{max}(\F_{\rho_m}[F_{\rho_{m-1}}])
\]
and
\[
C_F(\F) \lessapprox_M \prod_{m=1}^M C_F(\F_{\rho_m}[F_{\rho_{m-1}}]).
\]
\end{lem}

\begin{proof}
By Definition~\ref{Def:Frostman} and the telescoping identity 
\begin{align*}
    \Delta(\F, B_1) = \frac{|F||\F|}{|B_1|} \approx_M \prod_{m=1}^M \frac{|F_{\rho_m}||\F_{\rho_m}[F_{\rho_{m-1}}]|}{|F_{\rho_{m-1}}|} = \prod_{m=1}^M \Delta(\F_{\rho_m}, F_{\rho_{m-1}}),
\end{align*}
the second claim follows immediately from the first one. So it suffices to prove the first claim.

Suppose that \( K \subset B_1 \) is a convex set with $|\F[K]| \geq 1$ (thus the shortest axis has length $\geq \delta$). We have to estimate \( \Delta(\F, K) = \frac{|F| |\F[K]|}{|K|} \).

Let \( K_m \) denote the \( \rho_m \)-neighborhood of \( K \). For a flake \( F \in \F \), let \( F_{\rho_m} \) be the corresponding thickened flake in \( \F_{\rho_m} \), and notice that if \( F \subset K \) then \( F_{\rho_m} \subset K_m \). So if \( F \subset K \), then \( F_{\rho_1} \subset K_1 \), and \( F_{\rho_2} \subset F_{\rho_1} \cap K_2 \). In general, we have \( F_{\rho_m} \subset F_{\rho_{m-1}} \cap K_m \). Hence
\begin{align}\label{Eq:count}
|\F[K]| \lessapprox_M \prod_{m=1}^M |\F_{\rho_m}[F_{\rho_{m-1}} \cap K_m]|
\lessapprox_M \prod_{m=1}^M \Delta_{max}(\F_{\rho_m}[F_{\rho_{m-1}}]) \frac{|K_m \cap F_{\rho_{m-1}}|}{|F_{\rho_m}|}.
\end{align}

Now
\[
\frac{|K_m \cap F_{\rho_{m-1}}|}{|F_{\rho_{m}}|}
= \frac{|K_m \cap F_{\rho_{m-1}}|}{|F_{\rho_{m-1}}|} \cdot \frac{|F_{\rho_{m-1}}|}{|F_{\rho_m}|}
\sim \frac{|K_m|}{|K_{m-1}|} \cdot \frac{|F_{\rho_{m-1}}|}{|F_{\rho_m}|}.
\]

Plugging this into \eqref{Eq:count} and simplifying the telescoping product, we obtain
\[
|\F[K]| \lessapprox_M \left[ \prod_{m=1}^M \Delta_{max}(\F_{\rho_m}[F_{\rho_{m-1}}]) \right] \frac{|K|}{|F|}.
\]
This proves the first claim.
\end{proof}

\begin{rem}\label{Rem:KT_multi_relax}
    One can see from the proof of Lemma~\ref{Lem:KT_multi} that the conclusions still hold true if we only assume that $\F$ is uniform at a sequence of scales \( \delta = \rho_M \leq \rho_{M-1} \leq \dots \leq \rho_1 \leq \rho_0 = 1 \). In this case, we may replace $\lessapprox_M$ with $\lesssim_M$. Moreover, for $\F$ not necessarily uniform, we still have
    \[
    \Delta_{max}(\F) \lessapprox_M \prod_{m=1}^M \Big(\max_{F_{\rho_{m-1}} \in \F_{\rho_{m-1}}} \Delta_{max}(\F_{\rho_m}[F_{\rho_{m-1}}])\Big).
    \]
\end{rem}

\begin{lem}\label{Lem:random_(KT)->(F)}
Let \( \eps, \delta > 0 \) and let \( \F \) be a {\color{black} uniform} set of \( \delta \)-flakes in \( B_1 \subset \mathbb{R}^n \). Suppose that there is a sequence of scales \( \delta = \rho_M < \ldots < \rho_0 = 1 \) so that for each \( m \) we have
\begin{itemize}
    \item \( \frac{\rho_m}{\rho_{m-1}} \gtrapprox \delta^\eps \).
    \item \( \Delta_{max}(\F_{\rho_m}[F_{\rho_{m-1}}]) \lessapprox \left( \frac{\rho_{m-1}}{\rho_m} \right)^\eps \).
\end{itemize}
Then provided that $\delta > 0$ is sufficiently small depending on $n$ and $M$, there is a set of translations \( \{R_j\}_{j=1}^J \) of cardinality \( J \lessapprox_M \delta^{-\eps} (|F||\F|)^{-1} \) so that \( \bigsqcup_j R_j(\F) \) is {\color{black} $\lessapprox 1$ essentially distinct} and \( \lessapprox_M \delta^{-(2n-1)\eps} \) Frostman at every scale.
\end{lem}

The proof of Lemma~\ref{Lem:random_(KT)->(F)} is similar to that of \cite[Lemma~7.5]{GWZ26streamlined}, but we still present the details for completeness. We will use random translations \( R_j \). If \( v \in \mathbb{R}^n \), we let \( R_v(x) = x + v \). A random translation of \textit{size} \( \leq \rho \) is a translation \( R_v \) by a vector \( v \in B(0,\rho) \) selected randomly from the uniform measure on \( B(0,\rho) \). We will use the following probability lemma, which is a small variation on Lemma~\ref{Lem:random_rigid_F}.
\begin{lem}\label{Lem:random_translate}
Let \( \delta \leq \rho \leq 1 \). Suppose that \( \F \) is a set of {\color{black} essentially distinct} \( \delta \)-flakes in \( F_\rho \subset \mathbb{R}^n \) with \( \frac{|F_\rho|}{|\F||F_\delta|} \sim J \in \N^+ \). Let \( R_1, \ldots, R_J \) be a set of \( J \) {\color{black}independent} random translations of size \( \leq \rho \). Let \( \F' = \bigsqcup_{j=1}^J R_j(\F) \). Then with high probability {\rm(}at least $1-\delta^{100n^2}$ when $\delta$ is sufficiently small{\rm)} the flakes of \( \F' \) are essentially distinct up to multiplicity \( \lessapprox 1 \), and we have \( \Delta_{max}(\F') \lessapprox \Delta_{max}(\F) \). After refining \( \F' \) to make the flakes essentially distinct, we have with high probability \( |\F'| \approx J|\F| \) and \( \Delta_{max}(\F') \lessapprox \Delta_{max}(\F) \).
\end{lem}

We will prove Lemma~\ref{Lem:random_translate} in Appendix~\ref{Appendix:Probability}, where we collect a couple of related probability arguments. Now we are ready to prove Lemma~\ref{Lem:random_(KT)->(F)}. The proof consists of using Lemma~\ref{Lem:random_translate} at many scales.

\begin{proof}[Proof of Lemma~\ref{Lem:random_(KT)->(F)}]
For each \( m = 1, \dots, M \), let \( \mathcal{R}_m \) be a set of {\color{black}independent} random translations of size \(\leq \rho_{m-1}\), with cardinality
\begin{align}\label{Eq:R_k_number}
\left| \mathcal{R}_m \right| \sim \max \Big\{ 1, \frac{(\rho_m / \rho_{m-1})^{-2}}{| \F_{\rho_m} [F_{\rho_{m-1}}] |} \Big\}.
\end{align}

Define the set of translations \( \{ R_j \}_{j=1}^J \) as \( R_1' \circ R_2' \circ \cdots \circ R_M' \) for \( R_m' \in \mathcal{R}_m, \, m = 1, \dots, M \). Recall that \( \F' = \bigsqcup_j R_j(\F) \). Recall that $\Delta_{max}(\F_{\rho_m}[F_{\rho_{m-1}}]) \lessapprox \left( \frac{\rho_{m-1}}{\rho_m} \right)^\eps$, so we have $|\F_{\rho_m}[F_{\rho_{m-1}}]| \lessapprox \left( \frac{\rho_{m}}{\rho_{m-1}} \right)^{-2} \left( \frac{\rho_{m-1}}{\rho_m} \right)^\eps$, and 
\[
|\mathcal{R}_m| \lessapprox \Big( \frac{\rho_{m-1}}{\rho_m} \Big)^\eps \frac{(\rho_m / \rho_{m-1})^{-2}}{| \F_{\rho_m} [F_{\rho_{m-1}}] |}.
\]
This implies
\[
J \leq 
\prod_{m=1}^M |\mathcal{R}_m| \lessapprox_M \delta^{-\eps} \frac{1}{|F||\F|}.
\]

We will bound \( \Delta_{max}(\F') \) using Lemma~\ref{Lem:KT_multi}. To do that, we need to estimate \( \Delta_{max}(\F'_{\rho_m}[F'_{\rho_{m-1}}]) \). Suppose \( F'_{\rho_{m-1}} \in \F'_{\rho_{m-1}} \). Then we must have \( F'_{\rho_{m-1}} \sim R_1' \circ \cdots \circ R_{m-1}'(F_{\rho_{m-1}}) \) for some \( F_{\rho_{m-1}} \in \F_{\rho_{m-1}} \). By Lemma~\ref{Lem:random_translate}, with high probability, the choice of \( R_1', \dots, R_{m-1}' \) is \(\lessapprox 1\). Therefore with probability at least $3/4$ we have, for each \( m = 1, \dots, M \),
\begin{align}\label{Eq:to_R_m_trans}
\Delta_{max}(\F'_{\rho_m}[F'_{\rho_{m-1}}]) \lessapprox
\Delta_{max} \Big( \bigsqcup_{R_m' \in \mathcal{R}_m} R_m'(\F_{\rho_m}[F_{\rho_{m-1}}]) \Big).
\end{align}

Fix a flake \( F_{\rho_{m-1}} \). Applying Lemma~\ref{Lem:random_translate}, we have that with probability at least $1-\delta^{100n^2}$, 
\begin{align}\label{Eq:R_m_trans}
\Delta_{max} \Big( \bigsqcup_{R_m' \in \mathcal{R}_m} R_m'(\F_{\rho_m}[F_{\rho_{m-1}}]) \Big) \lessapprox \Delta_{max}(\F_{\rho_m}[F_{\rho_{m-1}}]) \lessapprox 
\Big( \frac{\rho_{m-1}}{\rho_m} \Big)^{\eps}.
\end{align}

There are \(\lesssim M \delta^{-50n}\) many flakes in total. Thus if \(\delta > 0\) is sufficiently small depending on $n$ and $M$, we conclude that with probability at least $3/4$, \eqref{Eq:R_m_trans} holds for every \( F_{\rho_{m-1}} \in \F_{\rho_{m-1}} \). Combining \eqref{Eq:to_R_m_trans} and \eqref{Eq:R_m_trans}, we see that for each \( m \),
\[
\Delta_{max}(\F_{\rho_m}'[F_{\rho_{m-1}}']) \lessapprox \Big( \frac{\rho_{m-1}}{\rho_m} \Big)^{\eps}.
\]

Now we can apply Lemma~\ref{Lem:KT_multi} to show that for each \( m \),
\begin{align}\label{Eq:apply_KT_multi}
\Delta_{max}(\F'[F'_{\rho_m}]) \lessapprox_M \prod_{m'=m+1}^M \Delta_{max}(\F'_{\rho_{m'}}[F_{\rho_{m'-1}}]) \lessapprox_M \delta^{-\eps}.
\end{align}

On the other hand, \eqref{Eq:R_k_number} shows that
\[
|\F'_{\rho_m}[F'_{\rho_{m-1}}]| \approx 
\Big| 
\bigsqcup_{R_m' \in \mathcal{R}_m} R_m'(\F_{\rho_m}[F_{\rho_{m-1}}]) \Big| 
= \sum_{R_m'\in\mathcal{R}_m} |
R_m'(\F_{\rho_m}[F_{\rho_{m-1}}])|
\gtrsim 
\Big( \frac{\rho_m}{\rho_{m-1}} \Big)^{-2},
\]
which implies that
\begin{align}\label{Eq:rho_m_number}
|\F'[F'_{\rho_m}]| \gtrapprox_M \Big( \frac{\delta}{\rho_m} \Big)^{-2} \quad \text{ and so } \quad \Delta_{max}(\F'[F'_{\rho_m}]) \gtrapprox_M 1.
\end{align}

Combining \eqref{Eq:apply_KT_multi} and \eqref{Eq:rho_m_number}, we see that for each \( m \),
\[
C_F(\F'[F'_{\rho_m}]) \lessapprox_M \delta^{-\eps}.
\]

Finally consider an arbitrary \( \rho \in [\delta, 1] \). We can choose \( m \) so that \( \rho_m \leq \rho \leq \rho_{m-1} \). Since \( \rho_m \gtrapprox \delta^{\eps}\rho_{m-1} \), using the last equation for \( m \) we can conclude that
\[
C_F(\F'[F'_{\rho}]) \lessapprox_M \delta^{-(2n-1)\eps},
\]
as desired.
\end{proof}

\begin{proof}[Proof of Part {\rm(B)} of Theorem~\ref{Thm:2-plane_sticky} assuming Part {\rm(A)}]

Fix \(\eps > 0\) and let \(\eta_1 > 0\) be the $\eta$ obtained from Theorem~\ref{Thm:2-plane_sticky}(A). Decreasing \(\eta_1\) if necessary, we can assume that $\eta_1 < \eps/3$, and \(2(2n-1)/\eta_1 = M\) is an integer. Let \(\eta = 1/M^2 \). Let \((\F, Y)\) be a set of essentially distinct \(\delta\)-flakes that is \(\delta^{-\eta}\) Katz--Tao at every scale and has \(\lambda(\F, Y) \geq \delta^\eta\). It suffices to prove that \(|U(\F, Y)| \geq \delta^{2\eps}|F| |\F|\). 

By Lemma~\ref{Lem:sticky_uniform}(B), there exist $\lambda \gtrapprox \delta^\eta$ and a uniform subset $\F' \subset \F$, such that $\lambda(\F', Y) \gtrapprox \lambda(\F, Y)$, $|Y(F)|/|F| \sim \lambda$ for all $F \in \F'$, and \( \F' \) is $\delta^{-\eta}$ Katz--Tao at every scale. Note that $|\F'| \gtrapprox \delta^\eta|\F|$.

Consider the sequence of scales \(\rho_m = \delta^{m/M}\), \(m = M, \ldots, 0\). We have \(\rho_m/\rho_{m-1} = \delta^{1/M} \geq \delta^{\eta_1}\), and 
\[
\Delta_{max}(\F'_{\rho_m}[{F_{\rho_{m-1}}}]) \leq \Delta_{max}(\F'_{\rho_m}) \leq \delta^{-\eta} = (\rho_{m-1}/\rho_m)^{\eta M} \leq (\rho_{m-1}/\rho_m)^{\eta_1}.
\] 
Now we can apply Lemma~\ref{Lem:random_(KT)->(F)} to $\F'$. We obtain a set of \(J\) translations \(\{R_j\}_{j=1}^J\), with \(J \lessapprox_M \delta^{-1/M} (|F||\F'|)^{-1}\) (thus $J \leq \delta^{-\eta_1} (|F||\F'|)^{-1}$), so that the set \(\F'' = \bigsqcup_{j=1}^J R_j(\F')\) is {\color{black} $\lessapprox 1$ essentially distinct} and \(\lessapprox_M \delta^{-(2n-1)/M}\) (thus $\delta^{-2\eta_1/3}$) Frostman at every scale. Define the shading \(Y'\) on \(\F''\) in the obvious way: For any \(F \in \F'\) and \(1 \leq j \leq J\), let \(Y'(R_j(F)) \coloneqq R_j(Y(F))\). Then \( \lambda(\F'', Y') = \lambda(\F', Y) \gtrapprox \delta^\eta \). {\color{black} Intuitively we want to apply Theorem~\ref{Thm:2-plane_sticky}(A) to $(\F'',Y')$, but since the flakes may not be genuinely essentially distinct, we need some technical steps.}

Since $\F''$ is $\lessapprox 1$ essentially distinct, there is $\F''' \subset \F''$ with $|\F'''| \gtrapprox |\F''|$ so that $\F'''$ is essentially distinct. By our construction of $\F''$, we have $|Y'(F)|/|F| \sim \lambda$ for all $F \in \F''$, so $(\F''',Y')$ is a $\gtrapprox 1$ refinement of $(\F'',Y')$, and $\lambda(\F''',Y') \sim \lambda(\F'',Y') \gtrapprox \delta^\eta$. For each $\delta \leq \rho \leq 1$, by pigeonholing, there is an $F_\rho$ such that $|\F'''[F_\rho]| \gtrapprox |\F''[F_\rho]|$, which means $C_F(\F'''[F_\rho]) \lessapprox C_F(\F''[F_\rho]) \leq \delta^{- 2\eta_1/3}$. Then by uniformity of $\F'''$, this holds for all $F_\rho \in \F'''_\rho$. We can conclude that $\F'''$ is $\delta^{-\eta_1}$ Frostman at every scale. 

Now \((\F''', Y')\) satisfies the hypotheses of Theorem~\ref{Thm:2-plane_sticky}(A) with $\eta_1$ in place of $\eta$. Thus we get 
\begin{align*}
|U(\F, Y)| &\geq |U(\F', Y)| \\
&= \frac{1}{J} \sum_{j=1}^J |R_j(U(\F', Y))| 
= \frac{1}{J} \sum_{j=1}^J | U(R_j(\F'), R_j(Y) )|\\
&\geq \frac{1}{J} \Big| \bigcup_{j=1}^J U(R_j(\F'), R_j(Y)) \Big| 
= \frac{1}{J} |U(\F'', Y')|\\
&\geq \frac{1}{J} |U(\F''', Y')|
\geq \frac{1}{J} \delta^{\eps} 
\geq \delta^{\eps+\eta_1}|F| |\F'| \gtrapprox \delta^{\eps+\eta_1+\eta}|F| |\F|.
\end{align*}
Since $\eta \leq \eta_1 \leq \eps/3$, we get $|U(\F, Y)| \geq \delta^{2\eps}|F| |\F|$. Finally, as \(\eps > 0\) was arbitrary, this proves Theorem~\ref{Thm:2-plane_sticky}(B).
\end{proof}

\subsection{Theorem~\ref{Thm:Main_Sticky} implies Theorem~\ref{Thm:2-plane_sticky}(A)} 
For this purpose, we need the following analogue of \cite[Lemma~7.2]{guth2025outline}:
\begin{lem}\label{Lem:reduce_delta_-2}
Let \( \eps, \delta > 0 \) and let \( \F \) be a {\color{black} uniform} set of \( \delta \)-flakes in \( B_1 \subset \mathbb{R}^n \). Suppose that there is a sequence of scales \( \delta = \rho_M < \ldots < \rho_0 = 1 \) so that for each \( m \) we have
\begin{itemize}
    \item \( \frac{\rho_m}{\rho_{m-1}} \gtrapprox \delta^\eps \).
    \item \( C_F(\F_{\rho_m}[F_{\rho_{m-1}}]) \lessapprox \left( \frac{\rho_{m-1}}{\rho_m} \right)^\eps \).
\end{itemize}
Then $\F$ contains a subset $\F'$ which is $\delta^{-3\eps}$-sticky.
\end{lem}

The proof of Lemma~\ref{Lem:reduce_delta_-2} is similar to that of \cite[Lemma~7.2]{guth2025outline}, but we still present
the details for completeness. The key ingredient is the following probability lemma, which should be compared with Lemma~\ref{Lem:random_rigid_F}.
\begin{lem}\label{Lem:random_sub_(F)}
Let $\F$ be a set of $\delta$-flakes in $B_1 \subset \R^n$. Then it contains a subset $\F'$ with $C_F(\F') \lessapprox C_F(\F)$ and $|\F'| \approx \delta^{-2} C_F(\F)^{-1}$.
\end{lem}

We will prove Lemma~\ref{Lem:random_sub_(F)} in Appendix~\ref{Appendix:Probability}, where we collect a couple of related
probabilistic arguments. Now we are ready to prove Lemma~\ref{Lem:reduce_delta_-2}. The proof consists of using Lemma~\ref{Lem:random_sub_(F)} at many scales.

\begin{proof}[Proof of Lemma~\ref{Lem:reduce_delta_-2}]
    We select $\F'$ in a random way one level of the tree at a time. We first select $\F_{\rho_1}' \subset \F_{\rho_1}$. Since $C_F(\F_{\rho_1}) \lessapprox \rho_1^{-\eps}$, applying Lemma~\ref{Lem:random_sub_(F)} to  $\F_{\rho_1}$ yields $\F_{\rho_1}' \subset \F_{\rho_1}$ such that $C_F(\F_{\rho_1}') \lessapprox C_F(\F_{\rho_1})\lessapprox \rho_1^{-\eps}$ and $\rho_1^{-2+\eps} \lessapprox |\F_{\rho_1}'| \lessapprox \rho_1^{-2}$. 
    Then for each $F_{\rho_1} \in \F_{\rho_1}'$, we select $\F_{\rho_2}'[F_{\rho_1}]$ by applying Lemma~\ref{Lem:random_sub_(F)} to (a rescaled version of) $\F_{\rho_2}[F_{\rho_1}]$. This yields $C_F(\F_{\rho_2}'[F_{\rho_1}]) \lessapprox C_F(\F_{\rho_2}[F_{\rho_1}]) \lessapprox (\rho_2/\rho_1)^{-\eps}$ and $(\rho_2/\rho_1)^{-2+\eps} \lessapprox |\F_{\rho_2}'[F_{\rho_1}]| \lessapprox (\rho_2/\rho_1)^{-2}$. Proceeding this way defines $\F'$. We have to check that $\F'$ has the desired properties.

    For each $m = 1, \dots, M$, we have
    \begin{align*}
        (\rho_m/\rho_{m-1})^{-2+\eps} \lessapprox |\F_{\rho_m}'[F_{\rho_{m-1}}]| \lessapprox (\rho_m/\rho_{m-1})^{-2}.
    \end{align*}
    This defines $\F'$, and implies that for every $\rho_m$,
    \begin{align*}
        (\delta/\rho_{m})^{-2+\eps} \lessapprox_M |\F'[F_{\rho_{m}}]| \lessapprox_M (\delta/\rho_{m})^{-2}.
    \end{align*}
    Then for any $\delta \leq \rho \leq 1$, there exists $\rho_m \leq \rho \leq \rho_{m-1}$ (recall $\rho_m \gtrapprox \delta^\eps \rho_{m-1}$), and so
    \begin{align*}
        \delta^{3\eps-\eps^2} (\delta/&\rho)^{-2} \leq \delta^{2\eps-\eps^2}(\delta/\rho)^{-2+\eps} \lessapprox (\delta/\rho_{m})^{-2+\eps} \lessapprox_M |\F'[F_{\rho_{m}}]| \leq |\F'[F_{\rho}]|; \\
        &|\F'[F_{\rho}]| \leq |\F'[F_{\rho_{m-1}}]| \lessapprox_M (\delta/\rho_{m-1})^{-2} \lessapprox \delta^{-2\eps}(\delta/\rho)^{-2}.
    \end{align*}
    On the other hand, for each $m = 1, \dots, M$, we have
    \begin{align*}
        C_F(\F_{\rho_m}'[F_{\rho_{m-1}}]) \lessapprox C_F(\F_{\rho_m}[F_{\rho_{m-1}}]) \lessapprox 
        (\rho_m/\rho_{m-1})^{-\eps}.
    \end{align*}
    Given our bound for $\F_{\rho_m}'[F_{\rho_{m-1}}]$, this implies
    \begin{align*}
        \Delta_{max}(\F_{\rho_m}'[F_{\rho_{m-1}}]) \lessapprox (\rho_m/\rho_{m-1})^{-\eps}.
    \end{align*}
    By the submultiplicative property of $\Delta_{max}$ in Lemma~\ref{Lem:KT_multi} (and also Remark~\ref{Rem:KT_multi_relax}), we get
    \begin{align*}
        \Delta_{max}(\F') \lessapprox_M \delta^{-\eps}.
    \end{align*}
    Therefore, we can conclude that $\F'$ which is $\delta^{-3\eps}$-sticky.
\end{proof}

\begin{proof}[Proof of Theorem~\ref{Thm:2-plane_sticky}{\rm(A)} assuming Theorem~\ref{Thm:Main_Sticky}]
Fix \(\eps > 0\) and let \(\eta_1 > 0\) be the $\eta$ obtained from Theorem~\ref{Thm:Main_Sticky}. Decreasing \(\eta_1\) if necessary, we can suppose that \(6/\eta_1 = M\) is an integer. Let \(\eta = 1/M^2\). Let \((\F, Y)\) be a set of essentially distinct \(\delta\)-flakes that is \(\delta^{-\eta}\) Frostman at every scale and has \(\lambda(\F, Y) \geq \delta^\eta\). It suffices to prove that \(|U(\F, Y)| \geq \delta^{\eps}\).

By Lemma~\ref{Lem:sticky_uniform}(A), there exist $\lambda \gtrapprox \delta^\eta$ and a uniform subset $\F' \subset \F$, such that $\lambda(\F', Y) \gtrapprox \lambda(\F, Y)$, $|Y(F)|/|F| \sim \lambda$ for all $F \in \F'$, and \( \F' \) is $\lessapprox\delta^{-2\eta}$ Frostman at every scale.

Consider the sequence of scales \(\rho_m = \delta^{m/M}\), \(m = M, \ldots, 0\). We have \(\rho_m/\rho_{m-1} = \delta^{1/M} \geq \delta^{\eta_1/3}\), and uniformity of $\F'$ implies 
\[
C_F(\F'_{\rho_m}[F_{\rho_{m-1}}]) \lessapprox C_F(\F'[F_{\rho_{m-1}}]) \lessapprox \delta^{-2\eta} = (\rho_{m-1}/\rho_m)^{2\eta M} = (\rho_{m-1}/\rho_m)^{\eta_1/3}.
\]
Now we can apply Lemma~\ref{Lem:reduce_delta_-2} to $\F'$, and obtain a subset $\F''$ which is $\delta^{-\eta_1}$-sticky. Since $|Y(F)|/|F| \sim \lambda$ for all $F \in \F'' \subset \F'$, we know that $\lambda(\F'',Y) \sim \lambda \gtrapprox \delta^\eta$ (thus $\lambda(\F'',Y) \geq \delta^{\eta_1}$). Now by applying Theorem~\ref{Thm:Main_Sticky} to $(\F'',Y)$, we get
\begin{align*}
    |U(\F,Y)| \geq |U(\F'',Y)| \geq \delta^{\eps},
\end{align*}
as desired. This proves Theorem~\ref{Thm:2-plane_sticky}(A).
\end{proof}

Therefore, Theorem~\ref{Thm:2-plane_sticky} holds true.

\section{General cases and two main lemmas}\label{Sec:partial_estimates}
Having solved the sticky case (Theorem~\ref{Thm:2-plane_sticky}), we now turn to the reduction of the general case. The idea follows tightly with the paper by Guth, Wang and Zahl \cite{GWZ26streamlined}. We will establish the statements $K_{KT}(\beta)$ and $K_F(\beta)$ and prove two self-improving lemmas that will establish Theorem~\ref{Thm:Main}.

In this section, we first introduce and analyze Kakeya-type estimates for $\delta$-flakes satisfying Katz--Tao and Frostman convex Wolff axioms. Then we will set up the general induction framework of proving Theorem~\ref{Thm:Main}.

\begin{definition}\label{Def:Katz_Tao_beta}
    We write $K_{KT}(\beta)$ for the following statement. For any $\eps>0$, there exists $\eta = \eta(\eps,\beta)$ so that the following holds for any $0< \delta \leq \delta(\eps,\beta)$. Let $\F$ be a set of $\delta$-flakes in $B_1\subset \R^n$ with the following properties.
    \begin{itemize}
        \item $\Delta_{max}(\F) \leq \delta^{-\eta}$.
        \item $Y$ is a shading on $\F$ with $\lambda(\F, Y) \geq \delta^{\eta}$.
    \end{itemize}
    Then
    \begin{align*}
        \mu(\F,Y) \leq \delta^{-\eps} |\F|^\beta.
    \end{align*}
\end{definition}

Note that $K_{KT}(1)$ is trivial. Our goal is to prove $K_{KT}(0)$. When $|\F|\sim \delta^{-2}$, $K_{KT}(\beta)$ gives $\mu(\F) \lesssim \delta^{-2\beta-\eps}$, which is equivalent to $|U(\F)| \gtrsim \delta^{2\beta+\eps}$.

\begin{definition}\label{Def:Frostman_beta}
    We write $K_{F}(\beta)$ for the following statement. For any $\eps>0$, there exists $\eta = \eta(\eps,\beta)$ so that the following holds for any $0< \delta \leq \delta(\eps,\beta)$. Let $\F$ be a set of {\color{black} essentially distinct} $\delta$-flakes in $B_1\subset \R^n$ with the following properties.
    \begin{itemize}
        \item $C_F(\F,B_1) \leq \delta^{-\eta}$.
        \item $Y$ is a shading on $\F$ with $\lambda(\F, Y) \geq \delta^{\eta}$.
    \end{itemize}
    Then
    \begin{align*}
        \mu(\F, Y) \leq \delta^{-\eps} \delta^{-2\beta} (|\F||F|)^{1-\beta/{\color{black} n}}.
    \end{align*}
    Or equivalently,
    \begin{align*}
        |U(\F, Y)| \geq \delta^{\eps} \delta^{2\beta} (|\F||F|)^{\beta/{\color{black} n}}.
    \end{align*}
\end{definition}
The above two definitions depend on the dimension $n$, which is nevertheless often clear from the context. So for clarity, we will not specify the dependence on $n$ in $K_{KT}(\beta)$ and $K_{F}(\beta)$.

Note that $C_F(\F,B_1) \leq \delta^{-\eta}$ implies $|\F| \geq \delta^{-2+\eta}$, and is in fact equivalent to $\Delta_{max}(\F) \lesssim \delta^{-\eta}$ when $|\F|\sim \delta^{-2}$. In this case, both $K_{KT}(\beta)$ and $K_F(\beta)$ imply $|U(\F)| \gtrsim \delta^{2\beta+\eps}$. However, when $|\F| \gg \delta^{-2}$ (and $\beta>0$), we have $(|\F||F|)^{\beta/n} \gg 1$, so $K_F(\beta)$ gives a stronger bound. The exact exponent $\beta/n$ in
this final factor is not crucial in our analysis, and the argument would work equally well with any positive exponent $c\beta$, $0<c<1/(n-2)$. There is one step (see {\color{black} Remark~\ref{Rem:numerology}}) where it is important
that $c<1/(n-2)$, which differs from the range in \cite{GWZ26streamlined}, and we will explain it at an appropriate time.

\begin{proposition}\label{Prop:tech_tau}
    If $K_{KT}(\beta)$ holds, then we get for free the following slightly more general version of $K_{KT}(\beta)$: 
    For any $\eps>0$, there exists $\eta = \eta(\eps,\beta)$ so that the following holds for any $0< \delta \leq \delta(\eps,\beta)$. Suppose $\F$ is a set of $\delta$-flakes in $B_1\subset \R^n$, $\tau\leq \delta$, and
    \begin{itemize}
        \item $\Delta_{max}(\F) \leq \tau^{-\eta}$.
        \item $Y$ is a shading on $\F$ with $\lambda(\F, Y) \geq \tau^{\eta}$.
    \end{itemize}
    Then
    \begin{align*}
        \mu(\F,Y) \leq \tau^{-\eps} |\F|^\beta.
    \end{align*}
    Similarly, $K_F(\beta)$ implies a slightly more general version of $K_F(\beta)$.
\end{proposition}
\begin{proof}
    Exactly the same as that of \cite[Remark~3.6]{GWZ26streamlined}, and so omitted.
\end{proof}
These slightly more general versions will be useful when we replace a set of $\delta$-flakes $(\F,Y)$ with
a thickening $(\F_\rho,Y_\rho)$, for which we only know that $\lambda(\F_\rho,Y_\rho) \geq \delta^\eta$. 

For $\beta \in [0,1]$, the two main lemmas in the whole proof are
\begin{Mlem}\label{Main_lem_1}
    $K_{KT}(\beta)$ implies $K_F(\beta)$.
\end{Mlem}

\begin{Mlem}\label{Main_lem_2}
    If $\beta>0$ and $K_{KT}(\beta)$ and $K_F(\beta)$ hold, then $K_{KT}(\beta-\nu)$ holds for $\nu = \nu(\beta) > 0$. Moreover $\nu(\beta)$ is monotone in $\beta$.
\end{Mlem}
Since the bound $K_{KT}(1)$ is trivial, the two main lemmas imply that $K_{KT}(\beta)$ and $K_{F}(\beta)$ hold
for all $\beta>0$, which implies Theorem~\ref{Thm:Main}.

Before we proceed with the proofs of these two main lemmas, let us first give a brief roadmap:
\begin{itemize}
    \item Firstly, the proofs require some flexibility in applications of the two statements $K_{KT}(\beta)$ and $K_F(\beta)$: In Section~\ref{subsec:large_CWA}, we will extend their usage to sets of flakes with larger convex Wolff constants; while in Section~\ref{Sec:slab_plank}, we will further extend their usage to planks instead of just flakes, as well as to flakes that are factored by planks.

    \item Next, in Section~\ref{Sec:reduce_to_sticky}, we will establish a multiscale decomposition (Lemma~\ref{Lem:multi-scale_decomposition}) that roughly provides a useful dichotomy between sticky and non-sticky cases. 
    
    \item Then, in Section~\ref{Sec:Main_Lemma_1}, we prove Main Lemma~\ref{Main_lem_1} by using the ``Frostman'' part of sticky Kakeya (Theorem~\ref{Thm:2-plane_sticky}(A)) and the multiscale decomposition (Lemma~\ref{Lem:multi-scale_decomposition}(A)).

    \item Finally, in Section~\ref{Sec:Main_Lemma_2}, we prove Main Lemma~\ref{Main_lem_2} by first reducing it to a ``very not sticky'' case (Lemma~\ref{Lem:NSAAS}) using the ``Katz--Tao'' part of sticky Kakeya (Theorem~\ref{Thm:2-plane_sticky}(B)) and the multiscale decomposition (Lemma~\ref{Lem:multi-scale_decomposition}(B)), and then prove Lemma~\ref{Lem:NSAAS} by using a carefully designed two-scale induction scheme.
\end{itemize}

\subsection{Multiplicity bounds for larger values of \texorpdfstring{$\Delta_{max}(\F)$}{∆max(P)} or \texorpdfstring{$C_F(\F)$}{CF(P)}}\label{subsec:large_CWA} 
The next two results
say that $K_{KT}(\beta)$ and $K_{F}(\beta)$ imply multiplicity bounds for families of flakes with any value of
$\Delta_{max}(\F)$ or $C_F(\F)$.

\begin{lem}\label{Lem:general_KT}
    Suppose that $K_{KT}(\beta)$ holds. Then for any $\eps > 0$, there exists $\eta = \eta(\eps, \beta) > 0$ such that the following holds for any $0 < \delta < \delta(\eps, \beta)$. Suppose that $(\F, Y)$ is a set of $\delta$-flakes in $B_1\subset \R^n$ with $\lambda(\F,Y) \geq \delta^\eta$. Then 
    \begin{align*}
        \mu(\F,Y) \leq \delta^{-\eps} \Delta_{max}(\F)^{1-\beta}|\F|^\beta.
    \end{align*}
\end{lem}
\begin{proof}
    Exactly the same as that for \cite[Lemma~3.7]{GWZ26streamlined}, and so omitted. 
\end{proof}

\begin{lem}\label{Lem:random_rigid_F}
    Suppose that $\F$ is a set of {\rm(}essentially distinct{\rm)} $\delta$-flakes in $B_1 \subset \R^n$. Let $R_1, \dots , R_J$ be {\color{black}independent} random rigid motions of size $\leq 1$, with $J = \lceil C_F(\F)\rceil$, and let $\F' = \bigsqcup_{j=1}^J R_j(\F)$. Then with high probability the flakes of $\F'$ are essentially distinct up to multiplicity $\lessapprox 1$ and also $C_F(\F') \lessapprox 1$. After refining $\F'$ to make the flakes essentially distinct, we have $|\F'| \approx C_F(\F)|\F|$ and $C_F(\F') \lessapprox 1$.
\end{lem} 
Here a random rigid motion \textit{of size $\leq 1$} is a random rotation composed with a
random translation of size $\leq 1$. We will prove Lemma~\ref{Lem:random_rigid_F} in Appendix~\ref{Appendix:Probability}, where we collect a couple of related probability
lemmas.

\begin{lem}\label{Lem:general_F}
    Suppose that $K_{F}(\beta)$ holds. Then for any $\eps > 0$, there exists $\eta = \eta(\eps, \beta) > 0$ such that the following holds for any $0 < \delta < \delta(\eps, \beta)$. Suppose that $(\F, Y)$ is a set of {\color{black} essentially distinct} $\delta$-flakes in $B_1\subset \R^n$ with $\lambda(\F,Y) \geq \delta^\eta$. Then 
    \begin{align*}
        \mu(\F,Y) \leq \delta^{-\eps} C_F(\F)^{1-\beta/n} \delta^{-2\beta}(\delta^2|\F|)^{1-\beta/n}.
    \end{align*}
\end{lem}
\begin{proof}
    Let $\F'$ be the union of $\lceil C_F(\F)\rceil$ {\color{black}independent} random rigid motions of $\F$, and $Y'$ be the shading on $\F'$ defined by the images of $Y$ under these motions. By Lemma~\ref{Lem:random_rigid_F}, with high probability $\F'$ is $\lessapprox 1$ essentially distinct, $C_F(\F') \lessapprox 1$, and $|\F'| \approx C_F(\F)|\F|$. Let $\F''$ be an essentially distinct subset of $\F'$ such that $(\F'', Y')$ is a $\gtrapprox 1$ refinement of $(\F', Y')$. Then $C_F(\F'') \lessapprox C_F(\F') \lessapprox 1$ and $|\F''| \lessapprox C_F(\F)|\F|$. Applying $K_F(\beta)$ to $(\F'', Y')$ with $\eps/2$ in place of $\eps$, we have 
    \begin{align*}
        \mu(\F, Y) \leq \mu(\F', Y') &\lessapprox 
        \mu(\F'', Y')
        \\ 
        &\leq \delta^{-\eps/2}  \delta^{-2\beta}(\delta^2|\F''|)^{1-\beta/n} \\
        &\lessapprox \delta^{-\eps/2} C_F(\F)^{1-\beta/n} \delta^{-2\beta}(\delta^2|\F|)^{1-\beta/n},
    \end{align*}
    which immediately yields the desired bound.
\end{proof}

We point out that the analogue of Proposition~\ref{Prop:tech_tau} holds for Lemma~\ref{Lem:general_KT} and Lemma~\ref{Lem:general_F}. For simplicity, we will often omit citing technical justifications like Proposition~\ref{Prop:tech_tau} when applying the lemmas.

\section{Multiplicity estimates for slabs and planks}\label{Sec:slab_plank}
In Section~\ref{Sec:organize} we reduced the analysis of a general collection of convex sets in $\R^n$ to that of two new families of convex sets: One is Katz--Tao, and the other is Frostman. The estimates $K_{KT}(\beta)$ and $K_F(\beta)$
control the intersection multiplicity of sets of flakes in $\R^n$ 
that are Katz--Tao and Frostman,
respectively. But since the maximal density factoring of flakes can be general planks, we also need to control the intersection multiplicity of sets of planks that are Frostman or Katz--Tao. In this section, we will show that incidence bounds for planks in $\R^n$ can be systematically reduced to incidence bounds for flakes.

{\bf Katz--Tao multiplicity estimates.} Lemma~\ref{Lem:plank_Katz_Tao}  below says that the Kakeya estimate $K_{KT}(\beta)$
implies an analogous estimate for planks of dimensions $a \times b \times \1$. Furthermore, a stronger estimate
holds when $a \ll b$, provided that the planks do not concentrate into slabs (i.e., sets of dimensions $\theta \times 1 \times \1$ for $\theta \ll 1$). To make this latter statement precise, we must define what it means for
planks to concentrate into slabs. 

If $\LL$ is a set of planks in $\R^n$ of dimensions $a \times b \times \1$, and $S$ is a $\theta \times 1 \times \1$ slab with $\theta \in [a/b,1]$, then we define
\begin{align}\label{Eq:L_S}
    \LL_S \coloneqq \{ L \in \LL: L\subset S \text{ and } \angle(TL, TS) \leq \theta\}.
\end{align}
In the above, $\angle(TL, TS)$ refers to the angle between the hyperplanes $TL, TS$, where $TL$ is the hyperplane
spanned by the $(n-1)$ longest axes of $L$, and similarly for $S$ (this is defined up to accuracy $\theta$).

\begin{lem}\label{Lem:plank_Katz_Tao}
    Suppose that \(K_{KT}(\beta)\) holds. Then for every \(\eps > 0\), there exists \(\eta, b_0 > 0\) such that the following holds for any \(0 < a \leq b \leq b_0\). Let \(\LL\) be a set of \(a \times b \times \1\) planks in \(B_1\), and let \(Y \) be a shading on \(\LL\) with density \(\lambda(\LL, Y) \geq a^{\eta}\). Suppose that the planks do not concentrate into slabs, in the following sense: There exists \(0 \leq \gamma \leq 1\) so that for each \(\theta \in [a/b, 1]\) and every \(\theta \times 1 \times \1\) slab \(S\), \(|\LL_S| \leq a^{-\eta} \theta^{\gamma} |\LL|\).

If \(\LL\) is Katz--Tao, then \(\mu(\LL) \leq a^{-\eps} (\frac{a}{b})^{\gamma\beta} |\LL|^{\beta}\). More generally, we have
\begin{align}\label{Eq:plank_Katz_Tao}
    \mu(\LL) \leq a^{-\eps} \Delta_{max}(\LL)^{1-\beta} \left(\frac{a}{b}\right)^{\gamma\beta} |\LL|^{\beta}.
\end{align}
\end{lem}

\begin{rem}
    We will only apply Lemma~\ref{Lem:plank_Katz_Tao} in the special cases $\gamma = 0$ (for which the nonconcentration condition is vacuous, and can be ignored), and $\gamma = 1$. When $\LL$ is Katz--Tao and $\gamma = 0$, the RHS of \eqref{Eq:plank_Katz_Tao} is $\lesssim a^{-\eps} |\LL|^\beta$, which matches the estimate $K_{KT}(\beta)$ for flakes in Definition~\ref{Def:Katz_Tao_beta}.
\end{rem}

\begin{rem}\label{Rem:trivial_case}
    The inequality automatically holds in the case when $\mu(\LL) \leq a^{-\eta}$. To see this, note that the nonconcentration hypothesis gives that $|\LL| \geq a^\eta (b/a)^\gamma$ (by taking $\theta = a/b$), and $\eta$ is much less than $\eps$ .
\end{rem}

{\bf Frostman multiplicity estimates.} Lemma~\ref{Lem:plank_Frostman} below says that the Kakeya estimate \( K_F(\beta) \) implies an analogous estimate for planks of dimensions \( a \times b \times \1\), provided we avoid the following bad scenario. Let \( \LL \) be a set of planks of dimensions \( a \times b \times \1\), \( L \in \LL \), and \( \theta \in [a/b, 1] \). We define the \textit{thickened plank} \( L_\theta \coloneqq N_{\theta b}(L) \), which has dimensions \( \theta b \times b \times \1\). If \( L_\theta \) contains a very large number of planks from \( \LL \), then this can cause an unacceptably large intersection multiplicity \( \mu(\LL) \). The next result says that this is the only obstruction to obtaining an analogue of \( K_F(\beta) \) for planks.

\begin{lem}\label{Lem:plank_Frostman}
Suppose that \( K_F(\beta) \) holds. Then for every \( \eps > 0 \), there exists \( \eta, b_0 > 0 \) such that the following holds for any \( 0 < a \leq b \leq b_0 \). Let \( \LL \) be a set of \( a \times b \times \1\) {\color{black} essentially distinct} planks in \( B_1 \), and let \( Y \) be a shading on \( \LL \) with density \( \lambda(\LL, Y) \geq a^\eta \). Suppose that the planks do not concentrate into thickened planks, in the following sense: There is a number \( M \geq 1 \) so that for each \( \theta \in [a/b, 1] \) and every \( L \in \LL \), we have \( |\LL[L_\theta]| \leq M\theta \). Then
\[
|U(\LL)| \geq a^{\eps} C_F(\LL)^{\beta/{\color{black} n} - 1} b^{2\beta} (M^{-1}b^2|\LL|)^{\beta/{\color{black} n}},
\]
and
\[
\mu(\LL) \leq a^{-\eps} C_F(\LL)^{1-\beta/{\color{black} n}} M^{\beta/{\color{black} n}} \frac{a}{b} b^{-2\beta} (b^2|\LL|)^{1-\beta/{\color{black} n}}.
\]
\end{lem}

\begin{rem}\label{Rem:trivial_case'}
    Both inequalities automatically hold in the case when \( \mu(\LL) \leq a^{-\eta} \). To see this, note that \( M \geq b/a \) (by taking \( \theta = a/b \)), \( C_F(\LL) \gtrsim \frac{1}{|\LL|ab} \), and $\eta$ is much less than $\eps$.
\end{rem}

{\bf Flakes that factor through planks.} Lemmas~\ref{Lem:plank_Katz_Tao}~and~\ref{Lem:plank_Frostman} allow us to obtain improved multiplicity bounds for families of flakes that factor through planks. The precise statement is as follows.

\begin{proposition}\label{Prop:plank_factor_improve}
Suppose that \( K_{KT}(\beta) \) and \( K_F(\beta) \) hold. Then for all \( \eps > 0 \), there exists \( \eta, \delta_0 > 0 \) so that the following holds for all \( \delta \leq \delta_0 \). Let \( (\F, Y) \) be uniform, with \( \lambda(\F, Y) \geq \delta^\eta \).

\begin{enumerate}[label={\rm (\Alph*)}]
    \item
    Suppose {\color{black} $\F$ is essentially distinct, and} there is a scale \( \rho \) so that for each \( F_\rho \in \F_\rho \), there is a family of \( a \times b \times \1\) planks that factors \( \F[F_\rho] \). Then \begin{align}\label{Eq:Part(A)_F}
    \mu(\F, Y) \leq \delta^{-\eps} C_F(\F)^{1-\beta/{\color{black} n}} \left( \frac{a}{b} \right)^{\color{black} (n+2)\beta/n} \delta^{-2\beta} (\delta^2 |\F|)^{1-\beta/{\color{black} n}}.
    \end{align}

    \item
    Suppose there is a scale \( \rho \) and a family of \( a \times b \times \1\) planks that factors \( \F_\rho \). Then
    \begin{align}\label{Eq:Part(B)_KT}
    \mu(\F, Y) \leq \delta^{-\eps} \Delta_{max}(\F)^{1-\beta} \left( \frac{a}{b} \right)^{\beta} |\F|^{\beta}. 
    \end{align}
\end{enumerate}
\end{proposition}

Note that the analogue of Proposition~\ref{Prop:tech_tau} applies to Lemmas~\ref{Lem:plank_Katz_Tao} and \ref{Lem:plank_Frostman}, and to Proposition~\ref{Prop:plank_factor_improve}. The remainder of this section is devoted to proving Lemmas~\ref{Lem:plank_Katz_Tao} and \ref{Lem:plank_Frostman}, and Proposition~\ref{Prop:plank_factor_improve}.

\subsection{Slabs.} 
A standard $L^2$ method gives strong bounds for the incidences of slabs. For example, Lemma~\ref{Lem:slab_0} gives the nice $\Gr(n,n-1)$ version of Theorem~\ref{Thm:Main} in a fairly easy way, and has been used many times so far. In this subsection, we state another two simple variants of Lemma~\ref{Lem:slab_0}, which will be useful in the remaining parts of this paper.

A key character is the notion of ``typical intersection angle'', which is meaningfully defined up to accuracy $\delta$. Suppose that $\Sb$ is a set of $\delta \times 1 \times \1$ slabs with a shading $Y$. Define the set of incidence triples
\begin{align*}
    Tri(\Sb) = \{(x, S_1, S_2): x\in Y(S_1)\cap Y(S_2)\}.
\end{align*}
There is a natural measure on $Tri(\Sb)$ coming from the volume measure for $x$ and the counting
measure for $S_1, S_2$. By Definition~\ref{Def:major_angle} with $k = n-1$, we can define the \textit{angle} between two slabs $S_1, S_2$ as their major angle (well-defined up
to additive error $\delta$). Then we define
\[
Tri_{\theta}(\Sb) = \{(x, S_1, S_2) \in Tri(\Sb) : \angle(S_1, S_2) \sim \theta\}.
\]

\begin{definition}\label{Def:typical_angle_slab}
We say that \(\theta\) is a \textit{typical angle of intersection} for \(\mathbb{S}\) if \(|Tri_{\theta}(\mathbb{S})| \approx |Tri(\mathbb{S})|\) up to \(\log(1/\delta)\) factors. If \(\mu(\mathbb{S}) \geq 2\), then \(\mathbb{S}\) always has a typical intersection angle \(\theta \in [\delta, 1]\).
\end{definition}

\begin{lem}\label{Lem:slab_1}
Suppose \(\mathbb{S}\) is a set of \(\delta \times 1 \times \1 \) slabs contained in \(\tilde{S}\), a \(\theta \times 1 \times \1\) slab, with \(\delta \leq \theta \leq 1\). Suppose the slabs of \(\mathbb{S}\) are shaded with \(\lambda(\mathbb{S}, Y) \geq \delta^\eta\). If \(\theta\) is a typical intersection angle for \(\mathbb{S}\), then \(|U(\mathbb{S}, Y)| \gtrsim \delta^{2\eta}|\tilde{S}|\).
\end{lem}

\begin{proof}
Exactly the same as that of \cite[Lemma~6.8]{GWZ26streamlined}, and so omitted.
\end{proof}

We also have the following analogue of Lemma~\ref{Lem:slab_0}, which has some flexibility in applications. 
\begin{lem}\label{Lem:slab_2}
Suppose \( \Sb \) is a set of \( \delta \times b_1 \times \cdots \times b_{n-1}\) convex sets in \( B_1 \subset \R^n\), and $Y$ is a shading on $\Sb$ with \( \lambda(\mathbb{S},Y) \geq \delta^\eta \). Suppose also \( b_j \geq \delta^\eta \) for each $j$. If \( \Delta_{max}(\Sb) \leq \delta^{-\eta} \), then \( \mu(\Sb,Y) \lesssim \delta^{-C_n\eta} \). If \( C_F(\Sb) \leq \delta^{-\eta} \), then \( |U(\Sb,Y)| \gtrsim \delta^{C_n\eta} \). Here $C_n$ is a constant depending only on $n$.
\end{lem}
\begin{proof}
    Similar to Lemma~\ref{Lem:slab_0}, and so omitted.
\end{proof}

\subsection{Planks}
Next we suppose that \( \LL \) is a set of \( a \times b \times \1\) planks in \( B_1 \). We will study \( \LL \) by relating it to incidence problems involving slabs and flakes. The incidence problem for slabs is solved completely using Lemma~\ref{Lem:slab_1}. Assuming \( K_{KT}(\beta) \) or \( K_F(\beta) \), we get information about the incidence problem for flakes. In this way, we are able to prove estimates for \( \mu(\LL) \) based on various geometric assumptions about \( \LL \).

The heuristic picture behind the whole process is very similar to that explained in \cite[Section~6.2]{GWZ26streamlined}. Given two intersecting planks, \( L_1, L_2 \), define the \textit{angle} between them as the angle between the hyperplanes \( TL_1, TL_2 \), where \( TL \) is the hyperplane spanned by the $(n-1)$ longest axes of \( L \). If the planks \( L_1 \) and \( L_2 \) have dimensions \( a \times b \times \1 \), then this angle is meaningfully defined up to accuracy \( a/b \). Let \( \theta \) denote a “typical” angle of intersection for \( \LL \).

First, we note that if the angle between \( L_1 \) and \( L_2 \) is \( \gg \theta \), then their intersection will not play an important role. Therefore, we can subdivide \( \LL \) into smaller collections of planks as follows. If \( S \) is a \( \theta \times 1 \times \1\) slab in \( B_1 \), we define \( \LL_S \) as in \eqref{Eq:L_S}. We can now decompose \( \LL \) as \( \LL = \bigsqcup_S \LL_S \), where \( S \) varies over some essentially distinct slabs. If \( \theta \) is a typical angle of intersection of \( (\LL, Y) \), then morally the sets \( U(\LL_S, Y) \) are disjoint and we can study them separately. Fix a single \( S \) and we will analyze \( \LL_S \).

Next, let \( Q \) denote a \( \theta b \times b \times b\1\) box, and let \( \LL_Q \) be the set of planks \( L \in \LL_S \) so that \( L \cap Q \) has dimensions \( a \times b \times b\1\). In other words, \( \LL_Q \) are the planks that intersect \( Q \) “tangentially”, making the intersection as large as possible. Assuming that \( \mu(\LL) \gg 1 \) and assuming that \( \LL_Q \) is not empty, then morally Lemma~\ref{Lem:slab_1} tells us that \( |U(\LL_S, Y) \cap Q| \approx |Q| \). For each plank \( L \in \LL_S \), the \( \theta b \)-neighborhood of \( L \) can be covered by boxes \( Q \) as above, and so we see that morally most of this neighborhood is contained in \( U(\LL_S, Y) \). This motivates the definition of the thickened version of \( \LL_S \) as follows.

\begin{definition}\label{Def:L_theta_S}
For each \( L \in \LL \), define \( L_\theta \) to be the \( \theta b \)-neighborhood of \( L \). Define \( \LL_\theta \) to be a set of essentially distinct prisms of the form \( L_\theta \), so that each \( L' \in \LL \) is contained in a set \( L_\theta \in \LL_\theta \) for which \( L'_\theta \) is comparable to \( L_\theta \). Note that if \( L_\theta \) and \( L'_\theta \) are comparable, then they are contained in comparable slabs \( S \) and \( S' \). Thus \( \LL_\theta \) decomposes into a boundedly overlapping union \( \bigcup_S \LL_{\theta,S} \), where the union is taken over incomparable slabs of dimensions \( \theta \times 1 \times \1\).
\end{definition}

In the previous paragraph, we saw that if \( \theta \) is a typical angle of intersection for \( \LL \), then morally \( U(\LL_S) \) is almost the same as \( U(\LL_{\theta,S}) \).

It remains to study \( U(\LL_{\theta,S}) \). There is a change of variables that takes \( S \) to the unit ball and takes each \( L_\theta \in \LL_{\theta,S} \) to a $b$-flake. Therefore, the problem of estimating \( |U(\LL_{\theta,S})| \) is equivalent to a problem about $b$-flakes which we can bound using \( K_{KT}(\beta) \) and/or \( K_F(\beta) \).

The following two lemmas give a precise version of the moral statements above. First we carefully define the typical angle of intersection.

\begin{lem}\label{Lem:typical_angle}
Let \((\LL, Y)\) be a family of \(a \times b \times \1\) planks and their associated shading. Let \(\eta > 0\) and suppose that \(\mu(\LL, Y) \geq a^{-\eta}\). Then there is an \(\approx 1\) refinement \((\LL, Y')\) of \((\LL, Y)\) and a number \(\theta \in [a/b, 1]\) with the following properties. In what follows, let \(A = e^{\sqrt{\log a^{-1}}}\).

\begin{itemize}
    \item \((\LL, Y')\) has constant multiplicity on \(U(\LL, Y')\) up to a factor \(\sim 1\).
    \item \(\max_{L_1, L_2 \in \LL_{Y'}(x)} \angle(TL_1, TL_2) \sim \theta\).
    \item If \(\LL_{Y''}(x) \subset \LL_{Y'}(x)\) with \(|\LL_{Y''}(x)| \geq \frac{1}{A}|\LL_{Y'}(x)|\), then
    \[
    \max_{L_1, L_2 \in \LL_{Y''}(x)} \max\{a/b, \angle(TL_1, TL_2)\} \approx \theta.
    \]
\end{itemize}
\end{lem}
\begin{proof}
    Similar to \cite[Lemma~6.11]{GWZ26streamlined}, and so omitted.
\end{proof}

\begin{definition}
We say that \( \theta \) is a \textit{typical angle of intersection} for \( (\LL, Y) \) if \( (\LL, Y) \) satisfies the conclusions of Lemma~\ref{Lem:typical_angle} for some \( A \geq e^{(\log(a^{-1}))^{1/4}} \).

We made the definition a little bit looser than the conclusion of Lemma~\ref{Lem:typical_angle} so that some mild refinements of \( (\LL, Y') \) will still have typical angle of intersection \( \theta \). In other words, at every point \( x \in U(\LL, Y) \), the intersecting planks have angle at most \( \theta \) and there are many pairs of planks with angle \( \sim \theta \). Furthermore, this remains true if the shading \( Y \) is replaced by a mild refinement. This is stronger than the definition that we used for typical intersection angle of slabs, Definition~\ref{Def:typical_angle_slab}. 
\end{definition}

\begin{lem}\label{Lem:plank_thicken}
Let \( \eta > 0 \) and let \( (\LL, Y) \) be a set of \( a \times b \times \1 \) planks in $B_1$ with \( \lambda(\LL, Y) \geq a^\eta \) and \( \mu(\LL, Y) \geq a^{-\eta} \).

Then there is a $\approx 1$ refinement \( (\LL', Y') \) of \( (\LL, Y) \) with \( |\LL'| \gtrapprox |\LL| \) and \( \lambda(\LL', Y') \gtrapprox \lambda(\LL, Y) \), a number \( N \geq 1 \), a number \( \theta \in [a/b, 1] \), and a set of \( \theta \times 1 \times \1 \) slabs \( \Sb \) so that the following holds. For each \( S \in \Sb \), let \( \LL_S' \) and \( \LL_{\theta,S}' \) be as defined in \eqref{Eq:L_S} and Definition~\ref{Def:L_theta_S}. \( \LL' = \bigcup_{S \in \Sb} \LL_S' \). Then each \( L_\theta \in \LL_{\theta,S}' \) contains \( \sim N \) planks \( L \in \LL_S' \). {\color{black} And $|\LL_{\theta,S}'|$ is constant over $S \in \Sb$ up to a factor $\sim1$.} Furthermore, there is a shading \( Y_{\theta,S} \) on \( \LL_{\theta,S}' \) so that

\begin{enumerate}[label={\rm (\arabic*)}]
    \item If \( U(\LL', Y') \cap B_{\theta b} \) is non-empty, then \( |U(\LL', Y') \cap B_{\theta b}| \gtrapprox a^{4\eta} |B_{\theta b}| \).
    \item For each \( S \in \Sb \), \( \lambda(\LL_{\theta,S}', Y_{\theta,S}) \gtrapprox a^{4\eta} \) and for each \( L_\theta \in \LL_{\theta,S}', Y_{\theta,S}(L_\theta) \subset U(\LL_S', Y') \).
    \item \( |U(\LL, Y)| \gtrapprox \sum_{S \in \Sb} |U(\LL_{\theta,S}', Y_{\theta,S})| \).
    \item \( \mu(\LL, Y) \lessapprox \max_{S \in \Sb} a^{-4\eta} \frac{aN}{b\theta} \mu(\LL_{\theta,S}', Y_{\theta,S}) \).
\end{enumerate}
\end{lem}
\begin{proof}
    Similar to \cite[Lemma~6.13]{GWZ26streamlined}, and so omitted. We only point out that we apply Lemma~\ref{Lem:typical_angle} instead of \cite[Lemma~6.11]{GWZ26streamlined} to find the typical angle of intersection $\theta$, and apply Lemma~\ref{Lem:slab_1} instead of \cite[Lemma~6.8]{GWZ26streamlined} to deal with intersection of local $a \times b \times b\1$ slabs.
\end{proof}
\begin{rem}\label{Rem:angle->thicken}
If \( (\LL, Y) \) has a typical angle of intersection \( \theta \), then Lemma~\ref{Lem:plank_thicken} holds for this value of \( \theta \).
\end{rem}

By a linear change of variables, we can convert $\LL_{\theta,S}$ from a set of $\theta b \times b \times \1$ planks in $S$ into a
set $\F$ of $b$-flakes in $B_1$. To see this more carefully, recall that $\LL_{\theta,S}$ consists of $\theta b \times b \times \1$ planks $L$
in the slab $S$ so that $TL$ is tangent to $TS$. The slab $S$ has dimensions $\theta \times 1 \times \1$. Let $ \phi: S \rightarrow B_1$
be a linear change of variables. Because of the tangency condition on $L$, $\phi(L)$ has dimensions
$b \times b \times \1$.

To summarize, we have reduced the original incidence problem for $\LL$ to an incidence problem
for $\LL_{\theta,S}$, which is equivalent to an incidence problem for $b$-flakes $\F$. If we assume $K_{KT}(\beta)$ or
$K_{F}(\beta)$, then we can get bounds for $\mu(\F)$ and hence bounds for $\mu(\LL)$. We are now ready to prove
Lemmas~\ref{Lem:plank_Katz_Tao} and \ref{Lem:plank_Frostman}.

\subsection{A Katz--Tao multiplicity estimate}
\begin{proof}[Proof of Lemma~\ref{Lem:plank_Katz_Tao}]
The arguments closely follow those in \cite[Section~6.3]{GWZ26streamlined}. We may suppose that \(\mu(\LL, Y) \geq a^{-\eta}\), since otherwise we are done by Remark~\ref{Rem:trivial_case}. Apply Lemma~\ref{Lem:plank_thicken}. Abusing notation, we will continue to refer to the resulting refinement by \((\LL, Y)\), and we write $\LL_S$ and $\LL_{\theta,S}$ in place of $\LL_S'$ and $\LL_{\theta,S}'$. By Lemma~\ref{Lem:plank_thicken}, there is a number \(\theta \in [a/b, 1]\), a number \(N\), and a slab \(S\) of dimensions \(\theta \times 1 \times \1\) so that we have
\begin{align}\label{Eq:to_plank_thicken}
\mu(\LL, Y) \lessapprox \max_{S \in \Sb} a^{-4\eta} \frac{aN}{b\theta} \mu(\LL_{\theta,S}, Y_{\theta,S}).
\end{align}

Lemma~\ref{Lem:plank_thicken} also gives a density lower bound
\[
\lambda(\LL_{\theta, S}, Y_{\theta, S}) \gtrapprox a^{4\eta}.
\]

By a linear change of variables, we can convert $\LL_{\theta,S}$ into a
set $\F$ of $b$-flakes in $B_1$. Since convex sets are invariant with respect to a linear change of variables, we have
\[
\Delta_{max}(\F) = \Delta_{max}(\LL_{\theta, S}) \leq \Delta_{max}(\LL_\theta) \lesssim \frac{b\theta}{aN} \Delta_{max}(\LL).
\]

Now apply Lemma~\ref{Lem:general_KT} (noting Proposition~\ref{Prop:tech_tau}) with \(b, \eps/2\) in place of \(\delta, \eps\). If \(\eta\) and \(b_0\) are chosen sufficiently small, then
\begin{align*}
\mu(\LL_{\theta, S}) 
&
= \mu(\F) \leq a^{-\eps/2} \Delta_{max}(\F)^{1-\beta} |\F|^\beta\\
& 
\lesssim a^{-\eps/2} \Delta_{max}(\LL)^{1-\beta} \left( \frac{b\theta}{aN} \right)^{1-\beta} (|\LL_{\theta, S}|)^\beta\\
&
\lesssim a^{-\eps/2} \Delta_{max}(\LL)^{1-\beta} \left( \frac{b\theta}{aN} \right)^{1-\beta} (N^{-1}|\LL_{S}|)^\beta.
\end{align*}

And so plugging in \eqref{Eq:to_plank_thicken}, we get
\begin{align*}
    \mu(\LL, Y) 
    &
    \lessapprox a^{-4\eta} \frac{aN}{b\theta} \mu(\LL_{\theta, S})\\
    &
    \lesssim a^{-\eps/2-4\eta} \Delta_{max}(\LL)^{1-\beta} \left( \frac{b\theta}{aN} \right)^{-\beta} (N^{-1}|\LL_{S}|)^\beta\\
    &
    = a^{-\eps/2-4\eta} \Delta_{max}(\LL)^{1-\beta} \left( \frac{a}{b\theta} \right)^\beta |\LL_S|^\beta.
\end{align*}

We know \(a \leq b\theta\) (and thus \(\frac{a}{b} \leq (\frac{a}{b})^{\gamma}\theta^{1-\gamma}\)), and \(|\LL_S| \leq a^{-\eta}\theta^\gamma|\LL|\). Plugging in, we get
\[
\mu(\LL) \leq a^{-\eps} \Delta_{max}(\LL)^{1-\beta} \left( \frac{a}{b} \right)^{\gamma\beta} |\LL|^\beta,
\]
as desired.
\end{proof}

\subsection{A Frostman multiplicity estimate}
\begin{proof}[Proof of Lemma~\ref{Lem:plank_Frostman}]
The arguments closely follow those in \cite[Section~6.4]{GWZ26streamlined}. An additional twist is that the exponents in $\beta$ need to be {\color{black} altered accordingly} to match those in Lemma~\ref{Lem:general_F}. We may suppose that \(\mu(\LL, Y) \geq a^{-\eta}\), since otherwise we are done by Remark~\ref{Rem:trivial_case'}. Apply Lemma~\ref{Lem:plank_thicken}. Abusing notation, we will continue to refer to the resulting refinement by \((\LL, Y)\), and we write $\LL_S$ and $\LL_{\theta,S}$ in place of $\LL_S'$ and $\LL_{\theta,S}'$.

By Lemma~\ref{Lem:plank_thicken}, there is a number \(\theta \in [a/b, 1]\) and a set \(\Sb\) of \(\theta \times 1 \times \1\) slabs so that
\[
|U(\LL, Y)| \gtrapprox |\Sb| \, |U(\LL_{\theta,S}, Y_{\theta,S})|
\]
for some $S \in \Sb$.

Lemma~\ref{Lem:plank_thicken} also gives a density lower bound
\[
\lambda(\LL_{\theta,S}, Y_{\theta,S}) \gtrapprox a^{4\eta}.
\]

Define \(C = C_F(\LL, B_1)\). This does not imply that \(C_F(\LL_{\theta,S}, S) \leq C\). So we have to work out what we know about \(\LL_{\theta,S}\). Let \(\LL_\theta\) be formed by thickening the planks of \(\LL\) to dimensions \(\theta b \times b \times \1\). By Lemma~\ref{Lem:plank_thicken} and our nonconcentration assumption, there is a number \(N \leq M\theta\) so that each plank \(L_\theta \in \LL_\theta\) contains \(\sim N\) planks of \(\LL\) (recall that we have replaced \(\LL\) by a \(\approx 1\) refinement, but the refinement preserves the upper bound \(|\LL[L_\theta]| \leq M\theta\)).

By Proposition~\ref{Prop:inherit}(A), \(C_F(\LL_\theta, B_1) \lessapprox C_F(\LL, B_1) = C\), and so
\[
|\LL_{\theta,S}| \leq |\LL_\theta[S]| \lessapprox C\theta |\LL_\theta|.
\]

The set \(\LL_{\theta,S}\) is not necessarily Frostman, but since \(\LL_\theta\) is Frostman, we do know that for any convex \(K \subset S\),
\[
\Delta(\LL_{\theta,S}, K) \leq \Delta(\LL_\theta, K) \lessapprox C \Delta(\LL_\theta, B_1) \sim C |L_\theta| |\LL_\theta|.
\]

Since \(|\LL_\theta| \approx |\Sb|\,|\LL_{\theta,S}|\),
\[
\Delta(\LL_{\theta,S}, K) \lessapprox C|\Sb||S| \frac{|L_\theta|\,|\LL_{\theta,S}|}{|S|} = C|\Sb||S|\,\Delta(\LL_{\theta,S}, S).
\]
Therefore \(C_F(\LL_{\theta,S}, S) \lessapprox C|\Sb||S|\).

Now by a linear change of variables, we can convert \(\LL_{\theta,S}\) to a set of \(b\)-flakes \(\F\) in \(B_1\). Let \(Y_{\F}\) be the image of the shading \(Y_{\theta,S}\) under this change of variables. Since convex sets are invariant with respect to a linear change of variables, \(C_F(\F) \lessapprox C|\Sb||S|\). Also, \(\lambda(\F, Y_{\F}) = \lambda(\LL_{\theta,S}, Y_{\theta,S}) \gtrapprox a^{4\eta}\).

If \(\eta\) and \(b_0\) are chosen sufficiently small depending on \(\eps\), then by Lemma~\ref{Lem:general_F} (noting Proposition~\ref{Prop:tech_tau}) to $(\F, Y_\F)$ with $\eps/2$ in place of $\eps$, we have
\begin{align*}
    |U(\F)| &= \frac{\sum_\F |Y(F)|}{\mu(\F,Y)}\\
    &\gtrapprox 
    a^{\eps/2+4\eta} (C|\Sb||S|)^{\beta/n-1} b^{2\beta} (b^2|\F|)^{\beta/n}\\
    &\approx a^{\eps/2 + 4\eta} C^{\beta/n-1} (|\Sb||S|)^{-1} b^{2\beta} (b^2 |S| |\LL_\theta|)^{\beta/n}\\
    &= a^{\eps/2 + 4\eta} C^{\beta/n-1} (|\Sb||S|)^{-1} b^{2\beta} (b^2 \theta |\LL|/N)^{\beta/n}\\
    &\geq a^{\eps/2 + 4\eta} C^{\beta/n-1} (|\Sb||S|)^{-1} b^{2\beta} (M^{-1} b^2 |\LL|)^{\beta/n}.
\end{align*}

Therefore, if \(b_0 > 0\) is sufficiently small, then we have
\[
|U(\LL)| \gtrapprox |\Sb||U(\LL_{\theta,S})| \sim |\Sb| |S| |U(\F)|
\gtrapprox
a^{\eps/2 + 4\eta} C^{\beta/n-1} b^{2\beta} (M^{-1} b^2 |\LL|)^{\beta/n}.
\]
This gives the desired lower bound for \(|U(\LL)|\). To get an upper bound for \(\mu(\LL)\), we relate \(\mu(\LL)\) to \(|U(\LL)|\). This part is just algebra and we leave the details to the reader.
\end{proof}

\subsection{Flakes that factor through planks}
In this section we will prove Proposition~\ref{Prop:plank_factor_improve}. The arguments closely follow those in \cite[Section~6.5]{GWZ26streamlined}. An additional twist is that for Part (A), the exponents in $\beta$ need to be {\color{black} altered accordingly} to match those in Lemma~\ref{Lem:plank_Frostman}, {\color{black} and the upper bound of one $M$ during the proof is also different}.

{\bf Part (A).} 
Let \( C = C_F(\F) \). After refining \( \F_\rho \) by a \( O(1) \) factor (this induces a harmless refinement on \( \F \)), we may suppose that each \( F \in \F \) is contained in exactly one flake \( F_\rho \in \F_\rho \). By the hypothesis, for each \( F_\rho \in \F_\rho \), there is a set \( \W_{F_\rho} \) of \( a \times b \times \1 \) planks that factors \( \F[F_\rho] \), i.e.
\[
\F[F_\rho] = \bigsqcup_{\W_{F_\rho}} \F[F_\rho]_W,
\]
where \( \F[F_\rho]_W \) is \( \sim 1 \) Frostman in \( W \). Define \( \W = \bigsqcup_{\F_\rho} \W_{F_\rho} \). There is a corresponding decomposition \( \F = \bigsqcup_\W \F_W \), where \( \F_W \) is \( \sim 1 \) Frostman in \( W \). 

Apply Proposition~\ref{Prop:fact_mult} to the pair \( (\F, Y) \) and \( \W \); we obtain a set \( \W' \subset \W \), and a shading \( Y_{\W'} \) on \( \W' \) with \( \lambda(\W', Y_{\W'}) \gtrapprox \delta^{2\eta} \). We also obtain a shading \( Y'(F) \subset Y(F) \) on \( \F' = \bigcup_{\W'} \F_W' \), so that \( (\F', Y') \) is an \( \approx 1 \) refinement of \( (\F, Y) \), and
\begin{align}\label{Eq:W_factor_P}
\mu(\F, Y) \lessapprox \mu(\W', Y_{\W'})\,\mu(\F_W', Y') \quad \text{for all } W \in \W'.
\end{align}
After pigeonholing and refining $\F'$ and \( \W' \), we may suppose that \( |\F_W'| \) is approximately the same for each \( W \in \W' \). After this step we also have \( C_F(\F') \lessapprox C \), \( \lambda(\F', Y) \gtrapprox \delta^{\eta} \), and \( C_F(\F_W', W) \lessapprox 1 \) for each \( W \in \W' \). Since \( C_F(\F') \lessapprox C \) and \( |\F_W'| \) is approximately the same for all \( W \in \W' \), by Proposition~\ref{Prop:inherit}(A) we have that \( C_F(\W') \lessapprox C \). Thus we will use Lemma~\ref{Lem:plank_Frostman} to bound each of the terms on the RHS of \eqref{Eq:W_factor_P}. In order to apply Lemma~\ref{Lem:plank_Frostman}, we need to control to what extent the planks in \( \W' \) can concentrate into thicker planks, and similarly for the (rescaled) flakes in \( \F_W' \) (this is the quantity “\( M \)” from Lemma~\ref{Lem:plank_Frostman}).

We begin with \( \mathbb{W}' \). For each \( \theta \in [a/b, 1] \) we need to bound the maximum number of \( a \times b \times \1 \) planks from \( \mathbb{W}' \) that are contained in a \( \theta b \times b \times \1\) plank. Denote this number by \( N(\theta) \). Fix a plank \( W \in \mathbb{W}' \). Recall that \( W \in \mathbb{W}_{F_\rho} \) for some \( F_\rho \in \F_\rho \). Let \( W_\theta \) be the corresponding \( \theta b \times b \times \1 \) plank. Let \( W' \in \mathbb{W}'[W_\theta] \). Since \( W' \) contains at least one $\delta$-flake from \( \F \), and each $\delta$-flake in \( \F \) is contained in exactly one $\rho$-flake from \( \F_\rho \), we have \( F_\rho \) is the unique tube containing \( W' \), and thus \( W' \in \mathbb{W}_{F_\rho} \). We conclude that \( \mathbb{W}'[W_\theta] \subset  \mathbb{W}_{F_\rho}[W_\theta] \). Recall that \( \Delta_{max}(\mathbb{W}_{F_\rho}) \lesssim 1 \), and thus
\[
|\mathbb{W}_{F_\rho}[W_\theta]| \lesssim \frac{|W_\theta|}{|W|} = \frac{\theta b \times b \times \1}{a \times b \times \1} = \frac{\theta b}{a}.
\]
Thus \( |\mathbb{W}'[W_\theta]| \leq M\theta \) for \( M \lesssim b/a \). Applying Lemma~\ref{Lem:plank_Frostman} (and recalling Proposition~\ref{Prop:tech_tau}) with $\eps/3$ in place of $\eps$, we have
\begin{align}\label{Eq:(A)_coarse}
\mu(\mathbb{W}', Y_{\mathbb{W}'}) \lesssim
\delta^{-\eps/3} C_F(\F)^{1-\beta/n} (a/b)^{1-\beta/n} b^{-2\beta} (b^2|\mathbb{W}'|)^{1-\beta/n}.
\end{align}

Next we consider \( \F_W' \). Fix a plank \( W \in \mathbb{W}' \). We first have to make a change of coordinates that maps \( W \) to \( B_1 \) and each \( \delta \)-flake \( F \in \F_W' \) to a \( b^{-1}\delta \times a^{-1}\delta \times \1 \) plank. Next we need to bound the number of \( b^{-1}\delta \times a^{-1}\delta \times \1 \) planks in a \( \theta a^{-1}\delta \times a^{-1}\delta \times \1 \) thick plank. Undoing the coordinate change, this is the number of $\delta$-flakes of \( \F_W' \) in a \( \delta \times \frac{\theta b}{a}\delta \times \1 \) plank. Since \( \F \) is essentially distinct, there are 
\(
\lesssim (\frac{\theta b}{a} )^{\color{black} n-1} \leq (b/a)^{n-1} \theta
\)
such $\delta$-flakes. Applying Lemma~\ref{Lem:plank_Frostman} and remembering that the dimensions of each plank are  
\(
b^{-1}\delta \times a^{-1}\delta \times \1
\)
(and using our bound \( M \lesssim (b/a)^{\color{black} n-1} \)), we get
\begin{align}\label{Eq:(A)_fine}
\mu(\F_W', Y') \lesssim \delta^{-\eps/3}(a/b)^{1- \color{black}{(n-1)\beta/n}}(a^{-1}\delta)^{-2\beta}(a^{-2}\delta^2|\F_W'|)^{1-\beta/n}.
\end{align}

Combining \eqref{Eq:W_factor_P}, \eqref{Eq:(A)_coarse}, and \eqref{Eq:(A)_fine}, we get \eqref{Eq:Part(A)_F}.

{\bf Part (B).}
By hypothesis, there is a set \( \W \) of \( a \times b \times \1 \) planks that factors \( \F_\rho \), i.e., 
\[
\F_\rho = \bigsqcup_\W \F_{\rho,W},
\]
and \( \F_{\rho,W} \) is \( \sim 1 \) Frostman in \( W \). This induces a decomposition \( \F = \bigsqcup_\W \F_W \). Apply Proposition~\ref{Prop:fact_mult} to \( (\F, Y) \) and \( \W \); while \( \F_W \) might not be Frostman in \( W \), we have that \( \F_{\rho,W} \) is Frostman in \( W \). By Remark~\ref{Rem:weaker_fact_mult}, this is enough to apply the proposition. We obtain a set \( \W' \subset \W \), a shading \( Y_{\W'} \) on \( \W' \) with
$\lambda(\W', Y_{\W'}) \gtrapprox \delta^{2\eta}$,
and a subshading $Y'(F) \subset Y(F)$ 
on \( \F' = \bigsqcup_{\W'} \F_W' \), so that \( (\F', Y') \) is an \( \approx 1 \) refinement of \( (\F, Y) \), and
\begin{align}\label{Eq:W_factor_P'}
\mu(\F, Y) \lessapprox \mu(\W', Y_{\W'})\,\mu(\F_W', Y') \quad \text{for all } W \in \W'.
\end{align}
{\color{black} After pigeonholing and refining $\F'$ and \( \W' \), we may suppose that $|Y'(F)|$ is approximately the same for each $F \in \F$, 
and \( |\F_W'| \) is approximately the same for each \( W \in \W' \).} Thus $|\F_W'| \gtrapprox |\F_W|$, and \( \lambda(\F_W', Y') \gtrapprox \delta^{\eta} \) for each \( W \in \W' \). 
We also have $\Delta_{max}(\W') \leq \Delta_{max}(\W) \lesssim 1$ and 
$\Delta_{max}(\F_W') \leq \Delta_{max}(\F)$.
We will bound the two terms on the RHS of \eqref{Eq:W_factor_P'} using Lemma~\ref{Lem:plank_Katz_Tao} with \( \gamma = 0 \) and \( \gamma = 1 \), respectively.

We begin with $\W'$. 
Applying Lemma~\ref{Lem:plank_Katz_Tao} with \( \gamma = 0 \) (and recalling Proposition~\ref{Prop:tech_tau}), we obtain
\begin{align}\label{Eq:(B)_coarse}
\mu(\W', Y_{\W'}) \lesssim \delta^{-\eps/3}\,|\W'|^\beta.
\end{align}

Next we consider $\F_W'$. To apply Lemma~\ref{Lem:plank_Katz_Tao} with \( \gamma = 1 \), we make a change of coordinates that maps \( W \) to \( B_1 \). This change of coordinates converts \( \F_W' \), a set of \( \delta \)-flakes \( F \subset W \), to \( \LL \), a set of \( b^{-1}\delta \times a^{-1}\delta \times \1 \) planks \( L \subset B_1 \). Notice that the eccentricity of the planks is \( \frac{b^{-1}\delta}{a^{-1}\delta} = \frac{a}{b} \). Let \( Y_{\LL} \) be the image of the shading \( Y'(F) \) after this transformation.

Next we need to check that \( \LL \) obeys the assumptions of Lemma~\ref{Lem:plank_Katz_Tao} with \( \gamma = 1 \). The nonconcentration assumption of Lemma~\ref{Lem:plank_Katz_Tao} says that if \( \theta \in [a/b, 1] \) and \( S \) is a \( \theta \times 1 \times \1 \) slab, then \( |\LL_S| \lesssim (b^{-1}\delta)^{-\eta} \theta |\LL| \). If we change coordinates back from \( B_1 \) to \( W \), then \( S \) is converted to a convex set \( S' \subset W \) and our desired bound becomes
\[
|\F_W'[S']| \lessapprox (b^{-1}\delta)^{-\eta} \frac{|S'|}{|W|} |\F_W'|.
\]
Suppose that \( S' \) has dimensions \( a_1 \times a_2 \times \cdots \times a_n \). We claim that \( a_1 \gtrsim a \geq \rho \). We will come back to this claim below. Assuming this claim, we have
\[
\frac{|\F_W'[S']|}{|\F_W'|} \lessapprox \frac{|\F[S']|}{|\F_W|} \sim \frac{|\F_\rho[S']|}{|\F_{\rho,W}|}
\lesssim \frac{|S'|}{|W|} C_F(\F_{\rho,W}, W) \lesssim \frac{|S'|}{|W|}.
\]
This gives the desired bound on \( |\F_W'[S']| \).

To finish we have to check that \( a_1 \gtrsim a \). Choose coordinates \( x_1, x_2, \dots, x_n \) on \( W \) so that \( W \) is the rectangular box
\[
0 \leq x_1 \leq a,\quad 0 \leq x_2 \leq b,\quad 0 \leq x_j \leq 1 \quad (3 \leq j \leq n).
\]
The change of coordinates \( \phi \) mapping \( W \) to \( B_1 \) is
\[
\phi(x_1, x_2, x_3, \dots, x_n) = \left( \frac{x_1}{a}, \frac{x_2}{b}, x_3, \dots, x_n \right).
\]

We can assume that \( \LL_S \) is non‑empty, and so there is at least one \( L \in \LL \) with \( TL = TS \) (up to angle \( \theta \)). Suppose \( L = \phi(F) \) where \( F \) is a $\delta$-flake from \( \F[W] \).

Recall that \( TS \) is an affine $(n-1)$‑plane. Let \( TS' = \phi^{-1}(TS) \), another affine $(n-1)$‑plane. Choose unit vectors \( v_1, v_2, \dots, v_{n-1} \) to make a basis for \( TS' \) with \( v_1 \) in the \( (x_1, x_2) \)-plane and \( v_2, \dots, v_{n-1} \) almost parallel to \( F \). We have
\[
|v_1 \wedge v_2 \wedge \cdots \wedge v_{n-1}| \sim 1,
\]
since \( F \) is transverse to the \( (x_1, x_2) \) plane. Because \( \phi \) maps \( TS' \) to \( TS = TL \),
\[
|\phi(v_1) \wedge \phi(v_2) \wedge \cdots \wedge \phi(v_{n-1})| \sim a^{-1}.
\]
Since \( |\phi(v_j)| \sim 1 \) for all $2\leq j \leq n-1$, we see \( |\phi(v_1)| \sim a^{-1} \), and so \( v_1 \) is transverse to the \( x_2 \)-direction. Let \( \ell_{S'} \) be the line normal to \( TS' \). Since \( v_2, \dots, v_{n-1} \) are almost in the \( x_3, \dots, x_n \) directions and \( v_1 \) is transverse to the \( x_2 \)-direction, we have \( \ell_{S'} \) transverse to the \( x_1 \)-direction. Now \( W \) can be covered by \( \theta^{-1} \) parallel copies of \( S' \) and so \( a_1 \theta^{-1} \gtrsim b \) and consequently \( a_1 \gtrsim \theta b \gtrsim a \). This confirms the nonconcentration assumption in Lemma~\ref{Lem:plank_Katz_Tao} with \( \gamma = 1 \).

We can now apply Lemma~\ref{Lem:plank_Katz_Tao} with \(\gamma = 1\). Recall that \(\Delta_{max}(\LL) = \Delta_{max}(\F_W') \leq \Delta_{max}(\F)\). We have
\begin{align}\label{Eq:(B)_fine}
\mu(\F_W', Y') \lesssim \delta^{-\eps/3} \Delta_{max}(\F)^{1-\beta} \left( \frac{a}{b} \right)^{\beta} |\F_W'|^{\beta}.
\end{align}

Combining \eqref{Eq:W_factor_P'}, \eqref{Eq:(B)_coarse}, and \eqref{Eq:(B)_fine}, we obtain \eqref{Eq:Part(B)_KT}.

\section{A multiscale decomposition}\label{Sec:reduce_to_sticky}

In order to apply Theorem~\ref{Thm:2-plane_sticky}, our collection of flakes must be Frostman at every scale or Katz--Tao at every scale, respectively. Our main result in this section says that either this is true, or else there is a long sequence of scales that are far from being Frostman (resp. Katz--Tao). It is reminiscent of \cite[Lemma~7.7]{GWZ26streamlined} but is much more complicated. The reason is that to make the application of multiscale decompositions in Section~\ref{Sec:Main_Lemma_1}~and~\ref{Sec:Main_Lemma_2} legitimate, one needs to preserve lower bounds of some convex Wolff constants under the maximal density factoring (Lemma~\ref{Lem:factoring}). Also, we want certain ``strongly non-sticky'' property to be preserved under the multiplicity factoring (Lemma~\ref{Lem:fact_mult}). These are not guaranteed by a version like \cite[Lemma~7.7]{GWZ26streamlined}.\footnote{This technicality is nontrivial and also appeared in \cite{GWZ26streamlined}, but was neglected therein.}

Before everything, let us first state a simple variant of Lemma~\ref{Lem:factoring}:
\begin{lem}\label{Lem:factoring'}
    Let $\F$ be a set of essentially distinct $\delta$-flakes in $B_1 \subset \R^n$. Then there is a set $\F' \subset \F$ with $|\F'| \gtrapprox |\F|$; a family $\W$ of convex sets in $\R^n$, each of approximately the same dimensions; and a partition
    \begin{align*}
        \F' = \bigsqcup_{W\in\W} \F_W,
    \end{align*}
    where each $\F_W$ is a uniform set contained in $W$. Moreover, 
    \begin{enumerate}[label={\rm (\roman*)}]
        \item $\Delta_{max}(\W) \lesssim 1$.
        \item For each $W \in \W$, we have $C_F(\F_W, W) \lessapprox 1$ and $\Delta(\F_W, W) \approx \Delta_{max}(\F')$.
        \item For any two scales $\delta \leq \tau \leq \theta \leq 1$, there exist constants $U_1$, $U_2$, $U_3$, such that for any $W \in \W$, we have 
        $$\quad |(\F_W)_\tau[F_\theta]| \approx U_1,\, C_F((\F_W)_\tau[F_\theta]) \approx U_2,\, \Delta_{max}((\F_W)_\tau[F_\theta]) \approx U_3$$ for all $F_\theta \in \F_\theta$ with $(\F_W)_\tau[F_\theta] \neq \varnothing$.
    \end{enumerate}
\end{lem}
Note that Item~(iii) above is stronger than the mere uniformity of each $\F_W$: It requires the \textit{uniformized data} 
$$|(\F_W)_\tau[F_\theta]|, \, C_F((\F_W)_\tau[F_\theta]), \, \Delta_{max}((\F_W)_\tau[F_\theta])$$ 
to be essentially the same for all $W \in \W$. 
The proof of Lemma~\ref{Lem:factoring'} goes as follows. We first use Lemma~\ref{Lem:factoring} to get an initial partition, and then apply Proposition~\ref{Prop:CWA_uniform} to each $\F_W$ in the output. These give the uniformity of each $\F_W$ and Items~(i)~and~(ii). Item~(iii) can be easily obtained by applying another uniformization procedure to the uniformized data of different $\F_W$, which also preserves Items~(i)~and~(ii). We omit the details. 

\begin{definition}\label{Def:weak_uniform}
    We call a set of $\delta$-flakes $\F$ \textit{weakly uniform} at scales $[\tau_0, \theta_0] \subset [\delta,1]$ up to a factor $C$, if for any $\tau_0 \leq \tau \leq \theta \leq \theta_0$, there exist constants $U_1$, $U_2$, $U_3$, such that $$|\F_\tau[F_\theta]| \in [C^{-1}U_1, CU_1], \, C_F(\F_\tau[F_\theta]) \in [C^{-1}U_2, CU_2], \, \Delta_{max}(\F_\tau[F_\theta]) \in [C^{-1}U_3, CU_3]$$ 
    for all $F_\theta \in \F_\theta$.
\end{definition}

The benefit of Lemma~\ref{Lem:factoring'} is that we have strong uniformity for each $\F_W$, while keeping the factoring properties essentially unchanged. For convenience, we introduce the following definition, which is a variant of Definition~\ref{Def:factor}.

\begin{definition}\label{Def:factor'}
    Let $\F$ be a set of essentially distinct $\delta$-flakes in $B_1 \subset \R^n$, and let $\W$ be a family of convex sets in $\R^n$. Suppose that there is a partition $\F = \bigsqcup_{W\in\W} \F_W$. We say that \textit{$\W$ uniformly factors $\F$} if the following is true: 
    \begin{itemize}
        \item $\Delta_{max}(\W) \lesssim 1$.
        \item Each set $W \in \W$ has approximately the same dimensions. 
        \item For each $W \in \W$, we have $C_F(\F_W, W) \lessapprox 1$ and $\Delta(\F_W, W) \approx \Delta_{max}(\F)$.
        \item Item~(iii) in Lemma~\ref{Lem:factoring'} holds.
    \end{itemize}
\end{definition}
Thus Lemma~\ref{Lem:factoring'} says that $\W$ uniformly factors $\F'$.

For any set $A$, let $2^A$ be the power set of $A$, i.e., the collection of all subsets of $A$. 
\begin{definition}\label{Def:seleciton_map}
    Let $\rho \in [\delta
    ,1]$, $\lambda_* \in (0,1]$, and $\F$ be a set of $\delta$-flakes. Then we call a map $\Phi_{\F, \rho}: 2^\F \rightarrow 2^{\F_\rho}$ a \textit{$\lambda_*$-dense selection map at scale $\rho$}, if for any $\underline{\F} \subset \F$, we have $\Phi_{\F, \rho}(\underline{\F}) \subseteq \underline{\F}_\rho$, and $|\bigsqcup_{F_\rho \in \Phi_{\F, \rho}(\underline{\F})} \underline{\F}[F_\rho]| \geq \lambda_* |\underline{\F}|$.
\end{definition}

Here is the main result in this section:
\begin{lem}\label{Lem:multi-scale_decomposition}
Let \( N \geq 10^5 \) be an integer. Define \( \eps = (10/N)^{1/5} < 1/3 \), and let \( \eta \leq \eta_0 \leq \eta_1 \leq \eta_2 \leq \ldots \leq \eta_N \leq \eps \). Let \( \delta > 0 \) and let \( \F \) be a set of {\color{black} essentially distinct} \( \delta \)-flakes in $B_1 \subset \R^n$. Let $\lambda_* \in (0,1]$. For each $\rho \in [\delta, 1]$, fix a $\lambda_*$-dense selection map $\Phi_{\F, \rho}$ at scale $\rho$. 

\begin{enumerate}[label={\rm (\Alph*)}]
\item 
Fix $\gamma>0$, and let $\eta_{j-1}' \coloneqq \frac{14\eta_{j-1}}{\eps^2\gamma}$ for each $1 \leq j \leq N$. Suppose that $\eta_0 \geq \eta/\eps$, and $\eta_j \geq \eta_{j-1}'/\eps$ for each $1 \leq j \leq N$. If \( C_F(\F) \leq \delta^{-\eta} \), then there exists a $\gtrapprox_\eps \lambda_*^{2/\eps^3}$ refinement $\F' \subset \F$ such that at least one of the following holds.
\begin{enumerate}[label={\rm (\Roman*)}]
    \item
    \( \F' \) is \( \lessapprox_\eps \lambda_*^{-2/\eps^3} \delta^{-{\color{black} 2n}\eps} \) Frostman at every scale.
    \item
    There are scales \( \delta \leq \tau \leq \theta \leq 1 \) with \( \tau \leq \delta^{\eps^2}\theta \), an integer \( 1 \leq J \leq N \), and $\tau \leq a \leq b \leq \theta(\tau/\theta)^\eps$, such that the following holds:
    \begin{itemize} 
        \item $\F'$ is uniform at scales $\tau$, $b$, and $\theta$.
        \item $\F'_\tau \subset \Phi_{\F, \tau}(\underline{\F})$ for some $\F' \subset \underline{\F} \subset \F$, and $\F' = \bigsqcup_{\F'_\tau} \underline{\F}[F_\tau]$.
        
        \item For each $F_\tau \in \F'_\tau$, \( C_F(\F'[F_\tau]) \lessapprox_\eps (\tau/\delta)^{2\eta_{J-1}} \).
        \item For each $F_\theta \in \F'_\theta$, \( C_F(\F'_\tau[F_\theta]) \lessapprox_\eps \lambda_*^{-2/\eps^3}(\theta/\tau)^{2\eta_{J-1}} \). 
        
        \item For each $F_{\theta(\tau/\theta)^\eps} \in \F'_{\theta(\tau/\theta)^\eps}$, there exists a family of convex sets $\W_{F_{\theta(\tau/\theta)^\eps}}$ of dimensions $a \times b \times \1$ in $F_{\theta(\tau/\theta)^\eps}$, such that $\W_{F_{\theta(\tau/\theta)^\eps}}$ uniformly factors $\F'_\tau[F_{\theta(\tau/\theta)^\eps}]$. At least one of the following holds.        \begin{enumerate}[label={\rm (\roman*)}]
        \item $a \leq (\tau/\theta)^{\eta_{J-1}'}b$.
        
        \item $a \geq (\tau/\theta)^{\eta_{J-1}'}b$. For all \( \rho \in [\tau(\theta/\tau)^\eps, \theta(\tau/\theta)^\eps] \) and each $F_\rho \in \F'_\rho$, we have \( C_F(\F'_\tau[F_\rho]) \geq (\rho/\tau)^{\eta_J} \). Moreover, $\F'$ is weakly uniform at scales $[\tau, \theta]$ up to a factor $\approx (\theta/\tau)^{\eta_{J-1}'}$. \end{enumerate}
    \end{itemize}
\end{enumerate}

\item
Fix $\beta>0$, and let $\eta_{j-1}' \coloneqq \frac{14\eta_{j-1}}{\eps^2\beta}$ for each $1 \leq j \leq N$. Suppose that $\eta_0 \geq \eta/\eps$, and $\eta_j \geq \eta_{j-1}'/\eps$ for each $1 \leq j \leq N$. If \( \Delta_{max}(\F) \leq \delta^{-\eta} \), then there exists a $\gtrapprox_\eps \lambda_*^{2/\eps^3}$ refinement $\F' \subset \F$ such that at least one of the following holds.
\begin{enumerate}[label={\rm (\Roman*)}]
    \item
    \( \F' \) is \( \lessapprox_\eps \delta^{-{\color{black} 4}\eps} \) Katz--Tao at every scale.
    \item
    There are scales \( \delta \leq \tau \leq \theta \leq 1 \) with \( \tau \leq \delta^{\eps^2}\theta \), an integer \( 1 \leq J \leq N \), and $\tau(\theta/\tau)^\eps \leq a \leq b \leq \theta$, such that the following holds:
    \begin{itemize}
        \item $\F'$ is uniform at scales $\tau$, $b$, and $\theta$.
        
        \item $\F'_\tau \subset \Phi_{\F, \tau}(\underline{\F})$ for some $\F' \subset \underline{\F} \subset \F$, and $\F' = \bigsqcup_{\F'_\tau} \underline{\F}[F_\tau]$.
        
        \item \( \Delta_{max}(\F'_\theta) \lessapprox_\eps \theta^{-2\eta_{J-1}} \).
        
        \item For each $F_\theta \in \F'_\theta$, \( \Delta_{max}(\F_\tau'[F_\theta]) \lessapprox_\eps (\theta/\tau)^{2\eta_{J-1}} \).

        \item For each $F_{\theta} \in \F'_{\theta}$, there exists a family of convex sets $\W_{F_{\theta}}$ of dimensions $a \times b \times \1$ in $F_{\theta}$, such that $\W_{F_{\theta}}$ uniformly factors $\F'_{\tau(\theta/\tau)^\eps}[F_{\theta}]$. At least one of the following holds.        \begin{enumerate}[label={\rm (\roman*)}]
        \item $a \leq (\tau/\theta)^{\eta_{J-1}'}b$.
        
        \item $a \geq (\tau/\theta)^{\eta_{J-1}'}b$. For all \( \rho \in [\tau(\theta/\tau)^\eps, \theta(\tau/\theta)^\eps] \) and each $F_\theta \in \F'_\theta$, we have \( \Delta_{max}(\F'_\rho[F_\theta]) \geq (\theta/\rho)^{\eta_J} \). Moreover, $\F'$ is weakly uniform at scales $[\tau, \theta]$ up to a factor $\approx (\theta/\tau)^{\eta_{J-1}'}$. 
        \end{enumerate}
    \end{itemize}
\end{enumerate}
\end{enumerate}
\end{lem}

\begin{proof}[Proof of Part~{\rm(}A{\rm)}]
To begin, let $\F^{(0)} \coloneqq \Phi_{\F, \delta}(\F)$. 

We apply (a rescaled version of) Lemma~\ref{Lem:factoring'} to each $F_{\delta^\eps} \in (\F^{(0)})_{\delta^\eps}$, which yields a subset $\F^{(0)}[F_{\delta^\eps}]' \subset \F^{(0)}[F_{\delta^\eps}]$ with $|\F^{(0)}[F_{\delta^\eps}]'| \gtrapprox |\F^{(0)}[F_{\delta^\eps}]|$ and a family of convex sets $\W_{F_{\delta^\eps}}$ of dimensions $a_0 \times b_0 \times \1$ ($\delta \leq a_0 \leq b_0 \leq \delta^\eps$) in $F_{\delta^\eps}$, such that $\W_{F_{\delta^\eps}}$ uniformly factors $\F^{(0)}[F_{\delta^\eps}]'$. By pigeonholing, we may assume that $a_0$ and $b_0$ are the same, and the uniformized data are essentially the same up to a factor $\approx 1$ for all $F_{\delta^\eps} \in (\F^{(0)})_{\delta^\eps}$. We may also assume that $\F^{(1)} \coloneqq \cup_{F_{\delta^\eps} \in (\F^{(0)})_{\delta^\eps} } \F^{(0)}[F_{\delta^\eps}]'$ is uniform at scales $[b_0, 1]$. We clearly have $\F^{(1)} \subset \Phi_{\F, \delta}(\F)$ and $\F^{(1)} = \bigsqcup_{\F^{(1)}_\delta} \F[F_\delta]$. Note that $|\F^{(1)}| \gtrapprox |\F^{(0)}| = |\Phi_{\F, \delta}(\F)| \geq \lambda_* |\F|$, and so \( C_F(\F^{(1)}) \lessapprox \lambda_*^{-1} C_F(\F) \leq \lambda_*^{-1} \delta^{-\eta} \leq \lambda_*^{-1/\eps^3} (1/\delta)^{2\eta_{J-1}} \). 
\begin{itemize}
    \item If $a_0 \leq \delta^{\eta_0'} b_0$, then we are in Conclusion~(II)(i) with $\F' = \F^{(1)}$, $\underline{\F} = \F$, $\tau = \delta$, $\theta = 1$, $J=1$, $a = a_0$, and $b=b_0$.
    \item Otherwise, we have $a_0 \geq \delta^{\eta_0'} b_0$, and the key observation is that for each $\rho \in [\delta, b_0]$ and $F_\rho \in (\F^{(0)}[F_{\delta^\eps}]')_\rho$, the convex sets in $\W_{F_{\delta^\eps}}$ containing a $\delta$-flake in $F_\rho$ are all contained in some $Cb_0$-flake, and so has cardinality $\in [1, Cb_0/a_0]$ by using $\Delta_{max}(\W_{F_{\delta^\eps}}) \lesssim 1$. This means that $\F^{(1)}$ is weakly uniform at scales $[\delta,b_0]$ (and thus $[\delta,1]$) up to a factor $\approx \delta^{-\eta_0'}$. And now we can check whether $C_F(\F^{(1)}[F_\rho]) \leq (\rho/\delta)^{\eta_1}$ holds for some $\rho \in [\delta^{1-\eps}, \delta^\eps]$ and $F_\rho \in (\F^{(1)})_\rho$.
    \begin{itemize}
        \item If this is not true, then we have $C_F(\F^{(1)}[F_\rho]) \geq (\rho/\delta)^{\eta_1}$ for all $\rho \in [\delta^{1-\eps}, \delta^\eps]$ and each $F_\rho \in (\F^{(1)})_\rho$. Thus we are in Conclusion~(II)(ii) with $\F' = \F^{(1)}$, $\underline{\F} = \F$, $\tau = \delta$, $\theta = 1$, $J=1$, $a = a_0$, and $b=b_0$.
        
        \item If this is true, then we pick one such $\rho = \rho_1$. By the weak uniformity, we have 
        $$\qquad C_F(\F^{(1)}[F_{\rho_1}]) \lessapprox \delta^{-\eta_0'}(\rho_1/\delta)^{\eta_1} \leq (\rho_1/\delta)^{\eta_1 + \eta_0'/\eps} \leq (\rho_1/\delta)^{2\eta_1}$$ 
        for each $F_{\rho_1} \in (\F^{(1)})_{\rho_1}$. This is the only unsolved case, and we continue with the remainder of the proof.
    \end{itemize}
\end{itemize} 

In general, we run the following iterative algorithm. Let $j$ be the loop index, with $1$ as its initial value. Suppose we have $|\F^{(j)}| \gtrapprox_\eps \lambda_*|\F|$, and $C_F(\F^{(j)}[F_{\rho_j}]) \lessapprox (\rho_j/\delta)^{2\eta_j}$ for each $F_{\rho_j} \in (\F^{(j)})_{\rho_j}$.

Applying (a rescaled version of) Lemma~\ref{Lem:factoring'} to each \( \F^{(j)}[F_{\rho_j(\delta/\rho_j)^\eps}] \) where $F_{\rho_j(\delta/\rho_j)^\eps} \in (\F^{(j)})_{\rho_j(\delta/\rho_j)^\eps}$, we get a subset $\F^{(j)}[F_{\rho_j(\delta/\rho_j)^\eps}]' \subset \F^{(j)}[F_{\rho_j(\delta/\rho_j)^\eps}]$ with $|\F^{(j)}[F_{\rho_j(\delta/\rho_j)^\eps}]'| \gtrapprox |\F^{(j)}[F_{\rho_j(\delta/\rho_j)^\eps}]|$, and a family of convex sets $\W_{F_{\rho_j(\delta/\rho_j)^\eps}}$ of dimensions $a_j \times b_j \times \1$ ($\delta \leq a_j \leq b_j \leq \rho_j(\delta/\rho_j)^\eps$) in $F_{\rho_j(\delta/\rho_j)^\eps}$, such that $\W_{F_{\rho_j(\delta/\rho_j)^\eps}}$ uniformly factors $\F^{(j)}[F_{\rho_j(\delta/\rho_j)^\eps}]'$. By pigeonholing, we may assume that $a_j$ and $b_j$ are the same, and the uniformized data are essentially the same up to a factor $\approx 1$ for all $F_{\rho_j(\delta/\rho_j)^\eps} \in (\F^{(j)})_{\rho_j(\delta/\rho_j)^\eps}$. We may also assume that $$\F^{(j+1)} \coloneqq \bigcup_{F_{\rho_j(\delta/\rho_j)^\eps} \in (\F^{(j)})_{\rho_j(\delta/\rho_j)^\eps}} \F^{(j)}[F_{\rho_j(\delta/\rho_j)^\eps}]'$$ is uniform at scales $[b_j, 1]$. Note that $|\F^{(j+1)}| \gtrapprox |\F^{(j)}| \gtrapprox_\eps \lambda_* |\F|$. So we may only retain those $F_{\rho_{j}} \in (\F^{(j+1)})_{\rho_{j}}$ with $|\F^{(j+1)}[F_{\rho_{j}}]| \gtrapprox |\F^{(j)}[F_{\rho_{j}}]|$, which gives $$C_F(\F^{(j+1)}[F_{\rho_j}]) \lessapprox (\rho_j/\delta)^{2\eta_j}$$
for each $F_{\rho_{j}} \in (\F^{(j+1)})_{\rho_{j}}$.
\begin{itemize}
    \item If $a_j \leq (\delta/\rho_j)^{\eta_j'} b_j$, then we stop. 
    \item Otherwise, we have $a_j \geq (\delta/\rho_j)^{\eta_j'} b_j$, and the key observation is that for each $\rho \in [\delta, b_j]$ and $F_\rho \in (\F[F_{\rho_j(\delta/\rho_j)^\eps}]')_\rho$, the convex sets in $\W_{F_{\rho_j(\delta/\rho_j)^\eps}}$ containing a $\delta$-flake in $F_\rho$ are all contained in some $Cb_j$-flake, and so has cardinality $\in [1, Cb_j/a_j]$ by using $\Delta_{max}(\W_{F_{\rho_j(\delta/\rho_j)^\eps}}) \lesssim 1$. This means that $\F^{(j+1)}$ is weakly uniform at scales $[\delta,b_j]$ (and thus $[\delta,\rho_j]$) up to a factor $\approx (\rho_j/\delta)^{\eta_j'}$. And now we can check whether $C_F(\F^{(j+1)}[F_\rho]) \leq (\rho/\delta)^{\eta_{j+1}}$ holds for some $\rho \in [\delta(\rho_j/\delta)^\eps, \rho_j(\delta/\rho_j)^\eps]$ and $F_\rho \in (\F^{(j+1)})_\rho$.
    \begin{itemize}
        \item If this is not true, then we have $C_F(\F^{(j+1)}[F_\rho]) \geq (\rho/\delta)^{\eta_{j+1}}$ for all $\rho \in [\delta(\rho_j/\delta)^\eps, \rho_j(\delta/\rho_j)^\eps]$ and each $F_\rho \in (\F^{(j+1)})_\rho$. And we stop. 
        
        \item If this is true, then we pick one such $\rho = \rho_{j+1}$. By the weak uniformity, we have $$\qquad C_F(\F^{(j+1)}[F_{\rho_{j+1}}]) \lessapprox (\rho_j/\delta)^{\eta_j'} (\rho_{j+1}/\delta)^{\eta_{j+1}} \leq (\rho_{j+1}/\delta)^{\eta_{j+1} + \eta_j'/\eps} \leq (\rho_{j+1}/\delta)^{2\eta_{j+1}}$$ 
        for each $F_{\rho_{j+1}} \in (\F^{(j+1)})_{\rho_{j+1}}$ (as $\rho_j/\delta \leq (\rho_{j+1}/\delta)^{1/\eps}$). If $\rho_{j+1} < \delta^{1-\eps}$, then we stop. Otherwise, we go back to the start of the loop and proceed with $j+1$ in place of $j$.
    \end{itemize}
\end{itemize} 

As long as $\rho_j \geq \delta^{1-\eps}$, we have $\rho_{j+1} \leq \rho_j(\delta/\rho_j)^\eps \leq \rho_j \delta^{\eps^2}$ for each $j$, and so $\rho_{j+1} \leq \rho_1 \delta^{j\eps^2} \leq \delta^{j\eps^2}$ by induction. As $\rho_{j+1} \geq \delta$, we must have $j \leq \lfloor 1/\eps^2 \rfloor$. In other words, the algorithm must stop in at most $\lfloor 1/\eps^2 \rfloor$ steps. Let $j = j_1$ be the loop index at which we stop. Then $\rho_{j_1} \geq \delta^{1-\eps}$, and the output of the above algorithm is $\F_1 \coloneqq \F^{(j_1 + 1)} \subset \F$, such that $|\F_1| \gtrapprox_\eps |\F^{(1)}| \gtrapprox \lambda_*|\F|$ and $\F_1$ is uniform at scale $b_{j_1}$. We clearly have $\F_1 \subset \Phi_{\F, \delta}(\F)$ and $\F_1 = \bigsqcup_{(\F_1)_\delta} \F[F_\delta]$. Note that $\rho_{j_1 + 1} \geq \delta(\rho_{j_1}/\delta)^\eps \geq \delta^{1-\eps^2}$, and \( C_F(\F_1) \lessapprox_\eps C_F(\F^{(1)}) \lessapprox \lambda_*^{-1/\eps^3} (1/\delta)^{2\eta_{J-1}} \). 

\begin{itemize}
    \item Case $1$: We have $\rho_{j_1 + 1} < \delta^{1-\eps}$, and $$C_F(\F_1[F_{\rho_{j_1 + 1}}]) \lessapprox (\rho_{j_1 + 1}/\delta)^{2\eta_{j_1 + 1}}$$ 
    for each $F_{\rho_{j_1 + 1}} \in (\F_1)_{\rho_{j_1 + 1}}$. 
    
    \item Case $2$: Otherwise, recall that we have $$C_F(\F_1[F_{\rho_{j_1}}]) \lessapprox (\rho_{j_1}/\delta)^{2\eta_{j_1}}$$ 
    for each $F_{\rho_{j_1}} \in (\F_1)_{\rho_{j_1}}$. For each $F_{\rho_{j_1}(\delta/\rho_{j_1})^\eps} \in (\F_1)_{\rho_{j_1}(\delta/\rho_{j_1})^\eps}$, there exists a family of convex sets $\W_{F_{\rho_{j_1}(\delta/\rho_{j_1})^\eps}}$ of dimensions $a_{j_1} \times b_{j_1} \times \1$ in $F_{\rho_{j_1}(\delta/\rho_{j_1})^\eps}$, such that $\W_{\rho_{j_1}(\delta/\rho_{j_1})^\eps}$ uniformly factors $\F_1[F_{\rho_{j_1}(\delta/\rho_{j_1})^\eps}]$. At least one of the following holds.
    \begin{enumerate}[label={\rm (\roman*)}]
        \item $a_{j_1} \leq (\delta/\rho_{j_1})^{\eta_{j_1}'}b_{j_1}$.
        
        \item $a_{j_1} \geq (\delta/\rho_{j_1})^{\eta_{j_1}'} b_{j_1}$. For all $\rho \in [\delta(\rho_{j_1}/\delta)^\eps, \rho_{j_1}(\delta/\rho_{j_1})^\eps]$ and each $F_\rho \in (\F_1)_\rho$, we have \( C_F(\F_1[F_\rho]) \geq (\rho/\delta)^{\eta_{j_1 + 1}} \). Moreover, $\F_1$ is weakly uniform at scales $[\delta, \rho_{j_1}]$ up to a factor $\approx (\rho_{j_1}/\delta)^{\eta_{j_1}'}$.
    \end{enumerate}
\end{itemize}
In Case $1$, we set $j_1' \coloneqq j_1 + 1$; while in Case $2$, we set $j_1' \coloneqq j_1$.\footnote{If both cases holds, then we just choose one of them to define the properties at scale $\rho_{j_1'}$.} Then in both cases, we have $C_F(\F_1[F_{\rho_{j_1'}}]) \lessapprox (\rho_{j_1'}/\delta)^{2\eta_{j_1'}}$ for each $F_{\rho_{j_1'}} \in (\F_1)_{\rho_{j_1'}}$. By pigeonholing, we may further assume that $\F_1$ is uniform at scale $\rho_{j_1'}$. This does not affect the properties in either case.

In Case $2$, we are actually done, as Conclusion~(II) holds with 
$$\F' = \F_1, \, \underline{\F} = \F, \, \tau = \delta, \, \theta = \rho_{j_1}, \, J=j_1+1, \, a = a_{j_1}, \, b = b_{j_1}.$$

Therefore, it remains to consider Case $1$. Recall that $\delta^{1-\eps^2} \leq \rho_{j_1'} < \delta^{1-\eps}$. The key point is that we can iterate the process described so far until we arrive at a good place. More precisely, we run the following iterative algorithm. Let $k \leq \lfloor 1/\eps^3 \rfloor$ be the loop index, with $1$ as its initial value. Let $\rho_{j_0'} \coloneqq \delta$.

Suppose we have got $\F_k \subset \F$ that satisfies $|\F_k| \gtrapprox_\eps \lambda_*^k|\F|$, and is uniform at scales $\{\rho_{j_i'}\}_{i=1}^k$ with $\rho_{j_i'} \geq (\rho_{j_{i-1}'})^{1-\eps^2}$ for each $1 \leq i \leq k$. Suppose also that 
$\rho_{j_{k}'} < \delta^\eps$, and 
\begin{align}\label{Eq:C_F_upper_bound_1}
C_F(\F_k[F_{\rho_{j_k'}}]) \lessapprox_\eps (\rho_{j_k'}/\delta)^{2\eta_{j_k'}}
\end{align} 
for each $F_{\rho_{j_k'}} \in (\F_k)_{\rho_{j_k'}}$. Define 
$$\F_k^* \coloneqq \bigsqcup_{F_{\rho_{j_k'}} \in \Phi_{\F, \rho_{j_k'}}(\F_k)} \F_k[F_{\rho_{j_k'}}].$$
Then $\F_k^* \subset \F_k$ is still uniform at scales $\{\rho_{j_i'}\}_{i=1}^k$, and $|\F_k^*| \geq \lambda_* |\F_k| \gtrapprox_\eps \lambda_*^{k+1}|\F|$ by our assumption on $\Phi_{\F, \rho}$. 

Repeat all the previous arguments with $\rho_{j_k'}$ in place of $\delta$, $(\F_k^*)_{\rho_{j_k'}}$ in place of $\F^{(0)}$, and $\eta_j$ starting at $\eta_{j_k'+1}$. By regarding $\rho_{j_k'}$ as the smallest scale, whenever we take a $\gtrapprox_\eps 1$ refinement of $(\F_k^*)_{\rho_{j_k'}}$, all the previous properties are preserved for every surviving $F_{\rho_{j_k'}}$. This will induce a $\gtrapprox_\eps 1$ refinement of $\F_k^*$ due to the uniformity of $\F_k^*$ at scale $\rho_{j_k'}$.

To begin, similarly to the initial $\F^{(1)} \subset \F^{(0)}$, we get some $\F_k^{*,(1)} \subset \F_k^*$ with $|\F_k^{*,(1)}| \gtrapprox |\F_k^*| \gtrapprox_\eps \lambda_*^{k+1}|\F|$, and so 
$$C_F(\F_k^{*,(1)}) \lessapprox_\eps \lambda_*^{-(k+1)} C_F(\F) \leq \lambda_*^{-2/\eps^3} \delta^{-\eta} \leq \lambda_*^{-2/\eps^3} \rho_{j_k'}^{-\eta/\eps} \leq \lambda_*^{-2/\eps^3} (1/\rho_{j_k'})^{2\eta_0}.$$
This accounts for the $\lambda_*^{-2/\eps^3}$ loss for ``$C_F(\F_\tau'[F_\theta])$''. We omit the other details, as they are essentially the same.

In general, we run an iterative algorithm starting from $\F_k^{*,(1)}$. The outcome can be summarized as: There exist $\F_{k+1} \subset \F_k^*$ with $|\F_{k+1}| \gtrapprox_\eps |\F_k^*| \gtrapprox_\eps \lambda_*^{k+1}|\F|$, some $j_k' < j_{k+1}' \leq j_k' + 2 + \lfloor 1/\eps^2 \rfloor$, and $\rho_{j_{k+1}'} \geq (\rho_{j_k'})^{1-\eps^2}$, such that 
\begin{align}\label{Eq:C_F_upper_bound_2}    C_F((\F_{k+1})_{\rho_{j_k'}}[F_{\rho_{j_{k+1}'}}]) \lessapprox (\rho_{j_{k+1}'}/\rho_{j_k'})^{2\eta_{j_{k+1}'}}
\end{align} 
for each $F_{\rho_{j_{k+1}'}} \in (\F_{k+1})_{\rho_{j_{k+1}'}}$. Moreover, $\F_{k+1}$ is uniform at scale $b_{j_{k+1}'}$, $(\F_{k+1})_{\rho_{j_k'}} \subset (\F_k^*)_{\rho_{j_k'}} \subset \Phi_{\F, _{\rho_{j_k'}}}(\F_k)$, $\F_{k+1} = \bigsqcup_{ (\F_{k+1})_{\rho_{j_k'}} } \F_k[F_{\rho_{j_k'}}]$, and one of the following two possibilities holds:
\begin{itemize}
    \item Case $1$: $\rho_{j_{k+1}'} < \rho_{j_k'}^{1-\eps}$.  
    
    \item Case $2$: $\rho_{j_{k+1}'} \geq \rho_{j_k'}^{1-\eps}$. For each $F_{\rho_{j_{k+1}'}(\rho_{j_{k}'}/\rho_{j_{k+1}'})^\eps} \in (\F_{k+1})_{\rho_{j_{k+1}'}(\rho_{j_{k}'}/\rho_{j_{k+1}'})^\eps}$, there exists a family of convex sets $\W_{F_{\rho_{j_{k+1}'}(\rho_{j_{k}'}/\rho_{j_{k+1}'})^\eps}}$ of dimensions $a_{j_{k+1}'} \times b_{j_{k+1}'} \times \1$ in $F_{\rho_{j_{k+1}'}(\rho_{j_{k}'}/\rho_{j_{k+1}'})^\eps}$, that uniformly factors $(\F_{k+1})_{\rho_{j_k'}}[F_{\rho_{j_{k+1}}(\rho_{j_{k}}/\rho_{j_{k+1}})^\eps}]$. At least one of the following holds.
    \begin{enumerate}[label={\rm (\roman*)}]
        \item $a_{j_{k+1}'} \leq (\rho_{j_{k}'}/\rho_{j_{k+1}'})^{\eta_{j_{k+1}'}'}b_{j_{k+1}'}$.
        
        \item $a_{j_{k+1}'} \geq (\rho_{j_{k}'}/\rho_{j_{k+1}'})^{\eta_{j_{k+1}'}'}b_{j_{k+1}'}$. For all $\rho \in [\rho_{j_k'}(\rho_{j_{k+1}'}/\rho_{j_k'})^\eps, \rho_{j_{k+1}'}(\rho_{j_k'}/\rho_{j_{k+1}'})^\eps]$ and each $F_\rho \in (\F_{k+1})_\rho$, we have \( C_F((\F_{k+1})_{\rho_{j_k'}}[F_\rho]) \geq (\rho/\rho_{j_k'})^{\eta_{j_{k+1}' + 1}} \). Moreover, $\F_{k+1}$ is weakly uniform at scales $[\rho_{j_{k}'}, \rho_{j_{k+1}'}]$ up to a factor $\approx (\rho_{j_{k+1}'}/\rho_{j_{k}'})^{\eta_{j_{k+1}'}'}$.
    \end{enumerate}
\end{itemize}
By pigeonholing, we may further assume that $\F_{k+1}$ is uniform at scales $\{\rho_{j_i'}\}_{i=1}^{k+1}$. This does not affect the properties in either case. By construction of $\F_{k+1}$, \eqref{Eq:C_F_upper_bound_1} implies 
\begin{align}\label{Eq:C_F_upper_bound_1'}
C_F(\F_{k+1}[F_{\rho_{j_k'}}]) \lessapprox_\eps (\rho_{j_k'}/\delta)^{2\eta_{j_k'}} \leq (\rho_{j_k'}/\delta)^{2\eta_{j_{k+1}'}}
\end{align}
for each $F_{\rho_{j_k'}} \in (\F_{k+1})_{\rho_{j_k'}}$.

In Case $2$, we are actually done, as Conclusion~(II) holds with 
$$\F' = \F_{k+1}, \, \underline{\F} = \F_k, \, \tau = \rho_{j_{k}'}, \, \theta = \rho_{j_{k+1}'}, \, J=j_{k+1}'+1, \, a = a_{j_{k+1}'}, \, b = b_{j_{k+1}'}.$$ 
To see this, note that $\rho_{j_k'} \leq \delta^\eps$ and $\rho_{j_{k+1}'} \geq \rho_{j_k'}^{1-\eps}$ imply $\rho_{j_{k+1}'} \geq \rho_{j_k'}\delta^{-\eps^2}$.

Therefore, without loss of generality, we may always assume that Case $1$ holds. By using (a rescaled version of) Lemma~\ref{Lem:KT_multi} (Remark~\ref{Rem:KT_multi_relax}), we have 
\begin{align*}
C_F(\F_{k+1}[F_{\rho_{j_{k+1}'}}]) &\lessapprox 
C_F((\F_{k+1})_{\rho_{j_k'}}[F_{\rho_{j_{k+1}'}}]) C_F(\F_{k+1}[F_{\rho_{j_k'}}]) 
\\
&\lessapprox_\eps 
(\rho_{j_{k+1}'}/\rho_{j_k'})^{2\eta_{j_{k+1}'}} 
(\rho_{j_k'}/\delta)^{2\eta_{j_k'}} 
\\
&\leq 
(\rho_{j_{k+1}'}/\delta)^{2\eta_{j_{k+1}'}}
\end{align*}
for each $F_{\rho_{j_{k+1}'}} \in (\F_{k+1})_{\rho_{j_{k+1}'}}$. In the second line above, we used \eqref{Eq:C_F_upper_bound_2} and \eqref{Eq:C_F_upper_bound_1'}.

If $\rho_{j_{k+1}'} \geq \delta^\eps$, then we stop. Otherwise, we go back to the start of the loop and proceed with $k+1$ in place of $k$.

As long as $\rho_{j_k'} \leq \delta^\eps$, we have $\rho_{j_{k+1}'} \geq (\rho_{j_k'})^{1-\eps^2} \geq \rho_{j_k'}\delta^{-\eps^3}$ for each $k$, and so $\rho_{j_{k+1}'} \geq \rho_{j_{1}'}\delta^{-k\eps^3}$ by induction. As $\rho_{j_{k+1}'} \leq 1$ and $\rho_{j_{1}'} \geq \delta$, we must have $k \leq \lfloor 1/\eps^3 \rfloor$, which matches the initial constraint on $k$. In other words, the algorithm must stop in at most $\lfloor 1/\eps^3 \rfloor$ steps. Let $k_0$ be the loop index at which we stop, then $\rho_{j_{k_0 + 1}'} \geq \delta^{\eps}$, and the output of the above algorithm is $\F' \coloneqq \F_{k_0 + 1} \subset \F$ with $|\F'| \gtrapprox_\eps \lambda_*^{k_0 + 1}|\F| \geq \lambda_*^{2/\eps^3}|\F|$. Note that $j_{k_0 + 1}' \leq (k_0 + 1) (2 + \lfloor 1/\eps^2 \rfloor) \leq 6/\eps^5 = 6N/10 \leq N-1$.

We are in the following case: For each $0 \leq k \leq k_0$\footnote{The $k=0$ case follows from the discussions about $\rho_{j_1'}$ before the iteration on $k$.}, $\F'$ satisfies Case~$1$ between scales $\rho_{j_k'}$ and $\rho_{j_{k+1}'}$ (with $\F_{k+1}$ replaced with $\F'$). 
For notational consistency, let $j_{k_0 +2}' \coloneqq j_{k_0 +1}' + 1 \leq N$ and $\rho_{j_{k_0 +2}'} \coloneqq 1$. For each \(1 \leq k \leq k_0 + 2\) and $F_{\rho_{j_k'}} \in \F'_{\rho_{j_k'}}$, we have 
\begin{align*}
C_F(\F'[F_{\rho_{j_k'}}]) \lessapprox_\eps \lambda_*^{-2/\eps^3} (\rho_{j_k'}/\delta)^{2\eta_{j_{k}'}},
\end{align*}
where the $\lambda_*^{-2/\eps^3}$ loss comes from the case when $k = k_0 + 2$.

Recall that $\rho_{j_0'} = \delta$, $\rho_{j_1'} < \delta^{1-\eps}$, $\rho_{j_{k_0 + 1}'} \geq \delta^{\eps}$, $\rho_{j_{k_0 +2}'} = 1$, and we are always in Case~$1$. Thus for each \(0 \leq k \leq k_0 + 1\), the ratio \( \rho_{j_{k+1}'} / \rho_{j_k'} \) is at most \( \rho_{j_k'}^{-\eps} \leq \delta^{-\eps} \). We claim that \( \F' \) is \( \lessapprox \lambda_*^{-2/\eps^3} \delta^{-{\color{black} 2n}\eps} \) Frostman at every scale, and hence Conclusion (I) holds. To verify this, let \( \rho \in [\delta, 1] \) and choose \( 0 \leq k \leq k_0 + 1\) so that \( \rho_{j_k'} \leq \rho \leq \rho_{j_{k+1}'} \). Then for each $F_\rho \in \F'_\rho$, we have 
\begin{align*}
C_F(\F'[F_\rho]) 
&\lesssim C_F(\F'[F_{\rho_{j_{k+1}'}}]) \frac{|\F'[F_{\rho_{j_{k+1}'}}]|}{|\F'[F_\rho]|}\\
&\lesssim C_F(\F'[F_{\rho_{j_{k+1}'}}]) \frac{|\F'[F_{\rho_{j_{k+1}'}}]|}{|\F'[F_{\rho_{j_{k}'}}]|}\\
&\lessapprox_\eps \lambda_*^{-2/\eps^3} (\rho_{j_{k+1}'}/\delta)^{2\eta_{j_{k+1}'}} \left( \frac{\rho_{j_{k+1}'}}{\rho_{j_k'}} \right)^{\color{black} 2(n-1)} 
\leq 
\lambda_*^{-2/\eps^3} \delta^{-{\color{black} 2n}\eps}.
\end{align*}
So the proof is completed.
\end{proof}

\begin{proof}[Proof of Part~{\rm(}B{\rm)}]
The strategy is similar to that for Part~(A). To begin, by pigeonholing, there exists $\F^{(0)} \subset \Phi_{\F, \delta}(\F)$ with $|\F^{(0)}| \gtrapprox |\Phi_{\F, \delta}(\F)|$ such that $\F^{(0)}$ is uniform at scales $[\delta, \delta^{1-\eps}]$. 

We apply Lemma~\ref{Lem:factoring'} to $(\F^{(0)})_{\delta^{1-\eps}}$, which yields a subset $(\F^{(0)})_{\delta^{1-\eps}}' \subset (\F^{(0)})_{\delta^{1-\eps}}$ with $|(\F^{(0)})_{\delta^{1-\eps}}'| \gtrapprox |(\F^{(0)})_{\delta^{1-\eps}}|$ and a family of convex sets $\W_{F_1}$ of dimensions $a_0 \times b_0 \times \1$ ($\delta^{1-\eps} \leq a_0 \leq b_0 \leq 1$) in $B_1$, such that $\W_{F_1}$ uniformly factors $(\F^{(0)})_{\delta^{1-\eps}}'$\footnote{The uniformity is originally for $\delta^{1-\eps}$-flakes at scales $[\delta^{1-\eps},1]$ by Definition~\ref{Def:factor'}, but can actually be made for $\delta$-flakes at scales $[\delta^{1-\eps},1]$ by modifying the uniformization procedure and using the existing uniformity at scales $[\delta, \delta^{1-\eps}]$. Such a convention is also adopted later on in the proof.}. By pigeonholing, we may assume that $\F^{(1)} \coloneqq \cup_{F_{\delta^{1-\eps}} \in (\F^{(0)})_{\delta^{1-\eps}}'} \F^{(0)}[F_{\delta^{1-\eps}}]$ is uniform at scales $[b_0, 1]$. We clearly have $\F^{(1)} \subset \Phi_{\F, \delta}(\F)$, $\F^{(1)} = \bigsqcup_{\F^{(1)}_\delta} \F[F_\delta]$, and $\F^{(1)}$ is still uniform at scales $[\delta, \delta^{1-\eps}]$. Note that $|\F^{(1)}| \gtrapprox |\F^{(0)}| \gtrapprox |\Phi_{\F, \delta}(\F)| \geq \lambda_* |\F|$, and \( \Delta_{max}(\F^{(1)}) \leq \Delta_{max}(\F) \leq \delta^{-\eta} \leq (1/\delta)^{2\eta_0} \). 
\begin{itemize}
    \item If $a_0 \leq \delta^{\eta_0'} b_0$, then we are in Conclusion~(II)(i) with $\F' = \F^{(1)}$, $\underline{\F} = \F$, $\tau = \delta$, $\theta = 1$, $J=1$, $a = a_0$, and $b=b_0$.
    \item Otherwise, we have $a_0 \geq \delta^{\eta_0'} b_0$, and the key observation is that for each $\rho \in [\delta^{1-\eps}, b_0]$ and $F_\rho \in (\F^{(1)})_\rho$, the convex sets in $\W_{F_1}$ containing a $\delta$-flake in $F_\rho$ are all contained in some $Cb_0$-flake, and so has cardinality $\in [1, Cb_0/a_0]$ by using $\Delta_{max}(\W_{F_1}) \lesssim 1$. This means that $\F^{(1)}$ is weakly uniform at scales $[\delta^{1-\eps},b_0]$ (and thus $[\delta,1]$) up to a factor $\approx \delta^{-\eta_0'}$. And now we can check whether $\Delta_{max}((\F^{(1)})_\rho) \leq \rho^{-\eta_1}$ holds for some $\rho \in [\delta^{1-\eps}, \delta^\eps]$.
    \begin{itemize}
        \item If this is not true, then we have $\Delta_{max}((\F^{(1)})_\rho) \geq \rho^{-\eta_1}$ for all $\rho \in [\delta^{1-\eps}, \delta^\eps]$. Thus we are in Conclusion~(II)(ii) with $\F' = \F^{(1)}$, $\underline{\F} = \F$, $\tau = \delta$, $\theta = 1$, $J=1$, $a = a_0$, and $b=b_0$.
        
        \item If this is true, then we pick one such $\rho = \rho_1$. By the weak uniformity, we have 
        $$\qquad \Delta_{max}((\F^{(1)})_{\rho_1}) \lessapprox \delta^{-\eta_0'} \rho_1^{-\eta_1} \leq \rho_1^{-(\eta_1 + \eta_0'/\eps)} \leq \rho_1^{-2\eta_1}.$$ 
        This is the only unsolved case, and we continue with the remainder of the proof.
    \end{itemize}
\end{itemize} 

In general, we run the following iterative algorithm. Let $j$ be the loop index, with $1$ as its initial value. Suppose $\rho_j \geq \delta^{1-\eps}$, $|\F^{(j)}| \gtrapprox_\eps \lambda_*|\F|$, and $\F^{(j)}$ is uniform at scales $[\delta, \delta(\rho_j/\delta)^\eps]$. Suppose also that $\Delta_{max}((\F^{(j)})_{\rho_{j}}) \lessapprox_\eps \rho_j^{-2\eta_j}$.

Applying (a rescaled version of) Lemma~\ref{Lem:factoring'} to each \( (\F^{(j)})_{\delta(\rho_j/\delta)^\eps}[F_{\rho_j}] \) where $F_{\rho_j} \in (\F^{(j)})_{\rho_j}$, we obtain $(\F^{(j)})_{\delta(\rho_j/\delta)^\eps}[F_{\rho_j}]' \subset (\F^{(j)})_{\delta(\rho_j/\delta)^\eps}[F_{\rho_j}]$ with $|(\F^{(j)})_{\delta(\rho_j/\delta)^\eps}[F_{\rho_j}]'| \gtrapprox |(\F^{(j)})_{\delta(\rho_j/\delta)^\eps}[F_{\rho_j}]|$ 
and a family of convex sets $\W_{F_{\rho_j}}$ of dimensions $a_j \times b_j \times \1$ ($\delta(\rho_j/\delta)^\eps \leq a_j \leq b_j \leq \rho_j$) in $F_{\rho_j}$, such that $\W_{F_{\rho_j}}$ uniformly factors $(\F^{(j)})_{\delta(\rho_j/\delta)^\eps}[F_{\rho_j}]'$. By pigeonholing, we may assume that $a_j$ and $b_j$ are the same, and the uniformized data are essentially the same up to a factor $\approx 1$ for all $F_{\rho_j} \in (\F^{(j)})_{\rho_j}$. We may also assume that 
$$\F^{(j+1)} \coloneqq \cup_{F_{\rho_j} \in (\F^{(j)})_{\rho_j}} \F^{(j)}[(\F^{(j)})_{\delta(\rho_j/\delta)^\eps}[F_{\rho_j}]']$$ 
is uniform at scales $[b_j, \rho_j]$. Note that $|\F^{(j+1)}| \gtrapprox |\F^{(j)}| \gtrapprox_\eps \lambda_* |\F|$, and $\F^{(j+1)}$ is still uniform at scales $[\delta, \delta(\rho_j/\delta)^\eps]$. 
\begin{itemize}
    \item If $a_j \leq (\delta/\rho_j)^{\eta_j'} b_j$, then we stop. 
    \item Otherwise, we have $a_j \geq (\delta/\rho_j)^{\eta_j'} b_j$, and the key observation is that for each $\rho \in [\delta(\rho_j/\delta)^\eps, b_j]$ and $F_\rho \in (\F^{(j)}[(\F^{(j)})_{\delta(\rho_j/\delta)^\eps}[F_{\rho_j}]'])_\rho$, the convex sets in $\W_{F_{\rho_j}}$ containing a $\delta$-flake in $F_\rho$ are all contained in some $Cb_j$-flake, and so has cardinality $\in [1, Cb_j/a_j]$ by using $\Delta_{max}(\W_{F_{\rho_j}}) \lesssim 1$. This means that $\F^{(j+1)}$ is weakly uniform at scales $[\delta(\rho_j/\delta)^\eps,b_j]$ (and thus $[\delta,\rho_j]$) up to a factor $\approx (\rho_j/\delta)^{\eta_j'}$. And now we can check whether $\Delta_{max}((\F^{(j+1)})_{\rho}[F_{\rho_{j}}]) \leq (\rho_j/\rho)^{\eta_{j+1}}$ holds for some $\rho \in [\delta(\rho_j/\delta)^\eps, \rho_j(\delta/\rho_j)^\eps]$ and $F_{\rho_j} \in (\F^{(j+1)})_{\rho_j}$.
    \begin{itemize}
        \item If this is not true, then we have $\Delta_{max}((\F^{(j+1)})_{\rho}[F_{\rho_{j}}]) \geq (\rho_j/\rho)^{\eta_{j+1}}$ for all $\rho \in [\delta(\rho_j/\delta)^\eps, \rho_j(\delta/\rho_j)^\eps]$ and each $F_{\rho_j} \in (\F^{(j+1)})_{\rho_j}$. And we stop. 
        
        \item If this is true, then we pick one such $\rho = \rho_{j+1}$. By the weak uniformity, we have 
        $$\Delta_{max}((\F^{(j+1)})_{\rho_{j+1}}[F_{\rho_{j}}]) \lessapprox (\rho_j/\delta)^{\eta_j'} (\rho_{j}/\rho_{j+1})^{\eta_{j+1}} \leq (\rho_{j}/\rho_{j+1})^{\eta_{j+1} + \eta_j'/\eps} \leq (\rho_{j}/\rho_{j+1})^{2\eta_{j+1}}$$ 
        for each $F_{\rho_{j}} \in (\F^{(j+1)})_{\rho_{j}}$ (as $\rho_j/\delta \leq (\rho_{j}/\rho_{j+1})^{1/\eps}$). Since $\F^{(j+1)} \subset \F^{(j)}$, by Proposition~\ref{Prop:inherit}(B), we still have $\Delta_{max}((\F^{(j+1)})_{\rho_{j}}) \lessapprox_\eps \rho_j^{-2\eta_j}$. Therefore, by applying Lemma~\ref{Lem:KT_multi} (Remark~\ref{Rem:KT_multi_relax}) to $(\F^{(j+1)})_{\rho_{j+1}}$ with $\delta = \rho_{j+1}$, we have 
        \begin{align*}
        \Delta_{max}((\F^{(j+1)})_{\rho_{j+1}}) &\lessapprox
        \Delta_{max}((\F^{(j+1)})_{\rho_{j}}) \cdot 
        \max_{F_{\rho_j} \in (\F^{(j+1)})_{\rho_{j}}} \Delta_{max}((\F^{(j+1)})_{\rho_{j+1}}[F_{\rho_j}]) 
        \\        &\lessapprox_\eps \rho_j^{-2\eta_j} (\rho_{j}/\rho_{j+1})^{2\eta_{j+1}}
        \\
        &\leq         \rho_{j+1}^{-2\eta_{j+1}}.
        \end{align*}
        
        If $\rho_{j+1} < \delta^{1-\eps}$, then we stop. Otherwise, we go back to the start of the loop and proceed with $j+1$ in place of $j$.
    \end{itemize}
\end{itemize} 

As long as $\rho_j \geq \delta^{1-\eps}$, we have $\rho_{j+1} \leq \rho_j(\delta/\rho_j)^\eps \leq \rho_j \delta^{\eps^2}$ for each $j$, and so $\rho_{j+1} \leq \rho_1 \delta^{j\eps^2} \leq \delta^{j\eps^2}$ by induction. As $\rho_{j+1} \geq \delta$, we must have $j \leq \lfloor 1/\eps^2 \rfloor$. In other words, the algorithm must stop in at most $\lfloor 1/\eps^2 \rfloor$ steps. Let $j = j_1$ be the loop index at which we stop. Then $\rho_{j_1} \geq \delta^{1-\eps}$, and the output of the above algorithm is $\F_1 \coloneqq \F^{(j_1 + 1)} \subset \F$, such that $|\F_1| \gtrapprox_\eps \lambda_*|\F|$ and $\F_1$ is uniform at scale $b_{j_1}$. We clearly have $\F_1 \subset \Phi_{\F, \delta}(\F)$, $\F_1 = \bigsqcup_{(\F_1)_\delta} \F[F_\delta]$, and $\Delta_{max}((\F_1)_\delta
[F_{\rho_{j_1}}]) \leq \Delta_{max}(\F) \leq \delta^{-\eta} \leq (\rho_{j_1}/\delta)^{\eta/\eps} \leq (\rho_{j_1}/\delta)^{\eta_{j_1}}$ (as $\eta_0 \geq \eta/\eps$). Note that $\rho_{j_1 + 1} \geq \delta(\rho_{j_1}/\delta)^\eps \geq \delta^{1-\eps^2}$. 

\begin{itemize}
    \item Case $1$: We have $\rho_{j_1 + 1} < \delta^{1-\eps}$, and 
    \begin{align*}
        \Delta_{max}((\F_1)_{\rho_{j_1+1}}) \lessapprox_\eps         \rho_{j_1+1}^{-2\eta_{j_1+1}}.
    \end{align*}
    
    \item Case $2$: Otherwise, by Proposition~\ref{Prop:inherit}(B), we have 
    \begin{align*}
        \Delta_{max}((\F_1)_{\rho_{j_1}}) \leq \Delta_{max}((\F^{(j_1)})_{\rho_{j_1}}) \lessapprox_\eps \rho_{j_1}^{-2\eta_{j_1}}.
    \end{align*} 

    For each $F_{\rho_{j_1}} \in (\F_1)_{\rho_{j_1}}$, there exists a family of convex sets $\W_{F_{\rho_{j_1}}}$ of dimensions $a_{j_1} \times b_{j_1} \times \1$ in $F_{\rho_{j_1}}$, such that $\W_{\rho_{j_1}}$ uniformly factors $(\F_1)_{\delta(\rho_{j_1}/\delta)^\eps}[F_{\rho_{j_1}}]$. At least one of the following holds.
    \begin{enumerate}[label={\rm (\roman*)}]
        \item $a_{j_1} \leq (\delta/\rho_{j_1})^{\eta_{j_1}'}b_{j_1}$.
        
        \item $a_{j_1} \geq (\delta/\rho_{j_1})^{\eta_{j_1}'} b_{j_1}$. For all $\rho \in [\delta(\rho_{j_1}/\delta)^\eps, \rho_{j_1}(\delta/\rho_{j_1})^\eps]$ and each $F_{\rho_{j_1}} \in (\F_1)_{\rho_{j_1}}$, we have $\Delta_{max}((\F_1)_{\rho}[F_{\rho_{j_1}}]) \geq (\rho_{j_1}/\rho)^{\eta_{j_1+1}}$. Moreover, $\F_1$ is weakly uniform at scales $[\delta, \rho_{j_1}]$ up to a factor $\approx (\rho_{j_1}/\delta)^{\eta_{j_1}'}$.
    \end{enumerate}
\end{itemize}
In Case $1$, we set $j_1' \coloneqq j_1 + 1$; while in Case $2$, we set $j_1' \coloneqq j_1$.\footnote{If both cases hold, then we just choose one of them to define the properties at scale $\rho_{j_1'}$.} Then in both cases, we have $\Delta_{max}((\F_1)_{\rho_{j_1'}}) \lessapprox_\eps \rho_{j_1'}^{-2\eta_{j_1'}}$. By pigeonholing, we may further assume that $\F_1$ is uniform at scale $\rho_{j_1'}$. This does not affect the properties in either case.

In Case $2$, we are actually done, as Conclusion~(II) holds with 
$$\F' = \F_1, \, \underline{\F} = \F, \, \tau = \delta, \, \theta = \rho_{j_1}, \, J=j_1+1, \, a = a_{j_1}, \, b = b_{j_1}.$$

Therefore, it remains to consider Case $1$. Recall that $\delta^{1-\eps^2} \leq \rho_{j_1'} < \delta^{1-\eps}$. The key point is that we can iterate the process described so far until we arrive at a good place. More precisely, we run the following iterative algorithm. Let $k \leq \lfloor 1/\eps^3 \rfloor$ be the loop index, with $1$ as its initial value. Let $\rho_{j_0'} \coloneqq \delta$.

Suppose we have got $\F_k \subset \F$ that satisfies $|\F_k| \gtrapprox_\eps \lambda_*^k|\F|$, and is uniform at scales $\{\rho_{j_i'}\}_{i=1}^k$ with $\rho_{j_i'} \geq (\rho_{j_{i-1}'})^{1-\eps^2}$ for each $1 \leq i \leq k$. Suppose also that 
$\rho_{j_{k}'} < \delta^\eps$, and 
\begin{align}\label{Eq:C_KT_upper_bound}
\Delta_{max}((\F_k)_{\rho_{j_k'}}) \lessapprox_\eps \rho_{j_k'}^{-2\eta_{j_k'}}.
\end{align}  
Define 
$$\F_k^* \coloneqq \bigsqcup_{F_{\rho_{j_k'}} \in \Phi_{\F, \rho_{j_k'}}(\F_k)} \F_k[F_{\rho_{j_k'}}].$$
Then $\F_k^* \subset \F_k$ is still uniform at scales $\{\rho_{j_i'}\}_{i=1}^k$, and $|\F_k^*| \geq \lambda_* |\F_k| \gtrapprox_\eps \lambda_*^{k+1}|\F|$ by our assumption on $\Phi_{\F, \rho}$. 

Repeat all the previous arguments with $\rho_{j_k'}$ in place of $\delta$, $(\F_k^*)_{\rho_{j_k'}}$ in place of $\Phi_{\F, \delta}(\F)$, and $\eta_j$ starting at $\eta_{j_k'+1}$. By Proposition~\ref{Prop:inherit}(B), all the previous bounds for ``$\Delta_{max}$'' are preserved. This will induce a $\gtrapprox_\eps 1$ refinement of $\F_k^*$ due to the uniformity of $\F_k^*$ at scale $\rho_{j_k'}$.

To begin, similarly to the initial $\F^{(1)} \subset \Phi_{\F, \delta}(\F)$, we get some $\F_k^{*,(1)} \subset \F_k^*$ with $|\F_k^{*,(1)}| \gtrapprox |\F_k^*| \gtrapprox_\eps \lambda_*^{k+1}|\F|$, and 
$$\Delta_{max}(\F_k^{*,(1)}) 
\leq \Delta_{max}(\F) \leq \delta^{-\eta} \leq \rho_{j_k'}^{-\eta/\eps} 
\leq  (1/\rho_{j_k'})^{2\eta_0}.$$
We omit the other details, as they are essentially the same.

In general, we run an iterative algorithm starting from $\F_k^{*,(1)}$. The outcome can be summarized as: There exist $\F_{k+1} \subset \F_k^*$ with $|\F_{k+1}| \gtrapprox_\eps |\F_k^*| \gtrapprox_\eps \lambda_*^{k+1}|\F|$, some $j_k' < j_{k+1}' \leq j_k' + 2 + \lfloor 1/\eps^2 \rfloor$, and $\rho_{j_{k+1}'} \geq (\rho_{j_k'})^{1-\eps^2}$, such that 
\begin{align*}
\Delta_{max}((\F_{k+1})_{\rho_{j_{k+1}'}}) \lessapprox_\eps \rho_{j_{k+1}'}^{-2\eta_{j_{k+1}'}}. 
\end{align*} 
Moreover, one of the following two possibilities holds:
\begin{itemize}
    \item Case $1$: $\rho_{j_{k+1}'} < \rho_{j_k'}^{1-\eps}$.  
    
    \item Case $2$: $\rho_{j_{k+1}'} \geq \rho_{j_k'}^{1-\eps}$. For each $F_{\rho_{j_{k+1}'}} \in (\F_{k+1})_{\rho_{j_{k+1}'}}$, there exists a family of convex sets $\W_{F_{\rho_{j_{k+1}'}}}$ of dimensions $a_{j_{k+1}'} \times b_{j_{k+1}'} \times \1$ in $F_{\rho_{j_{k+1}'}}$, such that $\W_{\rho_{j_{k+1}'}}$ uniformly factors $(\F_{k+1})_{\rho_{j_k'}(\rho_{j_{k+1}'}/\rho_{j_k'})^\eps}[F_{\rho_{j_{k+1}'}}]$. At least one of the following holds.
    \begin{enumerate}[label={\rm (\roman*)}]
        \item $a_{j_{k+1}'} \leq (\rho_{j_{k}'}/\rho_{j_{k+1}'})^{\eta_{j_{k+1}'}'}b_{j_{k+1}'}$.
        
        \item $a_{j_{k+1}'} \geq (\rho_{j_{k}'}/\rho_{j_{k+1}'})^{\eta_{j_{k+1}'}'}b_{j_{k+1}'}$. For all $\rho \in [\rho_{j_k'}(\rho_{j_{k+1}'}/\rho_{j_k'})^\eps, \rho_{j_{k+1}'}(\rho_{j_k'}/\rho_{j_{k+1}'})^\eps]$ and each $F_{\rho_{j_{k+1}'}} \in (\F_{k+1})_{\rho_{j_{k+1}'}}$, we have $$\Delta_{max}((\F_{k+1})_{\rho}[F_{\rho_{j_{k+1}'}}]) \geq (\rho_{j_{k+1}'}/\rho)^{\eta_{j_{k+1}'+1}}.$$
        Moreover, $\F_{k+1}$ is weakly uniform at scales $[\rho_{j_{k}'}, \rho_{j_{k+1}'}]$ up to a factor $\approx (\rho_{j_{k+1}'}/\rho_{j_{k}'})^{\eta_{j_{k+1}'}'}$.
    \end{enumerate}
\end{itemize}
By pigeonholing, we may further assume that $\F_{k+1}$ is uniform at scales $\{\rho_{j_i'}\}_{i=1}^{k+1}$. This does not affect the properties in either case. By construction of $\F_{k+1}$, we have $(\F_{k+1})_{\rho_{j_k'}} \subset (\F_k^*)_{\rho_{j_k'}} \subset \Phi_{\F, _{\rho_{j_k'}}}(\F_k)$, and $\F_{k+1} = \bigsqcup_{ (\F_{k+1})_{\rho_{j_k'}} } \F_k[F_{\rho_{j_k'}}]$.

In Case $2$, we are actually done, as Conclusion~(II) holds with 
$$\F' = \F_{k+1}, \, \underline{\F} = \F_k, \, \tau = \rho_{j_{k}'}, \, \theta = \rho_{j_{k+1}'}, \, J=j_{k+1}'+1, \, a = a_{j_{k+1}'}, \, b = b_{j_{k+1}'}.$$ 
To see this, note that $\rho_{j_k'} \leq \delta^\eps$ and $\rho_{j_{k+1}'} \geq \rho_{j_k'}^{1-\eps}$ imply $\rho_{j_{k+1}'} \geq \rho_{j_k'}\delta^{-\eps^2}$. Also, by Proposition~\ref{Prop:inherit}(B) and \eqref{Eq:C_KT_upper_bound}, we have 
\begin{align*}
    \Delta_{max}((\F_{k+1})_{\rho_{j_k'}}[F_{\rho_{j_{k+1}'}}]) 
    &\leq 
    \Delta_{max}((\F_{k+1})_{\rho_{j_k'}}) \\
    &\leq 
    \Delta_{max}((\F_{k})_{\rho_{j_k'}}) \\
    &\lessapprox_\eps \rho_{j_{k}'}^{-2\eta_{j_{k}'}} 
    \leq     (\rho_{j_{k+1}}'/\rho_{j_{k}}')^{2\eta_{j_{k}'}/\eps} \leq     (\rho_{j_{k+1}}'/\rho_{j_{k}}')^{2\eta_{j_{k+1}'}}
\end{align*}
for each $F_{\rho_{j_{k+1}'}} \in (\F_{k+1})_{\rho_{j_{k+1}'}}$.

Therefore, without loss of generality, we may always assume that Case $1$ holds.

If $\rho_{j_{k+1}'} \geq \delta^\eps$, then we stop. Otherwise, we go back to the start of the loop and proceed with $k+1$ in place of $k$.

As long as $\rho_{j_k'} \leq \delta^\eps$, we have $\rho_{j_{k+1}'} \geq (\rho_{j_k'})^{1-\eps^2} \geq \rho_{j_k'}\delta^{-\eps^3}$ for each $k$, and so $\rho_{j_{k+1}'} \geq \rho_{j_{1}'}\delta^{-k\eps^3}$ by induction. As $\rho_{j_{k+1}'} \leq 1$ and $\rho_{j_{1}'} \geq \delta$, we must have $k \leq \lfloor 1/\eps^3 \rfloor$, which matches the initial constraint on $k$. In other words, the algorithm must stop in at most $\lfloor 1/\eps^3 \rfloor$ steps. Let $k_0$ be the loop index at which we stop, then $\rho_{j_{k_0 + 1}'} \geq \delta^{\eps}$, and the output of the above algorithm is $\F' \coloneqq \F_{k_0 + 1} \subset \F$ with $|\F'| \gtrapprox_\eps \lambda_*^{k_0 + 1}|\F| \geq \lambda_*^{2/\eps^3}|\F|$. Note that $j_{k_0 + 1}' \leq (k_0 + 1) (2 + \lfloor 1/\eps^2 \rfloor) \leq 6/\eps^5 = 6N/10 \leq N-1$.

We are in the following case: For each $0 \leq k \leq k_0$\footnote{The $k=0$ case follows from the discussions about $\rho_{j_1'}$ before the iteration on $k$.}, $\F'$ satisfies Case~$1$ between scales $\rho_{j_k'}$ and $\rho_{j_{k+1}'}$ (with $\F_{k+1}$ replaced with $\F'$). 
For notational consistency, let $j_{k_0 +2}' \coloneqq j_{k_0 +1}' + 1 \leq N$ and $\rho_{j_{k_0 +2}'} \coloneqq 1$. For each \(0 \leq k \leq k_0 + 1\) and $F_{\rho_{j_k'}} \in \F'_{\rho_{j_k'}}$, we have 
\begin{align*}
\Delta_{max}(\F'_{\rho_{j_k'}}) \lessapprox_\eps \rho_{j_k'}^{-2\eta_{j_k'}},
\end{align*} 

Recall that $\rho_{j_0'} = \delta$, $\rho_{j_1'} < \delta^{1-\eps}$, $\rho_{j_{k_0 + 1}'} \geq \delta^{\eps}$, $\rho_{j_{k_0 +2}'} = 1$, and we are always in Case~$1$. Thus for each \(0 \leq k \leq k_0 + 1\), the ratio \( \rho_{j_{k+1}'} / \rho_{j_k'} \) is at most \( \rho_{j_k'}^{-\eps} \leq \delta^{-\eps} \). We claim that \( \F' \) is \( \lessapprox \delta^{-4\eps} \) Katz--Tao at every scale, and hence Conclusion (I) holds. To verify this, let \( \rho \in [\delta, 1] \) and choose \( 0 \leq k \leq k_0 + 1\) so that \( \rho_{j_k'} \leq \rho \leq \rho_{j_{k+1}'} \). Then we have 
\begin{align*}
\Delta_{max}(\F'_\rho) 
\lesssim \Delta_{max}(\F'_{\rho_{j_k'}}) \frac{|F_{\rho_{j_{k}'}}|}{|F_\rho|}
\lessapprox_\eps \rho_{j_k'}^{-2\eta_{j_k'}}  \left( \frac{\rho_{j_{k+1}'}}{\rho_{j_k'}} \right)^{2}  
\leq 
\delta^{-4\eps}.
\end{align*}
So the proof is completed.
\end{proof}

Lemma~\ref{Lem:multi-scale_decomposition} is often called the ``multiscale decomposition'', and is the key tool for the reduction we want. Lemma \ref{Lem:multi-scale_decomposition} (A) will be used for the proof of Main Lemma \ref{Main_lem_1}, and Lemma \ref{Lem:multi-scale_decomposition} (B) will be used for the proof of Main Lemma \ref{Main_lem_2}, as we will see in the next two sections.

\section{The proof of Main Lemma 1}\label{Sec:Main_Lemma_1}

In this section, we sketch the proof of Main Lemma~\ref{Main_lem_1}. Recall the statement:

\noindent{\bf Main Lemma 1.}
$K_{KT}(\beta)$ implies $K_F(\beta)$.

Let us first digest what the lemma is saying. Suppose that \( \F \) is a set of \( \delta \)-flakes in \( B_1 \) with \( C_F(\F) \leq \delta^{-\eta} \). We have to prove the bound
\begin{align}\label{Eq:Frostman_Goal}
\mu(\F) \leq \delta^{-\eps} (\delta^{-2})^\beta (\delta^2 |\F|)^{1-{\color{black} \beta/n}}. 
\end{align}

Since \( C_F(\F, B_1) \leq \delta^{-\eta} \), we have \( |\F| \geq \delta^{-2+\eta} \). If \( |\F| \approx \delta^{-2} \), then \( C_F(\F, B_1) \leq \delta^{-\eta} \) is equivalent to \( \Delta_{max}(\F) \lessapprox \delta^{-\eta} \), and then \( K_{KT}(\beta) \) gives \( \mu(\F) \leq \delta^{-\eps} \delta^{-2\beta} \), which is the desired bound. So the case \( |\F| \approx \delta^{-2} \) is trivial, and the content of the lemma is in the case when \( |\F| \gg \delta^{-2} \).

If \( |\F| \gg \delta^{-2} \), then \( \Delta_{max}(\F) \gg 1 \), so we cannot immediately apply \( K_{KT}(\beta) \) to \( \F \). If \( \tilde{\F} \subset \F \) is a random subset with \( |\tilde{\F}| \sim \delta^{-2} \), then it is not hard to check that \( \Delta_{max}(\tilde{\F}) \lessapprox 1 \). We can apply \( K_{KT}(\beta) \) to \( \tilde{\F} \), and we get \( \mu(\tilde{\F}) \lessapprox \delta^{-2\beta} \) which is equivalent to \( |U(\tilde{\F})| \gtrapprox \delta^{2\beta} \). We have \( |U(\F)| \geq |U(\tilde{\F})| \), and this translates into
\begin{align}\label{Eq:trivial_try}
    \mu(\F) \lessapprox (\delta^{-2})^\beta (\delta^2 |\F|).
\end{align}

Comparing \eqref{Eq:trivial_try} with our goal \eqref{Eq:Frostman_Goal}, we see that we need to reduce the exponent of the last factor from $1$ to \(1 - {\color{black} \beta/n}\). In our proof of $(n-2)$-plane Kakeya in $\R^n$, any exponent strictly less than 1 will work. 

\begin{rem}\label{Rem:difference}
    A similar observation has been pointed out in \cite{GWZ26streamlined} for line Kakeya in $\R^3$, in which the exponent can be reduced to $1- \beta/2$. In contrast, here we can only expect a worse improvement, which makes our modified formulation of $K_F(\beta)$ different from that in \cite[Definition~3.5]{GWZ26streamlined}. This ultimately comes down to the distinction between the dimension of sets of directional spaces: There are at most $\sim \delta^{-2}$ essentially distinct $\delta$-tubes passing through a fixed point in $\R^3$, while there are at most $\sim \delta^{-2(n-2)}$ essentially distinct $\delta$-flakes passing through a fixed point in $\R^n$. In other words, pointwise the flakes can be much more concentrated than tubes. Our Theorem~\ref{Thm:Main} shows that such much worse behavior still does not occur on average. The technical reason is that any exponent less than $1$ suffices for our purpose.
\end{rem} 

We can also describe the lemma in terms of \( |U(\F)| \), which may be more intuitive. Since we assume \( K_{KT}(\beta) \), we already know that if \( C_F(\F) \lessapprox 1 \) and \( |\F| \sim \delta^{-2} \), then \( |U(\F)| \gtrapprox \delta^{2\beta} \). Now we need to show that if \( C_F(\F) \lessapprox 1 \) and if \( |\F| \) is far larger than \( \delta^{-2} \), then \( |U(\F)| \) is a little larger than \( \delta^{2\beta} \). While this may sound intuitive, it is not at all trivial to prove, and the proof will depend crucially on $(n-2)$-plane sticky Kakeya.

We organize the proof as a bootstrapping or self-improving argument.

\begin{lem}\label{Lem:KT->F_bootstrap}
Suppose that \( K_{KT}(\beta) \) holds and that \( K_F(\gamma) \) holds for some \( \gamma > \beta \). Then \( K_F(\gamma - \nu) \) holds for \( \nu = \nu(\gamma, \beta) > 0 \). For \( \beta \) fixed, \( \nu = \nu(\gamma, \beta) > 0 \) is monotone in \( \gamma \).
\end{lem}

Since $K_{KT}(1)$ and $K_F(n-1)$ trivially hold, Lemma~\ref{Lem:KT->F_bootstrap} clearly implies Main Lemma~\ref{Main_lem_1}. 

\begin{proof}[Proof of Lemma~\ref{Lem:KT->F_bootstrap}]
By Definition~\ref{Def:Frostman_beta}, to prove $K_F(\gamma-\nu)$, we let $\F$ be a set of essentially distinct $\delta$-flakes in $B_1\subset \R^n$ with $C_F(\F,B_1) \leq \delta^{-\eta}$, and $Y$ be a shading on $\F$ with $\lambda(\F, Y) \geq \delta^{\eta}$. We need to bound $\mu(\F, Y)$ or $|U(\F, Y)|$. To begin with, let $(\F', Y)$ be a $\gtrapprox 1$ refinement of $(\F, Y)$ so that there is $\lambda > 0$ such that $|Y(F)|/|F| \sim \lambda$ for all $F \in \F'$. Then $|\F'| \gtrapprox \delta^{\eta} |\F|$, $C_F(\F') \lessapprox \delta^{-\eta} C_F(\F) \leq \delta^{-2\eta}$, and $\lambda \sim \lambda(\F',Y) \gtrapprox \lambda(\F,Y) \geq \delta^\eta$.

Let \(\eps_0\) and \(\delta_0\) be the output of Theorem~\ref{Thm:2-plane_sticky} with \(\gamma/n\) in place of \(\eps\). Let \(N \geq \max\{ 4(4n/\eps_0)^5, 10^5 \}\) be an integer to be chosen later, and define \(\eps = (10/N)^{1/5} \leq \eps_0/4n\). Let \( \eta \leq \eta_0 \leq \eta_1 \leq \ldots \leq \eta_N \leq \eps\) to be chosen later. Each \(\eta_{j-1}\) will be selected very small depending on \(\eta_{j}\), $\eps$, and \(\beta\).

Let $\eta_{j-1}' \coloneqq \frac{14\eta_{j-1}}{\eps^2\gamma}$ for each $1 \leq j \leq N$. 
By taking $3\eta \leq \eta_0 \eps$ and $\eta_j \geq \eta_{j-1}'/\eps$ for each $1\leq j \leq N$, we can apply Lemma~\ref{Lem:multi-scale_decomposition}(A) to $\F'$, with $\gamma$, $N$, $\eps$, $\{\eta_j\}_{j=0}^N$ as given, and $3\eta$ in place of $\eta$. The only subtle point is the choice of selection maps $\Phi_{\F',\rho}$, which is given by the following procedure. 

For any $\rho \in [\delta, 1]$ and $\underline{\F} \subset \F'$, by applying Lemma~\ref{Lem:fact_mult} to \((\underline{\F}, Y)\) and \(\underline{\F}_\rho\), we obtain a subset \( \underline{\F}'_\rho \subset \underline{\F}_\rho \), a shading \( Y_{\underline{\F}_\rho'}\) on \( \underline{\F}_\rho' \), and a shading \( Y' \) on \( \underline{\F}' = \bigsqcup_{F_\rho \in \underline{\F}_\rho'} \underline{\F}[F_\rho] \), such that $(\underline{\F}', Y')$ is a \(\gtrapprox 1\) refinement of \((\underline{\F}, Y)\), $|\underline{\F}[F_\rho]|$ are comparable for all $F_\rho \in \underline{\F}'_\rho$, and 
\begin{align*}
\mu(\underline{\F}, Y) \lessapprox 
\mu(\underline{\F}[F_\rho], Y') \mu(\underline{\F}'_\rho, Y_{\underline{\F}'_\rho})
\end{align*}
for all $F_\rho \in \underline{\F}'_\rho$. By pigeonholing, there is $\underline{\F}''_\rho\subset \underline{\F}'_\rho$ such that $\sum_{\sqcup_{\underline{\F}''_\rho} \underline{\F}[F_\rho]} |Y'(F)| \gtrapprox \sum_{\underline{\F}'} |Y'(F)|$, with $\lambda(\underline{\F}[F_\rho], Y') \sim \underline{\lambda}$ for all $F_\rho \in \underline{\F}''_\rho$. Note that 
\begin{align}\label{Eq:fine_scale_density}
\underline{\lambda} \sim \lambda\big(\bigsqcup_{\underline{\F}''_\rho} \underline{\F}[F_\rho], Y'\big) \gtrapprox \lambda(\underline{\F}', Y') \gtrapprox  \lambda(\underline{\F}, Y) \sim \lambda \gtrapprox \delta^\eta,
\end{align} 
so for each $F_\rho \in \underline{\F}''_\rho$, there must exist some $F \in \underline{\F}[F_\rho]$ such that $|Y'(F)|/|F| \gtrapprox \delta^\eta$. As $Y_{\underline{\F}_\rho'}$ is induced by $Y'$, this implies 
\begin{align}\label{Eq:coarse_scale_density}
|Y_{\underline{\F}_\rho'}(F_\rho)|/|F_\rho| \gtrapprox \delta^\eta
\end{align}
for all $F_\rho \in \underline{\F}''_\rho$. In particular, we have $\lambda(\underline{\F}_\rho'', Y_{\underline{\F}_\rho'}) \gtrapprox \delta^\eta$. On the other hand, note that 
\begin{align*}
|F||\underline{\F}''_\rho||\underline{\F}[F_\rho]| \geq \sum_{\sqcup_{\underline{\F}''_\rho} \underline{\F}[F_\rho]} |Y'(F)| \gtrapprox \sum_{\underline{\F}'} |Y'(F)| \gtrapprox \sum_{\underline{\F}} |Y(F)| \sim \lambda|F||\underline{\F}| \geq \lambda|F||\underline{\F}'_\rho||\underline{\F}[F_\rho]|,
\end{align*}
which implies $|\underline{\F}''_\rho| \gtrapprox \lambda |\underline{\F}'_\rho| \gtrapprox \delta^\eta |\underline{\F}'_\rho|$. In particular, 
\begin{align*}
\sum_{\underline{\F}''_\rho} |Y_{\underline{\F}'_\rho}(F_\rho)| = \lambda(\underline{\F}_\rho'', Y_{\underline{\F}_\rho'}) |F_\rho| |\underline{\F}''_\rho| \gtrapprox \delta^{2\eta} |F_\rho| |\underline{\F}'_\rho| \geq \delta^{2\eta} \sum_{\underline{\F}'_\rho} |Y_{\underline{\F}'_\rho}(F_\rho)|,
\end{align*}
which implies $\mu(\underline{\F}'_\rho, Y_{\underline{\F}'_\rho}) \lessapprox \delta^{-2\eta} \mu(\underline{\F}''_\rho, Y_{\underline{\F}'_\rho})$. Thus we have 
\begin{align}\label{Eq:mult_Phi_selection}
\mu(\underline{\F}, Y) \lessapprox \delta^{-2\eta} 
\mu(\underline{\F}[F_\rho], Y') \, \mu(\underline{\F}''_\rho, Y_{\underline{\F}'_\rho}).
\end{align}
for all $F_\rho \in \underline{\F}'_\rho$.

Now we can define the map $\Phi_{\F', \rho}: 2^{\F'} \rightarrow 2^{\F_\rho}$ by setting $\Phi_{\F', \rho}(\underline{\F}) \coloneqq \underline{\F}''_\rho$. Clearly $\Phi_{\F', \rho}(\underline{\F}) \subseteq \underline{\F}_\rho$, and  $$\big|\bigsqcup_{F_\rho \in \Phi_{\F', \rho}(\underline{\F})} \underline{\F}[F_\rho]\big| \sim |\underline{\F}''_\rho||\underline{\F}[F_\rho]| \gtrapprox \delta^\eta |\underline{\F}'_\rho||\underline{\F}[F_\rho]| \sim \delta^\eta |\underline{\F}'|.$$
As $(\underline{\F}', Y')$ is a \(\gtrapprox 1\) refinement of \((\underline{\F}, Y)\), we have 
$$|F||\underline{\F}'| \geq \sum_{\underline{\F}'}|Y'(F)| \gtrapprox \sum_{\underline{\F}}|Y(F)| \sim \delta^\eta |F||\underline{\F}|,$$ which implies $|\underline{\F}'| \gtrapprox \delta^\eta |\underline{\F}|$. So we can conclude that 
\begin{align*}
\big|\bigsqcup_{F_\rho \in \Phi_{\F', \rho}(\underline{\F})} \underline{\F}[F_\rho]\big| \geq \lambda_* |\underline{\F}|
\end{align*}
for some $\lambda_* \gtrapprox \delta^{2\eta}$. Apply Lemma~\ref{Lem:multi-scale_decomposition}(A) to this $\Phi_{\F', \rho}$ and $\lambda_*$.  

If Conclusion (I) holds, then there exists a $\gtrapprox_\eps \lambda_*^{2/\eps^3}$ refinement $\F'' \subset \F'$ that is \(\lessapprox_\eps \lambda_*^{-2/\eps^3} \delta^{-2n\eps} \) Frostman at every scale. As $\lambda_*^{-2/\eps^3} \delta^{-2n\eps} \lessapprox \delta^{-4\eta/\eps^3 - 2n\eps}$, by taking $4\eta \leq \eps^3\eps_0/3$ and recalling that $2n\eps \leq \eps_0/2$, we get that $\F''$ is $\delta^{-\eps_0}$ Frostman at every scale. Also, $\lambda(\F'', Y) \sim \lambda \gtrapprox \delta^\eta$. Hence by applying Theorem~\ref{Thm:2-plane_sticky}(A) to $(\F'', Y)$, we have
\[ 
|U(\F, Y)| \geq 
|U(\F'', Y)| \geq \delta^{(\gamma/2)/(n/2)} \gtrsim \delta^{2(\gamma/2)} (|\F| |F|)^{(\gamma/2)/{\color{black}n}},
\]
where we used the fact that \(|\F| \lesssim \delta^{-2(n-1)}\). In particular, \(|U(\F, Y)|\) satisfies the volume estimate corresponding to \(K_F(\gamma/2)\). In this case, Lemma~\ref{Lem:KT->F_bootstrap} holds with \(\nu = \gamma/2\).

\begin{rem}\label{Rem:numerology}
Recall that the statement \(K_F(\beta)\) roughly asserts \(|U(\F)| \gtrapprox \delta^{2\beta} (|\F| |F|)^{u}\) with \(u = \beta/n\). In the argument above, we actually only need \(u < \beta/(n-2)\).
\end{rem}

Henceforth we shall assume that Conclusion (II) of Lemma~\ref{Lem:multi-scale_decomposition}(A) holds. Let $\F'' \subset \underline{\F} \subset \F'$, $\delta \leq \tau \leq \theta \leq 1$ (with $\tau \leq \delta^{\eps^2}\theta$), \(1 \leq J \leq N\), and $\tau \leq a \leq b \leq \theta(\tau/\theta)^\eps$ be the output of that lemma. As $|\underline{\F}| \geq |\F''| \gtrapprox_\eps \lambda_*^{2/\eps^3} |\F'| \gtrapprox \delta^{4\eta/\eps^3 + \eta} |\F|$ and $|Y(F)|/|F| \sim \lambda$ for all $F \in \F'$, we know that $(\underline{\F}, Y)$ is a $\gtrapprox_\eps \delta^{5\eta/\eps^3}$ refinement of $(\F, Y)$. Thus 
\begin{align*}
    \mu(\F, Y) \lessapprox_\eps \delta^{-5\eta/\eps^3} \mu(\underline{\F}, Y).
\end{align*}
By our construction of $\Phi_{\F', \tau}$, the bound \eqref{Eq:mult_Phi_selection} says that \begin{align*}
\mu(\underline{\F}, Y) \lessapprox \delta^{-2\eta} 
\mu(\underline{\F}[F_\tau], Y') \, 
\mu(\Phi_{\F', \tau}(\underline{\F}), Y_{\tau}).
\end{align*}
for all $F_\tau \in \Phi_{\F', \tau}(\underline{\F})$. Here $Y_\tau$ is a shading on $\Phi_{\F', \tau}(\underline{\F})$ so that $$|Y_\tau(F_\tau)|/|F_\tau| \gtrapprox \delta^\eta$$ 
for all $F_\tau \in \Phi_{\F', \tau}(\underline{\F})$ by \eqref{Eq:coarse_scale_density}. Since $\F''_\tau \subset \Phi_{\F', \tau}(\underline{\F})$ and $\F'' = \bigsqcup_{\F''_\tau} \underline{\F}[F_\tau]$, for all $F_\tau \in \F''_\tau$, we have $\F''[F_\tau] = \underline{\F}[F_\tau]$, which yields $$\mu(\F''[F_\tau], Y') = \mu(\underline{\F}[F_\tau], Y'),$$
and 
$$\lambda(\F''[F_\tau], Y') = \lambda(\underline{\F}[F_\tau], Y') \gtrapprox \delta^\eta$$ 
by \eqref{Eq:fine_scale_density}. On the other hand, since $|\underline{\F}[F_\tau]|$ are comparable for all $F_\tau \in \Phi_{\F', \tau}(\underline{\F})$, we have 
$$|\F''_\tau|/|\Phi_{\F', \tau}(\underline{\F})| \sim |\F''|/\big|\bigsqcup_{\Phi_{\F', \tau}(\underline{\F})} \underline{\F}[F_\tau]\big| \geq |\F''|/|\underline{\F}| \geq |\F''|/|\F'| \gtrapprox_\eps \lambda_*^{2/\eps^3} \gtrapprox \delta^{4\eta/\eps^3},$$ 
and so 
$$\sum_{\F''_\tau} |Y_\tau(F_\tau)| \gtrapprox \delta^\eta|F_\tau||\F''_\tau| \gtrapprox_\eps \delta^{4\eta/\eps^3 + \eta} |F_\tau||\Phi_{\F', \tau}(\underline{\F})| \geq \delta^{5\eta/\eps^3} \sum_{\Phi_{\F', \tau}(\underline{\F})} |Y_\tau(F_\tau)|.$$ 
In other words, $(\F''_\tau, Y_\tau)$ is a $\gtrapprox_\eps \delta^{5\eta/\eps^3}$ refinement of $(\Phi_{\F', \tau}(\underline{\F}), Y_\tau)$, which means 
\begin{align*}
    \mu(\Phi_{\F', \tau}(\underline{\F}), Y_{\tau}) \lessapprox_\eps \delta^{-5\eta/\eps^3} 
    \mu(\F''_\tau, Y_\tau).
\end{align*}
Therefore, by putting all multiplicity bounds together, we have 
\begin{align}\label{Eq:mult_Phi_selection'}
    \mu(\F, Y) \lessapprox_\eps \delta^{-12\eta/\eps^3} 
    \mu(\F''[F_\tau], Y') \, 
    \mu(\F''_\tau, Y_{\F''_\tau})
\end{align}
for all $F_\tau \in \F''_\tau$. Here $Y_{\F''_\tau}$ is the restriction of $Y_\tau$ to $\F''_\tau$. To summarize, for all $F_\tau \in \F''_\tau$, we have $|Y_{\F''_\tau}(F_\tau)|/|F_\tau| \gtrapprox \delta^\eta$ and $\lambda(\F''[F_\tau], Y') \gtrapprox \delta^\eta$.

Next, as $\F''$ is uniform at scales $\tau$ and $\theta$, by applying Lemma~\ref{Lem:fact_mult} (Remark~\ref{Rem:factor_only_shading}) to \((\F_\tau'', Y_{\F_\tau''})\) and \(\F_\theta''\), we obtain a \(\gtrapprox 1\) refinement $(\F_\tau'', Y_{\F_\tau''}')$ of \((\F_\tau'', Y_{\F_\tau''})\) and a shading \(Y_{\F_\theta''}\) on \(\F_\theta''\) with $\lambda(\F_\theta'', Y_{\F_\theta''}) \gtrapprox \delta^\eta$, such that 
\begin{align*}
    \mu(\F_\tau'', Y_{\F_\tau''}) \lessapprox     \mu(\F_\tau''[F_\theta], Y_{\F_\tau''}') \, \mu(\F_\theta'', Y_{\F_\theta''})
\end{align*}
for all $F_\theta \in \F_\theta''$. Combining this multiplicity bound with
\eqref{Eq:mult_Phi_selection'} and taking $12\eta \leq \eta_0\eps^3$, we conclude that 
\begin{align}\label{Eq:3_scale}
\mu(\F, Y) \lessapprox_\eps \delta^{-\eta_0} 
\mu(\F''[F_\tau], Y') \, \mu(\F_\tau''[F_\theta], Y_{\F_\tau''}') \, \mu(\F_\theta'', Y_{\F_\theta''})
\end{align}
for all \(F_\tau \in \F_\tau''\) and \(F_\theta \in \F_\theta''\).

Our goal is to estimate the three quantities on the RHS of \eqref{Eq:3_scale}. Properties in conclusion (II) of Lemma~\ref{Lem:multi-scale_decomposition}(A) related to convex Wolff axioms and factoring structure (for $\F''$) would be very helpful.

For each \(F_\tau \in \F_\tau''\), we have \(C_F(\F''[F_\tau]) \lessapprox_\eps \delta^{-2\eta_{J-1}}\), and recall that $\lambda(\F''[F_\tau],Y') \gtrapprox \delta^\eta$. Then by taking $\eta \ll \eta_0$ and applying Lemma~\ref{Lem:general_F} to the rescaled sets \( (\F''[F_\tau], Y') \) (with $\gamma$ in place of $\beta$ and \(\eta_0\) in place of \(\eps\)), we have 
\begin{equation}\label{Eq:delta_tau}
\begin{split}
\mu(\F''[F_\tau], Y') 
&\lessapprox_\eps (\delta/\tau)^{-\eta_0} \delta^{-(1-\gamma/n)2\eta_{J-1}} ((\delta/\tau)^{-2})^{\gamma} ((\delta/\tau)^2 |\F''[F_\tau]|)^{1-\gamma/n} \\
\Longrightarrow \mu(\F''[F_\tau], Y')&\leq \delta^{-3\eta_{J-1}} ((\delta/\tau)^{-2})^{\gamma} ((\delta/\tau)^2 |\F''[F_\tau]|)^{1-\gamma/n}.
\end{split}
\end{equation}

On the other hand, by Proposition~\ref{Prop:inherit}(A) and the uniformity of $\F''$ at scale $\theta$, we have $C_F(\F_\theta'') \lesssim C_F(\F'') \lessapprox_\eps \delta^{-5\eta/\eps^3} C_F(\F) \leq \delta^{-6\eta/\eps^3}$ (recall that $|\F''| \gtrapprox_\eps \delta^{5\eta/\eps^3}|\F|$). 
Thus $C_F(\F_\theta'') \lessapprox_\eps \delta^{-6\eta/\eps^3}$. Recall that $\lambda(\F_\theta'', Y_{\F_\theta''}) \gtrapprox \delta^\eta$. Then by taking $\eta \ll \eta_0\eps^3$ and applying \(K_F(\gamma)\) to $(\F_\theta'', Y_{\F_\theta''})$ (Definition~\ref{Def:Frostman_beta}, with $\gamma$ in place of $\beta$ and \(\eta_0\) in place of \(\eps\)), we have
\begin{align}\label{Eq:theta_unit}
    \mu(\F_\theta'', Y_{\F_\theta''}) \leq \delta^{-\eta_{J-1}} ( \theta^{-2})^{\gamma} (\theta^2 |\F_\theta''|)^{1-\gamma/n}.
\end{align}

Therefore, by combining \eqref{Eq:3_scale}, \eqref{Eq:delta_tau}, and \eqref{Eq:theta_unit}, one can easily see that to prove Lemma~\ref{Lem:KT->F_bootstrap} with \( \nu = \eta_0 \), it suffices to find some \( F_\theta \in \F_\theta'' \) with 
\begin{align}\label{Eq:tau_theta}    \mu(\F_\tau''[F_\theta], Y_{\F_\tau''}') \leq \delta^{10\eta_{J-1}}((\tau/\theta)^{-2})^{\gamma}((\tau/\theta)^2|\F_\tau''[F_\theta]|)^{1-\gamma/n}.
\end{align}
We will show that this holds, provided each \( \eta_{j-1} \) is selected sufficiently small compared to \( \eta_j \). 

As $\lambda(\F_\tau'',Y_{\F_\tau''}') \gtrapprox \lambda(\F_\tau'',Y_{\F_\tau''}) \gtrapprox \delta^\eta$, by pigeonholing, there exists some $F_\theta \in \F_\theta''$ such that $\lambda(\F_\tau''[F_\theta],Y_{\F_\tau''}') \gtrapprox \delta^\eta$. Fix one such $F_\theta$. Let \( (\tilde{\F}, \tilde{Y}) \) be the image of $(\F_\tau''[F_\theta],Y_{\F_\tau''}')$ under the rescaling taking \( F_\theta \) to $B_1$. Hence \( \tilde{\F} \) is a family of \( \tau/\theta \)-flakes, and $\lambda(\tilde{\F}, \tilde{Y}) \gtrapprox \delta^\eta$. {\color{black} Let $\tilde{\delta} \coloneqq \tau/\theta$.} Recall that \( \tilde{\delta} \leq \delta^{\eps^2} \). Thus to prove \eqref{Eq:tau_theta}, we need to show that
\begin{align}\label{Eq:tau_theta_rescaled}
    \mu(\tilde{\F}, \tilde{Y}) \leq \tilde{\delta}^{10\eta_{J-1}/\eps^2}\tilde{\delta}^{-2\gamma}(\tilde{\delta}^2|\tilde{\F}|)^{1-\gamma/n}.
\end{align}

Recall that $\F''$, $\delta \leq \tau \leq \theta \leq 1$, \(1 \leq J \leq N\), and $\tau \leq a \leq b \leq \theta(\tau/\theta)^\eps$ are the output of Conclusion (II) of Lemma~\ref{Lem:multi-scale_decomposition}(A). So we have 
$$C_F(\tilde{\F}) \lessapprox_\eps \lambda_*^{-2/\eps^3} \tilde{\delta}^{-2\eta_{J-1}} \lessapprox \delta^{-4\eta/\eps^3}\tilde{\delta}^{-2\eta_{J-1}} \leq \tilde{\delta}^{-4\eta/\eps^5 - 2\eta_{J-1}} \leq \tilde{\delta}^{- 3\eta_{J-1}}$$
by taking $4\eta \leq \eta_0\eps^5$. For each $F_{\tilde{\delta}^{\eps}} \in \tilde{\F}_{\tilde{\delta}^{\eps}}$, there exists a family of convex sets $\tilde{\W}_{F_{\tilde{\delta}^\eps}}$ of dimensions $\tilde{a} \times \tilde{b} \times \1$ ($\tilde{a} = a/\theta$, $\tilde{b} = b/\theta$) in $F_{\tilde{\delta}^{\eps}}$, such that $\tilde{\W}_{F_{\tilde{\delta}^\eps}}$ uniformly factors $\tilde{\F}[F_{\tilde{\delta}^{\eps}}]$. Note that $\tilde{\delta} \leq \tilde{a} \leq \tilde{b} \leq \tilde{\delta}^\eps$. 
Moreover, $\tilde{\F}$ is uniform at scale $\tilde{b}$, which implies $C_F(\tilde{\F}_{\tilde{b}}) \lesssim C_F(\tilde{\F}) \lessapprox_\eps \tilde{\delta}^{- 3\eta_{J-1}}$ by Proposition~\ref{Prop:inherit}(A). At least one of the following holds.
\begin{enumerate}[label={\rm (\roman*)}]
\item $\tilde{a} \leq \tilde{\delta}^{\eta_{J-1}'} \tilde{b}$.

\item $\tilde{a} \geq \tilde{\delta}^{\eta_{J-1}'} \tilde{b}$. In this case, we additionally have 
\begin{align}\label{Eq:C_F_P_rho}
C_F(\tilde{\F}[F_\rho]) 
\geq 
(\rho/\tilde{\delta})^{\eta_J} \quad \text{for all } \rho \in [\tilde{\delta}^{1-\eps}, \tilde{\delta}^\eps] \text{ and } F_\rho \in \tilde{\F}_\rho.
\end{align}
\end{enumerate}

Since $\lambda(\tilde{\F}, \tilde{Y}) \gtrapprox \delta^\eta \geq \tilde{\delta}^{\eta/\eps^2}$ and $C_F(\tilde{\F}) \lessapprox_\eps \tilde{\delta}^{-3\eta_{J-1}}$, by taking $\eta/\eps^2 \ll \eta_0$ and applying Proposition~\ref{Prop:plank_factor_improve}(A) to $(\tilde{\F}, \tilde{Y})$ (with \( \gamma \) in place of \( \beta \), $\tilde{\delta}^\eps$ in place of $\rho$, and \( \eta_0 \) in place of \( \eps \)), we have
\[
\mu(\tilde{\F}, \tilde{Y}) \leq \tilde{\delta}^{-4\eta_{J-1}/\eps^2}\left(\frac{\tilde{a}}{\tilde{b}}\right)^{(n+2)\gamma/n}(\tilde{\delta}^{-2})^{\gamma}(\tilde{\delta}^2|\tilde{\F}|)^{1-\gamma/n}.
\]
In particular, either \eqref{Eq:tau_theta_rescaled} holds (and hence we are done), or else \( \tilde{a} \) is not too much smaller than \( \tilde{b} \), i.e., 
\[
\tilde{a} \geq \tilde{\delta}^{\frac{{\color{black}14}\eta_{J-1}}{\eps^2\gamma}} \tilde{b} = \tilde{\delta}^{\eta_{J-1}'} \tilde{b}.
\]

Pick a $W \in \tilde{\W}_{F_{\tilde{\delta}^\eps}}$. Since $W$ must be contained in some $F_{\tilde{b}} \in \tilde{\F}_{\tilde{b}}[F_{\tilde{\delta}^\eps}]$, we have 
\[
\Delta(\tilde{\F}, F_{\tilde{b}}) \gtrsim (\tilde{a}/\tilde{b})\Delta(\tilde{\F}, W) \approx (\tilde{a}/\tilde{b})\Delta_{max}(\tilde{\F}[F_{\tilde{\delta}^{\eps}}]) \geq (\tilde{a}/\tilde{b})\Delta_{max}(\tilde{\F}[F_{\tilde{b}}]).
\]
Therefore,
\begin{align}\label{Eq:C_F(P_b)}
C_F(\tilde{\F}[F_{\tilde{b}}]) \lessapprox \tilde{b}/\tilde{a} \leq \tilde{\delta}^{-\eta_{J-1}'}.
\end{align}
Since $\tilde{\F}$ is uniform at scale $\tilde{b}$, \eqref{Eq:C_F(P_b)} actually holds for all $F_{\tilde{b}} \in \tilde{\F}_{\tilde{b}}$.

Due to \( \tilde{b} \leq \tilde{\delta}^{\eps} \), if we select each \( \eta_{j-1} \) sufficiently small depending on \( \eta_j \) and \( \gamma \) (and depending on \( \eps \), which in turn depends on \( \gamma \)), then we must have \( \tilde{b} \leq \tilde{\delta}^{1-\eps} \) -- if not, then \eqref{Eq:C_F(P_b)} would violate \eqref{Eq:C_F_P_rho} (with \( \tilde{b} \) in place of \( \rho \)).

Let \( (\tilde{\F}, \bar{Y}) \) be the image of $(\F_\tau''[F_\theta],Y_{\F_\tau''})$ under the rescaling taking \( F_\theta \) to $B_1$. 
By pigeonholing, there is a \(\gtrapprox 1\) refinement $(\tilde{\F}', \tilde{Y})$ of \( (\tilde{\F}, \tilde{Y}) \), such that $|\tilde{Y}(F_{\tilde{\delta}})|/|F_{\tilde{\delta}}| \gtrapprox \delta^\eta$ for all $F_{\tilde{\delta}} \in \tilde{\F}'$. Then by applying Lemma~\ref{Lem:fact_mult} to \( (\tilde{\F}', \tilde{Y}) \) and \( \tilde{\F}_{\tilde{b}} \), we obtain a \(\gtrapprox 1\) refinement $(\tilde{\F}'', \tilde{Y}')$ of \( (\tilde{\F}', \tilde{Y}) \) and a shading \(Y_{\tilde{\F}_{\tilde{b}}'}\) on \( \tilde{\F}_{\tilde{b}}' \subset \tilde{\F}_{\tilde{b}} \) with $\lambda(\tilde{\F}_{\tilde{b}}', Y_{\tilde{\F}_{\tilde{b}}'}) \gtrapprox \delta^\eta$, such that 
\begin{align*}
    \mu(\tilde{\F}, \tilde{Y}) 
    \lessapprox
    \mu(\tilde{\F}', \tilde{Y})
    \lessapprox     \mu(\tilde{\F}''[F_{\tilde{b}}], \tilde{Y}') \, \mu(\tilde{\F}_{\tilde{b}}', Y_{\tilde{\F}_{\tilde{b}}'})
\end{align*}
for all $F_{\tilde{b}} \in \tilde{\F}_{\tilde{b}}'$. Note that $(\tilde{\F}'', \tilde{Y}')$ is a $\gtrapprox 1$ refinement of $(\tilde{\F}, \bar{Y})$. Thus by pigeonholing, there exists some $F_{\tilde{b}} \in \tilde{\F}_{\tilde{b}}''$ such that $(\tilde{\F}''[F_{\tilde{b}}], \tilde{Y}')$ is a $\gtrapprox 1$ refinement of $(\tilde{\F}[F_{\tilde{b}}], \bar{Y})$ (thus of $(\tilde{\F}'[F_{\tilde{b}}], \tilde{Y})$). Fix one such $F_{\tilde{b}}$. 

Note that $\lambda(\tilde{\F}''[F_{\tilde{b}}], \tilde{Y}') \gtrapprox \lambda(\tilde{\F}'[F_{\tilde{b}}], \tilde{Y}) \gtrapprox \delta^\eta$. Also, by 
$$|F_{\tilde{\delta}}||\tilde{\F}''[F_{\tilde{b}}]| \geq \sum_{\tilde{\F}''[F_{\tilde{b}}]}|\tilde{Y}'(F_{\tilde{\delta}})| \gtrapprox \sum_{\tilde{\F}[F_{\tilde{b}}]}|\bar{Y}(F_{\tilde{\delta}})| \gtrapprox 
\delta^\eta |F_{\tilde{\delta}}| |\tilde{\F}[F_{\tilde{b}}]|,
$$  
we have $|\tilde{\F}''[F_{\tilde{b}}]| \gtrapprox \delta^\eta |\tilde{\F}[F_{\tilde{b}}]|$, which implies that $C_F(\tilde{\F}''[F_{\tilde{b}}]) \lessapprox \delta^{-\eta} C_F(\tilde{\F}[F_{\tilde{b}}]) \lessapprox \tilde{\delta}^{-2\eta_{J-1}'}$ by using \eqref{Eq:C_F(P_b)}. By taking $\eta \ll \eta_0$ and applying Lemma~\ref{Lem:general_F} to the rescaled sets \( (\tilde{\F}''[F_{\tilde{b}}], \tilde{Y}') \) (with $\gamma$ in place of $\beta$ and \(\eta_0\) in place of \(\eps\)), we get
\begin{align*}
    \mu(\tilde{\F}''[F_{\tilde{b}}], \tilde{Y}') \lesssim \tilde{\delta}^{-3\eta_{J-1}'} (\tilde{\delta}/\tilde{b})^{-2\gamma}((\tilde{\delta}/\tilde{b})^2|\tilde{\F}[F_{\tilde{b}}]|)^{1-\gamma/n}.
\end{align*}
Therefore, in order to establish \eqref{Eq:tau_theta_rescaled}, it suffices to prove that
\begin{align}\label{Eq:tilde_P_b}
\mu(\tilde{\F}_{\tilde{b}}', Y_{\tilde{\F}_{\tilde{b}}'}) \leq \tilde{\delta}^{{\color{black} (3+\gamma)}\eta_{J-1}'} \tilde{b}^{-2\gamma}(\tilde{b}^2 |\tilde{\F}_{\tilde{b}}|)^{1-\gamma/n}.
\end{align}

On the other hand, since each flake \( F_{\tilde{b}} \in \tilde{\F}_{\tilde{b}}[F_{\tilde{\delta}^\eps}] \) essentially contains at least one \( W \in \tilde{\W}_{F_{\tilde{\delta}^\eps}} \), we have
\[
\Delta_{max}(\tilde{\F}_{\tilde{b}}[F_{\tilde{\delta}^\eps}]) \lesssim \frac{|F_{\tilde{b}}|}{|W|} \Delta_{max}(\tilde{\W}_{F_{\tilde{\delta}^\eps}}) \lesssim \frac{\tilde{b}}{\tilde{a}} \leq \tilde{\delta}^{-\eta_{J-1}'}.
\]
Therefore, 
\[
\Delta_{max}(\tilde{\F}_{\tilde{b}}')
\leq
\Delta_{max}(\tilde{\F}_{\tilde{b}})
\lesssim (\tilde{\delta}^\eps)^{-2(n-1)} \Delta_{max}(\tilde{\F}_{\tilde{b}}[F_{\tilde{\delta}^\eps}]) \lesssim \tilde{\delta}^{-2(n-1)\eps - \eta_{J-1}'}.
\]
Recall that $\lambda(\tilde{\F}_{\tilde{b}}', Y_{\tilde{\F}_{\tilde{b}}'}) \gtrapprox \delta^\eta$. 
By taking $\eta \ll \eta_0$ and applying Lemma~\ref{Lem:general_KT} to $(\tilde{\F}_{\tilde{b}}', Y_{\tilde{\F}_{\tilde{b}}'})$ (with \(\eta_0\) in place of \(\eps\)), we get 
\begin{align*}
\mu(\tilde{\F}_{\tilde{b}}', Y_{\tilde{\F}_{\tilde{b}}'}) \leq \tilde{\delta}^{-(2n-1)\eps - \eta_{J-1}'} |\tilde{\F}_{\tilde{b}}|^{\beta}
= \left( \tilde{\delta}^{-(2n-1)\eps - \eta_{J-1}'} (\tilde{b}^2 |\tilde{\F}_{\tilde{b}}|)^{\gamma/n + \beta - 1} \right) \tilde{b}^{-2\beta} (\tilde{b}^2 |\tilde{\F}_{\tilde{b}}|)^{1-\gamma/n}.
\end{align*}
We wish to bound the first bracketed term. The key point is that \(\tilde{b}^2 |\tilde{\F}_{\tilde{b}}|\) is close to $1$. To see this, we note that
\[
\tilde{\delta}^{3\eta_{J-1}} \lessapprox_\eps C_F (\tilde{\F}_{\tilde{b}})^{-1} \leq \tilde{b}^2 |\tilde{\F}_{\tilde{b}}| \lesssim \Delta_{max}(\tilde{\F}_{\tilde{b}}) \lesssim \tilde{\delta}^{-2(n-1)\eps - \eta_{J-1}'}.
\]

We do not know the sign of \(\gamma/n + \beta - 1\), but by using the above upper and lower bounds for \(\tilde{b}^2 |\tilde{\F}_{\tilde{b}}|\) and recalling $\tilde{b} \leq \tilde{\delta}^{1-\eps}$, we get
\[
\mu(\tilde{\F}_{\tilde{b}}', Y_{\tilde{\F}_{\tilde{b}}'}) \leq \tilde{\delta}^{-(4n-3)\eps - 2\eta_{J-1}' + 2(1-\eps)(\gamma - \beta)} \tilde{b}^{-2\gamma} (\tilde{b}^2 |\tilde{\F}_{\tilde{b}}|)^{1-\gamma/n}.
\]
Provided we select \(N\) sufficiently large (recall that \(\eps = (10/N)^{1/5}\)) and \(\eta_N\) sufficiently small (recall \(\eta_{J-1} \ll \eta_N\)) so that \[-(4n-3)\eps - 2\eta_{J-1}' + 2(1-\eps)(\gamma - \beta) \geq (\gamma - \beta),\] 
we have
\[
\mu(\tilde{\F}_{\tilde{b}}', Y_{\tilde{\F}_{\tilde{b}}'}) \leq \tilde{\delta}^{(\gamma - \beta)} \tilde{b}^{-2\gamma} (\tilde{b}^2 |\tilde{\F}_{\tilde{b}}|)^{1-\gamma/n}.
\]
This implies our desired bound \eqref{Eq:tilde_P_b}, provided we select each \(\eta_{j-1}\) sufficiently small so that \((3+\gamma)\eta_{j-1}' \leq \gamma - \beta\).
\end{proof}

\section{The proof of Main Lemma 2}\label{Sec:Main_Lemma_2}

Our last goal is to prove Main Lemma~\ref{Main_lem_2}. Recall the statement:

{\bf Main Lemma 2.}
If \( K_{KT}(\beta) \) and \( K_F(\beta) \) hold, then \( K_{KT}(\beta - \nu) \) holds, where \( \nu = \nu(\beta) > 0 \) is monotone.

\subsection{Reduction to the strongly non-sticky case.}
The $(n-2)$-plane sticky Kakeya theorem plays a crucial role in this proof. By using Theorem~\ref{Thm:2-plane_sticky}(B), we are able to reduce the problem to a situation which is not sticky in a strong sense.

\begin{lem}[Very not sticky case of Main Lemma~\ref{Main_lem_2}]\label{Lem:NSAAS}
For all \( \beta > 0 \), there is some \( \eps_{\rm{scal}} = \eps_{\rm{scal}}(\beta) \in (0,1/2) \) so that for all \( \zeta > 0 \), there is some \( \nu = \nu(\beta, \zeta) > 0 \) and \( \eta = \eta(\beta, \zeta) > 0 \) so that the following is true for small $\delta$:

If \( K_{KT}(\beta) \) and \( K_F(\beta) \) hold and if \( (\F, Y) \) is a uniform set of \( \delta \)-flakes in \( B_1 \subset \R^n \) obeying
\begin{itemize}
    \item \( \Delta_{max}(\F) \leq \delta^{-\eta} \) and \( \lambda(\F, Y) \geq \delta^{\eta} \).
    \item For each \( \rho \in [\delta^{1-\eps_{\rm{scal}}}, \delta^{\eps_{\rm{scal}}}], \, |\F_\rho| \geq \rho^{-2-\zeta} \).
\end{itemize}
Then
\[
\mu(\F, Y) \leq \delta^\nu |\F|^\beta.
\]
\end{lem}

For technical reasons, let us first show that Lemma~\ref{Lem:NSAAS} immediately implies a seemingly stronger version of it, which only requires weak uniformity up to a small factor:
\begin{lem}\label{Lem:NSAAS'}
For all \( \beta > 0 \), there is some \( \eps_{\rm{scal}} = \eps_{\rm{scal}}(\beta) \in (0,1/2) \) so that for all \( \zeta > 0 \), there is some \( \nu = \nu(\beta, \zeta) > 0 \) and \( \eta = \eta(\beta, \zeta) > 0 \) so that the following is true for small $\delta$:

If \( K_{KT}(\beta) \) and \( K_F(\beta) \) hold and if \( (\F, Y) \) is a set of {\color{black} essentially distinct} \( \delta \)-flakes in \( B_1 \subset \R^n \) obeying
\begin{itemize}
    \item For each $\rho \in [\delta, 1]$, there is $N_\rho > 0$ such that $$\delta^\eta N_\rho \leq |\F[F_\rho]| \leq \delta^{-\eta} N_\rho$$ 
    for all $F_\rho \in \F_\rho$.
    \item \( \Delta_{max}(\F) \leq \delta^{-\eta} \) and \( \lambda(\F, Y) \geq \delta^{\eta} \).
    \item For each  $\rho \in [\delta^{1-\eps_{\rm{scal}}}, \delta^{\eps_{\rm{scal}}}]$, $|\F_\rho| \geq \rho^{-2-\zeta}$.
\end{itemize}
Then
\[
\mu(\F, Y) \leq \delta^\nu |\F|^\beta.
\]
\end{lem}
\begin{proof}
    Note that our assumptions imply $|\F| \geq |\F_\rho| \delta^\eta N_\rho$, and so $N_\rho \leq \delta^{-\eta} |\F|/|\F_\rho|.$
    
    Let $(\F', Y)$ be a $\gtrapprox 1$ refinement of $(\F, Y)$ so that $|Y(F)|/|F|$ is constant for all $F \in \F'$. By Proposition~\ref{Prop:CWA_uniform}, there is a uniform subset $\F'' \subset \F'$, with $|\F''| \gtrapprox |\F'|$. Then $(\F'', Y)$ is a $\gtrapprox 1$ refinement of $(\F', Y)$ (thus of $(\F, Y)$). As in Section~\ref{subsec:uniform}, by applying dyadic pigeonholing to the shading $Y$, we find a further uniform $\gtrapprox$ refinement $(\F'', Y')$ of $(\F, Y)$. Note that $\lambda(\F'', Y') \gtrapprox \lambda(\F, Y) \geq \delta^\eta$, and by Proposition~\ref{Prop:inherit}(B), we have $\Delta_{max}(\F'') \leq \Delta_{max}(\F) \leq \delta^{-\eta}$.
    
    Since 
    $$|\F''||F| \geq \sum_{\F''}|Y(F)| \gtrapprox \sum_{\F}|Y(F)| \geq \delta^\eta |\F||F|,$$ 
    we have $|\F''| \gtrapprox \delta^{\eta} |\F|$. Now the key observation is that for each $\rho \in [\delta,1]$ and $F_\rho \in \F''_\rho$, we have 
    $$\F''[F_\rho] \leq |\F[F_\rho]| \leq \delta^{-\eta}N_\rho \leq \delta^{-2\eta}|\F|/|\F_\rho|,$$ 
    which by a simple counting argument implies 
    \begin{align*}
        |\F''_\rho| \geq |\F''|/(\delta^{-2\eta}|\F|/|\F_\rho|) 
        \gtrapprox \delta^{3\eta} |\F_\rho| \geq 
        \delta^{3\eta}
        \rho^{-2-\zeta}.
    \end{align*}

    For each $\beta > 0$, let $\eps_{\rm scal} = \eps_{\rm scal}(\beta)$ be as in Lemma~\ref{Lem:NSAAS}. For each $\rho \in [\delta^{1-\eps_{\rm{scal}}},\delta^{\eps_{\rm{scal}}}]$, we have 
    \begin{align*}
        |\F''_\rho| \geq \delta^{4\eta} \rho^{-2-\zeta} 
        \geq 
        \rho^{-2-\zeta + 4\eta/\eps_{\rm{scal}}}.
    \end{align*}
    Then for any $\zeta >0$, take $\nu = \nu(\beta, \zeta/2)/2 > 0$ and $\eta = \min\{ \eps_{\rm{scal}}\zeta/8, \eta(\beta, \zeta/2)/2 \}>0$, where the binary functions $\nu(\cdot, \cdot)$ and $\eta(\cdot, \cdot)$ are as in Lemma~\ref{Lem:NSAAS}. One can easily verify that $\Delta_{max}(\F'') \leq \delta^{-\eta(\beta, \zeta/2)}$, $\lambda(\F'', Y') \geq \delta^{\eta(\beta, \zeta/2)}$, and $|\F_\rho''| \geq \rho^{-2-\zeta/2}$ for each $\rho \in [\delta^{1-\eps_{\rm{scal}}},\delta^{\eps_{\rm{scal}}}]$. As $(\F'', Y')$ is uniform, by applying Lemma~\ref{Lem:NSAAS} to it, we get 
    \begin{align*}
        \mu(\F, Y) \lessapprox \mu(\F'', Y') \leq 
        \delta^{\nu(\beta, \zeta/2)} |\F''|^\beta \leq 
        \delta^{\nu(\beta, \zeta/2)} |\F|^\beta.
    \end{align*}
    So $\mu(\F, Y) \leq \delta^{\nu}|\F|^\beta$, as desired.  
\end{proof}

Our goal in this subsection is to show that Lemma~\ref{Lem:NSAAS'} implies Main Lemma~\ref{Main_lem_2}. The proof is similar to the proof of Main Lemma~\ref{Main_lem_1}, using Lemma~\ref{Lem:multi-scale_decomposition}(B), Proposition~\ref{Prop:plank_factor_improve}(B), and Theorem~\ref{Thm:2-plane_sticky}(B) instead of their Part (A) counterparts.

Before formally beginning the proof, let us first do some two-ends reduction, which allows us to only consider the case when $\F$ is broad (thus $|\F|$ is large). For this purpose, we need the following analogues of \cite[Definition~7.6~and~7.7]{WZ25volume}:
\begin{definition}\label{Def:direction_broad}
    Let $\eps, \delta>0$ and $C\geq 1$. We say that a multiset $\VV \subset \Gr(n,k)$ is \textit{$\eps$-broad with error $C$ at scales $\geq \delta$}, if for all $V_0 \in \VV$ and all $r\in[\delta,1]$, we have 
    \begin{align*}
        |\{ V \in \VV: \angle(V, V_0) \leq r \}| \leq Cr^\eps|\VV|.
    \end{align*}

    More generally, let $\VV$ be contained in a metric ball\footnote{Recall from Remark~\ref{Rem:major_metric} that the metric is induced by the major angle.} $B \subset \Gr(n,k)$ of radius $\rho$. We say that $\VV$ is \textit{$\eps$-broad with error $C$ at scales $\geq \delta$ inside $B$}, if for all $V_0 \in \VV$ and all $r\in[\delta,\rho]$, we have 
    \begin{align*}
        |\{ V \in \VV: \angle(V, V_0) \leq r \}| \leq C(r/\rho)^\eps|\VV|.
    \end{align*}
\end{definition}

\begin{definition}\label{Def:flakes_broad}
    Let $(\F, Y)$ be a set of $(\delta,k)$-flakes in $B_1 \subset \R^n$ and their associated shading. 
    \begin{enumerate}[label={\rm (\Alph*)}]
    \item We say that $(\F,Y)$ is \textit{$\eps$-broad with error $C$} if for each $x \in U(\F,Y)$, the set $\{\dir(F): x\in Y(F)\}$ is $\eps$-broad with error $C$ at scales $\geq \delta$.

    \item We say that $(\F,Y)$ is \textit{$\eps$-broad with error $C$ relative to $\F_\rho$}\footnote{A set of $(\rho,k)$-flakes with the smallest cardinality with $\F = \bigcup_{F_\rho \in \F_\rho}\F[F_\rho]$.} ($\rho \in [\delta,1]$) if for each $F_\rho \in \F_\rho$, the image of $(\F[F_\rho], Y)$ under the rescaling taking $F_\rho$ to $B_1$, as a set of $(\delta/\rho)$-flakes in $B_1$, is $\eps$-broad with error $C$.
    \end{enumerate}
\end{definition}

The key observation here is that by Remark~\ref{Rem:major_metric} and Lemma~\ref{Lem:flake_contain}, the metric on $\Gr(n,k)$ has the doubling property and is compatible with coarse coverings of flakes, and so all the discussions in \cite[Section~7.1]{WZ25volume} can be carried out in parallel with $\Gr(n,k)$ in place of $\mathbb{S}^{n-1}$ and $\F$ in place of $\T$. In particular, we have the following analogue of \cite[Corollary~7.10]{WZ25volume}:
\begin{lem}\label{Lem:broad_decomposition}
    Let $\eps>0$, and let $(\F,Y)$ be a set of $(\delta,k)$-flakes in $B_1 \subset \R^n$. Then there exists a scale $\rho \in [\delta,1]$ and a $\gtrapprox 1$ refinement $(\F',Y')$ of $(\F,Y)$, so that the following holds:
    \begin{enumerate}[label={\rm (\roman*)}]
        \item 
        Each point $x \in \R^n$ is contained in $\lessapprox \rho^{-\eps}$ sets of the form $U(\F'[F_\rho], Y')$ where $F_\rho \in \F_\rho$.

        \item 
        $\F'$ is $\eps$-broad with error $\approx 1$ relative to $\F'_\rho$.

        \item $\lambda(\F'[F_\rho], Y') \gtrapprox \lambda(\F,Y)$ for all $F_\rho \in \F'_\rho$. 
    \end{enumerate}
\end{lem}
\begin{proof}
    The proof is exactly the same as that of \cite[Corollary~7.10]{WZ25volume} -- only harmless constants will be changed -- and so omitted for simplicity. We only point out Item~(iii) comes from dyadic pigeonholing on the average density inside each $F_\rho$ in the middle of the proof. 
\end{proof}

Now let us turn to the two-ends reduction. By Definition~\ref{Def:Katz_Tao_beta}, to prove $K_{KT}(\beta - \nu)$, we take any set of $\delta$-flakes $(\F, Y)$ in $B_1 \subset \R^n$ with $\Delta_{max}(\F) \leq \delta^{-\eta}$ and $\lambda(\F, Y) \geq \delta^\eta$, and need to find for any $\eps>0$ some $\nu = \nu(\beta) > 0$ such that 
$$\mu(\F, Y) \leq \delta^{-\eps}|\F|^{\beta-\nu}$$ 
holds when $\eta$ and $\delta$ are small enough. Since $\Delta_{max}(\F) \leq \delta^{-\eta}$ implies that $\F$ is essentially distinct up to multiplicity $\lesssim \delta^{-\eta}$, by a harmless pigeonholing argument, we may only consider the case that $\F$ is essentially distinct. Also, we may assume $\eps < 1/4$.

Fix $0 < \eps < 1/4$. By applying Lemma~\ref{Lem:broad_decomposition} to $(\F,Y)$ with $\eps/2$ in place of $\eps$ therein, we get a scale $\rho \in [\delta,1]$ and a $\gtrapprox 1$ refinement $(\F',Y')$ of $(\F,Y)$, so that 
\begin{enumerate}[label={\rm (\roman*)}]
    \item 
    Each point $x \in \R^n$ is contained in $\lessapprox \rho^{-\eps/2}$ sets of the form $U(\F'[F_\rho], Y')$ where $F_\rho \in \F_\rho$.

    \item 
    $(\F',Y')$ is $(\eps/2)$-broad with error $\approx 1$ relative to $\F'_\rho$.

    \item $\lambda(\F'[F_\rho], Y') \gtrapprox \lambda(\F,Y) \geq \delta^\eta$ for all $F_\rho \in \F'_\rho$. 
\end{enumerate}
Note that Item~(i) implies that 
\begin{align*}
    \sum_{F_\rho \in \F'_\rho} |U(\F'[F_\rho], Y')| 
    \lessapprox \rho^{-\eps/2} |U(\F', Y')|,
\end{align*}
and so 
\begin{align*}
    \mu(\F', Y') = \frac{\sum_{F\in\F'}|Y'(F)|}{|U(\F', Y')|} 
    &\lessapprox 
    \frac{\sum_{F_\rho\in \F_\rho'} \sum_{F\in\F'[F_\rho]} |Y'(F)|}{\rho^{\eps/2} \sum_{F_\rho \in \F'_\rho} |U(\F'[F_\rho], Y')|} \\
    & \leq 
    \rho^{-\eps/2} 
    \max_{F_\rho\in\F_\rho'} \frac{\sum_{F\in\F'[F_\rho]} |Y'(F)|}{|U(\F'[F_\rho], Y')|}
    =
    \rho^{-\eps/2}     \max_{F_\rho\in\F_\rho'} \mu(\F'[F_\rho], Y').
\end{align*}

If $\rho \leq \delta^{1-\eps/4n}$, then we have the trivial bound 
\begin{align*}
    \mu(\F', Y') \lessapprox 
    \rho^{-\eps/2}     \max_{F_\rho\in\F_\rho'} |\F'[F_\rho]| \lesssim 
    \rho^{-\eps/2} (\rho/\delta)^{2n-2} 
    \leq 
    \delta^{-\eps/2 - (n-1)\eps/2n},
\end{align*}
which gives $\mu(\F', Y') \leq \delta^{-\eps} \leq \delta^{-\eps}|\F'|^{\beta-\nu}$ (for any $\nu \leq \beta$), as desired. 

If $\rho \geq \delta^{1-\eps/4n}$, then let $F_\rho^* \in \F'_\rho$ maximize $\mu(\F'[F_\rho], Y')$, and let $(\tilde{\F}, \tilde{Y})$ be the image of $(\F'[F_\rho^*], Y')$ under the rescaling taking $F_\rho^*$ to $B_1$. Let $\tilde{\delta} \coloneqq \delta/\rho \leq \delta^{\eps/4n}$. Note that Proposition~\ref{Prop:inherit}(B) yields 
\begin{align*}
    \Delta_{max}(\tilde{\F}) = \Delta_{max}(\F'[F_\rho^*]) \leq \Delta_{max}(\F) \leq \delta^{-\eta} \leq \tilde{\delta}^{-4n\eta/\eps},
\end{align*}
while Item~(iii) yields 
\begin{align*}
\lambda(\tilde{\F}, \tilde{Y}) = \lambda(\F'[F_\rho^*], Y') \gtrapprox \delta^\eta \geq \tilde{\delta}^{4n\eta/\eps}.
\end{align*} 
On the other hand, we claim that Item~(ii) implies $|\tilde{\F}| \geq \tilde{\delta}^{-1/2}$.
\begin{proof}[Proof of the claim]
By pigeonholing, there exists $F_0 \in \tilde{\F}$ with $|\tilde{Y}(F_0)| \gtrapprox \delta^\eta |F_0| \gtrsim \tilde{\delta}^{2+4n\eta/\eps}$. By Definition~\ref{Def:flakes_broad}, Item~(ii) means that $(\tilde{\F}, \tilde{Y})$ is $(\eps/2)$-broad with error $\approx 1$, i.e., for each $x \in U(\tilde{\F}, \tilde{Y})$, the set $\{\dir(F): x\in \tilde{Y}(F)\}$ is $(\eps/2)$-broad with error $\approx 1$ at scales $\geq \tilde{\delta}$. In particular, by Definition~\ref{Def:direction_broad} with $V_0 = \dir(F_0)$, we know that for each $x \in \tilde{Y}(F_0)$, we have 
$$|\{ F \in \tilde{F}: x \in \tilde{Y}(F), \angle(F, F_0) \leq r \}| \lessapprox r^{\eps/2}|\tilde{\F}_{\tilde{Y}}(x)|$$ 
for all $r \in [\tilde{\delta}, 1]$. By taking $r = \tilde{\delta}^\eps$, we see that for each $x \in \tilde{Y}(F_0)$, there exists some $F_x \in \tilde{\F}$ such that $\angle(F_x, F_0) \geq \tilde{\delta}^\eps$, and so $|(10 F_x) \cap F_0|\lesssim \tilde{\delta}^3/\angle(F_x, F_0) \leq \tilde{\delta}^{3-\eps}$ by Lemma~\ref{Lem:flake_volume}. Now we can iteratively pick $F_x$ as follows until we cannot continue any more: Fix any $x_1 \in \tilde{Y}(F_0)$, then for each $j \geq 2$, we pick $x_j \in \tilde{Y}(F_0) \setminus \bigcup_{i=1}^{j-1} (10 F_{x_i})$. Note that $\{F_{x_j}\}_j \subset \tilde{\F}$ are essentially distinct and has cardinality at least $|\tilde{Y}(F_0)|/\tilde{\delta}^{3-\eps} \gtrapprox \tilde{\delta}^{-1 + 4n\eta/\eps + \eps} \geq \tilde{\delta}^{-3/4 + 4n\eta/\eps}$. Thus $|\tilde{\F}| \geq \tilde{\delta}^{-1/2}$ for $\eta \leq \eps/(32n)$. 
\end{proof}
Suppose we know that $K_{KT}(\beta - \nu)$ holds for some $\nu = \nu(\beta) > 0$ and all $\delta$-flakes $(\F,Y)$ that additionally satisfy $|\F| \geq \delta^{-1/2}$, then by taking $\eta \ll \eps$ and applying $K_{KT}(\beta - \nu)$ to $(\tilde{\F}, \tilde{Y})$ (Definition~\ref{Def:Katz_Tao_beta}, with $\tilde{\delta}$ in place of $\delta$ and $\eps/4$ in place of $\eps$), we have 
\begin{align*}    \mu(\F'[F_\rho^*], Y') 
= \mu(\tilde{\F}, \tilde{Y}) 
\leq \tilde{\delta}^{-\eps/4} |\tilde{\F}|^{\beta-\nu} 
= \tilde{\delta}^{-\eps/4} |\F'[F_\rho^*]|^{\beta-\nu} 
\leq \delta^{-\eps/4} |\F'|^{\beta-\nu},
\end{align*}
which implies 
\begin{align*}
    \mu(\F,Y) \lessapprox \mu(\F',Y') \lessapprox \rho^{-\eps/2} \mu(\F'[F_\rho^*], Y') \leq \delta^{-3\eps/4} |\F'|^{\beta-\nu} \leq \delta^{-3\eps/4} |\F|^{\beta-\nu},
\end{align*}
and so $\mu(\F,Y) \leq \delta^{-\eps}|\F|^{\beta-\nu}$, as desired. 

Therefore, from now on, we shall only focus on essentially distinct $\F$ with $|\F| \geq \delta^{-1/2}$, which is the most interesting case and requires using Lemma~\ref{Lem:NSAAS'}.
\begin{rem}\label{Rem:broad_reduction}
    This largeness reduction is technically necessary for the forthcoming proof to be entirely rigorous. If we do not have any additional assumption on $|\F|$, then via similar arguments we can still get estimates of the form     \begin{align}\label{Eq:KP_beta-nu_0}
    \mu(\F, Y) \leq \delta^{-\eps} |\F|^{\beta-\nu}.
    \end{align}
    This is superficially sufficient, but actually not, as $\nu$ depends on $\eps$ if one carefully tracks the dependence of the parameters. This issue also happened in their streamlined proof paper \cite[Section~9.1]{GWZ26streamlined}, but was neglected. However, we borrow the idea of the two-end reduction from Wang and Zahl's original paper \cite[Section 7.1]{WZ25volume} to resolve this subtle problem.
\end{rem}

\begin{proof}[Proof of Main Lemma~\ref{Main_lem_2} using Lemma~\ref{Lem:NSAAS'}]

Let \( \F \) be a set of {\color{black} essentially distinct} \( \delta \)-flakes in \( B_1 \subset \mathbb{R}^n \) with \( \Delta_{max}(\F) \leq \delta^{-\eta} \) and $|\F| \geq \delta^{-1/2}$, and $Y$ be a shading on $\F$ with \( \lambda(\F, Y) \geq \delta^\eta \). {\color{black} Our goal is to show that} if \( \eta \) and \( \delta \) are sufficiently small, then
\begin{align}\label{Eq:KP_beta-nu}
\mu(\F, Y) \leq |\F|^{\beta-\nu}
\end{align}
for some $\nu = \nu(\beta) >0$. Note that if this could be done, then we immediately get $K_{KT}(\beta-\nu)$ for $\F$ with $|\F| \geq \delta^{-1/2}$, as desired.

To begin with, let $(\F', Y)$ be a $\gtrapprox 1$ refinement of $(\F, Y)$ so that there is $\lambda > 0$ such that $|Y(F)|/|F| \sim \lambda$ for all $F \in \F'$. Then $|\F'| \gtrapprox \delta^{\eta} |\F|$, $\Delta_{max}(\F') \leq \Delta_{max}(\F) \leq \delta^{-\eta}$, and $\lambda \sim \lambda(\F',Y) \gtrapprox \lambda(\F,Y) \geq \delta^\eta$.

Let \( \eps_0 \) and \( \delta_0 \) be the output of Theorem~\ref{Thm:2-plane_sticky}(B) with \( \beta/8 \) in place of $\eps$. Let \( \eps_{\rm{scal}} = \eps_{\rm{scal}}(\beta) \) be the quantity from Lemma~\ref{Lem:NSAAS} with \( \beta \) as above. Let \( \eps = \eps(\eps_0, \eps_{\rm{scal}}, \beta) \) be a number to be chosen later (since \( \eps_0 \) and \( \eps_{\rm{scal}} \) depend only on \( \beta \), \( \eps \) depends only on \( \beta \)). Let \(N \geq \max\{ 10(5/\eps_0)^5, 10^5 \}\) be an integer to be chosen later, and define \(\eps = (10/N)^{1/5} \leq \eps_0/5\). Let \( \eta \leq \eta_0 \leq \eta_1 \leq \ldots \leq \eta_N \leq \eps \) be a sequence of numbers to be chosen later. Each \( \eta_{j-1} \) will be selected very small depending on \( \eta_{j}, \eps \), and \( \beta \). The quantity \( \nu \) from the conclusion of Main Lemma~\ref{Main_lem_2} will be selected small compared to \( \eta_0 \).

Let $\eta_{j-1}' \coloneqq \frac{14\eta_{j-1}}{\eps^2\beta}$ for each $1 \leq j \leq N$. 
By taking $\eta \leq \eta_0 \eps$ and $\eta_j \geq \eta_{j-1}'/\eps$ for each $1\leq j \leq N$, we can apply Lemma~\ref{Lem:multi-scale_decomposition}(B) to $\F'$, with $\beta$, $N$, $\eps$ (in place of $\eps$), $\{\eta_j\}_{j=0}^N$ as given. The only subtle point is the choice of selection maps $\Phi_{\F',\rho}$. Fortunately, we can define it in the same way as in the proof of Lemma~\ref{Lem:KT->F_bootstrap}. In particular, for any $\rho \in [\delta,1]$ and $\underline{\F} \subset \F'$, we have  $\Phi_{\F', \rho}(\underline{\F}) \subset \underline{\F}_\rho$, and 
\begin{align*}
\big|\bigsqcup_{F_\rho \in \Phi_{\F', \rho}(\underline{\F})} \underline{\F}[F_\rho]\big| \geq \lambda_* |\underline{\F}|
\end{align*}
for some $\lambda_* \gtrapprox \delta^{2\eta}$. Apply Lemma~\ref{Lem:multi-scale_decomposition}(B) to this $\Phi_{\F', \rho}$ and $\lambda_*$.  

If Conclusion (I) holds, then there exists a $\gtrapprox_{\eps} \lambda_*^{2/\eps^3}$ refinement $\F'' \subset \F'$ that is \( \lessapprox_{\eps} \delta^{-4\eps} \) Katz--Tao (thus $\delta^{-\eps_0}$ Katz--Tao by $\eps \leq \eps_0/5$) at every scale. Also, $\lambda(\F'', Y) \sim \lambda \gtrapprox \delta^\eta$. As $|Y(F)|/|F| \sim \lambda$ for all $F \in \F'$, $(\F'', Y)$ is a $\gtrapprox_{\eps} \lambda_*^{2/\eps^3}$ refinement of $(\F', Y)$ (thus of $(\F, Y)$). Hence by taking $4\eta \leq \eps^3\beta/8$ and applying Theorem~\ref{Thm:2-plane_sticky}(B) to $(\F'', Y)$, we have
\[
\mu(\F, Y) \lessapprox_{\eps} 
\lambda_*^{-2/\eps^3} \mu(\F'', Y) 
\lessapprox
\delta^{-4\eta/\eps^3}
\mu(\F'', Y)
\leq 
\delta^{-\beta/8}\cdot
\delta^{-\beta/8} = \delta^{-\beta/4}, 
\]
which implies 
\begin{align*}
    \mu(\F, Y) \leq \delta^{-\beta/3} \leq |\F|^{2\beta/3}
\end{align*}
by using $|\F| \geq \delta^{-1/2}$.\footnote{This is the only place in the proof that we need the additional assumption that $|\F|$ is large.} 
{\color{black} In particular, \eqref{Eq:KP_beta-nu} holds with \(\nu = \beta/3\).}

Henceforth we shall assume that Conclusion (II) of Lemma~\ref{Lem:multi-scale_decomposition}(B) holds. Let $\F'' \subset \underline{\F} \subset \F'$, $\delta \leq \tau \leq \theta \leq 1$ (with $\tau \leq \delta^{\eps^2}\theta$), \(1 \leq J \leq N\), and $\tau(\theta/\tau)^{\eps} \leq a \leq b \leq \theta$ be the output of that lemma. As in the proof of Lemma~\ref{Lem:KT->F_bootstrap}, by taking $12\eta \leq \eta_0\eps^3$ and using $\underline{\F}$ and $\Phi_{\F',\tau}(\underline{\F})$ as intermediate bridges to handle $\mu(\F, Y)$, we have 
\begin{equation}
\begin{split}
    \mu(\F, Y) \lessapprox_{\eps} \delta^{-\eta_{J-1}} 
    \mu(\F''[F_\tau], Y') \, 
    \mu(\F''_\tau, Y_{\F''_\tau})
\end{split}\label{Eq:mult_Phi_selection_0}
\end{equation}
for all $F_\tau \in \F''_\tau$. Here $Y'$ is a shading on $\F''$, and $Y_{\F''_\tau}$ is a shading on $\F''_\tau$. Moreover, for all $F_\tau \in \F''_\tau$, we have $|Y_{\F''_\tau}(F_\tau)|/|F_\tau| \gtrapprox \delta^\eta$ and $\lambda(\F''[F_\tau], Y') \gtrapprox \delta^\eta$.  

By Proposition~\ref{Prop:inherit}(B), for each \( F_\tau \in \F_\tau'' \), we have \( \Delta_{max}(\F''[F_\tau]) \leq \Delta_{max}(\F) \leq \delta^{-\eta} \), and recall that $\lambda(\F''[F_\tau],Y') \gtrapprox \delta^\eta$. Thus by taking $\eta \ll \eta_0$ and applying \( K_{KT}(\beta) \) to $(\F''[F_\tau],Y')$ (Definition~\ref{Def:Katz_Tao_beta}, with \( \eta_0 \) in place of \( \eps \)), we have
\begin{align}\label{Eq:delta_tau_0}
\mu(\F''[F_\tau], Y') \leq \delta^{-\eta_{J-1}} |\F''[F_\tau]|^{\beta}.
\end{align}

Our goal is to estimate $\mu(\F''_\tau, Y_{\F''_\tau})$ on the RHS of \eqref{Eq:mult_Phi_selection_0}. Properties in conclusion (II) of Lemma~\ref{Lem:multi-scale_decomposition}(B) related to convex Wolff axioms and factoring structure (for $\F''$) would be very helpful. 

{\color{black} Let $\tilde{\delta} \coloneqq \tau/\theta$.} Then $\tilde{\delta} \leq \delta^{\eps^2}$. We set \( \rho_0 \coloneqq \tilde{\delta}^{1-\eps} \).

{\bf The eccentric case.} Let us first consider the case when \[
b \geq \tilde{\delta}^{-\eta_{J-1}'} a,
\]
i.e., Conclusion~(II)(i) holds.

As $\F''$ is uniform at scales $\tau$ and $\theta$, by applying Lemma~\ref{Lem:fact_mult} (Remark~\ref{Rem:factor_only_shading}) to \((\F_\tau'', Y_{\F_\tau''})\) and \(\F_\theta''\), we obtain a \(\gtrapprox 1\) refinement $(\F_\tau'', Y_{\F_\tau''}')$ of \((\F_\tau'', Y_{\F_\tau''})\) and a shading \(Y_{\F_\theta''}\) on \(\F_\theta''\) with $\lambda(\F_\theta'', Y_{\F_\theta''}) \gtrapprox \delta^\eta$, such that 
\begin{align}\label{Eq:theta_sep}
    \mu(\F_\tau'', Y_{\F_\tau''}) \lessapprox     \mu(\F_\tau''[F_\theta], Y_{\F_\tau''}') \, \mu(\F_\theta'', Y_{\F_\theta''})
\end{align}
for all $F_\theta \in \F_\theta''$. 

Since $\Delta_{max}(\F_\theta'') \lessapprox_{\eps} \delta^{-2\eta_{J-1}}$ and $\lambda(\F_\theta'', Y_{\F_\theta''}) \gtrapprox \delta^\eta$, by taking $\eta \ll \eta_0$ and applying Lemma~\ref{Lem:general_KT} to $(\F_\theta'', Y_{\F_\theta}'')$ (again, with \( \eta_0 \) in place of \( \eps \)), we conclude that
\begin{align}\label{Eq:theta_unit_0}
    \mu(\F_\theta'', Y_{\F_\theta}'') \leq \delta^{-3\eta_{J-1}} |\F_\theta''|^\beta.
\end{align}

Therefore, by combining \eqref{Eq:mult_Phi_selection_0}, \eqref{Eq:delta_tau_0}, \eqref{Eq:theta_sep}, and \eqref{Eq:theta_unit_0}, one can easily see that {\color{black} to prove \eqref{Eq:KP_beta-nu} with \( \nu = \eta_0 \)}, it suffices to find some \( F_\theta \in \F_\theta'' \) with 
\begin{align}\label{Eq:tau_theta_0}    \mu(\F_\tau''[F_\theta], Y_{\F_\tau''}') \leq \delta^{10\eta_{J-1}} |\F_\tau''[F_\theta]|^\beta.
\end{align}

As $\lambda(\F_\tau'',Y_{\F_\tau''}') \gtrapprox \lambda(\F_\tau'',Y_{\F_\tau''}) \gtrapprox \delta^\eta$, by pigeonholing, there exists some $F_\theta \in \F_\theta''$ such that $\lambda(\F_\tau''[F_\theta],Y_{\F_\tau''}') \gtrapprox \delta^\eta$. Fix one such $F_\theta$. Let \( (\tilde{\F}, \tilde{Y}) \) be the image of $(\F_\tau''[F_\theta],Y_{\F_\tau''}')$ under the rescaling taking \( F_\theta \) to $B_1$. Hence \( \tilde{\F} \) is a family of \( \tilde{\delta} \)-flakes, and $\lambda(\tilde{\F}, \tilde{Y}) \gtrapprox \delta^\eta$. Recall that \( \tilde{\delta} \leq \delta^{\eps^2} \). Thus to prove \eqref{Eq:tau_theta_0}, we need to show that
\begin{align}\label{Eq:tau_theta_rescaled_0}
    \mu(\tilde{\F}, \tilde{Y}) \leq \tilde{\delta}^{10\eta_{J-1}/\eps^2} |\tilde{\F}|^\beta.
\end{align}

Recall that $\F''$, $\delta \leq \tau \leq \theta \leq 1$, \(1 \leq J \leq N\), and $\tau(\theta/\tau)^{\eps} \leq a \leq b \leq \theta$ are the output of Conclusion (II) of Lemma~\ref{Lem:multi-scale_decomposition}(B). So we have 
\begin{align*}
    \Delta_{max}(\tilde{\F}) \lessapprox_{\eps} \tilde{\delta}^{-2\eta_{J-1}}.
\end{align*}
There exists a family of convex sets $\tilde{\W}$ of dimensions $\tilde{a} \times \tilde{b} \times \1$ ($\tilde{a} = a/\theta$, $\tilde{b} = b/\theta$), such that $\tilde{\W}$ uniformly factors $\tilde{\F}_{\rho_0}$. Note that $\rho_0 \leq \tilde{a} \leq \tilde{b} \leq 1$. 
Moreover, we have $\tilde{a}/ \tilde{b} \leq \tilde{\delta}^{\eta_{J-1}'}$ in the eccentric case.

Since $\lambda(\tilde{\F}, \tilde{Y}) \gtrapprox \delta^\eta \geq \tilde{\delta}^{\eta/\eps^2}$ and $\Delta_{max}(\tilde{\F}) \lessapprox_{\eps} \tilde{\delta}^{-2\eta_{J-1}}$, by taking $\eta/\eps^2 \ll \eta_0$ and applying Proposition~\ref{Prop:plank_factor_improve}(B) to $(\tilde{\F}, \tilde{Y})$ (with $\rho_0$ in place of $\rho$ and \(\eta_0\) in place of \(\eps\)), we have
\[
\mu(\tilde{\F}, \tilde{Y}) \lessapprox_{\eps} \tilde{\delta}^{-3\eta_{J-1}/\eps^2} \left( \frac{a}{b} \right)^{\beta} |\tilde{\F}|^{\beta}
\]
Thus \eqref{Eq:tau_theta_rescaled_0} holds by using \(\tilde{a}/\tilde{b} \leq \tilde{\delta}^{\eta_{J-1}'} 
= \tilde{\delta}^{\frac{14\eta_{J-1}}{\eps^2 \beta}}\).

{\bf The non-eccentric case.} 
Henceforth we shall assume that 
\begin{align*}
    b < \tilde{\delta}^{-\eta_{J-1}'} a,
\end{align*}
i.e., Conclusion~(II)(ii) holds.

As $\F''$ is uniform at scales $\tau$ and $b$, by applying Lemma~\ref{Lem:fact_mult} (Remark~\ref{Rem:factor_only_shading}) to \((\F_\tau'', Y_{\F_\tau''})\) and \(\F_b''\), we obtain a \(\gtrapprox 1\) refinement $(\F_\tau'', Y_{\F_\tau''}')$ of \((\F_\tau'', Y_{\F_\tau''})\) and a shading \(Y_{\F_b''}\) on \(\F_b''\) with $\lambda(\F_b'', Y_{\F_b''}) \gtrapprox \delta^\eta$, such that 
\begin{align}\label{Eq:scale_b_sep}
    \mu(\F_\tau'', Y_{\F_\tau''}) \lessapprox     \mu(\F_\tau''[F_b], Y_{\F_\tau''}') \, \mu(\F_b'', Y_{\F_b''})
\end{align}
for all $F_b \in \F_b''$. 

Recall that $\F''$, $\delta \leq \tau \leq \theta \leq 1$, \(1 \leq J \leq N\), and $\tau(\theta/\tau)^{\eps} \leq a \leq b \leq \theta$ are the output of Conclusion (II) of Lemma~\ref{Lem:multi-scale_decomposition}(B). So we have 
\begin{align*}
    \Delta_{max}(\F_\theta'') \lessapprox_{\eps} \delta^{-2\eta_{J-1}}. 
\end{align*}
For each $F_{\theta} \in \F''_{\theta}$, there exists a family of convex sets $\W_{F_{\theta}}$ of dimensions $a \times b \times \1$ in $F_{\theta}$, such that $\W_{F_{\theta}}$ uniformly factors $\F''_{\tau(\theta/\tau)^\eps}[F_{\theta}]$. Since each flake $F_b \in \F''_b[F_\theta]$ essentially contains at least one $W \in \W_{F_{\theta}}$, we have 
\begin{align}\label{Eq:Delta_max_P_b}
\Delta_{max}(\F''_b[F_\theta]) \lesssim \frac{|F_b|}{|W|} \Delta_{max}(\W_{F_\theta}) \lesssim \frac{b}{a} \leq \tilde{\delta}^{-\eta'_{J-1}},
\end{align}

Since $\F''$ is uniform at scales $b$ and $\theta$, by Lemma~\ref{Lem:KT_multi} (Remark~\ref{Rem:KT_multi_relax}), we have 
\begin{align*}
    \Delta_{max}(\F''_b) 
    \lesssim 
    \Delta_{max}(\F_\theta'') \Delta_{max}(\F''_b[F_\theta])\lessapprox_{\eps} \delta^{-2\eta_{J-1}} \tilde{\delta}^{-\eta_{J-1}'}.
\end{align*}
Recall that $\lambda(\F_b'', Y_{\F_b''}) \gtrapprox \delta^\eta$. By taking $\eta \ll \eta_0$ and applying Lemma~\ref{Lem:general_KT} to $(\F_b'', Y_{\F_b''})$ (with \(\eta_0\) in place of \(\eps\)), we have
\begin{align}\label{Eq:mu_b_coarse}
\mu(\F_b'', Y_{\F_b''}) \leq \delta^{-3\eta_{J-1}} \tilde{\delta}^{-\eta_{J-1}'} |\F_b''|^\beta.
\end{align}

Therefore, by combining \eqref{Eq:mult_Phi_selection_0}, \eqref{Eq:delta_tau_0}, \eqref{Eq:scale_b_sep}, and \eqref{Eq:mu_b_coarse}, one can easily see that {\color{black} to prove \eqref{Eq:KP_beta-nu} with \( \nu = \eta_0 \)}, it suffices to find some \( F_b \in \F_b'' \) with 
\begin{align}\label{Eq:tau_b_0}    \mu(\F_\tau''[F_b], Y_{\F_\tau''}') \leq \delta^{10\eta_{J-1}} \tilde{\delta}^{\eta_{J-1}'} |\F_\tau''[F_b]|^\beta.
\end{align}
We will show that this holds, provided each \( \eta_{j-1} \) is selected sufficiently small compared to \( \eta_j \). 

As $\lambda(\F_\tau'',Y_{\F_\tau''}') \gtrapprox \lambda(\F_\tau'',Y_{\F_\tau''}) \gtrapprox \delta^\eta$, by pigeonholing, there exists some $F_b \in \F_b''$ such that $\lambda(\F_\tau''[F_b],Y_{\F_\tau''}') \gtrapprox \delta^\eta$. Fix one such $F_b$, and fix a $F_\theta \in \F''_\theta$ containing it. Let \( (\tilde{\F}, \tilde{Y}) \) be the image of $(\F_\tau''[F_\theta],Y_{\F_\tau''}')$ under the rescaling taking \( F_\theta \) to $B_1$. Let $\tilde{a} = a/\theta$, $\tilde{b} = b/\theta$, and so $\rho_0 \leq \tilde{a} \leq \tilde{b} \leq 1$. Let $F_{\tilde{b}}$ be the image of $F_b$ under this rescaling. Hence \( \tilde{\F} \) is a family of \( \tilde{\delta} \)-flakes, and $\lambda(\tilde{\F}[F_{\tilde{b}}], \tilde{Y}) \gtrapprox \delta^\eta$. Recall that \( \tilde{\delta} \leq \delta^{\eps^2} \), which means $\tilde{\delta}^{\beta \eta_{J-1}'} \leq \tilde{\delta}^{10\eta_{J-1}/\eps^2} \leq \delta^{10\eta_{J-1}}$. Thus to prove \eqref{Eq:tau_b_0}, we need to show that
\begin{align}\label{Eq:tau_b_rescaled_0}
    \mu(\tilde{\F}[F_{\tilde{b}}], \tilde{Y}) \leq \tilde{\delta}^{(\beta + 1)\eta'_{J-1}} |\tilde{\F}[F_{\tilde{b}}]|^\beta.
\end{align}

Note that in the rescaled setting, \eqref{Eq:Delta_max_P_b} becomes 
\begin{align}\label{Eq:Delta_max_P_b_0}
    \Delta_{max}(\tilde{\F}_{\tilde{b}}) 
    \lesssim 
    \tilde{\delta}^{-\eta'_{J-1}},
\end{align}
while $\Delta_{max}(\F_\tau''[F_\theta]) \lessapprox_{\eps} (\theta/\tau)^{2\eta_{J-1}}$ becomes 
\begin{align}\label{Eq:KT_tilde_P<=}
\Delta_{max}(\tilde{\F}) \lessapprox_{\eps} \tilde{\delta}^{-2\eta_{J-1}}.
\end{align} 
Also, Conclusion~(II)(ii) of Lemma~\ref{Lem:multi-scale_decomposition}(B) tells that 
\begin{align}\label{Eq:KT_tilde_P_rho>=}
\Delta_{max}(\tilde{\F}_\rho) \geq \rho^{-\eta_J} \quad \text{for all } \rho \in [\tilde{\delta}^{1-\eps}, \tilde{\delta}^{\eps}].
\end{align}
Moreover, $\tilde{\F}$ is weakly uniform at scales $[\tilde{\delta}, 1]$ up to a factor $\approx \tilde{\delta}^{-\eta_{J-1}'}$. 

Due to $\tilde{b} \geq \rho_0 = \tilde{\delta}^{1-\eps}$, if we select each \( \eta_{j-1} \) sufficiently small depending on \( \eta_j \) and \( \beta \) (and depending on \( \eps \), which in turn depends on \( \beta \)), then we must have \( \tilde{b} \geq \tilde{\delta}^{\eps}\) -- if not, then \eqref{Eq:Delta_max_P_b_0} would violate \eqref{Eq:KT_tilde_P_rho>=} (with \( \tilde{b} \) in place of \( \rho \)).

Recall that there exists a family of convex sets $\W_{F_{\theta}}$ of dimensions $a \times b \times \1$ in $F_{\theta}$, such that $\W_{F_{\theta}}$ uniformly factors $\F''_{\tau(\theta/\tau)^\eps}[F_{\theta}]$. Pick a $W \in \W_{F_{\theta}}[F_b]$, then we have 
\begin{align*}
    \Delta(\tilde{\F}_{\rho_0}, F_{\tilde{b}}) 
    = \Delta(\F''_{\tau(\theta/\tau)^\eps}[F_{\theta}], F_b) 
    &\geq \frac{|W|}{|F_b|} \Delta(\F''_{\tau(\theta/\tau)^\eps}[F_{\theta}], W) \\
    &\gtrapprox \frac{a}{b} \Delta_{max}(\F''_{\tau(\theta/\tau)^\eps}[F_{\theta}])
    \geq  \tilde{\delta}^{\eta'_{J-1}} \Delta_{max}(\tilde{\F}_{\rho_0}),
\end{align*}
i.e., \( C_F(\tilde{\F}_{\rho_0}[F_{\tilde{b}}]) \lessapprox \tilde{\delta}^{-\eta'_{J-1}} \). Thus by the weak uniformity of $\tilde{\F}$ and Proposition~\ref{Prop:inherit}(A) (Remark~\ref{Rem:inherit}), we have
\begin{align}\label{Eq:C_F_P_sigma}
C_F(\tilde{\F}_\sigma[F_{\tilde{b}}]) 
\lessapprox 
\tilde{\delta}^{-\eta'_{J-1}}
C_F(\tilde{\F}_{\rho_0}[F_{\tilde{b}}])
\lessapprox \tilde{\delta}^{-2\eta'_{J-1}} \quad \text{for all } \sigma \in [\rho_0, \tilde{b}].
\end{align}

On the other hand, using \eqref{Eq:KT_tilde_P_rho>=} (with $\sigma$ in place of $\rho$), the weak uniformity of $\tilde{\F}$, and Lemma~\ref{Lem:KT_multi} (Remark~\ref{Rem:KT_multi_relax}) at scale \(\sigma \in [\tilde{\delta}^{1-\eps}, \tilde{\delta}^{\eps}]\), we have that
\[
\sigma^{-\eta_J} 
\leq 
\Delta_{max}(\tilde{\F}_\sigma) \lessapprox
\tilde{\delta}^{-\eta'_{J-1}} 
\Delta_{max}(\tilde{\F}_\sigma[F_{\tilde{b}}]) \, \Delta_{max}(\tilde{\F}_{\tilde{b}})
\lesssim 
\tilde{\delta}^{-2\eta'_{J-1}}
\Delta_{max}(\tilde{\F}_\sigma[F_{\tilde{b}}]),
\]
i.e.,
\begin{align}\label{Eq:Delta_max_P_sigma}
\Delta_{max}(\tilde{\F}_\sigma[F_{\tilde{b}}]) 
\gtrapprox 
\tilde{\delta}^{2\eta'_{J-1}} \sigma^{-\eta_J}.
\end{align}
Comparing \eqref{Eq:C_F_P_sigma} and \eqref{Eq:Delta_max_P_sigma}, we conclude that
\begin{align}\label{Eq:|P_sigma|}
|\tilde{\F}_\sigma[F_{\tilde{b}}]| \gtrapprox \tilde{\delta}^{4\eta_{J-1}'} (\sigma/{\tilde{b}})^{-2}\sigma^{-\eta_J} \quad \text{for all } \sigma \in [\tilde{\delta}^{1-\eps}, \tilde{\delta}^{\eps}].
\end{align}

We will bound \(\mu(\tilde{\F}[F_{\tilde{b}}], \tilde{Y})\) {\color{black} using Lemma~\ref{Lem:NSAAS'}}, with \(\zeta = \eta_J / 2\). In order to apply Lemma~\ref{Lem:NSAAS'}, we need to first rescale \(F_{\tilde{b}}\) to \(B_1\), converting \( (\tilde{\F}[F_{\tilde{b}}], \tilde{Y}) \) to a new set of flakes \( (\overline{\F}, \Bar{Y}) \) of thickness \(\delta' = \tilde{\delta}/{\tilde{b}} \leq \tilde{\delta}^{1-\eps}\). We need to check that for every \(\sigma' \in [(\delta')^{1-\eps_{\rm{scal}}}, (\delta')^{\eps_{\rm{scal}}}]\), we have 
$$|\overline{\F}_{\sigma'}| \geq (\sigma')^{-2-\zeta}.$$ 
We choose \(\eps \leq \eps_{\rm{scal}}/2\), and $\eta'_{j-1} \leq \eta_j \eps/32$ for each $j$. Then this estimate follows from \eqref{Eq:|P_sigma|}.

Now let \(\eta'\) be the value of \(\eta\) from Lemma~\ref{Lem:NSAAS'}, with \(\beta\) as above and \(\zeta = \eta_J / 2\). By Proposition~\ref{Prop:inherit}(B), we have \( \Delta_{max}(\tilde{\F}[F_{\tilde{b}}]) \leq \Delta_{max}(\tilde{\F}) \lessapprox_{\eps} \tilde{\delta}^{-2\eta_{J-1}} \) by \eqref{Eq:KT_tilde_P<=}. Recall that $\lambda(\tilde{\F}[F_{\tilde{b}}], \tilde{Y}) \gtrapprox \delta^\eta \geq \tilde{\delta}^{\eta/\eps^2}$, and $\tilde{\F}[F_{\tilde{b}}] \subset \tilde{\F}$ is weakly uniform at scales $[\tilde{\delta}, \tilde{b}]$ up to a factor $\approx \tilde{\delta}^{-\eta_{J-1}'}$. To ensure that the hypotheses of Lemma~\ref{Lem:NSAAS'} are satisfied, we select each \(\eta_{j-1}\) sufficiently small compared to \(\eta_j\), \(\eps\), and \(\beta\), so that 
$$\Delta_{max}(\overline{\F}) = \Delta_{max}(\tilde{\F}[F_{\tilde{b}}]) \leq (\delta')^{-\eta'}, \quad \lambda(\overline{\F}, \Bar{Y}) = \lambda(\tilde{\F}[F_{\tilde{b}}], \tilde{Y}) \geq (\delta')^{\eta'},$$ 
and $\Bar{\F}$ is weakly uniform at scales $[\delta', 1]$ up to a factor $\leq (\delta')^{-\eta'}$. Applying Lemma~\ref{Lem:NSAAS'} to $(\overline{\F}, \Bar{Y})$ with $\delta'$ in place of $\delta$ and then changing coordinates back to $(\tilde{\F}[F_{\tilde{b}}], \tilde{Y})$, we conclude that there exists \(\nu(\beta, \eta_J) > 0\) so that
\begin{align}\label{Eq:mu_b_fine}
\mu(\tilde{\F}[F_{\tilde{b}}], \tilde{Y}) = \mu(\overline{\F}, \Bar{Y}) \leq (\delta')^{\nu(\beta, \eta_J)} |\overline{\F}|^\beta \leq \tilde{\delta}^{(1-\eps)\nu(\beta, \eta_J)} |\tilde{\F}[F_{\tilde{b}}]|^{\beta}. 
\end{align}

This implies our desired bound \eqref{Eq:tau_b_rescaled_0}, provided we select each \(\eta_{j-1}\) sufficiently small (depending on \(\eta_j\), \(\eps\), and \(\beta\)) so that \( (\beta + 1)\eta'_{j-1} \leq (1-\eps)\nu(\beta, \eta_j) \).
\end{proof}

\subsection{Maximal density factoring revisited.}
It remains to prove Lemma~\ref{Lem:NSAAS}. For technical reasons, we will need the following variation of Lemma~\ref{Lem:factoring} with biased choices:
\begin{lem}[{\cite[Lemma~9.2]{GWZ26streamlined}}, Biased maximal density factoring lemma]\label{Lem:biased_factoring}
Let \( \mathbb{V} \) be a finite family of convex subsets of \( \mathbb{R}^n \) contained in a convex set \( U \), and suppose \( \eta_{\rm{bias}} > 0 \). Then there is a set \( \mathbb{V}' \subset \mathbb{V} \) with \( \sum_{\mathbb{V}'} |V| \gtrapprox \sum_{\mathbb{V}} |V| \); a family \( \mathbb{W} \) of convex sets in \( \mathbb{R}^n \), each of approximately the same dimensions; and a partition
\[
\mathbb{V}' = \bigsqcup_{W \in \mathbb{W}} \mathbb{V}_W,
\]
where each set in \( \mathbb{V}_W \) is contained in \( W \). This partition has the following properties:
\begin{enumerate}[label={\rm (\arabic*)}]
    \item If \( K \subset W \in \mathbb{W} \),
    \begin{align}\label{Eq:(F)_V_W}
        \Delta(\mathbb{V}_W, K) \lesssim \left( \frac{|K|}{|W|} \right)^{\eta_{\rm{bias}}} \Delta(\mathbb{V}_W, W).
    \end{align}
    In particular, \( C_F(\mathbb{V}_W, W) \lesssim 1 \).
    \item For each \( W \in \mathbb{W} \),
    \[
    \Delta(\mathbb{V}_W, W) \gtrsim \left( \frac{|W|}{|U|} \right)^{\eta_{\rm{bias}}} \Delta_{max}(\mathbb{V}').
    \]
    \item Finally, 
    \begin{align*}
        \Delta_{max}(\mathbb{W}) \lesssim \left( \frac{|W|}{|U|} \right)^{-\eta_{\rm{bias}}}.
    \end{align*}
\end{enumerate}
\end{lem}

One may consult \cite[Section~9.2]{GWZ26streamlined} for the motivation behind Lemma~\ref{Lem:biased_factoring}. We will call the decomposition \( \mathbb{V}' = \bigsqcup_{W} \mathbb{V}_W \) a \textit{maximal density factoring with bias \( \eta_{\rm bias} \)}. Observe that we can still apply Proposition~\ref{Prop:fact_mult} to the factoring \( \bigsqcup_{W} \mathbb{V}_W \)—Proposition~\ref{Prop:fact_mult} requires that \( \mathbb{V}_W \) be \( C \)-Frostman in \( W \), and this remains the case with \( C \sim 1 \).

\subsection{Setup for the proof of Lemma~\ref{Lem:NSAAS}}\label{subsec:Setup_NSAAS}
To prove Lemma~\ref{Lem:NSAAS}, we will study how \( \F \) behaves in smaller balls at various scales. {\color{black} The proof is similar to that in \cite[Section~9]{GWZ26streamlined} and does not contain new insights, except for some modifications to handle flakes instead of tubes.} For simplicity, we only provide a sketch of the proof by omitting details about technical shading manipulations. 

Our goal is to prove that
\begin{align}\label{Eq:KT_Goal}
\mu(\F, Y) \lessapprox \delta^\nu |\F|^\beta,
\end{align}
for some \( \nu = \nu(\beta, \zeta) > 0 \). For this purpose, it also suffices to show that 
\begin{align}\label{Eq:KT_Goal'}
|U(\F, Y)| \gtrapprox \delta^{-\nu} |F| |\F|^{1-\beta},
\end{align}
for some \( \nu = \nu(\beta, \zeta) > 0 \).

The argument has several cases. In some cases, we will prove \eqref{Eq:KT_Goal} and in some cases we will prove \eqref{Eq:KT_Goal'}. In some of the cases, we will prove \eqref{Eq:KT_Goal'} by showing that \( U(\F, Y) \) fills a larger than expected fraction of a typical \( a \)-ball that intersects \( U(\F, Y) \), and we will pair this with an estimate on \( U(\F_a) \).

We have to do some pigeonholing to prepare for our argument. First we refine \( (\F, Y) \) to be $(\F,Y')$ so that for each \( x \in U(\F, Y') \), the set \( \F_{Y'}(x) \) is uniform. We identify \( \F_{Y'}(x) \) with a subset of {\color{black} \( \Gr(n,n-2) \)} corresponding to the direction subspaces of the flakes. Then for each \( x \), we can choose a uniform subset $U_x$ of \( \F_{Y'}(x) \) with \( |U_x| \gtrapprox |\F_{Y'}(x)| \). We define \( Y' \) so that \( x \in Y'(F) \) if and only if \( F \in U_x \). This is a \( \gtrapprox 1 \) refinement. Because of the uniformity, there is a function \( \mu(\rho) \) so that for any \( F_0 \in \F \) and \( x \in Y'(F_0) \),
\begin{align}\label{Eq:major_angle_uniform}
\#\{F \in \F_{Y'}(x) : \angle(F, F_0) \leq \rho\} \approx \mu(\rho),
\end{align}

In particular, we have \( \mu(\F, Y) \lessapprox \mu(\F, Y') \approx \mu(1) \). Moreover, for each \( F_\rho \in \F_\rho \), we have
\[
\mu(\F[F_\rho], Y') \approx \mu(\rho).
\]
Abusing notation, we will refer to \( Y' \) as \( Y \) from here on.

Next we apply Lemma~\ref{Lem:fact_mult} to \( (\F, Y) \) and \( \F_a \), for \( a = \delta^{\eta j}, j = 1, \ldots, \lfloor\eta^{-1}\rfloor \). We obtain a \( \approx 1 \) refinement of \( (\F, Y) \) (abusing notation, we will continue to refer to this as \( (\F, Y) \)) so that for each \( x \in U(\F_a, Y_{\F_a}) \), we have
\begin{align}\label{Eq:a_density}
|U(\F, Y)| \gtrapprox |U(\F_a, Y_{\F_a})| \frac{|U(\F, Y) \cap B(x, a)|}{|B(x, a)|}.
\end{align}

Recall that \( \lambda(\F_a, Y_{\F_a}) \gtrapprox \lambda(\F, Y) \geq \delta^\eta \). We also have
\[
\Delta_{max}(\F_a) \leq \delta^{-\eta} \frac{(a/\delta)^2}{|\F[F_a]|}.
\]
Applying Lemma~\ref{Lem:general_KT} to $(\F_a, Y_{\F_a})$ with \( \nu/10 \) in place of \( \eps \), we have (provided \( \eta > 0 \) is chosen sufficiently small)
\begin{align}\label{Eq:a_coarse}
|U(\F_a, Y_{\F_a})| \geq \delta^{\nu/9} \left( |\F_a| |F_a| \right) \Big/ \left( \Big( \frac{(a/\delta)^2}{|\F[F_a]|} \Big)^{1-\beta} |\F_a|^{\beta} \right)
\approx
\delta^{\nu/9} \left( \frac{a}{\delta} \right)^{2\beta} |F| |\F|^{1-\beta}.
\end{align}
Chaining \eqref{Eq:a_density} and \eqref{Eq:a_coarse}, we conclude that in order to obtain our desired estimate \eqref{Eq:KT_Goal'}, it suffices to prove that there exists an \( a \)-ball \( B_a \) with
\begin{align}\label{Eq:a_fine}
\frac{|U(\F) \cap B_a|}{|B_a|} \gtrapprox \delta^{-2\nu} \Big( \frac{\delta}{a} \Big)^{2\beta}.
\end{align}

With these preparations, we are ready to set up the main objects in the proof. Let \( r_1 = \delta^{\eps_{\rm scal}} \). Divide \( U(\F_{r_1}, Y) \) into boundedly overlapping balls \( B \) of radius \( r_1 \).

For each ball \( B \), we define \( \F_B \) to be a set of essentially distinct \( \delta \times \delta \times r_1\1 \) flake segments contained in \( B \) of the form \( F \cap B \) with \( F \in \F \). For each \( F_B \in \F_B \), we let \( \F(F_B) \) be the set of flakes \( F \in \F \) so that \( F \cap B \) is comparable to \( F_B \).

First, we define a shading on \( \F_B \). For each \( F_B \in \F_B \), we consider the multiplicity \( |Y_{\F(F_B)}(x)| \). We choose \( Y_B(F_B) \subset F_B \) by pigeonholing the multiplicity \( |Y_{\F(F_B)}(x)| \). Next, we refine \( Y \) so that if \( F \in \F(F_B) \), \( Y(F) \cap F_B \subset Y_B(F_B) \). This is a \( \gtrapprox 1 \) refinement. Finally, we refine \( (\F_B, Y_B) \) by pigeonholing \( |Y_B(F_B)| \). Again we refine \( Y \) so that \( Y(F) \cap F_B \subset Y_B(F_B) \) and this is still a \( \gtrapprox 1 \) refinement.

We let \( \mathbb{B} \) be the set of balls \( B \) with \( \lambda(\F_B, Y_B) \gtrapprox \lambda(\F, Y) \). Then we refine \( Y \) by removing the parts of \( Y(F) \) that are not in the balls of \( \mathbb{B} \). Again this is a \( \gtrapprox 1 \) refinement.

In the sequel, we will use the following properties of the shading \( (\F_B, Y_B) \):
\begin{itemize}
    \item \( \lambda(\F_B, Y_B) \gtrapprox \delta^\eta \) for each \( B \in \mathbb{B} \).
    \item \( |Y_B(F_B)| \) is roughly the same for all \( F_B \in \F_B \).
    \item If \( Y(F) \cap B \) is non-empty, then \( F \in \F(F_B) \) for some \( F_B \in \F_B \) and \( Y(F) \cap B \subset Y_B(F_B) \).
    \item The multiplicity \( |\F(F_B)_Y(x)| \) is roughly constant for \( x \in Y_B(F_B) \).
\end{itemize}

Now we apply the biased maximal density factoring, Lemma~\ref{Lem:biased_factoring}, to \( \F_B \), and we let \( \mathbb{W}_B \) be the resulting family of convex sets (the maximal density factoring replaces \( \F_B \) by a large subset, but we will abuse notation and continue to refer to this set as \( \F_B \); this is harmless because of all the uniformity in the bullet points above). By pigeonholing, we can suppose that for each \( B \), the convex sets \( W \in \mathbb{W}_B \) have dimensions \( a \times b \times r_1\1 \), with \( a \leq b \leq r_1 \). Thus for each ball \( B \) we obtain a decomposition
\[
\F_B = \bigsqcup_{W \in \mathbb{W}_B} \F_{B,W}.
\]

The way we proceed depends on the dimensions of \( W \). We let \( \tau > 0 \) be a new small parameter, to be chosen later. If \( a \geq \delta^{1-\tau} \), then we are in the \emph{thick case}. In this case we will show that \( U(\F, Y) \) is large by proving \eqref{Eq:a_fine}. If \( a \leq \delta^{1-\tau} \), then we are in the \emph{thin case}. This case divides into several subcases. In the most important subcase, we will estimate \( \mu(\mathbb{W}_B) \) by an induction on scales. We deal with these cases in the next two sections.

\subsection{The thick case}\label{subsec:thick_case}
In this section, we prove Lemma~\ref{Lem:NSAAS} in the thick case. Recall that we applied the maximal density factoring lemma with bias \(\eta_{\text{bias}}\). Thus by \eqref{Eq:(F)_V_W}, for any \(W \in \mathbb{W}_B\) and any convex set \(K \subset W\),
\begin{align}\label{Eq:(F)_V_W'}
\frac{\Delta(\F_{B,W}, W)}{\Delta(\F_{B,W}, K)} \gtrsim \left( \frac{|W|}{|K|} \right)^{\eta_{\text{bias}}}.
\end{align}

In particular, taking \(K = F_B\), we have
\[
\Delta(\F_{B,W}, W) \gtrsim \left( \frac{|W|}{|F_B|} \right)^{\eta_{\text{bias}}} \Delta(\F_B, F_B) = \left( \frac{|W|}{|F_B|} \right)^{\eta_{\text{bias}}},
\]
and so
\begin{align}\label{Eq:P_B,W_>=}
|\F_{B,W}| \gtrsim \left( \frac{|W|}{|F_B|} \right)^{1 + \eta_{\text{bias}}}.
\end{align}

In the thick case, we will prove \eqref{Eq:a_fine} for balls of radius \(a\):
\begin{align*}
\frac{|U(\F) \cap B_a|}{|B_a|} \gtrapprox \delta^{-2\nu} \Big( \frac{\delta}{a} \Big)^{2\beta}.
\end{align*}
This inequality says that \(U(\F, Y)\) fills a surprisingly large fraction of a typical \(a\)-ball, and it implies our goal as explained in Section~\ref{subsec:Setup_NSAAS}. For the remainder of Section~\ref{subsec:thick_case} we will fix an \(r_1\)-ball \(B\) for which \(\lambda(\F_B, Y_B) \gtrapprox \delta^\eta\). Since each \(W \in \mathbb{W}_B\) has thickness \(a\), in order to establish \eqref{Eq:a_fine} it suffices to prove that there exists at least one \(W \in \mathbb{W}_B\) for which
\begin{align}\label{Eq:W_fine}
\frac{|U(\F_B, Y_B) \cap W|}{|W|} \gtrapprox \delta^{-2\nu} \Big( \frac{\delta}{a} \Big)^{2\beta}.
\end{align}
We fix a choice of \(W\) for which \(\lambda(\F_{B,W}, Y_B) \gtrapprox \delta^\eta\). Recall that since we are in the thick case, \(W\) has dimensions \(a \times b \times r_1\1 \), with \(\delta^{1-\tau} \leq a \leq b \leq r_1\).

We need to prove \eqref{Eq:W_fine}. To do so, we make a linear change of variables to convert \(W\) to \(B_1\). This converts \(\F_{B,W}\) to a set of planks \(\LL\) in \(B_1\). Let \(Y_{\LL}\) denote the image of the shading \(Y_B\) under this transformation. We have
\[
\frac{|U(\F_{B,W}, Y_B)|}{|W|} \sim |U(\LL, Y_{\LL})|,
\]
so it suffices to prove
\begin{align}\label{Eq:rescaled_fine}
|U(\LL, Y_{\LL})| \gtrapprox \delta^{-2\nu} \Big( \frac{\delta}{a} \Big)^{2\beta}.
\end{align}

If \(W\) has dimensions \(a \times b \times  r_1\1 \), then a plank \(L \in \LL\) has dimensions
$\frac{\delta}{b} \times \frac{\delta}{a} \times \1$.
Note that 
$C_F(\LL, B_1) \sim C_F(\F_{B,W}, W) \lesssim 1$ 
by \eqref{Eq:(F)_V_W}, and 
$\lambda(\LL, Y_{\LL}) \sim \lambda(\F_{B,W}, Y_B) \gtrapprox \delta^{\eta}$.

Now we apply Lemma~\ref{Lem:plank_Frostman}. Before doing so, we must estimate the quantity “\(M\)” from that lemma. For any \(\theta \in [a/b, 1]\), let \(L \in \LL\), and let \(L_\theta\) be the corresponding \(\theta \frac{\delta}{a} \times \frac{\delta}{a} \times \1\) thickened plank, as described in Definition~\ref{Def:L_theta_S}. From the statement of Lemma~\ref{Lem:plank_Frostman} we can choose any \(M \geq 1\) with \(M \geq |\LL[L_\theta]|/\theta\) for each \(\theta \in [a/b, 1]\). So we need to estimate \(|\LL[L_\theta]|\).

Recall that \(\LL\) is related to \(\F_{B,W}\) by a linear change of variables. When we undo that change of variables, the planks \(\LL\) become \(\F_{B,W}\) and each plank \(L\) becomes a \(\delta \times \delta \times r\1\) flake segment. The thickened plank \(L_\theta\) becomes a rectangular solid \(E_\theta\) with dimensions \(\delta \times b' \times r\1\) with \(\delta \leq b' \leq b\) (the exact value of \(b'\) depends on $\theta$ and the orientation of the plank \(L\)).

It is here that we make crucial use of the modified maximal density factoring lemma. By \eqref{Eq:(F)_V_W'}, we have
\[
\frac{\Delta(\F_{B,W}, W)}{\Delta(\F_{B,W}, E_\theta)} \gtrsim \left( \frac{|W|}{|E_\theta|} \right)^{\eta_{\text{bias}}}.
\]

Now we are ready to estimate \(|\LL[L_\theta]|\). Note that \(|\LL[L_\theta]| = |\F_{B,W}[E_\theta]|\) and \(|\LL| = |\F_{B,W}|\), so
\begin{align*}
\frac{|\LL[L_\theta]|}{|\LL|} = \frac{|\F_{B,W}[E_\theta]|}{|\F_{B,W}|}
&= \frac{\Delta(\F_{B,W}, E_\theta)}{\Delta(\F_{B,W}, W)} \cdot \frac{|E_\theta|}{|W|}\\
&\lesssim \left( \frac{|E_\theta|}{|W|} \right)^{1+\eta_{\text{bias}}}
= |L_\theta|^{1+\eta_{\text{bias}}}
= \left( \theta \frac{\delta^2}{a^2} \right)^{1+\eta_{\text{bias}}}
\leq \theta \left( \frac{\delta}{a} \right)^{2(1+\eta_{\text{bias}})}.
\end{align*}
Rearranging gives
\[
\frac{|\LL[L_\theta]|}{\theta} \lesssim \Big( \frac{\delta}{a} \Big)^{2(1+\eta_{\text{bias}})} |\LL|.
\]

We need to ensure \((\delta/a)^{2(1+\eta_{\text{bias}})} |\LL| \gtrsim 1\), i.e., \(|\LL| \gtrsim (a/\delta)^{2(1+\eta_{\text{bias}})}\). This follows from \eqref{Eq:P_B,W_>=}, which says
\[
|\LL| = |\F_{B,W}| \gtrsim \left( \frac{|W|}{|F_B|} \right)^{1+\eta_{\text{bias}}}
\sim 
\left( \frac{a \times b \times r_1\1}{\delta \times \delta \times r_1\1} \right)^{1+\eta_{\text{bias}}}
\geq \left( \frac{a}{\delta} \right)^{2(1+\eta_{\text{bias}})}.
\]
Thus we can apply Lemma~\ref{Lem:plank_Frostman} with \(M \sim (\delta/a)^{2(1+\eta_{\text{bias}})} |\LL|\) and $M \geq 1$. We conclude that
\[
|U(\LL, Y_{\LL})| \gtrsim \delta^{\eps} \Big( \frac{\delta}{a} \Big)^{2\beta} \left( M^{-1} \Big( \frac{\delta}{a} \Big)^2 |\LL| \right)^{\beta/{\color{black}n}}
\sim \delta^{\eps} \Big( \frac{\delta}{a} \Big)^{2\beta - 2\eta_{\text{bias}}\beta/n}.
\]
Since \(a \geq \delta^{1-\tau}\), if we select \(\eps < \eta_{\text{bias}}\beta\tau/n\), then this yields \eqref{Eq:rescaled_fine} with \(\nu = \eta_{\text{bias}}\beta\tau/(2n)\). This finishes the proof in the thick case.

\begin{rem}
Notice that in the thick case, the exponent \(\nu\) obeys
\[
\nu \sim \eta_{\text{bias}}\beta\tau.
\]
This value of \(\nu\) does not directly depend on \(\zeta\). But in the thin case, we will have to choose \(\eta_{\text{bias}}\) and \(\tau\) much smaller than \(\zeta\), so indirectly \(\nu\) depends on \(\zeta\).
\end{rem}

\subsection{The thin case}
Now we consider the thin case when \( a \in [\delta, \delta^{1-\tau}] \). Fix a ball \( B \). Recall that each set \( W \in \mathbb{W}_B \) has dimensions \( a \times b \times r_1\1\), with \( r_1 = \delta^{\eps_{\text{scal}}} \). We also have that each flake segment \( F_B \in \F_B \) has roughly the same density of shading, and \( |Y_B(F_B)| \gtrapprox \lambda(\F_B, Y_B)|F_B| \gtrapprox \delta^\eta |F_B| \) for each \( F_B \in \F_B \). After a harmless refinement we can suppose that \( Y_B(F_B) \) is covered by a union of \( \delta \)-balls, each of which satisfies
\begin{align}\label{Eq:delta_dense}
|B_\delta \cap Y_B(F_B)| \gtrapprox \delta^\eta |B_\delta|.
\end{align}

Apply Proposition~\ref{Prop:fact_mult} to \( \F_B = \bigsqcup_{\W_B} \F_{B,W} \). We obtain a shading \( Y'_B(F_B) \subset Y_B(F_B) \) on \( \F_B \) so that \( (\F_B, Y'_B) \) is an \( \approx 1 \) refinement of \( (\F_B, Y_B) \)\footnote{Abusing notation, we still write $\F_B$ for $\F_B'$, and $\F_{B,W}$ for $\F_{B,W}'$.}; a subset \( \mathbb{W}'_B \) of \( \mathbb{W}_B \); and a shading \( Y_{\mathbb{W}'_B} \) on \( \mathbb{W}_B \) so that
\begin{align}\label{Eq:W_B_dense}
\lambda(\mathbb{W}'_B, Y_{\mathbb{W}'_B}) \gtrapprox \delta^{2\eta},
\end{align}
and
\begin{align}\label{Eq:pointwise_cover}
(F_B)_{Y'_B}(x) \subset \bigcup_{W \in \mathbb{W}'_B,\; x \in Y_{\mathbb{W}'_B}(W)} (\F_{B,W})_{Y'_B}(x).
\end{align}

Recall that \( Y_{\mathbb{W}_B'}(W) \) is essentially the \( a \)-neighborhood of \( U(\F_{B,W}, Y') \). Combining this with \eqref{Eq:delta_dense} we have that
\begin{align}\label{Eq:thin->W_B}
|U(\F_B, Y'_B)| \gtrapprox (\delta/a)^{\color{black} n} \delta^\eta |U(\mathbb{W}'_B, Y_{\mathbb{W}'_B})| \geq
\delta^{{\color{black} n}\tau+\eta} |U(\mathbb{W}'_B, Y_{\mathbb{W}'_B})|.
\end{align}

\subsubsection{Slab case.}
Recall that \( r_1 = \delta^{\eps_{\text{scal}}} \), where \( \eps_{\text{scal}} > 0 \) is the constant in Lemma~\ref{Lem:NSAAS}, which is a tiny constant that we get to choose. We define the slab case to be the case \( b \geq \delta^{\eps_{\text{scal}}} r_1 = \delta^{2\eps_{\text{scal}}} \). In this case, \( W \) is a giant slab of dimensions at least \( a \times \delta^{2\eps_{\text{scal}}} \times \delta^{\eps_{\text{scal}}} \times \delta^{\eps_{\text{scal}}}\).

Applying Lemma~\ref{Lem:slab_2} to \( (\mathbb{W}'_B, Y_{\mathbb{W}'_B}) \), we have \( \mu(\mathbb{W}'_B, Y_{\mathbb{W}'_B}) \lesssim \delta^{-C_n'(\eta+\eps_{\text{scal}}+\eta_{\rm bias})} \) for some constant $C_n'$ depending only on $n$, i.e.,
\[
|U(\mathbb{W}'_B, Y_{\mathbb{W}'_B})| \gtrsim \delta^{C_n'(\eta+\eps_{\text{scal}}+\eta_{\rm bias})} \delta^{2\eta} |\mathbb{W}'_B| |W|.
\]
Combining this with \eqref{Eq:thin->W_B}, we have
\[
|U(\F_B, Y'_B)| \gtrapprox \delta^{C_n'(\eta+\eps_{\text{scal}}+\eta_{\rm bias})} \delta^{n\tau+3\eta} |\mathbb{W}'_B| |W|.
\]

Next we estimate \( \Delta_{max}(\F_B) \). Each flake segment \( F_B \in \F_B \) is contained in a flake \( F \in \F \). If \( F_B \subset K \) for some convex set $K$, then the corresponding flake \( F \) is contained in \( r_1^{-1}K \), the dilation of \( K \) by a factor \( r_1^{-1} \). Therefore \( \Delta_{max}(\F_B) \lesssim r_1^{-2} \Delta_{max}(\F) \leq r_1^{-2}\delta^{-\eta}\). Using this, we can see that
\[
|\mathbb{W}'_B| |W| \gtrapprox
\Delta_{max}(\F_B)^{-1} |\F_B| |F_B| \gtrsim r_1^2 \, \delta^{\eta} |\F_B| |F_B|.
\]
Summing over the \( r_1 \)-balls \( B \), we conclude that
\[
|U(\F, Y)| \gtrapprox \delta^{C_n'(\eta+\eps_{\text{scal}}+\eta_{\rm bias})} \delta^{n\tau+3\eta} r_1^2 \delta^{2\eta}  |\F| \, |F| 
\geq \delta^{(C_n'+5)(\eta+\eps_{\text{scal}}+\eta_{\rm bias}) + n\tau} |\F| \, |F|.
\]
We can choose \( \tau, \eta, {\color{black} \eta_{\rm bias}} \), and \( \eps_{\text{scal}}\) to be very small compared to $\beta$, and so this estimate directly implies our goal.

This finishes the case \( b \geq \delta^{\eps_{\text{scal}}} r_1 \), and so from now on, we can assume that \( b \leq \delta^{\eps_{\text{scal}}} r_1 \).

\subsubsection{Transverse case.}
Fix a ball \( B \). We will study the typical angle of intersection of the planks of \( \mathbb{W}'_B \) by applying Lemma~\ref{Lem:typical_angle} to \( (\mathbb{W}'_B, Y_{\mathbb{W}'_B}) \). This gives an \( \approx 1 \) refinement \( (\mathbb{W}''_B, Y_{\mathbb{W}''_B}) \) and a typical angle of intersection \( \theta \in [a/b, 1] \). Let \( \tau' \) be a small parameter to be determined, with \( \zeta \gg \tau' \gg \tau \). We define the transverse case to be the case that \( \theta \geq \delta^{-\tau'} a/b \).

When \( \theta \) is the typical angle of intersection, Lemma~\ref{Lem:plank_thicken} tells us that \( U(\mathbb{W}''_B, Y_{\mathbb{W}''_B}) \) morally almost fills each ball of radius \( \theta b \). More precisely, Lemma~\ref{Lem:plank_thicken} says that \( (\mathbb{W}''_B, Y_{\mathbb{W}''_B}) \) has a \( \gtrapprox 1 \) refinement \( (\mathbb{W}'''_B, Y_{\mathbb{W}'''_B}) \) so that for every ball \( B_{\theta b} \) that intersects \( U(\mathbb{W}'''_B, Y_{\mathbb{W}'''_B}) \), we have
\[
|U(\mathbb{W}'_B, Y_{\mathbb{W}'_B}) \cap B_{\theta b}| \geq |U(\mathbb{W}'''_B, Y_{\mathbb{W}'''_B}) \cap B_{\theta b}| \gtrapprox
a^{4\eta} |B_{\theta b}| \geq \delta^{4\eta} |B_{\theta b}|.
\]
(Here \( \theta \) is the typical angle in the paragraph above -- see Remark~\ref{Rem:angle->thicken}.)

Recall that \( Y_{\mathbb{W}'_B} \) consists of a union of \( a \)-balls that intersect \( U(\F, Y) \). Therefore, using \eqref{Eq:delta_dense}, for any ball \( B_r \) of radius \( r \geq a \), we have
\[
|U(\F, Y) \cap B_r| \gtrapprox (\delta/a)^n \delta^\eta |U(\mathbb{W}'_B, Y_{\mathbb{W}'_B}) \cap B_r| \geq \delta^{n\tau} \delta^\eta |U(\mathbb{W}'_B, Y_{\mathbb{W}'_B}) \cap B_r|,
\]
where we used the assumption \( a \leq \delta^{1-\tau} \) in the thin case.

In particular, for the ball \( B_{\theta b} \) above, we have
\[
|U(\F, Y) \cap B_{\theta b}| \gtrapprox \delta^{n\tau+5\eta} |B_{\theta b}|.
\]

Informally, this shows that \( U(\F, Y) \) almost fills \( B_{\theta b} \). Formally, we choose a radius \( r \) of the form \( r = \delta^{\eta j} \) with \( \theta b \leq r \leq \delta^{-\eta} \theta b \). Since \( \theta b \geq \delta^{1-\tau'} \), we have
\[
\frac{|U(\F, Y) \cap B_r|}{|B_r|} \gtrapprox \delta^{n\tau+(n+5)\eta} \geq \delta^{-2\beta\tau' + n\tau+(n+5)\eta} \Big( \frac{\delta}{r} \Big)^{2\beta}.
\]
Since we choose \( \tau' \gg \tau, \eta \), this gives \eqref{Eq:a_fine} with a value of \( \nu = \beta \tau'/2 \). Equation \eqref{Eq:a_fine} implies our conclusion, as explained in Section~\ref{subsec:Setup_NSAAS}.

This finishes the case \( \theta \geq \delta^{-\tau'} a/b \), and so from now on, we can assume that \( \theta \leq \delta^{-\tau'} a/b \).

\subsubsection{Tangential case}
Recall from the start of the transverse case that we fixed a ball \( B \) and studied the typical angle of intersection of the planks of \( \mathbb{W}'_B \) by applying Lemma~\ref{Lem:typical_angle} to \( (\mathbb{W}'_B, Y_{\mathbb{W}'_B}) \). This gave an \( \approx 1 \) refinement \( (\mathbb{W}''_B, Y_{\mathbb{W}''_B}) \) with typical angle of intersection \( \theta \in [a/b, 1] \). Let \( \tau' \) be a small parameter to be determined, with \( \zeta \gg \tau' \gg \tau \). We defined the transverse case when \( \theta \geq \delta^{-\tau'} a/b \). Now we consider the tangential case when \( \theta < \delta^{-\tau'} a/b \).

Recall from \eqref{Eq:pointwise_cover} that
\[
(F_B)_{Y'_B}(x) \subset \bigcup_{W \in \mathbb{W}'_B,\; x \in Y_{\mathbb{W}'_B}(W)} (\F_{B,W})_{Y'_B}(x),
\]

Next we refine \( (\F_B, Y'_B) \) to make it compatible with \( (\mathbb{W}''_B, Y_{\mathbb{W}''_B}) \). Suppose that \( F_B \in \F_{B,W} \). If \( W \in \mathbb{W}''_B \), we set \( Y''_B(F_B) = Y'_B(F_B) \cap Y_{\W_B''}(W) \), and if \( W \notin \mathbb{W}''_B \), we set \( Y''_B(F_B) \) to be empty. With this setup, it follows that
\[
(\F_B)_{Y''_B}(x) \subset \bigcup_{W \in \mathbb{W}''_B,\; x \in Y_{\mathbb{W}''_B}(W)} (\F_{B,W})_{Y''_B}(x).
\]
We claim that \( (\F_B, Y''_B) \) is a \( \gtrapprox 1 \) refinement of \( (\F_B, Y'_B) \), because \( (\mathbb{W}''_B, Y_{\mathbb{W}''_B}) \) is a \( \gtrapprox 1 \) refinement of \( (\mathbb{W}'_B, Y_{\mathbb{W}'_B}) \) and \( |U(\F_B, Y'_B) \cap B_a| \) is constant up to factors \( \approx 1 \) and \( (\F_{B,W}, Y'_B) \) has constant multiplicity. (The last two facts are part of Proposition~\ref{Prop:fact_mult} which was used to construct \( (\F_B, Y'_B) \) and \( \mathbb{W}'_B \).)

Next we refine \( (\F, Y) \) to make it compatible with \( (\F_B, Y''_B) \). Recall that \( |\F(F_B)_Y(x)| \) is roughly constant among \( x \in Y_B(F_B) \). We know that \( (\F_B, Y''_B) \) is a \( \gtrapprox 1 \) refinement of \( (\F_B, Y_B) \). 
If \( F \in \F(F_B) \), we define \( Y'(F) \cap B = Y(F) \cap Y''_B(F_B) \). Because of the uniformity above, \( (\F, Y') \) is a \( \gtrapprox 1 \) refinement of \( (\F, Y) \). From the definition of \( Y' \) we see that
\[
\F_{Y'}(x) \subset \bigcup_{F_B \in \F_B,\; x \in Y''_B(F_B)} \F(F_B)
\]
and so
\[
\F_{Y'}(x) \subset \bigcup_{W \in \mathbb{W}''_B,\; x \in Y_{\mathbb{W}''_B}(W)} 
\bigcup_{F_B \in \F_{B,W},\; x \in Y''_B(F_B)}
\F(F_B).
\]

By a rescaled version of Lemma~\ref{Lem:flake_contain}, the {\color{black} major angle} between any {\color{black} two flakes} in \( \F_{B,W} \) is \( \lesssim b/r_1 \). Define \( \rho_2 = b/r_1 \). This means that for each \( W \), there exists a $(n-2)$-plane {\color{black} \( v(W) \in \Gr(n,n-2)\)} so that each \( F_B \in \F_{B,W} \) makes a {\color{black} major angle} at most \( \rho_2 \) with \( v(W) \). Thus for each \( x \in U(\F, Y') \) we have
\begin{align}\label{Eq:narrow_cover}
\F_{Y'}(x) \subset \bigcup_{W \in \mathbb{W}''_B,\; x \in Y_{\mathbb{W}''_B}(W)} \{F \in \F_Y(x): \angle(\dir(F), v(W)) \leq \rho_2\}.
\end{align}

Now the cardinality of the inner set is \( \lessapprox \mu(\rho_2) \approx \mu(\F[F_{\rho_2}], Y) \) for any \( F_{\rho_2} \in \F_{\rho_2} \). Also, by Proposition~\ref{Prop:fact_mult}, we know that \( (\W', Y_{\W'}) \) has constant multiplicity \( \mu(\W', Y_{\W'}) \). Since \( (\W''_B, Y_{\W''_B}) \) is a refinement of \( (\W'_B, Y_{\W'_B}) \), we get that for any \( F_{\rho_2} \in \F_{\rho_2} \),
\begin{align}\label{Eq:two_scale}
\mu(\F, Y) \lessapprox \mu(\W''_B, Y_{\W''_B}) \mu(\F[F_{\rho_2}], Y).
\end{align}

Using \( K_{KT}(\beta) \) with \( \eta_{\text{bias}} \) in place of \( \eps \), we have
\begin{align}\label{Eq:fine_scale}
\mu(\F[F_{\rho_2}], Y) \leq \delta^{-\eta_{\text{bias}}} |\F[F_{\rho_2}]|^\beta.
\end{align}

Next we consider \( |\F_{\rho_2}| \). Since we are not in the slab case, \( b \leq \delta^{\eps_{\rm{scal}}} r_1 \), and so \( \rho_2 = b/r_1 \leq \delta^{\eps_{\rm{scal}}} \). On the other hand, \( b \geq \delta \) and so \( \rho_2 = b/r_1 \geq \delta/r_1 = \delta^{1-\eps_{\rm{scal}}} \). So \( \rho_2 \in [\delta^{1-\eps_{\rm{scal}}}, \delta^{\eps_{\rm{scal}}}], \) and so by the hypothesis of Lemma~\ref{Lem:NSAAS}, we have
\[|\F_{\rho_2}| \geq \rho_2^{-2-\zeta}.\]

Our next goal is to estimate \( \mu(\W'', Y_{\W''_B}) \). Recall that \( \theta \) is a typical angle of intersection for \( (\W''_B, Y_{\W_B''}) \). In particular, if \( W_1, W_2 \in \W''_B \) and if \( x \in Y_{\W''_B}(W_1) \cap Y_{\W''_B}(W_2) \), then
\begin{align}\label{Eq:narrow_intersect}
\angle(TW_1, TW_2) \lesssim \theta \leq \delta^{-\tau'} a/b.
\end{align}
Here \( TW \) is the hyperplane spanned by the $(n-1)$ longest axes of \( W \), which is well defined up to an accuracy $a/b$.

Let \( S \) denote a slab segment of dimensions \((a/b)r_1 \times r_1 \times r_1\1\), and let \( \mathbb{S} \) denote a maximal set of essentially distinct slabs \( S \subset B \). Recall from \eqref{Eq:L_S} that we defined
\[
\mathbb{W}_S'' = \{ W \in \mathbb{W}_B'' : W \subset S \text{ and } \angle(TW, TS) \leq \theta \},
\]
and that
\[
\mathbb{W}_B'' = \bigsqcup_{S \in \mathbb{S}} \mathbb{W}_S''.
\]

By \eqref{Eq:narrow_intersect}, \( \mathbb{W}_{Y_{\W_B''}}''(x) \) is contained in at most \( \delta^{-O(\tau')} \) sets \( \mathbb{W}_S'' \). Therefore,
\[
\mu(\mathbb{W}_B'', Y_{\mathbb{W}_B''}) \lesssim \delta^{-O(\tau')} \max_{S \in \mathbb{S}} \mu(\mathbb{W}_S'', Y_{\mathbb{W}_S''}).
\]

However, we would like to only include \( S \) with \( \lambda(\mathbb{W}_S'', Y_{\mathbb{W}_S''}) \) not too small. Recall from \eqref{Eq:W_B_dense} that \( \lambda(\mathbb{W}_B'', Y_{\mathbb{W}_B''}) \gtrapprox \lambda(\mathbb{W}_B', Y_{\mathbb{W}_B'}) \gtrapprox \delta^{2\eta} \). Define
\[
\mathbb{S}_{dens} = \{ S \in \mathbb{S} : \lambda(\mathbb{W}_S'', Y_{\mathbb{W}_S''}) \geq (1/100)\lambda(\mathbb{W}_B', Y_{\mathbb{W}_B'}) \}.
\]
Then we can define \( \mathbb{W}_B''' = \bigcup_{S \in \mathbb{S}_{dens}} \mathbb{W}_B'' \) and \( Y_{\mathbb{W}_B'''} = Y_{\mathbb{W}_B''} \), so that \( (\mathbb{W}_B''', Y_{\mathbb{W}_B'''}) \) is a \( \gtrapprox 1 \) refinement of \( (\mathbb{W}_B'', Y_{\mathbb{W}_B''}) \). Therefore
\[
\mu(\mathbb{W}_B'', Y_{\mathbb{W}_B''}) \lessapprox \mu(\mathbb{W}_B''', Y_{\mathbb{W}_B'''}) \lesssim \delta^{-O(\tau')} \max_{S \in \mathbb{S}_{dens}} \mu(\mathbb{W}_S'', Y_{\mathbb{W}_S''}).
\]

Select \( S \in \mathbb{S}_{dens} \) attaining the preceding maximum. We know that \( \lambda(\mathbb{W}_S'', Y_{\mathbb{W}_S''}) \gtrapprox \delta^{2\eta} \).

We can make a linear change of variables that converts \( S \) to \( B_1 \) and converts \( \mathbb{W}_S'' \) to a set \( \tilde{\F} \) of \( \rho_2 \)-flakes in \( B_1 \). By the maximal density factoring lemma, \( \Delta_{max}(\mathbb{W}_S'') \leq \Delta_{max}(\mathbb{W}_B) \) is small. With the original lemma, we would have \( \Delta_{max}(\mathbb{W}_B) \lesssim 1 \). With our modified lemma,  
Lemma~\ref{Lem:biased_factoring}, we have \( \Delta_{max}(\mathbb{W}_B) \lesssim \delta^{-O(\eta_{\text{bias}})} \). This implies \( \Delta_{max}(\mathbb{W}''_S) \lesssim \delta^{-O(\eta_{\text{bias}})} \), and so  
\[\Delta_{max}(\tilde{\F}) \lesssim \delta^{-O(\eta_{\text{bias}})}.\]
We also have  
\[\lambda(\tilde{\F}, Y_{\tilde{\F}}) = \lambda(\mathbb{W}''_S, Y_{\mathbb{W}''_S}) \gtrapprox \delta^{2\eta}.\]

Since \( \tilde{\F} \) is a set of \( \rho_2 \)-flakes in \( B_1 \) with \( \Delta_{max}(\tilde{\F}) \lesssim \delta^{-O(\eta_{\text{bias}})} \), we have \( |\tilde{\F}| \lesssim \delta^{-O(\eta_{\text{bias}})}\rho_2^{-2} \).  
If \( \eta > 0 \) is sufficiently small, then we can apply Lemma~\ref{Lem:general_KT} to \( \mu(\tilde{\F}, Y_{\tilde{\F}}) \) with \( \eta_{\text{bias}} \) in place of \( \eps \), which gives 
\begin{align}\label{Eq:rescaled_envelope}
\mu(\tilde{\F}, Y_{\tilde{\F}}) \leq \delta^{-O(\eta_{\text{bias}})}|\tilde{\F}|^{\beta} \lesssim \delta^{-O(\eta_{\text{bias}})}(\rho_2^{-2})^{\beta} \leq \delta^{-O(\eta_{\text{bias}})}\rho_2^{\beta\zeta}|\F_{\rho_2}|^{\beta}.
\end{align}
Therefore, we get 
\begin{align}\label{Eq:coarse_scale}
\mu(\mathbb{W}''_B, Y_{\mathbb{W}''_B}) \lessapprox \delta^{-O(\tau')}\mu(\tilde{\F}, Y_{\tilde{\F}}) \lesssim \delta^{-O(\eta_{\text{bias}}+\tau')}\rho_2^{\beta\zeta}|\F_{\rho_2}|^{\beta}.
\end{align}

Combining \eqref{Eq:two_scale}, \eqref{Eq:fine_scale}, and \eqref{Eq:coarse_scale}, we conclude that  
\[\mu(\F, Y) \lessapprox \delta^{-O(\eta_{\text{bias}}+\tau')}\rho_2^{\beta\zeta}|\F|^{\beta} \leq \delta^{-O(\eta_{\text{bias}}+\tau')+\eps_{\text{scal}}\beta\zeta}|\F|^{\beta},\]
where for the final inequality we used the fact that we are not in the slab case, and thus \( \rho_2 = b/r_1 \leq \delta^{\eps_{\text{scal}}} \). We choose \( \tau, \tau', \eta_{\text{bias}} \ll \eps_{\text{scal}}\beta\zeta \), and this gives the desired bound \eqref{Eq:KT_Goal} with a gain \( \nu = \eps_{\text{scal}}\beta\zeta/2 \).

Therefore, the proof of Main Lemma~\ref{Main_lem_2}, and thus Theorem~\ref{Thm:Main}, is completed.

\section{Applications}\label{Sec:applications}

\subsection{Rogers-type \texorpdfstring{$(n-2)$}{(n-2)}-plane Kakeya sets in \texorpdfstring{$\R^n$}{Rn}}
\label{subsec:apply_Rogers_R4}

In this section, we apply Theorem~\ref{Thm:Main} to a codimension-two analogue of the planar variant of the Kakeya problem introduced by Rogers~\cite{rogers2006}. This gives a full-dimensional result for suitable two-dimensional families of $(n-2)$-planes in $\mathbb R^n$. When $n=4$, this coincides exactly with the ``real'' $2$-plane problem considered by Rogers.

Let $\mathcal M\subset \Gr(n,n-2)$ be (a compact patch of) a smooth two-dimensional manifold, equipped with the normalized measure $\mu$ induced by the major-angle metric on $\Gr(n,n-2)$. Recall from Definition~\ref{Def:major_angle} that
$$\angle(V,W)=\|\pi_V-\pi_W\|.$$

For an $(n-1)$-dimensional linear subspace $\Pi\subset\mathbb R^n$ and
$a>0$, define $$\mathcal M_\Pi^a := \left\{V\in\mathcal M: \inf_{\substack{W\in \Gr(n,n-2)\\ W\subset\Pi}}
\angle(V,W)<a\right\}.$$

Also, for $V\in\mathcal M$ and $b>0$, write
$$B_b(V):=\{W\in\mathcal M:\angle(V,W)<b\}.$$

Recall that Rogers' curvature condition for $(n-2)$-dimensional families $\mathcal N$ of $2$-planes in $\mathbb R^n$ is
$$\mu\big(\mathcal N_{\Pi_m}^{\eps_1} \cap B_{\eps_2}(V) \big) \lesssim \eps_1^{n-m}\eps_2^{m-2},$$
for every $\Pi_m\in \Gr(n,m)$ with $3\leq m\leq n-1$, every $V\in\mathcal N$, and every $\eps_1,\eps_2>0$.

The natural analogue for families of $k$-planes is obtained by replacing the exponent $m-2$ by $m-k$:
$$\mu\big(\mathcal M_{\Pi_m}^{\eps_1}
\cap B_{\eps_2}(V) \big) \lesssim \eps_1^{n-m}\eps_2^{m-k},$$
for every $\Pi_m\in \Gr(n,m)$ with $k+1\leq m\leq n-1$, every $V\in\mathcal M$, and every $\eps_1,\eps_2>0$. In the codimension-two case $k=n-2$, the only intermediate dimension is $m=n-1$, and the corresponding condition reduces to
$$\mu\big(\mathcal M_\Pi^a\cap B_b(V)
\big) \lesssim ab,$$
for every $\Pi \in \Gr(n,n-1)$, every $V\in\mathcal M$, and every $a,b>0$.

We shall say that a smooth two-dimensional manifold
$\mathcal M\subset \Gr(n,n-2)$ is \textit{codimension-two curved in the sense of Rogers} if it satisfies this condition. When $n=4$, this is precisely Rogers' original curvature condition for the real $2$-plane problem in $\mathbb R^4$.

We now show that this curvature condition implies the convex nonconcentration condition used throughout this paper.

\begin{lem}
\label{lemma:rogers-implies-non-concentration}
Let $\mathcal M\subset \Gr(n,n-2)$ be a smooth two-dimensional manifold which is codimension-two curved in the sense of Rogers. Let $\Lambda\subset\mathcal M$ be a $\delta$-separated set, and
for each $V\in\Lambda$, let $F_V$ be an arbitrary $\delta$-flake with $\dir(F_V)=V$. Set
$$\F :=\{F_V:V\in\Lambda\}.$$
Then
$$\Delta_{\max}(\F )\lesssim 1,$$
where the implicit constant is allowed to depend on $\mathcal M$.
In particular,
$$\Delta_{\max}(\F )\leq \delta^{-\eta}$$
for every fixed $\eta>0$ and all sufficiently small $\delta$.
\end{lem}

\begin{proof}
Fix a convex set $K\subset\mathbb R^n$. We need to prove
$$ \Delta(\F ,K) =\frac{\sum_{F\in\F [K]}|F|}{|K|}
\lesssim 1.$$
This is trivial if $\F [K]=\varnothing$, so suppose that
$K$ contains at least one $\delta$-flake from $\F $.
Let $a_1\leq a_2\leq\cdots\leq a_n$ be the dimensions of $K$, in the sense of the axes of its John ellipsoid, and let $e_1,\ldots,e_n$ be the corresponding orthogonal principal directions. Since $K$ contains a $\delta$-flake, we have
$a_1,a_2\gtrsim\delta$ and $a_3,\ldots,a_n\gtrsim 1$. 

Set $a:=a_1$ and $b:=a_2$. Define 
$$ \Pi:=e_1^\perp = \operatorname{span}(e_2,\ldots,e_n),$$ and $$V_K:=\operatorname{span}(e_3,\ldots,e_n).$$
We claim that every $F_V\in\F [K]$ satisfies $V\in\mathcal M_\Pi^{Ca}$ and $\angle(V,V_K)\lesssim b$. To see this, first note that the unit ball in $V$ is contained in the difference set $F_V-F_V$, while $$F_V-F_V\subset K-K.$$ Since the widths of $K-K$ in the $e_1$ and $e_2$ directions are
$O(a)$ and $O(b)$, respectively, it follows that
$$
|v\cdot e_1|\lesssim a,
\qquad
|v\cdot e_2|\lesssim b,
$$
for every unit vector $v\in V$. The first estimate implies that $V$ is within major angle $O(a)$ of an $(n-2)$-plane contained in $\Pi$ (by Lemma~\ref{Lem:faithful_proj}, just consider the $(n-2)$-plane $\pi_{\Pi}(V)$). The two estimates together give $\|\pi_{V_K^\perp}|_V\|\lesssim b$, which by Item~(i) of Theorem~\ref{Thm:Principal_angle} gives
$$\angle(V,V_K)\lesssim b.$$

Note that the curvature condition requires the ball
$B_b(\cdot)$ to be centered at an element of $\mathcal M$, whereas $V_K$ need not belong to $\mathcal M$. This issue is easily fixed by choosing a flake $F_{V_0}\in\F [K]$. Then, $V_0\in\Lambda\subset\mathcal M$ and the preceding estimate gives
$$
\angle(V_0,V_K)\lesssim b.
$$

Consequently, for every $F_V\in\F [K]$, the triangle inequality for the major-angle metric gives
$$\angle(V,V_0) \le \angle(V,V_K)+\angle(V_K,V_0) \lesssim b.$$
Combining this with $V\in\mathcal M_\Pi^{Ca}$, we obtain 
$$\{\dir(F): F\in\F [K]\} \subset \mathcal M_\Pi^{Ca}\cap B_{Cb}(V_0).$$

We next discretize the codimension-two curvature condition. Since $\mathcal M$ is a smooth two-dimensional manifold and $\Lambda$ is $\delta$-separated, we obtain
$$ \#\big(\Lambda\cap\mathcal M_\Pi^{Ca}\cap B_{Cb}(V_0)
\big) \lesssim \delta^{-2} \mu\big( \mathcal M_\Pi^{C'a}\cap B_{C'b}(V_0)\big),$$
by a standard packing argument (since $a,b\gtrsim\delta$, the $O(\delta)$ losses arising from the packing argument can be absorbed into constant multiples of $a$ and $b$) for a slightly bigger constant $C' > 0$. Next, applying the curvature condition, we conclude that
$$\#\F [K] \lesssim ab\delta^{-2}.$$
Finally, $|F|\sim\delta^2$ for every $F\in\F $, while John's theorem gives $|K| \approx a_1a_2\ldots a_n \gtrsim ab.$ 
Therefore, we get 
$$\Delta(\F ,K) = \frac{\#\F [K]|F|}{|K|} \lesssim
\frac{ab\delta^{-2}\delta^2}{|K|} \lesssim 1.$$
Since $K$ was arbitrary, we can conclude that 
$$ \Delta_{\max}(\F )\lesssim1,$$
as desired.
\end{proof}

We shall now resolve the ``real'' $2$-plane problem posed by Rogers. In what follows, we will use $\mathcal H^s(E)$ to denote the $s$-dimensional Hausdorff measure of a set $E$.

\begin{theorem}[Rogers-type $(n-2)$-plane Kakeya problem in $\mathbb R^n$]
\label{thm:rogers-application}
Let ${\mathcal M\subset \Gr(n,n-2)}$ be a smooth two-dimensional manifold which is codimension-two curved in the sense of Rogers, and suppose that $E\subset\mathbb R^n$ contains a translate of every plane $V\in\mathcal M$. Then,
$$
\dim_H(E)=n.
$$
In particular, when $n=4$, every real $\mathcal K_2^4$-set in the sense of Rogers~\cite{rogers2006} has full Hausdorff dimension.
\end{theorem}

\begin{proof}
Suppose for a contradiction that
$$\dim_H(E)<s<n.$$ 
Fix $0<\beta<n-s,$ and let $\eta>0$ be given by Theorem~\ref{Thm:Main} with $\eps = \beta$ therein. 

Since $\dim_H(E)<s$, we have $\mathcal H^s(E)=0$. Hence, for every
$\varepsilon>0$ and every $\delta_0>0$, we may cover $E$ by balls
$$E\subset\bigcup_j B(x_j,r_j),$$
satisfying $r_j\leq\delta_0$, such that
$$\sum_j r_j^s<\varepsilon.$$

For $k\geq k_0$, where $k_0$ is chosen so that
$2^{-k_0-1}<\delta_0\leq 2^{-k_0}$, define
$$\mathcal B_k := \left\{B(x_j,r_j): 2^{-k-1}<r_j\leq 2^{-k} \right\}.$$
Set
$$ E_k:=\bigcup_{B\in\mathcal B_k}B,$$
and $\delta_k:=2^{-k}$.

For each $V\in\mathcal M$, choose a unit $(n-2)$-dimensional cube $S_V$ contained in the corresponding affine plane lying in $E$. Since the sets $E_k$ cover $S_V$, we have
$$1\le \sum_{k\geq k_0} \mathcal H^{n-2}(S_V\cap E_k).$$

Since $\sum_{k\geq k_0}k^{-2}<\infty$, it follows, after fixing a
sufficiently small absolute constant $c>0$, that for every $V\in\mathcal M$ there exists some $k\geq k_0$ such that
$$\mathcal H^{n-2}(S_V\cap E_k)\geq c k^{-2}.$$
Define
$$\mathcal M_k := \left\{ V\in\mathcal M: \mathcal H^{n-2}(S_V\cap E_k)\geq c k^{-2}\right\}.$$
Then,
$$\mathcal M=\bigcup_{k\geq k_0}\mathcal M_k.$$
Since $\mu(\mathcal M) = 1$, another summability argument shows that for some $k\geq k_0$,
$$\mu(\mathcal M_k)\gtrsim k^{-2}.$$

Fix such a $k$ and write $\delta:=2^{-k}.$ Let $\Lambda\subset\mathcal M_k$ be a maximal $\delta$-separated subset.
Since $\mathcal M$ is a smooth two-dimensional manifold, the
$\delta$-balls centered at points of $\Lambda$ cover $\mathcal M_k$, and each
such ball has $\mu$-measure $O(\delta^2)$. Hence,
$$|\Lambda|\gtrsim \mu(\mathcal M_k)\,\delta^{-2} \gtrsim k^{-2}\delta^{-2}.$$

For each $V\in\Lambda$, let $F_V$ be the $\delta$-flake obtained by taking the $\delta$-neighborhood of $S_V$, and set
$$\F :=\{F_V:V\in\Lambda\}.$$
By Lemma~\ref{lemma:rogers-implies-non-concentration},
$$\Delta_{\max}(\F )\lesssim1.$$
In particular, since $\eta>0$ is fixed, we have
$$\Delta_{\max}(\F )\leq\delta^{-\eta},$$
for all sufficiently small
$\delta$.

We now define a shading of $\F $. For $V\in\Lambda$, let $A_V:=S_V\cap E_k,$ and define
$$Y(F_V) :=\left\{x+u: x\in A_V,\ u\in V^\perp,\ 
\norm{u}\leq c_0\delta \right\},$$
where $c_0>0$ is a sufficiently small absolute constant. Then 
$$Y(F_V)\subset F_V\cap N_{c_0\delta}(E_k),$$
and 
$$|Y(F_V)| \sim \delta^2\mathcal H^{n-2}(A_V)
\gtrsim k^{-2}\delta^2,$$
by Fubini. Since $|F_V|\sim \delta^2$, it follows that
$$\lambda(\F ,Y)\gtrsim k^{-2}.$$
As $\delta=2^{-k}$, we have 
$$k^{-2}\gtrsim\delta^\eta,$$
for sufficiently large $k$.
Thus, $(\F ,Y)$ satisfies the hypotheses of Theorem~\ref{Thm:Main}. Therefore, we obtain
$$|U(\F ,Y)| \gtrsim \delta^\beta|\F ||F_V|.$$
Using $|\F |=|\Lambda|\gtrsim k^{-2}\delta^{-2}$ and $|F_V|\sim\delta^2$, this gives 
$$|U(\F ,Y)| \gtrsim k^{-2}\delta^\beta.$$
On the other hand,
$$U(\F ,Y)\subset N_{c_0\delta}(E_k).$$
If $N_k:=|\mathcal B_k|,$ then $E_k$ is a union of $N_k$ balls of radius comparable to $\delta$, and hence
$$|N_{c_0\delta}(E_k)| \lesssim N_k\delta^n.$$
Consequently,
$$N_k\delta^n \gtrsim k^{-2}\delta^\beta,$$
or equivalently,
$$N_k\delta^s \gtrsim k^{-2}\delta^{s+\beta-n}.$$
Since $\delta=2^{-k}$ and $s+\beta-n<0$, the right-hand side tends to infinity as $k\to\infty$.

However, every ball in $\mathcal B_k$ has radius comparable to $\delta$, so the choice of the cover of $E$ gives
$$ N_k\delta^s \lesssim \sum_{B(x_j,r_j)\in\mathcal B_k}r_j^s
<\varepsilon.$$

Since $k\geq k_0$ and $k_0\to\infty$ as $\delta_0\to0$, choosing
$\delta_0$ sufficiently small ensures that the preceding lower bound for $N_k\delta^s$ exceeds this upper bound $N_k\delta^s\lesssim\varepsilon,$ giving a contradiction. Therefore, $$\dim_H(E)\ge s.$$
Since this holds for every $s<n$, the conclusion follows.
\end{proof}

One can easily see from the proof of Theorem~\ref{thm:rogers-application} that we still have full Hausdorff dimension even if we only take a patch of each plane with positive $\mathcal{H}^{n-2}$-measure instead of the whole plane. Such a ``compact'' version is closer to the classical formulation of the Kakeya set conjecture.

It is also worth pointing out that in \cite{rogers2006}, Rogers proved that every real $\mathcal{K}_2^4$-set has Hausdorff dimension $\geq 7/2$ by adapting Wolff's hairbrush argument \cite{wolff95} to $2$-planes in $\R^4$. His proof actually yields something stronger, and may be used to obtain the corresponding maximal function estimates. Thus, Rogers' approach provides additional information beyond the Hausdorff dimension statement and remains particularly relevant from the perspective of maximal function estimates. It is a good question how to obtain the best possible such estimates, which may require the resolution of the Kakeya maximal conjecture for tubes in $\R^3$.

\subsection{Kakeya for \texorpdfstring{$a \times b \times \1$}{axbx1} planks in \texorpdfstring{$\R^n$}{Rn}}\label{subsec:apply_plank_Rn}
Before ending this paper, we point out that good Kakeya estimates for $a \times b \times \1$ planks in $B_1 \subset \R^n$ ($n \geq 3$) can also be obtained as by-products of the proof of Theorem~\ref{Thm:Main}. For example, since the assertion $K_{KT}(0)$ in Definition~\ref{Def:Katz_Tao_beta} holds, by Lemma~\ref{Lem:plank_Katz_Tao}, we immediately get the following result: 
\begin{theorem}\label{Thm:Rn_plank_Katz_Tao}
For every $\eps > 0$, there exists $\eta, a_0 > 0$ such that
the following holds whenever $0 < a \leq a_0$ and $a \leq b \leq 1$. Let $\LL$ be a set of $a \times b \times \1$ planks in $B_1 \subset \R^n$ with $\Delta_{max}(\LL) \leq a^{-\eta}$, and let $Y$ be a shading on $\LL$ with $\lambda(\LL, Y) \geq a^\eta$. Then $|U(\LL, Y)| \geq a^\eps|\LL||L|$.    
\end{theorem}
Note that Theorem~\ref{Thm:Rn_plank_Katz_Tao} includes Theorem~\ref{Thm:Main} as a special case by taking $a=b=\delta$. We can also require the planks to be Frostman instead of Katz--Tao. More precisely, since the assertion $K_F(0)$ in Definition~\ref{Def:Frostman_beta} holds, by Lemma~\ref{Lem:plank_Frostman}, we have:
\begin{theorem}\label{Thm:Rn_plank_Frostman}
For every $\eps > 0$, there exists $\eta, a_0 > 0$ such that
the following holds whenever $0 < a \leq a_0$ and $a \leq b \leq 1$. Let $\LL$ be a set of $a \times b \times \1$ planks in $B_1 \subset \R^n$ with $C_F(\LL) \leq a^{-\eta}$, and let $Y$ be a shading on $\LL$ with $\lambda(\LL, Y) \geq a^\eta$. Then $|U(\LL, Y)| \geq a^\eps$.    
\end{theorem}
Theorem~\ref{Thm:Rn_plank_Frostman} is a higher-dimensional generalization of Lemma~\ref{Lem:R3_plank_Frostman}. Note that if we take $a=b=\delta$, then Theorem~\ref{Thm:Rn_plank_Frostman} can be viewed as a ``Frostman'' version of Theorem~\ref{Thm:Main}. 

\begin{rem}
The ranges for $a$ and $b$ in Theorem~\ref{Thm:Rn_plank_Katz_Tao}~and~\ref{Thm:Rn_plank_Frostman} are initially $0 < a \leq b \leq b_0$ for some $b_0$ depending on $\eps$. However, similar arguments as in Remark~\ref{Rem:a_0_b_0} remove the smallness restriction on $b$ from both of them.
\end{rem}

In our proof of the $(n-2)$-plane sticky Kakeya in $\R^n$, Lemma~\ref{Lem:R3_plank_Frostman} is the key ingredient for growing medium grains in each slice. It is unclear to the authors whether Theorem~\ref{Thm:Rn_plank_Katz_Tao}~or~\ref{Thm:Rn_plank_Frostman} could be a useful input in other problems.

\appendix

\section{Probability lemmas}\label{Appendix:Probability}

In this appendix, we use simple probabilistic lemmas to prove Lemmas~\ref{Lem:random_translate},~\ref{Lem:random_sub_(F)}, and~\ref{Lem:random_rigid_F}. These 
lemmas are basically different versions of Chernoff’s inequality, and a detailed proof can be found in \cite{GWZ26streamlined} and \cite{Chernoff1952}. Here we only record the statements for the reader's convenience. We use $\mathcal{P}[\cdot]$ to denote the probability of an event.

\begin{lem}[{\cite[Lemma~A.1]{GWZ26streamlined}}]\label{Lem:Prob}
    Suppose that $X_1, \dots, X_N$ are i.i.d. random variables. Suppose that $0 \leq X_1 \leq M$ almost surely and that $m = \E X_1$. If $Nm\geq M$, then for every $S>0$, 
    \begin{align*}
        \mathcal{P}[X_1+\cdots+X_N>S] \leq e^{10} e^{-\frac{S}{Nm}}.
    \end{align*}
    If $Nm \leq M$, then for every $S>0$, 
    \begin{align*}
        \mathcal{P}[X_1+\cdots+X_N>S] \leq e^{10} e^{-\frac{S}{2M}}.
    \end{align*}
\end{lem}

\begin{lem}[Chernoff \cite{Chernoff1952}]\label{Lem:Chernoff}
    Let $X_1, \dots, X_N$ be independent random variables taking value $1$ with probability $p$ and value $0$ with probability $1 - p$. Let $A \geq 2$. Then 
    \begin{align*}
        \mathcal{P}[X_1+\cdots+X_N \leq pN/2] \leq e^{-pN/8}; \qquad
        \mathcal{P}[X_1+\cdots+X_N \geq ApN] \leq e^{-ApN/6}.
    \end{align*}
\end{lem}

We will consider random translations and rigid motions of $\R^n$. {\color{black} Recall that for $v \in \R^n$, we defined $R_v(x) = x + v$.
A random translation of size $\leq \rho$ is a translation $R_v$ by a vector $v \in B(0, \rho)$ selected randomly from
the uniform measure on $B(0, \rho)$. A random rigid motion of size $\leq 1$ is a random rotation composed with a
random translation of size $\leq 1$.}

First, we state and prove a more general version of Lemma~\ref{Lem:random_translate}.

\begin{lem}\label{Lem:random_translate'}
Let \( \delta \leq \rho \leq 1 \). Suppose that \( \F \) is a set of {\color{black} essentially distinct} \( (\delta ,k)\)-flakes in \( F_\rho \subset \mathbb{R}^n \) with \( \frac{|F_\rho|}{|\F||F_\delta|} \sim J \in \N^+ \). Let \( R_1, \ldots, R_J \) be a set of \( J \) {\color{black}independent} random translations of size \( \leq \rho \). Let \( \F' = \bigsqcup_{j=1}^J R_j(\F) \). Then with high probability {\rm(}at least $1-\delta^{100n^2}$ when $\delta$ is sufficiently small{\rm)} the flakes of \( \F' \) are essentially distinct up to multiplicity \( \lessapprox 1 \), and we have \( \Delta_{max}(\F') \lessapprox \Delta_{max}(\F) \). After refining \( \F' \) to make the flakes essentially distinct, we have with high probability \( |\F'| \approx J|\F| \) and \( \Delta_{max}(\F') \lessapprox \Delta_{max}(\F) \).
\end{lem}

\begin{proof}
First we check that with high probability, for any \( (\delta,k) \)-flake \( F_0 \subset F_\rho \), we have \( |\F'[100F_0]| \lessapprox 1 \). This will show that after refining \( \F' \) to be essentially distinct, we still have \( |\F'| \approx J|\F| \).

We define \( X_j = |R_j(\F)[100F_0]| \), which we think of as a random variable. Note that \( |\F'[100F_0]| = \sum_{j=1}^J X_j \). 

Since the flakes of \( \F \) are essentially distinct, we have \( X_j = |R_j(\F)[100F_0]| = |\F[R_j^{-1}(100F_0)]| \lesssim 1 \). We choose \( M \lesssim 1 \) so that \( 0 \leq X_j \leq M \) for each $j$.

Next we want to bound \( \mathbb{E} X_j \). By linearity of expectation, we have 
\[ \mathbb{E}X_j = \sum_{F \in \F} \mathcal{P}[R_j(F) \subset 100F_0]. \] 
Since \( R_j \) is a random translation of size \( \leq \rho \), we note that for each \( F \in \F \),
\[
\mathcal{P}[R_j(F) \subset 100F_0] \lesssim \frac{|F_\delta|}{|F_\rho|}.
\]
Therefore,
\[
m = \mathbb{E}X_j \lesssim |\F| |F_\delta| / |F_\rho|.
\]

We claim that \( Jm \lesssim 1 \). To see this, we note that \( Jm \lesssim J|\F| |F_\delta| / |F_\rho| \sim 1 \) by our assumption about \( J \). Thus by possibly enlarging the previous $M$, we also have $Jm \leq M$ for $M \lesssim 1$.

Therefore, if \( S \geq 1 \), then Lemma~\ref{Lem:Prob} (or the standard Chernoff bound) gives
\[
\mathcal{P}[X_1 + \cdots + X_J > S] \leq e^{10} e^{-\frac{S}{2M}}
\]
for some $M\lesssim 1$. Since the number of essentially distinct \( F_0 \) is at most \( \delta^{-O(1)} \), we can choose \( S \lessapprox 1 \) so that with high probability \( |\F'[100F_0]| \leq S \) for every \( F_0 \).

Next we estimate \( \Delta_{max}(\F') \).

Let \( K \subset 2F_\rho \) be a convex set of dimensions $a_1 \times \cdots \times a_{n-k} \times 1 \times \cdots \times 1$ with $\delta \leq a_1 \leq \cdots \leq a_{n-k} \leq \rho$. Let \( Z_j \coloneqq |R_j(\F)[K]| \), which we think of as a random variable. Note that \( |\F'[K]| = \sum_{j=1}^J Z_j \). We will show that (for each \( K \)),
\begin{align}\label{Eq:large_deviation'}
\mathcal{P}\left[\sum_{j=1}^J Z_j > A\Delta_{max}(\F)\frac{|K|}{|F_{\delta}|}\right] \leq e^{10} e^{-A/C}.
\end{align}

Suppose this has been established. Since there are only \( (\rho/\delta)^{O(1)} \) essentially different choices for the convex set \( K \subset F_\rho \), we have that with probability at least \( 1 - (\rho/\delta)^{O(1)} e^{10} e^{-A/C} \), for every \( K \) we have \( |\F'[K]| \lesssim A\Delta_{max}(\F)\frac{|K|}{|F_{\delta}|} \) and so \( \Delta_{max}(\F') \lesssim A\Delta_{max}(\F) \). We can choose $A \lessapprox 1$ so that the probability is at least $1-\delta^{100n^2}$ for $\delta$ sufficiently small.

We will check \eqref{Eq:large_deviation'} using Lemma~\ref{Lem:Prob}. Define \( M = \Delta_{max}(\F)|K|/|F_{\delta}| \). Note that
\[
M \geq |\F[R_j^{-1}(K)]| = Z_j.
\]

Clearly \( Z_j \geq 0 \), so \( 0 \leq Z_j \leq M \). For each \( F \in \F \), we claim that
\begin{align}\label{Eq:prob_claim}
\mathcal{P}[R_j(F) \subset K] \lesssim |K| / |F_\rho|.
\end{align}

To see this, pick a point \( x \in F \) (for example \( x \) could be the center of \( F \)) and bound
\[
\mathcal{P}[R_j(F) \subset K] \leq \mathcal{P}[R_j(x) \in K] \leq \frac{|B(x,\rho) \cap K|}{|B(x,\rho)|} \lesssim
\frac{a_1 \cdots a_{n-k} \cdot \rho^k}{\rho^n} \sim \frac{|K|}{|F_\rho|}.
\]

By linearity of expectation, \eqref{Eq:prob_claim} gives the bound
\[
\mathbb{E}Z_j \lesssim \frac{|\F||K|}{|F_\rho|}.
\]

Define \( m = \mathbb{E}Z_j \). By our definition of \( J \) we have \( Jm \sim \frac{|K|}{|F_\delta|} \lesssim \Delta_{max}(\F) \frac{|K|}{|F_\delta|} = M \). Therefore, Lemma~\ref{Lem:Prob} tells us that
\[
\mathbb{P}[Z_1 + \cdots + Z_J > S] \leq e^{10} e^{-\frac{S}{CM}}.
\]

By taking \( S = A \Delta_{max}(\F) \frac{|K|}{|F_{\delta}|} = AM \), we get \eqref{Eq:large_deviation'} as desired.
\end{proof}

Next, we state and prove a more general version of Lemma~\ref{Lem:random_sub_(F)}.
\begin{lem}\label{Lem:random_sub_(F)'}
     Let $\F$ be a set of $(\delta,k)$-flakes in $B_1 \subset \R^n$. Then it contains a subset $\F'$ with $C_F(\F') \lessapprox C_F(\F)$ and $|\F'| \approx \delta^{-(n-k)} C_F(\F)^{-1}$.
\end{lem}
\begin{proof}
    Note that $|\F| \gtrsim \delta^{-(n-k)} C_F(\F)^{-1}$. If $|\F| \leq (\log(1/\delta))^{2}\delta^{-(n-k)} C_F(\F)^{-1}$, then there is nothing to prove as we can simply take $\F' = \F$. So we may assume $|\F| \geq (\log(1/\delta))^{2} \delta^{-(n-k)} C_F(\F)^{-1}$. 
    
    Let us choose each flake in $\F$ independently with probability \begin{align*}
        p \coloneqq (\log(1/\delta))^{2} (\delta^{n-k}C_F(\F)|\F|)^{-1}.
    \end{align*} 
    Note that $0 \leq p \leq 1$ by our assumption. It suffices to show that with high probability, the resulting subset $\F'$ has the desired properties. 

    Let $N \coloneqq |\F|$ and $\F = \{F_j\}_{j=1}^N$. For $j = 1, \dots, |\F|$, let $X_j$ be the random variable that takes value $1$ if $F_j$ is chosen, and $0$ otherwise. Note that $|\F'| = \sum_{j=1}^N X_j$. Then Lemma~\ref{Lem:Chernoff} gives
    \begin{align*}
        \mathcal{P}[X_1+\cdots+X_N \leq (\log(1/\delta))^{2} \delta^{-(n-k)}C_F(\F)^{-1}/2] &\leq e^{-(\log(1/\delta))^{2} \delta^{-(n-k)}C_F(\F)^{-1}/8};\\
        \mathcal{P}[X_1+\cdots+X_N \geq 2(\log(1/\delta))^{2} \delta^{-(n-k)}C_F(\F)^{-1}] &\leq e^{-(\log(1/\delta))^{2} \delta^{-(n-k)}C_F(\F)^{-1}/3}.
    \end{align*}
    Thus with probability at least 
    $$1 - 2e^{-(\log(1/\delta))^{2} \delta^{-(n-k)}C_F(\F)^{-1}/10} \geq 1-2e^{-(\log(1/\delta))^{2}/C}$$ 
    (as $C_F(\F) \lesssim \delta^{-(n-k)}$), we have $|\F'| \sim (\log(1/\delta))^{2} \delta^{-(n-k)}C_F(\F)^{-1}$. In particular, $|\F'| \approx \delta^{-(n-k)}C_F(\F)^{-1}$. 

    Next we turn to bound $C_F(\F')$. Let $K \subset B_1$ be a convex set of dimensions $a_1 \times \cdots \times a_{n-k} \times 1 \times \cdots \times 1$ with $\delta \leq a_1 \leq \cdots \leq a_{n-k} \leq 1$. Applying Lemma~\ref{Lem:Prob} to those $X_j$'s associated to $\F[K]$, we know that
    \begin{align*}
        \mathcal{P}[|\F'[K]|>S] \leq e^{10} e^{-\frac{S}{2M}},
    \end{align*}
    where $M = \lceil |\F[K]|\,\E X_1 \rceil = \lceil |\F[K]|\,p \rceil$. In particular, by taking $S = AM$, we get
    \begin{align*}
        \mathcal{P}\big[|\F'[K]|>A\lceil |\F[K]|\,p \rceil\big] \leq e^{10} e^{-A/2}.
    \end{align*}    
    Since there are only $\delta^{-O(1)}$ essentially different choices for the convex set $K \subset B_1$, with probability at least \( 1 - \delta^{-O(1)} e^{10} e^{-A/2} \), we have that \( |\F'[K]| \leq A\lceil |\F[K]|\,p \rceil \) for every \( K \). We can choose $A \lessapprox 1$ so that the probability is at least $1-\delta^{100n^2}$ for $\delta$ sufficiently small. 
    
    Recall that $|\F[K]| \leq C_F(\F)|K||\F|$. Then with high probability, for every \( K \), we have 
    $$|\F'[K]| \lessapprox \lceil |\F[K]|\,p \rceil \leq \lceil C_F(\F)|K||\F|\,p \rceil = \lceil (\log(1/\delta))^2\delta^{-(n-k)}|K| \rceil.$$ If $\F'[K] = \varnothing$, then $|\F'[K]| \lessapprox C_F(\F)|K||\F'|$ trivially holds. Otherwise, we must have $|K| \geq |F| \sim \delta^{n-k}$, which means $(\log(1/\delta))^2 \delta^{-(n-k)}|K| \geq 1$, and so $|\F'[K]| \lessapprox (\log(1/\delta))^2 \delta^{-(n-k)}|K|$. Finally, recall that with high probability 
    $$|\F'| \sim (\log(1/\delta))^{2} \delta^{-(n-k)}C_F(\F)^{-1},$$ so we can conclude that $|\F'[K]| \lessapprox C_F(\F)|K||\F'|$. In particular, $C_F(\F') \lessapprox C_F(\F)$.
\end{proof}

Finally, we state and prove a more general version of Lemma~\ref{Lem:random_rigid_F}:
\begin{lem}\label{Lem:random_rigid_F'}
    Suppose that $\F$ is a set of {\color{black} essentially distinct} $(\delta,k)$-flakes in $B_1 \subset \R^n$. Let $R_1, \dots , R_J$ be {\color{black}independent} random rigid motions of size $\leq 1$, with $J = \lceil C_F(\F) \rceil$, and let $\F' = \bigsqcup_{j=1}^J R_j(\F)$. Then with high probability the flakes of $\F'$ are essentially distinct up to multiplicity $\lessapprox 1$ and also $C_F(\F') \lessapprox 1$. After refining $\F'$ to make the flakes essentially distinct, we have $|\F'| \approx C_F(\F)|\F|$ and $C_F(\F') \lessapprox 1$.
\end{lem}

The proof ideas are the same as for the previous Lemma~\ref{Lem:random_translate'}.

\begin{proof}
    We follow the arguments in \cite[Appendix~A.1]{GWZ26streamlined}. Although $\F$ is essentially distinct, $\F'$ may not be. However, we will show that with high probability, for any $(\delta,k)$-flake $F_0$, there are $\lessapprox 1$ flakes of $\F'$ contained in $100F_0$. 

    We define random variables $X_j \coloneqq |R_j(\F)[100F_0]|$, where $j=1,\dots,J$. Note that $|\F'[100F_0]| = \sum_{j=1}^J X_j$.

    Since the flakes of $\F$ are essentially distinct, we have $X_j = |\F[R_j^{-1}(100F_0)]| \lesssim 1$. In fact, there exists $M \lesssim 1$ so that $0 \leq X_j \leq M$  for each $j$.

    Next, we estimate $\E X_j$. By linearity of expectation, we have 
    $$\E X_j = \sum_{F \in \F} \mathcal{P}[R_j(F) \subset 100F_0].$$ 
    Since $R_j$ is a random rigid motion, {\color{black} we note that} for each $F \in \F$,
    \begin{align*}
        \mathcal{P}[R_j(F) \subset 100F_0] \lesssim |F_\delta|^{\color{black} k+1}.
    \end{align*}
    Therefore, 
    \begin{align*}
        m = \E X_j \lesssim |\F||F_\delta|^{k+1}.
    \end{align*}

    We claim that $Jm \lesssim 1$. To see this, first note that $Jm \lesssim C_F(\F)|\F||F_\delta|^{k+1}$. Suppose $\Delta_{max}(\F) \sim \Delta(\F,K)$ for some convex set $K \subset B_1$, then we have 
    $$C_F(\F) \sim \frac{\Delta(\F,K)}{\Delta(\F,B_1)} = \frac{|\F[K]|}{|\F||K|}.$$ 
    But since $\F$ is essentially distinct, we have $|\F[K]| \lesssim (|K|/|F_\delta|)^{\color{black} k+1}$. Therefore, we have
    \begin{align*}
        Jm \lesssim C_F(\F)|\F||F_\delta|^{k+1} \sim \frac{|\F[K]||F_\delta|^{k+1}}{|K|} \lesssim |K|^k \lesssim 1.
    \end{align*}
    Thus by possibly enlarging the previous $M$, we also have $Jm \leq M$ for $M \lesssim 1$.
    
    Therefore, if $S \geq 1$, then Lemma~\ref{Lem:Prob} gives
    \begin{align*}
        \mathcal{P}[X_1+\cdots+X_J>S] \leq e^{10} e^{-\frac{S}{2M}}
    \end{align*}
    for some $M \lesssim 1$.
    Since the number of essentially distinct $F_0$ is at most $\delta^{-O(1)}$, we can choose $S \lessapprox 1$ so that with high probability $|\F'[100F_0]|\leq S$ for every $F_0$.

    If we refine $\F'$ to be a set of essentially distinct flakes, then we still have $|\F'| \approx C_F(\F)|\F|$.

    Now we estimate $C_F(\F')$. Refining $\F'$ to make the flakes essentially distinct will only affect $C_F(\F')$ by a factor $\lessapprox 1$, so we can ignore the refinement.

    Let $K \subset B_2$ be a convex set of dimensions $a_1 \times \cdots \times a_{n-k} \times 1 \times \cdots \times 1$ with $\delta \leq a_1 \leq \cdots \leq a_{n-k} \leq 1$. We define random variables $Z_j \coloneqq |R_j(\F)[K]|$, where $j = 1, \dots, J$. Note that $|\F'[K]| = \sum_{j=1}^J Z_j$. We will show that for each $K$,
    \begin{align}\label{Eq:large_deviation}
        \mathcal{P}\Big[ \sum_{j=1}^J Z_j > A \Delta(\F',B_1)\frac{|K|}{|F_\delta|} \Big] \leq e^{10} e^{-A/C}.
    \end{align}
    
    Suppose this has been established. Since there are only $\delta^{-O(1)}$ essentially different choices for the convex set $K \subset B_2$, with probability at least $1 - \delta^{-O(1)}e^{10}e^{-A/C}$, we have $|\F'[K]| \lesssim A \Delta(\F',B_1)\frac{|K|}{|F_\delta|}$ for every $K$, which means $C_F(\F') \lesssim A$. We can choose $A \lessapprox 1$ so that the probability is at least $99/100$.

    We will justify \eqref{Eq:large_deviation} using Lemma~\ref{Lem:Prob}. Note that $\Delta_{max}(\F) = C_F(\F) \Delta(\F,B_1)\sim \Delta(\F',B_1) \sim |\F'||F_\delta|$. Thus
    \begin{align*}
        |\F'||K| \sim \frac{\Delta_{max}(\F)|K|}{|F_\delta|} \geq |\F[R_j^{-1}(K)]| = Z_j. 
    \end{align*}
    So there exists $M \sim |\F'||K|$ such that $0 \leq Z_j \leq M$. For each $F \in \F$, we claim that
    \begin{align*}
        \mathcal{P}[R_j(F) \subset K] \lesssim |K|.
    \end{align*}
    To see this, pick a point $x \in F$ (for example $x$ could be the center of $F$) and bound
    \begin{align*}
        \mathcal{P}[R_j(F) \subset K] \leq \mathcal{P}[R_j(x) \in K] \lesssim \frac{|K|}{|B_1(x)|} \lesssim |K|.
    \end{align*}
    
    By linearity of expectation, this immediately yields
    \begin{align*}
        \E Z_j \lesssim |\F||K|.
    \end{align*}
    Define $m = \E Z_j$. Since $J = \lceil C_F(\F) \rceil$, we have $Jm \lesssim C_F(\F) |\F||K| \sim |\F'||K| \sim M$.

    Therefore, Lemma~\ref{Lem:Prob} tells us that for every $S\geq 1$, 
    \begin{align*}
        \mathcal{P}[Z_1+\cdots+Z_J>S] \leq e^{10} e^{-\frac{S}{CM}}.
    \end{align*}
    By taking 
    $$S = A \Delta(\F',B_1)\frac{|K|}{|F_\delta|} \sim A|\F'||K| \sim AM,$$ 
    we get \eqref{Eq:large_deviation} as desired.
\end{proof}

\section*{Declaration of AI usage}

During the preparation of this manuscript, the authors made limited use of AI-assisted tools (OpenAI’s ChatGPT-5.5) for auxiliary purposes, namely to assist with sorting out preliminary literature and searching existing geometric and combinatorial results (Theorem~\ref{Thm:Principal_angle} and Lemma~\ref{Lem:booksize_graph}). All references identified through such tools were independently checked by the authors. All outputs were reviewed and verified by the author. Except for the uses described above, the formulation of the results, the proofs, the mathematical arguments, and the final exposition are entirely the work of the authors. OpenAI's ChatGPT-6.1 Sol was used at a later stage for language editing and manuscript consistency checking. The authors accept full responsibility for the originality, correctness, and integrity of all content contained in this manuscript.

\section*{Acknowledgments}
This work is partially supported by the AIM research community on the Fourier restriction conjecture. The authors would like to thank the organizers Po-Lam Yung, Alan Chang and Zane Li for bringing us together in an online research group. The authors would like to thank Joshua Zahl for insightful discussions at the beginning of this project, and Hong Wang for insightful conversations after an earlier version of this manuscript had been completed. The second author would like to thank Zhenbin Cao, Shengwen Gan, and Yumeng Ou for previous helpful discussions on the $2$-plane Kakeya problem. 

\bibliographystyle{plain}
\bibliography{references}

\end{document}